\documentclass{amsart}
\usepackage{amssymb}
\usepackage{mathrsfs}
\usepackage{stmaryrd} 
\usepackage{bbm}
\usepackage{oldgerm}
\usepackage[francais]{babel}
\usepackage[T1]{fontenc}
\usepackage[latin1]{inputenc}
\usepackage[all]{xy}
\usepackage{hyperref}

\makeindex

\newtheorem{teo}[subsection]{Théorème}
\newtheorem{prop}[subsection]{Proposition}
\newtheorem{cor}[subsection]{Corollaire}
\newtheorem{lem}[subsection]{Lemme}

\theoremstyle{definition}

\newtheorem{defi}[subsection]{Définition}
\newtheorem{rema}[subsection]{Remarque}
\newtheorem{remas}[subsection]{Remarques}

\numberwithin{equation}{subsection}

\newcommand{\gtimes}{\stackrel{\leftarrow}{\times}}

\newcommand{\mQ}{{\mathbb Q}}
\newcommand{\mN}{{\mathbb N}}

\newcommand{\mZ}{{\mathbb Z}}

\newcommand{\mK}{{\mathbb K}}
\newcommand{\mU}{{\mathbb U}}

\newcommand{\bA}{{\bf A}}
\newcommand{\bB}{{\bf B}}
\newcommand{\bC}{{\bf C}}
\newcommand{\bD}{{\bf D}}
\newcommand{\bE}{{\bf E}}
\newcommand{\bL}{{\bf L}}
\newcommand{\bT}{{\bf T}}
\newcommand{\bP}{{\bf P}}
\newcommand{\bQ}{{\bf Q}}
\newcommand{\Et}{{\bf \acute{E}t}}
\newcommand{\Sch}{{\bf Sch}}
\newcommand{\Ens}{{\bf Ens}}
\newcommand{\Pt}{{\bf Pt}}

\newcommand{\Top}{{\bf Top}}
\newcommand{\bHom}{{\bf Hom}}
\newcommand{\bExt}{{\bf Ext}}
\newcommand{\bRep}{{\bf Rep}}

\newcommand{\bMod}{{\bf Mod}}

\newcommand{\bpad}{{p{\text -{\bf ad}}}}

\newcommand{\crig}{{\cH^0_\rig}}

\newcommand{\et}{{\rm \acute{e}t}}
\newcommand{\fet}{{\rm f\acute{e}t}}
\newcommand{\rig}{{\rm rig}}
\newcommand{\ad}{{\rm ad}}

\newcommand{\MOD}{{\rm MOD}}

\newcommand{\zar}{{\rm zar}}

\newcommand{\coh}{{\rm coh}}
\newcommand{\cart}{{\rm cart}}
\newcommand{\scoh}{{\rm scoh}}
\newcommand{\Dolb}{{\rm Dolb}}
\newcommand{\sol}{{\rm sol}}

\newcommand{\tf}{{\rm tf}}

\newcommand{\aatf}{{\text {\rm -atf}}}
\newcommand{\atf}{{\rm atf}}

\newcommand{\rf}{{\rm f}}

\newcommand{\Spec}{{\rm Spec}}

\newcommand{\Spf}{{\rm Spf}}

\newcommand{\ob}{{\rm Ob}}

\newcommand{\tor}{{\rm tor}}

\newcommand{\disc}{{\rm disc}}
\newcommand{\coker}{{\rm coker}}

\newcommand{\tot}{{\rm tot}}

\newcommand{\cont}{{\rm cont}}

\newcommand{\gp}{{\rm gp}}
\newcommand{\id}{{\rm id}}

\newcommand{\rb}{{\rm b}}

\newcommand{\Hom}{{\rm Hom}}

\newcommand{\Ext}{{\rm Ext}}

\newcommand{\bMH}{{\bf MH}}
\newcommand{\bIH}{{\bf IH}}

\newcommand{\rE}{{\rm E}}

\newcommand{\rH}{{\rm H}}
\newcommand{\rT}{{\rm T}}

\newcommand{\rR}{{\rm R}}
\newcommand{\rS}{{\rm S}}

\newcommand{\rW}{{\rm W}}

\newcommand{\rp}{{\rm p}}

\newcommand{\oK}{{\overline{K}}}

\newcommand{\oR}{{\overline{R}}}
\newcommand{\oS}{{\overline{S}}}

\newcommand{\oU}{{\overline{U}}}

\newcommand{\oX}{{\overline{X}}}
\newcommand{\oY}{{\overline{Y}}}
\newcommand{\oZ}{{\overline{Z}}}
\newcommand{\oa}{{\overline{a}}}
\newcommand{\of}{{\overline{f}}}
\newcommand{\ogg}{{\overline{g}}}

\newcommand{\os}{{\overline{s}}}

\newcommand{\ox}{{\overline{x}}}
\newcommand{\oy}{{\overline{y}}}
\newcommand{\oz}{{\overline{z}}}

\newcommand{\oalpha}{{\overline{\alpha}}}

\newcommand{\oeta}{{\overline{\eta}}}

\newcommand{\otheta}{{\overline{\theta}}}

\newcommand{\oiota}{{\overline{\iota}}}

\newcommand{\ocB}{{\overline{\cB}}}

\newcommand{\ofp}{{\overline{\fp}}}

\newcommand{\un}{{\underline{n}}}
\newcommand{\upp}{{\underline{p}}}

\newcommand{\hu}{{\widehat{u}}}

\newcommand{\hA}{{\widehat{A}}}

\newcommand{\hE}{{\widehat{E}}}

\newcommand{\hoY}{\widehat{\oY}}
\newcommand{\hoZ}{\widehat{\oZ}}

\newcommand{\hRun}{{\widehat{R_1}}}

\newcommand{\hOmega}{\widehat{\Omega}}

\newcommand{\cA}{{\mathscr A}}
\newcommand{\cB}{{\mathscr B}}
\newcommand{\cC}{{\mathscr C}}
\newcommand{\cD}{{\mathscr D}}
\newcommand{\cE}{{\mathscr E}}
\newcommand{\cF}{{\mathscr F}}
\newcommand{\cG}{{\mathscr G}}

\newcommand{\cJ}{{\mathscr J}}
\newcommand{\cK}{{\mathscr K}}
\newcommand{\cL}{{\mathscr L}}
\newcommand{\cP}{{\mathscr P}}
\newcommand{\co}{{\mathscr O}}
\newcommand{\cR}{{\mathscr R}}
\newcommand{\cS}{{\mathscr S}}

\newcommand{\cH}{{\mathscr H}}
\newcommand{\cM}{{\mathscr M}}
\newcommand{\cN}{{\mathscr N}}
\newcommand{\cQ}{{\mathscr Q}}
\newcommand{\cU}{{\mathscr U}}
\newcommand{\cV}{{\mathscr V}}

\newcommand{\cHom}{{\mathscr Hom}}

\newcommand{\fF}{{\mathfrak F}}

\newcommand{\fM}{{\mathfrak M}}
\newcommand{\fN}{{\mathfrak N}}

\newcommand{\fS}{{\mathfrak S}}
\newcommand{\fV}{{\mathfrak V}}

\newcommand{\fX}{{\mathfrak X}}

\newcommand{\fgg}{{\mathfrak g}}
\newcommand{\fm}{{\mathfrak m}}

\newcommand{\fp}{{\mathfrak p}}
\newcommand{\fq}{{\mathfrak q}}

\newcommand{\ttT}{{\tt T}}
\newcommand{\tta}{{\tt a}}

\newcommand{\hM}{{\widehat{M}}}

\newcommand{\hR}{{\widehat{R}}}
\newcommand{\hX}{{\widehat{X}}}

\newcommand{\hoR}{{\widehat{\oR}}}

\newcommand{\hcC}{{\widehat{\cC}}}

\newcommand{\hiota}{{\widehat{\iota}}}

\newcommand{\hotimes}{{\widehat{\otimes}}}

\newcommand{\bvd}{{\breve{d}}}
\newcommand{\bvu}{{\breve{u}}}

\newcommand{\bvR}{{\breve{R}}}

\newcommand{\bvoS}{{\breve{\oS}}}

\newcommand{\bvoX}{{\breve{\oX}}}
\newcommand{\bvocB}{{\breve{\ocB}}}

\newcommand{\bvcC}{{\breve{\cC}}}

\newcommand{\bvcF}{{\breve{\cF}}}
\newcommand{\bvcH}{{\breve{\cH}}}
\newcommand{\bvcK}{{\breve{\cK}}}
\newcommand{\bvcM}{{\breve{\cM}}}
\newcommand{\bvcN}{{\breve{\cN}}}
\newcommand{\bvtau}{{\breve{\tau}}}
\newcommand{\bvnu}{{\breve{\nu}}}
\newcommand{\bvgamma}{{\breve{\gamma}}}
\newcommand{\bvalpha}{{\breve{\alpha}}}
\newcommand{\bvbeta}{{\breve{\beta}}}
\newcommand{\bvsigma}{{\breve{\sigma}}}

\newcommand{\bvrho}{{\breve{\rho}}}
\newcommand{\bvPhi}{{\breve{\Phi}}}
\newcommand{\bvpartial}{{\breve{\partial}}}

\newcommand{\bvtta}{{\breve{\tta}}}

\newcommand{\crH}{{\check{\rH}}}
\newcommand{\coS}{{\check{\oS}}}
\newcommand{\coU}{{\check{\oU}}}
\newcommand{\coX}{{\check{\oX}}}
\newcommand{\coY}{{\check{\oY}}}
\newcommand{\cog}{{\check{\ogg}}}
\newcommand{\calpha}{{\check{\alpha}}}
\newcommand{\cbeta}{{\check{\beta}}}

\newcommand{\tE}{{\widetilde{E}}}

\newcommand{\tS}{{\widetilde{S}}}

\newcommand{\tU}{{\widetilde{U}}}

\newcommand{\tX}{{\widetilde{X}}}
\newcommand{\tY}{{\widetilde{Y}}}

\newcommand{\te}{{\widetilde{e}}}
\newcommand{\ttf}{{\widetilde{f}}}
\newcommand{\tg}{{\widetilde{g}}}

\newcommand{\tcE}{{\widetilde{\cE}}}
\newcommand{\trT}{{\widetilde{\rT}}}
\newcommand{\tOmega}{{\widetilde{\Omega}}}

\newcommand{\tkappa}{{\widetilde{\kappa}}}
\newcommand{\ttheta}{{\widetilde{\theta}}}

\newcommand{\talpha}{{\widetilde{\alpha}}}

\newcommand{\tmK}{{\widetilde{\mK}}}

\begin{document}

\title[Sur la correspondance de Simpson $p$-adique. II]
{Sur la correspondance de Simpson $p$-adique.\\
II~: aspects globaux}
\author{Ahmed Abbes et Michel Gros}
\address{A.A. Laboratoire Alexander Grothendieck, ERL 9216 du CNRS, Institut des Hautes \'Etudes Scientifiques, 35 route de Chartres, 91440 Bures-sur-Yvette, France}
\address{M.G. CNRS UMR 6625, IRMAR, Université de Rennes 1,
Campus de Beaulieu, 35042 Rennes cedex, France}
\email{abbes@ihes.fr}
\email{michel.gros@univ-rennes1.fr}

\maketitle

\setcounter{tocdepth}{1}
\tableofcontents

\section{Introduction}
Nous poursuivons dans cet article la construction et l'étude de la correspondance de Simpson $p$-adique 
initiée dans \cite{ag1}, suivant l'approche générale résumée dans \cite{ag0}. 
Après avoir fixé les notations et conventions générales au § \ref{higgs3-not}, 
nous développons dans les sections \ref{higgs3-sli} à \ref{higgs3-spsa} quelques préliminaires utiles pour la suite.  
La section  \ref{higgs3-sli} contient des rappels et compléments sur la notion de schémas localement irréductibles. 
Nous précisons au § \ref{higgs3-slad} le cadre de géométrie logarithmique dans lequel seront placées nos constructions. 
L'intermède § \ref{higgs3-vck} regroupe quelques sorites utilisés à divers endroits du texte sur le complexe de Koszul. 
La section \ref{higgs3-caip} développe le formalisme des catégories 
additives à isogénie près. La section \ref{higgs3-spsa} est consacrée à l'étude des systèmes projectifs d'un topos annelé, 
en particulier à la notion de module adique et aux conditions de finitude appropriées à ce cadre.

Soit $K$ un corps de valuation discrète complet de 
caractéristique $0$, à corps résiduel algébriquement clos de caractéristique $p$ et soit $\oK$ une clôture algébrique de $K$.
On note $\co_K$ l'anneau de valuation de $K$, $\co_\oK$ la clôture intégrale de $\co_K$ dans $\oK$
et $\co_C$ le séparé complété $p$-adique de $\co_\oK$. On pose $S=\Spec(\co_K)$ et on le munit
de la structure logarithmique $\cM_S$ définie par son point fermé. 
On considère dans cet article un schéma logarithmique $(X,\cM_X)$ lisse et saturé au-dessus de $(S,\cM_S)$,
vérifiant une condition locale \eqref{higgs3-slad6} correspondant aux hypothèses faites dans la première partie 
(\cite{ag1} 6.2). On note $X^\circ$ le sous-schéma ouvert maximal de $X$ où la structure logarithmique 
$\cM_X$ est triviale. Le cadre topologique dans lequel se réalise la correspondance de Simpson $p$-adique 
est celui du topos de Faltings $\tE$ associé au morphisme canonique $X^\circ\otimes_K\oK\rightarrow X$, 
dont l'étude détaillée a été menée indépendamment dans \cite{ag2}. 
Nous l'équipons dans § \ref{higgs3-tfa} d'un anneau $\ocB$. Nous introduisons ensuite dans § \ref{higgs3-TFT}
le topos $\tE_s$ fibre spéciale de $\tE$ et le topos annelé $(\tE_s^{\mN^\circ},\bvocB)$ complété 
formel $p$-adique de $(\tE,\ocB)$, dans lequel auront lieu nos principales constructions. 

Comme il a été esquissé dans l'introduction générale \cite{ag0}, la correspondance de Simpson $p$-adique 
dépend d'une déformation (logarithmique) lisse de $X\otimes_{\co_K}\co_C$ au-dessus 
de l'épaississement infinitésimal $p$-adique universel de $\co_C$ d'ordre $\leq 1$, introduit par Fontaine (\cite{ag1} 9.8).  
Nous supposons dans la suite de cette introduction qu'il existe une telle déformation que nous fixons. 
Nous définissons dans la section \ref{higgs3-RGG}, pour tout nombre rationnel $r\geq 0$, 
la {\em $\bvocB$-algèbre de Higgs-Tate d'épaisseur $r$}, notée $\bvcC^{(r)}$. 
Ces algèbres forment un système inductif~: pour tous nombres rationnels $r\geq r'\geq 0$,
on a un homomorphisme canonique $\bvcC^{(r)}\rightarrow \bvcC^{(r')}$.
Ce sont des analogues faisceautiques des algèbres de Higgs-Tate introduites dans la première partie 
(\cite{ag1} 12.1). Elles sont naturellement munies de champs de Higgs. 
La section \ref{higgs3-coh} contient deux résultats d'acyclicité fondamentaux pour la suite. 
Nous démontrons \eqref{higgs3-coh22} que la limite inductive des complexes de Dolbeault de $\bvcC^{(r)}$, pour $r\in \mQ_{>0}$,
est une résolution de $\bvocB$, à isogénie près. D'autre part, notant $\fX$ le schéma formel complété $p$-adique de 
$X\otimes_{\co_K}\co_\oK$, on dispose d'un morphisme canonique de topos annelés 
\[
\top\colon (\tE_s^{\mN^\circ},\bvocB)\rightarrow (\fX_\zar,\co_{\fX}). 
\]
Nous démontrons \eqref{higgs3-coh17} que l'homomorphisme canonique 
\[
\co_{\fX}[\frac 1 p]\rightarrow \underset{\underset{r\in \mQ_{>0}}{\longrightarrow}}{\lim}\  \top_*(\bvcC^{(r)})[\frac 1 p]
\]
est un isomorphisme, et que pour tout entier $q\geq 1$,
\[
\underset{\underset{r\in \mQ_{>0}}{\longrightarrow}}{\lim}\ \rR^q\top_*(\bvcC^{(r)})[\frac 1 p] =0.
\]

La section \ref{higgs3-MF} est consacrée à la construction de la correspondance de Simpson $p$-adique. 
Nous introduisons les notions de module de Dolbeault et de fibré de Higgs soluble \eqref{higgs3-MF14} 
et nous montrons \eqref{higgs3-MF21} qu'elles donnent lieu à deux catégories équivalentes.
Nous construisons en fait deux équivalences de catégories explicites et quasi-inverses l'une de l'autre. 
Nous démontrons aussi la compatibilité de cette correspondance au passage à la cohomologie \eqref{higgs3-MF27}. 
Nous établissons dans §~\ref{higgs3-AFE} un lien avec les constructions 
développées dans la première partie \cite{ag1} pour des schémas affines d'un type particulier, qualifiés par Faltings de {\em petits}.
Nous étudions dans §~\ref{higgs3-FHT} la fonctorialité de la correspondance de Simpson $p$-adique par des morphismes étales. 
La section §~\ref{higgs3-PMH} est consacrée à l'étude de la catégorie fibrée des modules de Dolbeault 
au-dessus du site étale restreint de $X$.  
Nous montrons le caractère local pour la topologie étale de $X$ de la propriété de Dolbeault pour les modules \eqref{higgs3-CMD4}.
L'assertion analogue pour la propriété soluble pour les fibrés de Higgs est équivalente au fait que la catégorie fibrée 
des modules de Dolbeault  soit un champ \eqref{higgs3-CMD5}.
On dit qu'un fibré de Higgs est {\em petit} si son champ de Higgs vérifie une certaine condition de divisibilité 
relativement à un réseau, et qu'il est {\em localement petit} si cette condition est vérifiée localement \eqref{higgs3-PMH1}.
Nous montrons que tout fibré de Higgs soluble est localement petit \eqref{higgs3-PMH10}.
Inversement, si les modules de Dolbeault forment un champ, tout fibré de Higgs localement petit est soluble   \eqref{higgs3-PMH7}.
Pour un schéma affine petit, nous montrons inconditionnellement que tout fibré de Higgs petit est soluble \eqref{higgs3-PMH6}.

\section{Notations et conventions}\label{higgs3-not}

{\em Tous les anneaux considérés dans cet article possèdent un élément unité~;
les homomorphismes d'anneaux sont toujours supposés transformer l'élément unité en l'élément unité.
Nous considérons surtout des anneaux commutatifs, et lorsque nous parlons d'anneau 
sans préciser, il est sous-entendu qu'il s'agit d'un anneau commutatif~; en particulier, 
il est sous-entendu, lorsque nous parlons d'un topos annelé $(X,A)$ sans préciser, que $A$ est commutatif.}

\subsection{}\label{higgs3-not1}\index{1000201@$K$, $\oK$, $C$}
\index{1000203@$(S,\cM_S)$, $(\oS,\cM_\oS)$, $(\coS,\cM_\coS)$}\index{1000205@$\cS$}
\index{1000207@$\oX$, $\coX$, $X_n$ ($X$ un $S$-schéma)}
Dans cet article, $p$ désigne un nombre premier, $K$ un corps de valuation discrète complet de 
caractéristique $0$, à corps résiduel {\em algébriquement clos} $k$ de caractéristique $p$, 
et $\oK$ une clôture algébrique de $K$.
On note $\co_K$ l'anneau de valuation de $K$, $\co_\oK$ la clôture intégrale de $\co_K$ dans $\oK$, 
$\fm_\oK$ l'idéal maximal de $\co_\oK$ et $v$ la valuation de $\oK$ normalisée 
par $v(p)=1$. On désigne par $\co_C$ le séparé complété $p$-adique de $\co_\oK$, par $C$ son corps des fractions
et par $\fm_C$ son idéal maximal. Sauf mention explicite du contraire, on considère $\co_C$ comme un anneau adique, 
muni de la topologie $p$-adique (\cite{egr1} 1.8.7); c'est un anneau $1$-valuatif (\cite{egr1} 1.9.9). 

On choisit un système compatible $(\beta_n)_{n>0}$ de racines $n$-ièmes de $p$ dans $\co_\oK$. 
Pour tout nombre rationnel $\varepsilon>0$, 
on pose $p^\varepsilon=(\beta_n)^{\varepsilon n}$, où $n$ est un entier $>0$ tel que $\varepsilon n$ soit entier.

Pour tout groupe abélien $A$, on note $\hA$ son séparé complété $p$-adique.

On pose $S=\Spec(\co_K)$, $\oS=\Spec(\co_\oK)$ et $\coS=\Spec(\co_C)$. 
On note $s$ (resp. $\eta$, resp. $\oeta$) le point fermé de $S$ (resp. générique de $S$, resp. générique de $\oS$).
Pour tout entier $n\geq 1$, on pose $S_n=\Spec(\co_K/p^n\co_K)$. Pour tout $S$-schéma $X$, on pose
\begin{equation}\label{higgs3-not1c}
\oX=X\times_S\oS, \ \ \ \coX=X\times_S\coS\ \ \ {\rm et}\ \ \ X_n=X\times_SS_n. 
\end{equation}

On munit $S$ de la structure logarithmique $\cM_S$ définie par son point fermé, 
autrement dit, $\cM_S=j_*(\co_\eta^\times)\cap \co_S$, où $j\colon \eta\rightarrow S$ est l'injection canonique
(cf. \cite{ag1} 5.9). On notera qu'un homomorphisme de monoïdes $\iota\colon \mN\rightarrow \Gamma(S,\cM_S)$
est une carte pour $(S,\cM_S)$ (\cite{ag1} 5.13) si et seulement si $\iota(1)$ est une uniformisante de $\co_K$. 

On munit $\oS$ et $\coS$ des structures logarithmiques $\cM_\oS$ et $\cM_\coS$ 
images inverses de $\cM_S$ (cf. \cite{ag1} 5.10).

On désigne par $\cS=\Spf(\co_C)$ le schéma formel complété $p$-adique de $\oS$ ou, ce qui revient au même, de $\coS$.

\subsection{}\label{higgs3-not5}\index{1000210@$\cR_A$}\index{1000211@$\cA_2(A)$}\index{1000212@$\theta\colon \cA_2(A)\rightarrow \hA$}
Rappelons  (\cite{ag1} 9.3) que Fontaine associe fonctoriellement à toute $\mZ_{(p)}$-algèbre $A$ l'anneau 
\begin{equation}\label{higgs3-not5a}
\cR_A=\underset{\underset{x\mapsto x^p}{\longleftarrow}}{\lim}A/pA,
\end{equation} 
et un homomorphisme $\theta$ de l'anneau $\rW(\cR_A)$
des vecteurs de Witt de $\cR_A$ dans le séparé complété $p$-adique $\hA$ de $A$. On pose 
\begin{equation}\label{higgs3-not5b}
\cA_2(A)=\rW(\cR_A)/\ker(\theta)^2
\end{equation} 
et on note encore $\theta\colon \cA_2(A)\rightarrow \hA$ l'homomorphisme
induit par $\theta$. 

\subsection{}\label{higgs3-not4} \index{1000215@$\xi\in \rW(\cR_{\co_\oK})$}\index{1000217@$i_\oS\colon (\coS,\cM_\coS)\rightarrow (\cA_2(\oS),\cM_{\cA_2(\oS)})$}
Dans cet article, on se donne une suite $(p_n)_{n\geq 0}$ d'éléments de $\co_\oK$ telle que 
$p_0=p$ et $p_{n+1}^p=p_n$ pour tout $n\geq 0$.  
On désigne par $\upp$ l'élément de $\cR_{\co_\oK}$ \eqref{higgs3-not5a} induit par la suite $(p_n)_{n\geq 0}$ et on pose  
\begin{equation}\label{higgs3-not4a}
\xi=[\upp]-p \in \rW(\cR_{\co_\oK}),
\end{equation}
où $[\ ]$ est le représentant multiplicatif.  
L'homomorphisme $\theta\colon \rW(\cR_{\co_\oK})\rightarrow \co_C$ est surjectif et son noyau 
est engendré par $\xi$, qui n'est pas un diviseur de zéro dans $\rW(\cR_{\co_\oK})$ (\cite{ag1} 9.5). 
On a donc une suite exacte 
\begin{equation}\label{higgs3-not4b}
0\longrightarrow \co_C\stackrel{\cdot \xi}{\longrightarrow} \cA_2(\co_\oK)\stackrel{\theta}{\longrightarrow} 
\co_\oK\longrightarrow 0,
\end{equation}
où on a encore noté $\cdot \xi$ le morphisme induit par la multiplication par $\xi$ dans $\cA_2(\oR)$. 
L'idéal $\ker(\theta)$ de $\cA_2(\co_\oK)$ est de carré nul. 
C'est un $\co_C$-module libre de base $\xi$. Il sera noté $\xi\co_C$. 
On observera que contrairement à $\xi$, ce module ne dépend pas du choix de la suite $(p_n)_{n\geq 0}$. 

On note $\xi^{-1}\co_C$ le $\co_C$-module dual de $\xi\co_C$. 
Pour tout $\co_C$-module $M$, on désigne les $\co_C$-modules $M\otimes_{\co_C}(\xi \co_C)$ 
et $M\otimes_{\co_C}(\xi^{-1} \co_C)$ simplement par $\xi M$ et $\xi^{-1} M$, respectivement. 
On observera que contrairement à $\xi$, ces modules ne dépendent pas du choix de la suite $(p_n)_{n\geq 0}$. 
Il est donc important de ne pas les identifier à $M$. 

On pose 
\begin{equation}\label{higgs3-not4c}
\cA_2(\oS)=\Spec(\cA_2(\co_\oK))
\end{equation} 
que l'on munit de la structure logarithmique $\cM_{\cA_2(\oS)}$
définie dans  (\cite{ag1} 9.8). Le schéma logarithmique $(\cA_2(\oS),\cM_{\cA_2(\oS)})$ est 
alors fin et saturé, et $\theta$ induit une immersion fermée exacte 
\begin{equation}\label{higgs3-not4d}
i_\oS\colon (\coS,\cM_\coS)\rightarrow (\cA_2(\oS),\cM_{\cA_2(\oS)}).
\end{equation}

\subsection{}\label{higgs3-not7}
Pour toute catégorie abélienne $\bA$, on désigne par 
$\bD(\bA)$ sa catégorie dérivée et par $\bD^-(\bA)$, $\bD^+(\bA)$ et $\bD^\rb(\bA)$ les sous-catégories
pleines de $\bD(\bA)$ formées des complexes à cohomologie bornée
supérieurement, inférieurement et des deux côtés, respectivement. 
Sauf mention expresse du contraire, les complexes de $\bA$ sont à différentielles de degré $+1$,
le degré étant écrit en exposant. 

\subsection{}\label{higgs3-not0}
Dans tout cet article, on fixe un univers $\mU$ possédant un élément de cardinal infini. 
On appelle catégorie des $\mU$-ensembles et l'on note $\Ens$, 
la catégorie des ensembles qui se trouvent dans $\mU$. 
Sauf mention explicite du contraire, il sera sous-entendu que les schémas 
envisagés dans cet article sont éléments de l'univers $\mU$.
On désigne par $\Sch$ la catégorie des schémas éléments de $\mU$. 

\subsection{}\label{higgs3-not8}
Suivant les conventions de (\cite{sga4} VI), nous utilisons l'adjectif {\em cohérent} 
comme synonyme de quasi-compact et quasi-séparé. 

\subsection{}\label{higgs3-not9}\index{1000220@$\bMod(A,X)$ ou $\bMod(A)$}
\index{1000222@$\rS_A(M)$ ou $\rS(M)$, $\wedge_A(M)$ ou $\wedge(M)$}
Soit $(X,A)$ un topos annelé. On note $\bMod(A)$ ou $\bMod(A,X)$ 
la catégorie des $A$-modules de $X$. Si $M$ est un $A$-module, on désigne par $\rS_A(M)$ 
(resp. $\wedge_A(M)$, resp. $\Gamma_A(M)$) l'algèbre symétrique (resp. extérieure, resp. à puissances divisées) 
de $M$ (\cite{illusie1} I 4.2.2.6) et pour tout entier $n\geq 0$, par $\rS_A^n(M)$ (resp. $\wedge_A^n(M)$,
resp. $\Gamma_A^n(M)$) sa partie homogène de degré $n$. 
On omettra l'anneau $A$ des notations lorsqu'il n'y a aucun risque d'ambiguïté. 
Les formations de ces algèbres commutent à la localisation au-dessus d'un objet de $X$.

\begin{defi}[\cite{sga6} I 1.3.1]\label{higgs3-not3}\index{Module localement projectif de type fini}
Soit $(X,A)$ un topos annelé. 
On dit qu'un $A$-module $M$ de $X$ est {\em localement projectif de type fini} si les conditions 
équivalentes suivantes sont satisfaites :
\begin{itemize}
\item[{\rm (i)}] $M$ est de type fini et le foncteur $\cHom_A(M,\cdot)$ est exact~;
\item[{\rm (ii)}] $M$ est de type fini et tout épimorphisme de $A$-modules $N\rightarrow M$ admet localement une section~;
\item[{\rm (iii)}] $M$ est localement facteur direct d'un $A$-module libre de type fini. 
\end{itemize}
\end{defi}

Lorsque $X$ a suffisamment de points et que pour tout point $x$ de $X$, la fibre de $A$ en $x$
est un anneau local, les $A$-modules localement projectifs de type fini sont les $A$-modules 
localement libres de type fini (\cite{sga6} I 2.15.1).

\subsection{}\label{higgs3-not2}\index{1000230@$\Et_{/X}$, $X_\et$ (site et topos étales)}\index{1000232@$\Et_{\rf/X}$, $X_\fet$ (site et topos finis étales)}
\index{1000234@$\rho_X\colon X_\et\rightarrow X_\fet$}\index{1000236@$u_X\colon X_{\et}\rightarrow X_{\zar}$}
Pour tout schéma $X$, on note $\Et_{/X}$ (resp. $X_\et$) le site (resp. topos) étale de $X$. 
On désigne par $\Et_{\rf/X}$ le site fini étale de $X$, 
c'est-à-dire, la sous-catégorie pleine de $\Et_{/X}$ formée des schémas étales finis sur $X$, 
munie de la topologie induite par celle de $\Et_{/X}$,  
et par $X_\fet$ le topos fini étale de $X$, c'est-à-dire, le topos des faisceaux de $\mU$-ensembles sur $\Et_{\rf/X}$ 
(cf. \cite{ag2} 9.2). L'injection canonique $\Et_{\rf/X}\rightarrow \Et_{/X}$ induit un morphisme de topos 
\begin{equation}\label{higgs3-not2a}
\rho_X\colon X_\et\rightarrow X_\fet.
\end{equation}

On désigne par $X_\zar$ le topos de Zariski de $X$ et par 
\begin{equation}\label{higgs3-not2b}
u_X\colon X_{\et}\rightarrow X_{\zar}
\end{equation}
le morphisme canonique (\cite{sga4} VII 4.2.2). Si $F$ est un $\co_X$-module quasi-cohérent de $X_\zar$, 
on note encore $F$ le faisceau de $X_\et$ défini pour tout $X$-schéma étale $X'$ par (\cite{sga4} VII 2 c))
\begin{equation}\label{higgs3-not2c}
F(X')=\Gamma(X',F\otimes_{\co_X}\co_{X'}).
\end{equation}
Cet abus de notation n'induit aucune confusion. On a un isomorphisme canonique 
\begin{equation}\label{higgs3-not2d}
u_{X*}(F)\stackrel{\sim}{\rightarrow}F.
\end{equation}
Nous considérons donc $u_X$ comme un morphisme du topos annelé $(X_{\et},\co_X)$ vers le topos annelé  
$(X_{\zar},\co_X)$. Nous utilisons pour les modules la notation $u^{-1}_X$ pour désigner l'image
inverse au sens des faisceaux abéliens et nous réservons la notation 
$u^*_X$ pour l'image inverse au sens des modules. L'isomorphisme \eqref{higgs3-not2d} induit par adjonction un morphisme
\begin{equation}\label{higgs3-not2e}
u^*_X(F)\rightarrow F.
\end{equation}
Celui-ci est un isomorphisme si $F$ est un $\co_X$-module de présentation finie. En effet, la question étant locale,
on peut se borner au cas où il existe une suite exacte de $\co_X$-modules 
$\co_X^m\rightarrow \co_X^n\rightarrow F\rightarrow 0$ de $X_\zar$. Celle-ci induit une suite exacte 
de $\co_X$-modules $\co_X^m\rightarrow \co_X^n\rightarrow F\rightarrow 0$ de $X_\et$.
L'assertion s'ensuit compte tenu de l'exactitude à droite du foncteur $u^*_X$.

\subsection{}\label{higgs3-not6}\index{Revetement universel normalise@Revêtement universel normalisé}
\index{Foncteur fibre d'un topos fini etale@Foncteur fibre d'un topos fini étale}\index{1000240@$\nu_\ox\colon X_\fet\rightarrow \bB_{\pi_1(X,\ox)}$}
Soient $X$ un schéma connexe, $\ox$ un point géométrique de $X$.
On désigne par 
\begin{equation}\label{higgs3-not6a}
\omega_\ox\colon \Et_{\rf/X}\rightarrow \Ens
\end{equation}
le foncteur fibre en $\ox$, qui à tout revêtement étale $Y$ de $X$ associe l'ensemble des points géométriques de 
$Y$ au-dessus de $\ox$, par  $\pi_1(X,\ox)$ le groupe fondamental de $X$ en $\ox$ (c'est-à-dire
le groupe des automorphismes du foncteur $\omega_\ox$) et par $\bB_{\pi_1(X,\ox)}$ 
le topos classifiant du groupe profini $\pi_1(X,\ox)$, 
c'est-à-dire la catégorie des $\mU$-ensembles discrets munis d'une action continue à gauche de $\pi_1(X,\ox)$ 
(\cite{sga4} IV 2.7). Alors $\omega_\ox$ induit un foncteur pleinement fidèle 
\begin{equation}\label{higgs3-not6b}
\mu_\ox^+\colon \Et_{\rf/X}\rightarrow \bB_{\pi_1(X,\ox)}
\end{equation}
d'image essentielle la sous-catégorie pleine de $\bB_{\pi_1(X,\ox)}$ formée des ensembles finis
(\cite{sga1} V  §~4 et §~7).
Soit $(X_i)_{i\in I}$ un système projectif sur un ensemble ordonné filtrant $I$ dans $\Et_{\rf/X}$
qui pro-représente $\omega_\ox$, normalisé par le fait que les morphismes de transition $X_i\rightarrow X_j$
$(i\geq j)$ sont des épimorphismes et que tout épimorphisme $X_i\rightarrow X'$ de $\Et_{\rf/X}$ 
est équivalent à un épimorphisme $X_i\rightarrow X_j$ $(j\leq i)$ convenable. 
Un tel pro-objet est essentiellement unique. Il est appelé le {\em revêtement universel normalisé de $X$ en $\ox$} 
ou le {\em pro-objet fondamental normalisé de $\Et_{\rf/X}$ en $\ox$}. 
On notera que l'ensemble $I$ est $\mU$-petit. Le foncteur 
\begin{equation}\label{higgs3-not6c}
\nu_\ox\colon X_\fet\rightarrow \bB_{\pi_1(X,\ox)},
\ \ \ F\mapsto \underset{\underset{i\in I}{\longrightarrow}}{\lim}\ F(X_i)
\end{equation}
est une équivalence de catégories qui prolonge le foncteur $\mu_\ox^+$ (cf. \cite{ag2} 9.8). 
On l'appelle le {\em foncteur fibre} de $X_\fet$ en $\ox$.

\subsection{}\label{higgs3-not60}
Conservons les hypothèses de \ref{higgs3-not6}, soit de plus $R$ un anneau de $X_\fet$. 
Posons $R_\ox=\nu_\ox(R)$ qui est un anneau muni de la topologie discrète et d'une action 
continue de $\pi_1(X,\ox)$ par des homomorphismes d'anneaux. On désigne par 
$\bRep_{R_\oy}^{\disc}(\pi_1(Y,\oy))$ la catégorie 
des $R_\oy$-représentations continues de $\pi_1(Y,\oy)$ pour lesquelles la topologie est discrète (\cite{ag1} 3.1). 
En restreignant le foncteur $\nu_\ox$ aux $R$-modules, on obtient une équivalence 
de catégories que l'on note encore
\begin{equation}\label{higgs3-not60a}
\nu_\ox\colon \bMod(R)\stackrel{\sim}{\rightarrow} \bRep_{R_\ox}^{\disc}(\pi_1(X,\ox)). 
\end{equation}
Pour qu'un $R$-module $M$ de $X_\fet$ soit de type fini (resp. localement projectif de type fini \eqref{higgs3-not3}), 
il faut et il suffit que le $R_\ox$-module sous-jacent à $\nu_\ox(M)$ soit de type fini (resp. projectif de type fini).  
En effet, la condition est nécessaire en vertu de (\cite{ag2} 9.9). 
Montrons qu'elle est suffisante. Supposons d'abord le $R_\ox$-module $\nu_\ox(M)$ de type fini.
Soient $N$ un $R_\ox$-module libre de base $e_1,\dots,e_d$, $u\colon N\rightarrow \nu_\ox(M)$ 
un épimorphisme $R_\ox$-linéaire. L'assertion recherchée étant locale pour $X_\fet$, 
quitte à remplacer $X$ par un revêtement étale, on peut supposer que $\pi_1(X,\ox)$ 
fixe les éléments $u(e_1),\dots,u(e_d)$ de $\nu_\ox(M)$. Munissant $N$ de l'unique 
$R_\ox$-représentation de $\pi_1(X,\ox)$ telle que $e_1,\dots,e_d$ soient fixes, l'homomorphisme 
$u$ est alors $\pi_1(X,\ox)$-équivariant. Par suite, $M$ est un $R$-module de type fini.  
Supposons de plus le $R_\ox$-module $\nu_\ox(M)$ projectif de type fini.
Soit $v\colon \nu_\ox(M)\rightarrow N$ un scindage 
$R_\ox$-linéaire de $u$.  Quitte à remplacer de nouveau $X$ par un revêtement étale, on peut
supposer que $\pi_1(X,\ox)$ fixe les éléments $v(u(e_1)),\dots,v(u(e_d))$ de $N$. Par suite, $v$
est $\pi_1(X,\ox)$-équivariant. On en déduit que $M$ est facteur direct 
d'un $R$-module libre de type fini~; d'où l'assertion.

\section{Schémas localement irréductibles}\label{higgs3-sli}

\subsection{}\label{higgs3-sli1}\index{Schema@Schéma!localement irreductible@localement irréductible}
Soit $X$ un schéma dont l'ensemble des composantes irréductibles est localement fini. 
On rappelle que les conditions suivantes sont équivalentes (\cite{ega1n} 0.2.1.6)~:
\begin{itemize}
\item[{\rm (i)}] Les composantes irréductibles de $X$ sont ouvertes. 
\item[{\rm (ii)}]  Les composantes irréductibles de $X$ sont identiques à ses composantes connexes. 
\item[{\rm (iii)}] Les composantes connexes de $X$ sont irréductibles. 
\item[{\rm (iv)}] Deux composantes irréductibles distinctes de $X$ ne se rencontrent pas.
\end{itemize}
Le schéma $X$ est alors la somme des schémas induits sur ses composantes irréductibles. 
Lorsque ces conditions sont remplies, on dit que $X$ est {\em localement irréductible}.
Cette notion est clairement locale sur $X$, {\em i.e.}, si $(X_i)_{i\in I}$ est un recouvrement ouvert de $X$, 
pour que $X$ soit localement irréductible, il faut et il suffit que pour tout $i\in I$, 
il en soit de même pour $X_i$. 

\begin{remas}\label{higgs3-sli2}\index{Schema@Schéma!etale localement connexe@étale-localement connexe}
(i)\ L'ensemble des composantes irréductibles d'un schéma localement noethérien est localement fini 
(\cite{ega1n} 0.2.2.2).

(ii)\ Pour qu'un schéma normal soit localement irréductible, il faut et il suffit que  
l'ensemble de ses composantes irréductibles soit localement fini. En effet, tout schéma normal vérifie clairement la condition \ref{higgs3-sli1}(iv).

(iii)\ Un schéma localement irréductible $X$ est étale-localement connexe, {\em i.e.}, 
pour tout morphisme étale $X'\rightarrow X$, toute composante connexe de $X'$ est un ensemble ouvert dans $X'$ 
(\cite{ag2} 9.7.3). 
\end{remas}

\begin{lem}\label{higgs3-sli3}
Soient $X$ un schéma normal et localement irréductible, $f\colon Y\rightarrow X$ 
un morphisme étale. Alors $Y$ est normal et localement irréductible. 
\end{lem}
En effet, $Y$ est normal en vertu de (\cite{raynaud1} VII prop.~2). 
Il suffit donc de montrer que l'ensemble de ses composantes irréductibles est localement fini d'après \ref{higgs3-sli2}(ii).
La question étant locale sur $X$ et $Y$, on peut se borner au cas où ils sont affines,
de sorte que $f$ est quasi-compact et par suite quasi-fini. L'assertion recherchée résulte alors de (\cite{ega4} 2.3.6(iii)).

\begin{lem}\label{higgs3-sli4}
Soient $X$ un schéma normal, $j\colon U\rightarrow X$ une immersion ouverte dense et quasi-compacte. 
Pour que $X$ soit localement irréductible, il faut et il suffit qu'il en soit de même de~$U$.
\end{lem}

En effet, si $X$ est localement irréductible, il en est de même de $U$. Inversement, supposons $U$ 
localement irréductible et montrons qu'il en est de même de $X$. On peut clairement se borner au cas où $X$ 
est quasi-compact. Par suite, $U$ est quasi-compact et n'a donc qu'un nombre fini de composantes irréductibles
d'après \ref{higgs3-sli1}(i). Il en est alors de même de $X$ puisque $X$ et $U$ ont mêmes points génériques.
Comme $X$ est normal, il est localement irréductible en vertu de \ref{higgs3-sli2}(ii).

\begin{lem}\label{higgs3-sli5}
Soit $A$ un anneau local hensélien, $B$ une $A$-algèbre intègre et entière sur $A$. 
Alors $B$ est un anneau local hensélien. Si, de plus, $A$ est strictement local, il en est de même de $B$. 
\end{lem}
Considérons $B$ comme une limite inductive filtrante de sous-$A$-algèbres 
de type fini $(B_i)_{i\in I}$. Pour tout $i\in I$, $B_i$ étant intègre et fini sur $A$, il est local et hensélien.
Pour tous $(i,j)\in I^2$ tel que $i\leq j$, le morphisme de transition $B_i\rightarrow B_j$ étant fini, il est local.
Par suite, $B$ est local et hensélien (\cite{raynaud1} I § 3 prop.~1).
Supposons $A$ strictement local. Comme l'homomorphisme $A\rightarrow B$ est local,
le corps résiduel de $B$ est une extension algébrique de celui de $A$. 
Il est donc séparablement clos. Par suite, $B$ est strictement local.

\begin{lem}\label{higgs3-sli7}
Soient $X$ un schéma, $\ox$ un point géométrique de $X$, $X'$ le localisé strict de $X$
en $\ox$, $Y$ un $X$-schéma, $Y'=Y\times_XX'$, $f\colon Y'\rightarrow Y$ la projection canonique. 
Alors l'homomorphisme canonique $f^{-1}(\co_Y)\rightarrow \co_{Y'}$ est un isomorphisme de $Y'_\et$.
\end{lem}
Considérons $X'$ comme une limite projective cofiltrante de voisinages étales affines 
$(X_i)_{\in I}$ de $\ox$ dans $X$ (cf. \cite{sga4} VIII 4.5) et posons $Y_i=Y\times_XX_i$ pour tout $i\in I$.  
Le schéma $Y'$ est alors la limite projective des schémas $(Y_i)_{\in I}$ (\cite{ega4} 8.2.5). 
Soit $\oy$ un point géométrique de $Y'$. Pour tout $i\in I$, la projection canonique $Y_i\rightarrow Y$ induit 
un isomorphisme entre les localisés stricts de $Y_i$ et $Y$ en les images canoniques de $\oy$. D'autre part, il résulte de  
(\cite{ega1n} 0.6.1.6) et (\cite{raynaud1} I § 3 prop.~1) que le localisé strict de $Y'$ en $\oy$ 
est la limite projective des localisés 
stricts des $Y_i$ en les images canoniques de $\oy$. Par conséquent, la projection canonique $Y'\rightarrow Y$ induit 
un isomorphisme entre les localisés stricts de $Y'$ en $\oy$ et $Y$ en $f(\oy)$; autrement dit, 
l'homomorphisme canonique 
$\co_{Y,f(\oy)}\rightarrow \co_{Y',\oy}$ est un isomorphisme~; d'où la proposition.

\begin{lem}\label{higgs3-sli6}
Soient $f\colon Y\rightarrow X$ un morphisme entier de schémas, 
$\ox$ un point géométrique de $X$, $X'$ le localisé strict de $X$ en $\ox$.
Supposons $Y$ normal et $f_\ox\colon Y_\ox\rightarrow \ox$ un homéomorphisme universel. 
Alors, $X'\times_XY$ est normal et strictement local.
\end{lem}

Considérant $X'$ comme une limite projective cofiltrante de voisinages étales affines $(X_i)_{\in I}$ de $\ox$ 
dans $X$ (\cite{sga4} VIII 4.5), $X'\times_XY$ est canoniquement isomorphe à la limite projective des schémas 
$(X_i\times_XY)_{i\in I}$ (\cite{ega4} 8.2.5). Pour tout $i\in I$, $X_i\times_XY$ est normal (\cite{raynaud1} VII prop.~2). 
Pour tout $(i,j)\in I^2$ tel que $i\leq j$, le morphisme $X_j\rightarrow X_i$ 
étant étale, chaque composante irréductible de $X_j\times_XY$ domine une composante irréductible 
de $X_i\times_XY$  (\cite{ega4} 2.3.5(ii)). On en déduit que $X'\times_XY$ est normal (\cite{ega1n} 0.6.5.12(ii)). 
Par ailleurs, comme $X'\times_XY$ est entier sur $X'$ et que $f_\ox$ est un homéomorphisme universel, 
$X'\times_XY$ est strictement local (\cite{raynaud1} I § 3 prop.~2).

\section{Schémas logarithmiques adéquats}\label{higgs3-slad}

\subsection{}\label{higgs3-slad0}
Les conventions et notations de géométrie logarithmique introduites dans (\cite{ag1} § 5)  sont en vigueur dans cet article. 
On désigne par $(S,\cM_S)$ le trait logarithmique fixé dans \eqref{higgs3-not1}. 

\begin{prop}\label{higgs3-slad4}
Soient $f\colon (X,\cM_X)\rightarrow (S,\cM_S)$ un morphisme lisse et saturé {\rm (\cite{ag1} 5.18)} 
de schémas logarithmiques fins, 
$X^\circ$ le sous-schéma ouvert maximal de $X$ où la structure logarithmique $\cM_X$ est triviale, 
$j\colon X^\circ\rightarrow X$ l'injection canonique. Alors, 
\begin{itemize}
\item[{\rm (i)}] Le schéma $X$ est $S$-plat et le schéma $X_\os$ est réduit. 
\item[{\rm (ii)}] Le schéma logarithmique $(X,\cM_X)$ est régulier {\rm (\cite{kato2} 2.1 et \cite{niziol} 2.2)}.
\item[{\rm (iii)}] Les schémas $X$ et $X\times_S\oS$ sont normaux et localement irréductibles \eqref{higgs3-sli1}. 
\item[{\rm (iv)}]  L'immersion $j\colon X^\circ\rightarrow X$ 
est schématiquement dominante, et on a des isomorphismes canoniques
\begin{eqnarray}
\cM^\gp_X&\stackrel{\sim}{\rightarrow}&j_*(\co_{X^\circ}^\times),\label{higgs3-slad4a}\\
\cM_X&\stackrel{\sim}{\rightarrow}&j_*(\co_{X^\circ}^\times)\cap \co_X.\label{higgs3-slad4b}
\end{eqnarray}
En particulier, l'homomorphisme canonique $\cM_X\rightarrow \co_X$ est un monomorphisme.
\end{itemize}
\end{prop}

(i) Cela résulte de (\cite{kato1} 4.5) et de (\cite{tsuji4} II 4.2).

(ii) Comme $(X,\cM_X)$ est saturé d'après (\cite{tsuji4} II 2.12),
il est régulier en vertu de (\cite{kato2} 8.2); cf. aussi (\cite{niziol} 2.3) et la preuve de (\cite{tsuji1} 1.5.1). 

(iii) Le schéma $X$ est normal en vertu de (ii) et (\cite{kato2} 4.1); cf. aussi (\cite{tsuji1} 1.5.1). 
On en déduit par passage à la limite inductive 
que le schéma $X\times_S\oS$ est normal (cf. la preuve de \cite{ag1} 6.3(iii)). 
Comme $X$ est localement noethérien, il est localement irréductible d'après \ref{higgs3-sli2}(ii). 
Par ailleurs, comme $X\times_S\oS$ est $\oS$-plat, ses points génériques 
sont les points génériques du schéma $X\times_S\oeta$, qui est localement noethérien. Par suite, 
l'ensemble des  points génériques de $X\times_S\oS$ est localement fini, et 
$X\times_S\oS$ est localement irréductible en vertu de \ref{higgs3-sli2}(ii). 

(iv) Cela résulte de (ii) et (\cite{kato2} 11.6); cf. aussi (\cite{niziol} 2.6).

\begin{lem}\label{higgs3-slad1}
Soient $f\colon (X,\cM_X)\rightarrow (S,\cM_S)$ un morphisme lisse et saturé de schémas logarithmiques fins, 
$(\mN,\iota)$ une carte pour $(S,\cM_S)$,
$\ox$ un point géométrique de $X$ au-dessus de $\os$. Alors, il existe un voisinage étale $U$
de $\ox$ dans $X$, une carte $(P,\gamma)$ pour $(U,\cM_X|U)$ et un homomorphisme $\vartheta\colon \mN\rightarrow P$ 
de monoïdes tels que les conditions suivantes soient remplies~:  
\begin{itemize}
\item[{\rm (i)}]  $((P,\gamma),(\mN,\iota),\vartheta)$ est une carte pour la restriction
$f_U\colon (U,\cM|U)\rightarrow (S,\cM_S)$ de $f$ à $U$ {\rm (\cite{ag1} 5.14)}, autrement dit,
le diagramme d'homomorphismes de monoïdes 
\begin{equation}\label{higgs3-slad1a}
\xymatrix{
P\ar[r]^-(0.5)\gamma&{\Gamma(U,\cM_X)}\\
\mN\ar[r]^-(0.5){\iota}\ar[u]^\vartheta&{\Gamma(S,\cM_S)}\ar[u]_{f^\flat_U}}
\end{equation}
est commutatif, ou ce qui revient au même le diagramme associé de morphismes de schémas logarithmiques 
\begin{equation}\label{higgs3-slad1b}
\xymatrix{
{(U,\cM_X|U)}\ar[r]^-(0.5){\gamma^a}\ar[d]_{f_U}&{\bA_P}\ar[d]^{\bA_\vartheta}\\
{(S,\cM_S)}\ar[r]^-(0.5){\iota^a}&{\bA_\mN}}
\end{equation}
est commutatif.
\item[{\rm (ii)}]  $P$ est torique, i.e., $P$ est fin et saturé et $P^\gp$ est libre sur $\mZ$ {\rm (\cite{ag1} 5.1)}. 
\item[{\rm (iii)}] L'homomorphisme $\vartheta$ est saturé {\rm (\cite{ag1} 5.2)}. 
\item[{\rm (iv)}] L'homomorphisme $\vartheta^\gp\colon \mZ\rightarrow P^\gp$ est injectif, 
le sous-groupe de torsion de $\coker(\vartheta^\gp)$ est d'ordre premier à $p$ et le morphisme de schémas usuels
\begin{equation}\label{higgs3-slad1c}
U\rightarrow S\times_{\bA_\mN}\bA_P
\end{equation}
déduit de \eqref{higgs3-slad1b} est étale.  
\item[{\rm (v)}] Il existe un sous-groupe $A$ de $P$ tel que $\gamma$ induise un isomorphisme
\begin{equation}\label{higgs3-slad1d}
P/A\stackrel{\sim}{\rightarrow} \cM_{X,\ox}/\co_{X,\ox}^\times.
\end{equation}
\end{itemize}
\end{lem}

Nous reprenons la preuve de (\cite{fkato} 4.1; cf. § 6) en l'adaptant. Posons 
\begin{equation}
\tOmega^1_{\tX/\tS}=\Omega^1_{(X,\cM_X)/(S,\cM_S)},
\end{equation}
et notons $\lambda\in \Gamma(X,\cM_X)$ l'image canonique de $\iota(1)=\pi\in \Gamma(S,\cM_S)=\co_K-\{0\}$. 
Soient $t_1,\dots,t_r\in \cM_{X,\ox}$ tels que $d\log(t_1),\dots,d\log(t_r)$ forment une base du $\co_{X,\ox}$-module 
$\tOmega^1_{\tX/\tS,\ox}$. Posons $H=\mN^{r+1}$ et considérons l'homomorphisme
\begin{equation}
\varphi\colon H\rightarrow \cM_{X,\ox}, \ \ \ (n_1,\dots,n_{r+1})\mapsto \prod_{i=1}^r t_i^{n_i} \cdot \lambda^{n_{r+1}}.
\end{equation}
Notons $\alpha\colon \cM_{X,\ox}^\gp\rightarrow \cM_{X,\ox}^\gp/\co_{X,\ox}^\times$ la projection 
canonique et $L$ l'image de l'homomorphisme
\begin{equation}
\alpha\circ \varphi^\gp \colon H^\gp\rightarrow \cM_{X,\ox}^\gp/\co_{X,\ox}^\times.
\end{equation}
Comme $\cM_X$ est fin et saturé (\cite{tsuji4} II 2.12), 
$\cM_{X,\ox}^\gp/\co_{X,\ox}^\times$ est un $\mZ$-module libre de type fini. 
D'après l'étape 2 de (\cite{fkato} page 331), on voit que le conoyau de $\alpha\circ \varphi^\gp$ 
est annulé par un entier inversible dans $\co_{X,\ox}$. 
Il existe donc une $\mZ$-base $e_1,\dots,e_d$ de $\cM_{X,\ox}^\gp/\co_{X,\ox}^\times$ et des entiers $f_1,\dots,f_d$ tels que 
$e_1^{f_1},\dots,e_d^{f_d}$ forment une $\mZ$-base de $L$, que $f_i$ divise $f_{i+1}$ pour tout $1\leq i\leq d-1$ et que $f_d$ soit inversible dans $\co_{X,\ox}$. Soient $F_1,\dots,F_d\in H^\gp$ et $\te_1,\dots,\te_d\in \cM^\gp_{X,\ox}$ 
tels que $\alpha(\varphi^\gp(F_i))=e_i^{f_i}$  et $\alpha(\te_i)=e_i$ pour tout $1\leq i\leq d$. 
Il existe alors $u_i\in \co_{X,\ox}^\times$ tel que 
$\varphi^\gp(F_i)=u_i\te_i^{f_i}$. Comme $\co_{X,\ox}^\times$ est $f_i$-divisible, il existe $v_i\in \co_{X,\ox}^\times$
tel que $v_i^{f_i}=u_i$. Remplaçant $\te_i$ par $v_i\te_i$, on peut supposer que $\varphi^\gp(F_i)=\te_i^{f_i}$.
On désigne par $\beta\colon \cM_{X,\ox}^\gp/\co_{X,\ox}^\times\rightarrow \cM_{X,\ox}^\gp$ le scindage de $\alpha$
défini par $\beta(e_i)=\te_i$ pour tout $1\leq i\leq d$,  
par $\rho\colon H^\gp\rightarrow L$ l'homomorphisme surjectif induit par $\alpha\circ \varphi^\gp$, 
par $\sigma\colon L\rightarrow H^\gp$ le scindage de $\rho$ défini par 
$\sigma(e_i^{f_i})=F_i$ pour tout $1\leq i\leq d$, par $M$ le noyau de $\rho$, 
et par $\tau\colon H^\gp\rightarrow M$ l'homomorphisme qui à tout $h\in H^\gp$ associe $h-\sigma(\rho(h))$. 
On pose $G=M\oplus \cM_{X,\ox}^\gp/\co_{X,\ox}^\times$ et
\begin{equation}
\phi\colon G=M\oplus \cM_{X,\ox}^\gp/\co_{X,\ox}^\times \rightarrow \cM^\gp_{X,\ox}, 
\ \ \ (m,t)\mapsto \phi(m,t)=\varphi^\gp(m)\cdot \beta(t). 
\end{equation} 
On a $\phi\circ (\tau\oplus \alpha\circ \varphi^\gp)=\varphi^\gp \colon H^\gp\rightarrow \cM_{X,\ox}^\gp$. On le vérifie immédiatement
pour les éléments de $M$ et pour les éléments $(F_i)_{1\leq i\leq d}$. On pose $P=\phi^{-1}(\cM_{X,\ox})$. D'après (\cite{kato1}
2.10), il existe un voisinage étale $U$ de $\ox$ dans $X$ et un homomorphisme $\gamma\colon P\rightarrow \Gamma(U,\cM_X)$
qui est une carte pour $(U,\cM_X|U)$ et dont la fibre $\gamma_\ox\colon P\rightarrow \cM_{X,\ox}$ en $\ox$ est induite par $\phi$. 
L'homomorphisme $\tau\oplus \alpha\circ \varphi^\gp\colon H^\gp\rightarrow G$ 
induit un homomorphisme $H\rightarrow P$ et par suite un homomorphisme
$\vartheta\colon \mN\rightarrow P$ qui rend commutatif le diagramme \eqref{higgs3-slad1a}. 

Comme l'homomorphisme $G\rightarrow \cM^\gp_{X,\ox}/\co_{X,\ox}^\times$ induit par $\phi$ est surjectif, 
$P^\gp=G$ et $P$ est intègre.  
Il résulte aussitôt de la définition qu'il existe un sous-groupe $A$ de $P$ tel que $\gamma_\ox$ induise un isomorphisme  
\begin{equation}\label{higgs3-slad1e}
P/A\stackrel{\sim}{\rightarrow} \cM_{X,\ox}/\co_{X,\ox}^\times.
\end{equation}
Par suite, $A=P^\times$. 
Comme $\cM_{X,\ox}$ est saturé, $\cM_{X,\ox}/\co_{X,\ox}^\times$ est saturé et donc $P$ est saturé (\cite{ag1} 5.1). 
On en déduit que $P$ est torique. 
L'homomorphisme $\vartheta$ est saturé en vertu de \eqref{higgs3-slad1e} et (\cite{tsuji4} I 3.16). 

L'homomorphisme $\tau\oplus \alpha\circ \varphi^\gp\colon H^\gp\rightarrow G$ est injectif et 
son conoyau est isomorphe à celui de $\alpha\circ \varphi^\gp$. 
On en déduit que l'homomorphisme $\vartheta^\gp\colon \mZ\rightarrow P^\gp$ est injectif et que le sous-groupe de torsion de 
$\coker(\vartheta^\gp)$ est d'ordre premier à $p$. L'étape 4 de (\cite{fkato} § 6 page 332) montre que 
le morphisme de schémas usuels $U\rightarrow S\times_{\bA_\mN}\bA_P$ déduit de \eqref{higgs3-slad1b} est étale.

\begin{defi}\label{higgs3-slad3}\index{Carte logarithmique adequate@Carte (logarithmique) adéquate}
Soient $f\colon (X,\cM_X)\rightarrow (S,\cM_S)$ un morphisme de schémas logarithmiques,  
$(P,\gamma)$ une carte pour $(X,\cM_X)$, 
$(\mN,\iota)$ une carte pour $(S,\cM_S)$, $\vartheta\colon \mN\rightarrow P$ un homomorphisme 
de monoïdes tels que $((P,\gamma),(\mN,\iota),\vartheta)$ soit une carte pour $f$ (\cite{ag1} 5.14), autrement dit,  
que le diagramme d'homomorphismes de monoïdes 
\begin{equation}\label{higgs3-slad3a}
\xymatrix{
P\ar[r]^-(0.5)\gamma&{\Gamma(X,\cM_X)}\\
\mN\ar[r]^-(0.5){\iota}\ar[u]^\vartheta&{\Gamma(S,\cM_S)}\ar[u]_{f^\flat}}
\end{equation}
soit commutatif, ou ce qui revient au même que le diagramme associé de morphismes de schémas logarithmiques 
\begin{equation}\label{higgs3-slad3b}
\xymatrix{
{(X,\cM_X)}\ar[r]^-(0.5){\gamma^a}\ar[d]_f&{\bA_P}\ar[d]^{\bA_\vartheta}\\
{(S,\cM_S)}\ar[r]^-(0.5){\iota^a}&{\bA_\mN}}
\end{equation}
soit commutatif. On dit que la carte $((P,\gamma),(\mN,\iota),\vartheta)$ 
est {\em adéquate} si les conditions suivantes sont remplies~:
\begin{itemize}
\item[(i)] Le monoïde $P$ est torique, {\em i.e.}, $P$ est fin et saturé et $P^\gp$ est libre sur $\mZ$.
\item[(ii)] L'homomorphisme $\vartheta$ est saturé. 
\item[(iii)] L'homomorphisme $\vartheta^\gp\colon \mZ\rightarrow P^\gp$ est injectif, 
le sous-groupe de torsion de $\coker(\vartheta^\gp)$ est d'ordre premier à $p$ et le morphisme de schémas usuels
\begin{equation}\label{higgs3-slad3c}
X\rightarrow S\times_{\bA_\mN}\bA_P
\end{equation}
déduit de \eqref{higgs3-slad3b} est étale.  
\item[(iv)] Posons $\lambda=\vartheta(1)\in P$, 
\begin{eqnarray}
L&=&\Hom_{\mZ}(P^\gp,\mZ),\label{higgs3-slad3d}\\
\rH(P)&=&\Hom(P,\mN).\label{higgs3-slad3e}
\end{eqnarray} 
On notera que $\rH(P)$ est un monoïde fin, saturé et affûté et que l'homomorphisme canonique 
$\rH(P)^\gp\rightarrow \Hom((P^\sharp)^\gp,\mZ)$ est un isomorphisme (\cite{ogus} I 2.2.3), 
où $P^\sharp$ désigne le quotient $P/P^\times$ (\cite{ag1} 5.1).
On suppose qu'il existe $h_1,\dots,h_r\in \rH(P)$, qui sont $\mZ$-linéairement indépendants dans $L$, tels que  
\begin{equation}\label{higgs3-slad3f}
\ker(\lambda)\cap \rH(P)=\{\sum_{i=1}^ra_ih_i | \ (a_1,\dots,a_r)\in \mN^r\},
\end{equation}
où l'on considère $\lambda$ comme un homomorphisme $L\rightarrow \mZ$.
\end{itemize}
\end{defi}

On notera que tout morphisme de schémas logarithmiques qui admet une carte adéquate est lisse et saturé (\cite{ag1} 5.25
et \cite{tsuji4} chap.~II 3.5). 

\begin{lem}\label{higgs3-slad2}
Soient $X$ un schéma affine, $U$ un ouvert schématiquement dense de $X$, 
$j\colon U\rightarrow X$ l'injection canonique,
$\cM$ le sous-monoïde multiplicatif $j_*(\co_U^\times)\cap \co_X$ de $j_*(\co_U)$,
$\lambda\in \Gamma(X,\cM)$. Notons $\Gamma(X,\cM)_\lambda$ la localisation du monoïde $\Gamma(X,\cM)$ en $\lambda$ 
{\em (\cite{ogus} I 1.4.4)} et $X_\lambda$ l'ouvert
de $X$ où l'image canonique de $\lambda$ dans $\Gamma(X,\co_X)$ est inversible. Alors, l'homomorphisme canonique 
\begin{equation}
\Gamma(X,\cM)_\lambda\rightarrow \Gamma(X_\lambda,\cM)
\end{equation}
est un isomorphisme. 
\end{lem}

En effet, posons $A=\Gamma(X,\co_X)$ et identifions $\lambda$ à un élément de $A$
et $\Gamma(X,\cM)$ (resp. $\Gamma(X_\lambda,\cM)$) au sous-monoïde multiplicatif 
de $A$ (resp. $A_\lambda$) formé des éléments $f$ tels que $f|U$ soit une unité.  
On vérifie aussitôt que pour tout monoïde $P$, tout homomorphisme $u\colon \Gamma(X,\cM)\rightarrow P$
tel que $u(\lambda)$ soit inversible, se factorise uniquement à travers $\Gamma(X_\lambda,\cM)$; d'où la proposition.

\begin{prop}\label{higgs3-slad5}
Soient $f\colon (X,\cM_X)\rightarrow (S,\cM_S)$ un morphisme lisse et saturé
de schémas logarithmiques fins, $\ox$ un point géométrique de $X$ au-dessus de $s$, 
$X'$ le localisé strict de $X$ en $\ox$, $\cM_{X'}$ la structure logarithmique image inverse de $\cM_X$ sur $X'$ 
{\rm (\cite{ag1} 5.10)}. Alors, les conditions suivantes sont équivalentes~:
\begin{itemize}
\item[{\rm (a)}] Il existe un voisinage étale $U$ de $\ox$ dans $X$ tel que la restriction
$f_U\colon (U,\cM_X|U)\rightarrow (S,\cM_S)$ de $f$ à $U$ admette une carte adéquate \eqref{higgs3-slad3}.
\item[{\rm (b)}] Il existe un voisinage étale $U$ de $\ox$ dans $X$ tel que le schéma $U_\eta$ 
soit lisse sur $\eta$ et que la structure logarithmique $\cM_X|U_\eta$ sur $U_\eta$
soit définie par un diviseur à croisements normaux stricts sur $U_\eta$.
\item[{\rm (c)}] Le schéma $X'_\eta$ est régulier et la structure logarithmique $\cM_{X'}|X'_\eta$ est 
définie par un diviseur à croisements normaux stricts sur $X'_\eta$.
\end{itemize}
\end{prop}

L'implication (a)$\Rightarrow$(b) résulte de (\cite{ag1} 6.3(v)); on notera que
les conditions (C$_1$) et (C$_2$) de (\cite{ag1} 6.2) ne jouent aucun rôle dans la preuve de (\cite{ag1} 6.3(v)). 
On désigne par $\nu\colon X'\rightarrow X$ le morphisme canonique et par $X^\circ$ le sous-schéma ouvert maximal de $X$ 
où la structure logarithmique $\cM_X$ est triviale. On pose $X'^\circ=X'\times_XX^\circ$
et on note $j\colon X^\circ\rightarrow X$ et $j'\colon X'^\circ\rightarrow X'$ les injections canoniques.
L'homomorphisme canonique $\nu^{-1}(\co_X)\rightarrow \co_{X'}$ est un isomorphisme de $X'_\et$ \eqref{higgs3-sli7}. 
Par suite, l'homomorphisme canonique 
\begin{equation}\label{higgs3-slad5d}
\nu^{-1}(\cM_X)\rightarrow \cM_{X'}
\end{equation} 
est aussi un isomorphisme. On a $\cM_X=j_*(\co_{X^\circ}^\times)\cap \co_X$ en vertu de \ref{higgs3-slad4}(iv).
Comme $\nu$ est universellement localement acyclique, on en déduit que
\begin{equation}\label{higgs3-slad5e}
\cM_{X'}=j'_*(\co_{X'^\circ}^\times)\cap \co_{X'}.
\end{equation} 
L'implication (b)$\Rightarrow$(c) s'ensuit. 
Montrons l'implication (c)$\Rightarrow$(a). D'après \ref{higgs3-slad1}, 
quitte à remplacer $X$ par un voisinage étale de $\ox$ dans $X$, on peut supposer que 
$f$ admet une carte $((P,\gamma),(\mN,\iota),\vartheta)$ vérifiant les conditions (i), (ii) et (iii) de \ref{higgs3-slad3} et 
qu'il existe un sous-groupe $A$ de $P$ tel que $\gamma$ induise un isomorphisme
\begin{equation}\label{higgs3-slad5a}
P/A\stackrel{\sim}{\rightarrow} \cM_{X,\ox}/\co_{X,\ox}^\times.
\end{equation}
Posons $\lambda=\vartheta(1)\in P$, 
\begin{eqnarray}
L&=&\Hom_{\mZ}(P^\gp,\mZ),\label{higgs3-slad5b}\\
\rH(P)&=&\Hom(P,\mN).\label{higgs3-slad5c}
\end{eqnarray} 
On notera que $\rH(P)$ est un monoïde fin, saturé et affûté et que l'homomorphisme canonique 
$\rH(P)^\gp\rightarrow \Hom((P^\sharp)^\gp,\mZ)$ est un isomorphisme (\cite{ogus} I 2.2.3) 
où $P^\sharp$ désigne le quotient $P/P^\times$. 
Soit $F$ la face de $P$ engendrée par $\lambda$, c'est-à-dire l'ensemble des éléments 
$\alpha\in P$ tels qu'il existe $\beta\in P$ et $n\in \mN$ tels que $\alpha+\beta=n\lambda$ (\cite{ogus} I 1.4.2). 
Notons $P/F$ le quotient de $P$ par $F$ (cf. \cite{ogus} I 1.1.6). L'homomorphisme canonique 
\begin{equation}
\Hom(P/F,\mN)\rightarrow \ker(\lambda)\cap \rH(P),
\end{equation}
où l'on considère $\lambda$ comme un homomorphisme $L\rightarrow \mZ$, est un isomorphisme. 
Il suffit donc de montrer que le monoïde $P/F$ est libre de type fini.  
Soit $G$ la face de $\cM_{X,\ox}$ engendrée par $\gamma(\lambda)$. Compte tenu de \eqref{higgs3-slad5a}, 
il suffit encore de montrer que le monoïde $\cM_{X,\ox}/G$ est libre de type fini. On a un isomorphisme canonique
\begin{equation}
\cM_{X,\ox}/G\stackrel{\sim}{\rightarrow}(G^{-1}\cM_{X,\ox})^\sharp.
\end{equation}
Par ailleurs, on a $\cM_{X,\ox}=\Gamma(X',\cM_{X'})$ \eqref{higgs3-slad5d}, et l'homomorphisme canonique
\begin{equation}
G^{-1}\Gamma(X',\cM_{X'})\rightarrow \Gamma(X'_\eta,\cM_{X'})
\end{equation}
est un isomorphisme d'après \ref{higgs3-slad2} et \eqref{higgs3-slad5e}. La proposition recherchée résulte alors du fait que 
le monoïde $\Gamma(X'_\eta,\cM_{X'}^\sharp)$ est libre de type fini compte tenu de (c).

\begin{defi}\label{higgs3-slad6}\index{Morphisme adequat de schemas logarithmiques@Morphisme adéquat de schémas logarithmiques}
On dit qu'un morphisme $f\colon (X,\cM_X)\rightarrow (S,\cM_S)$ de schémas logarithmiques fins  
est {\em adéquat} s'il est lisse et saturé, si le morphisme de schémas sous-jacents $X\rightarrow S$ est de type fini
et si pour tout point géométrique $\ox$ de $X$ au-dessus de $s$, 
les conditions équivalentes de \ref{higgs3-slad5} sont remplies.
\end{defi} 

La notion de morphisme adéquat de schémas logarithmiques correspond à la notion 
de morphisme {\em à singularités toroïdales} de Faltings. Dans sa terminologie, 
le morphisme $f$ est qualifié de {\em petit} s'il admet une carte adéquate, si $X$ est affine 
et connexe et si $X_s$ est non-vide, autrement dit, si les conditions de (\cite{ag1} 6.2) sont remplies. 

\begin{prop}
Soient $f\colon (X,\cM_X)\rightarrow (S,\cM_S)$ un morphisme adéquat de schémas logarithmiques fins,   
$K'$ une extension finie de $K$, $\co_{K'}$ la fermeture intégrale de $\co_K$ dans $K'$, 
$S'=\Spec(\co_{K'})$. On munit $S'$ de la structure logarithmique $\cM_{S'}$ 
définie par son point fermé. On a un morphisme canonique $(S',\cM_{S'})\rightarrow (S,\cM_S)$. On pose 
\begin{equation}\label{higgs3-TFSLA12a}
(X',\cM_{X'})=(X,\cM_X)\times_{(S,\cM_S)}(S',\cM_{S'}),
\end{equation}
le produit étant pris dans la catégorie des schémas logarithmiques.
Alors, la projection canonique $f'\colon (X',\cM_{X'})\rightarrow (S',\cM_{S'})$ est adéquate.
\end{prop}

En effet, $f'$ est lisse et saturé (\cite{tsuji4} II 2.11). Comme $X'=X\times_SS'$, il est de type fini sur $S'$.
Par ailleurs, la condition \ref{higgs3-slad5}(b) est clairement satisfaite en tout point géométrique de $X'$ au-dessus 
du point fermé de $S'$.

\section{Variation sur le complexe de Koszul}\label{higgs3-vck}

\subsection{}\label{higgs3-vck1}
Dans cette section, $(X,A)$ désigne un topos annelé.
Pour tout morphisme de $A$-modules $u\colon E\rightarrow F$, 
il existe sur l'algèbre bigraduée $\rS(E)\otimes_A\Lambda(F)$ \eqref{higgs3-not9}
(anti-commutative pour le second degré~; \cite{illusie1} I 4.3.1.1) une unique $A$-dérivation 
\begin{equation}\label{higgs3-vck1a}
d_u\colon \rS(E)\otimes_A\Lambda(F)\rightarrow \rS(E)\otimes_A\Lambda(F)
\end{equation} 
de bidegré $(-1,1)$ telle que pour toutes sections locales $x_1,\dots,x_n$ de $E$ $(n\geq 1)$ et $y$ de $\wedge(F)$, on ait 
\begin{eqnarray}
d_u([x_1\otimes\dots \otimes x_n] \otimes y)&=&\sum_{i=1}^n [x_1\otimes\dots \otimes x_{i-1} \otimes x_{i+1} \otimes
\dots \otimes x_n]\otimes (u(x_i)\wedge y),\label{higgs3-vck1b}\\
d_u(1\otimes y)&=& 0.\label{higgs3-vck1c}
\end{eqnarray}
Elle satisfait de plus $d_u\circ d_u=0$. Pour tout entier $n\geq 0$, la partie homogène de degré $n$ de 
$\rS(E)\otimes_A\Lambda(F)$ donne un complexe 
\begin{equation}\label{higgs3-vck1d}
0\rightarrow \rS^n(E)\rightarrow \rS^{n-1}(E)\otimes_AF\rightarrow \dots \rightarrow E\otimes_A \wedge^{n-1}(F)
\rightarrow \wedge^n(F)\rightarrow 0.
\end{equation}
L'algèbre $\rS(E)\otimes_A\Lambda(F)$ munie de la dérivation $d_u$ dépend fonctoriellement de $u$. 

Si $A$ est une $\mQ$-algèbre, identifiant $\rS(E)$ à l'algèbre à puissances divisées $\Gamma(E)$ 
de $E$, $d_u$ s'identifie à la $A$-dérivation de $\Gamma(E)\otimes_A\Lambda(F)$ définie dans 
(\cite{illusie1} I 4.3.1.2(b)). Il résulte alors de (\cite{illusie1} I 4.3.1.6) que si $u$ est un isomorphisme de $A$-modules plats et si $n>0$, 
la suite \eqref{higgs3-vck1d} est exacte. 

\subsection{}\label{higgs3-vck2}
Soit 
\begin{equation}\label{higgs3-vck2a}
0\rightarrow A \rightarrow E\rightarrow F\rightarrow 0
\end{equation}
une suite exacte de $A$-modules plats. 
D'après (\cite{illusie1} I 4.3.1.7), celle-ci induit pour tout entier $n\geq 1$, une suite exacte localement scindée \eqref{higgs3-not9}
\begin{equation}\label{higgs3-vck2b}
0\rightarrow \rS^{n-1}(E)\rightarrow \rS^{n}(E)\rightarrow \rS^n(F)\rightarrow 0.
\end{equation}
On en déduit une suite exacte 
\begin{equation}\label{higgs3-vck2c}
0\rightarrow \rS^{n-1}(F)\rightarrow \rS^n(E)/\rS^{n-2}(E)\rightarrow \rS^n(F)\rightarrow 0.
\end{equation}
On définit ainsi un foncteur $\rS^n$ 
de la catégorie $\bExt(F,A)$ des extensions de $F$ par $A$ dans la catégorie $\bExt(\rS^n(F),\rS^{n-1}(F))$
des extensions de $\rS^n(F)$ par $\rS^{n-1}(F)$. Passant aux groupes des classes d'isomorphismes 
des objets de ces catégories, on obtient un homomorphisme, que l'on note encore
\begin{equation}\label{higgs3-vck2d}
\rS^n\colon \Ext^1_A(F,A)\rightarrow \Ext^1_A(\rS^n(F),\rS^{n-1}(F)).
\end{equation} 

\subsection{}\label{higgs3-vck3}
Soient $F$ un $A$-module localement projectif de type fini \eqref{higgs3-not3}, $n$ un entier $\geq 1$. La suite spectrale
qui relie les Ext locaux et globaux (\cite{sga4} V 6.1) fournit un isomorphisme 
\begin{equation}\label{higgs3-vck3a}
\Ext^1_A(F,A)\stackrel{\sim}{\rightarrow}\rH^1(X,\cHom_A(F,A)).
\end{equation}
De même, comme le $A$-module $\rS^n(F)$ est localement libre de type fini, on a un isomorphisme  canonique
\begin{equation}\label{higgs3-vck3b}
\Ext^1_A(\rS^n(F),\rS^{n-1}(F))\stackrel{\sim}{\rightarrow}\rH^1(X,\cHom_A(\rS^n(F),\rS^{n-1}(F))).
\end{equation}
Par ailleurs, on désigne par  
\begin{equation}\label{higgs3-vck3c}
J_n\colon \cHom_A(F,A)\rightarrow \cHom_A(\rS^n(F),\rS^{n-1}(F))
\end{equation}
le morphisme qui pour tout $U\in \ob(X)$, associe à tout morphisme $u\colon F|U\rightarrow A|U$ 
la restriction à $\rS^n(F)|U$ de la dérivation $d_u$ de $\rS(F|U)$ définie dans \eqref{higgs3-vck1a}. 
Celui-ci induit un accouplement
\begin{equation}\label{higgs3-vck3d}
\cHom_A(F,A)\otimes_A \rS^n(F) \rightarrow \rS^{n-1}(F).
\end{equation}

\begin{prop}\label{higgs3-vck4}
Pour tout $A$-module localement projectif de type fini $F$ et tout entier $n\geq 1$, le diagramme 
\begin{equation}\label{higgs3-vck4a}
\xymatrix{
{\Ext^1_A(F,A)}\ar[d]\ar[r]^-(0.5){\rS^n}&{\Ext^1_A(\rS^n(F),\rS^{n-1}(F))}\ar[d]\\
\rH^1(X,\cHom_A(F,A))\ar[r]^-(0.5){J_n}&{\rH^1(X,\cHom_A(\rS^n(F),\rS^{n-1}(F)))}}
\end{equation}
où $\rS^n$ est le morphisme \eqref{higgs3-vck2d}, $J_n$ est le morphisme \eqref{higgs3-vck3c}
et les flèches verticales sont les isomorphismes \eqref{higgs3-vck3a} et \eqref{higgs3-vck3b}, est commutatif. 
\end{prop}

On rappelle d'abord (\cite{sga4} V 3.4) que la suite spectrale de Cartan-Leray relative 
aux recouvrements de l'objet final de $X$ induit un isomorphisme 
\begin{equation}\label{higgs3-vck4b}
\crH^1(X,\cHom_A(F,A))\stackrel{\sim}{\rightarrow}\rH^1(X,\cHom_A(F,A)), 
\end{equation}
où la source désigne le groupe de cohomologie de Cech (\cite{sga4} V (2.4.5.4)).
On peut décrire explicitement l'isomorphisme 
\begin{equation}\label{higgs3-vck4c}
\Ext^1_A(F,A)\stackrel{\sim}{\rightarrow}\crH^1(X,\cHom_A(F,A))
\end{equation}
composé de \eqref{higgs3-vck3a} et l'inverse de \eqref{higgs3-vck4b} comme suit. Soit 
\begin{equation}\label{higgs3-vck4d}
0\rightarrow A \rightarrow E\stackrel{\nu}{\rightarrow} F\rightarrow 0
\end{equation}
une suite exacte de $A$-modules. Comme $F$ est localement projectif de type fini, 
il existe une famille $\cU=(U_i)_{i\in I}$ d'objets de $X$, épimorphique au-dessus de l'objet final, telle que
pour tout $i\in I$, il existe une section $(A|U_i)$-linéaire $\varphi_i\colon F|U_i\rightarrow E|U_i$ de $\nu|U_i$.  
Pour tout $(i,j)\in I^2$, posant $U_{i,j}=U_i\times U_j$, la différence $\varphi_{i,j}=\varphi_i|U_{i,j}-\varphi_j|U_{i,j}$
définit un morphisme de $F|U_{i,j}$ dans $A|U_{i,j}$. La collection $(\varphi_{i,j})$ est un $1$-cocycle pour 
le recouvrement $\cU$ à coefficients dans $\cHom_A(F,A)$ dont la classe dans $\crH^1(X,\cHom_A(F,A))$ 
est l'image canonique de l'extension \eqref{higgs3-vck4d} ({\em i.e.}, son image par l'isomorphisme \eqref{higgs3-vck4c}). Pour tout 
$i\in I$, on note $\psi^n_i$ 
le morphisme composé 
\begin{equation}\label{higgs3-vck4e}
\xymatrix{
{\rS^n(F)|U_i}\ar[r]^{\rS^n(\varphi_i)}&{\rS^n(E)|U_i}\ar[r]&{(\rS^n(E)/\rS^{n-2}(E))|U_i}},
\end{equation}
où la seconde flèche est la projection canonique. C'est clairement un scindage au-dessus de $U_i$ 
de la suite exacte \eqref{higgs3-vck2c}
\begin{equation}\label{higgs3-vck4f}
0\rightarrow \rS^{n-1}(F)\rightarrow \rS^n(E)/\rS^{n-2}(E)\rightarrow \rS^n(F)\rightarrow 0
\end{equation}
déduite de \eqref{higgs3-vck4d}. Pour tout $(i,j)\in I^2$, la différence 
$\psi^n_{i,j}=\psi^n_i|U_{i,j}-\psi^n_j|U_{i,j}$ définit un morphisme de $\rS^n(F)|U_{i,j}$ dans 
$\rS^{n-1}(F)|U_{i,j}$. La collection $(\psi^n_{i,j})$ est un $1$-cocycle pour 
le recouvrement $\cU$ à coefficients dans $\cHom_A(\rS^n(F),\rS^{n-1}(F))$ dont la classe dans 
$\crH^1(X,\cHom_A(\rS^n(F),\rS^{n-1}(F)))$ est l'image canonique de l'extension \eqref{higgs3-vck4f}. 

Pour tout $(i,j)\in I^2$ et toutes sections locales $x_1,\dots,x_n$ de $F|U_{i,j}$, 
\begin{eqnarray}
\psi^n_{i,j}([x_1\otimes\dots\otimes x_n])&=&\psi_i^n([x_1\otimes\dots\otimes x_n])-
\psi_j^n([x_1\otimes\dots\otimes x_n])\label{higgs3-vck4g}\\
&=&[(\varphi_j(x_1)+\varphi_{i,j}(x_1))\otimes\dots\otimes (\varphi_j(x_n)+\varphi_{i,j}(x_n))]\nonumber\\
&&-[\varphi_j(x_1)\otimes\dots\otimes \varphi_j(x_n)] \mod (\rS^{n-2}(E))\nonumber\\
&=&\sum_{\alpha=1}^n \varphi_{i,j}(x_\alpha) [x_1\otimes \dots \otimes x_{\alpha-1}\otimes x_{\alpha+1}\otimes\dots \otimes x_n]\nonumber\\
&=&J_n(\varphi_{i,j}) ([x_1\otimes\dots\otimes x_n]).\nonumber
\end{eqnarray}
La proposition s'ensuit. 

\subsection{}\label{higgs3-vck8}
Soient 
\begin{equation}\label{higgs3-vck8c}
0\rightarrow A \rightarrow E\rightarrow F\rightarrow 0
\end{equation}
une suite exacte de $A$-modules localement projectifs de type fini, $n,q$ deux entiers $\geq 0$.
D'après \ref{higgs3-vck2}, la suite exacte \eqref{higgs3-vck8c} induit une suite exacte 
\begin{equation}\label{higgs3-vck8d}
0\rightarrow \rS^n(F)\rightarrow \rS^{n+1}(E)/\rS^{n-1}(E)\rightarrow \rS^{n+1}(F)\rightarrow 0.
\end{equation}
Par ailleurs, l'accouplement \eqref{higgs3-vck3d} induit un accouplement 
\begin{equation}\label{higgs3-vck8b}
\rH^1(X,\cHom_A(F,A))\otimes_{A(X)} \rH^q(X,\rS^{n+1}(F))\rightarrow  \rH^{q+1}(X,\rS^n(F)).
\end{equation}
Il résulte aussitôt de \ref{higgs3-vck4} que  le morphisme 
\begin{equation}\label{higgs3-vck8e}
\rH^q(X,\rS^{n+1}(F))\rightarrow\rH^{q+1}(X,\rS^n(F))
\end{equation}
bord de la suite exacte longue de cohomologie déduite de la suite exacte courte \eqref{higgs3-vck8d}, 
est induit par le cup-produit avec la classe de l'extension \eqref{higgs3-vck8c} par l'accouplement \eqref{higgs3-vck8b}. 

Notons 
\begin{equation}\label{higgs3-vck8f}
\partial\colon \Gamma(X,F)\rightarrow \rH^1(X,A)
\end{equation}
le morphisme bord de la suite exacte longue de cohomologie déduite de la suite exacte courte \eqref{higgs3-vck8c}. 
Il résulte encore de \ref{higgs3-vck4} (plus précisément de \eqref{higgs3-vck4g}) qu'on a un diagramme commutatif 
\begin{equation}\label{higgs3-vck8h}
\xymatrix{
{\rS^{n+1}(\Gamma(X,F))}\ar[r]\ar[d]_{\alpha}&{\Gamma(X,\rS^{n+1}(F))}\ar[d]\\
{\rS^n(\Gamma(X,F))\otimes_{A(X)} \rH^1(X,A)}\ar[r]&{\rH^1(X,\rS^n(F))}}
\end{equation}
où $\alpha$ est la restriction à $\rS^{n+1}(\Gamma(X,F))$
de la $A(X)$-dérivation $d_\partial$ de $\rS(\Gamma(X,F))\otimes_{A(X)}\wedge(\rH^1(X,A))$ définie dans \eqref{higgs3-vck1a} 
relativement au morphisme $\partial$, la flèche horizontale supérieure (resp. inférieure) est le morphisme canonique
(resp. est induite par le cup-produit) et la flèche verticale
de droite est le morphisme \eqref{higgs3-vck8e} pour $q=0$. 
Par associativité du cup-produit, on en déduit que le diagramme 
\begin{equation}\label{higgs3-vck8g}
\xymatrix{
{\rS^{n+1}(\Gamma(X,F))\otimes_{A(X)} \rH^q(X,A)}\ar[r]\ar[d]_{\alpha\otimes \id}&{\rH^q(X,\rS^{n+1}(F))}\ar[dd]\\
{\rS^n(\Gamma(X,F))\otimes_{A(X)} \rH^1(X,A)\otimes_{A(X)} \rH^q(X,A)}\ar[d]_{\id\otimes \cup}&\\
{\rS^n(\Gamma(X,F))\otimes_{A(X)} \rH^{q+1}(X,A)}\ar[r]&{\rH^{q+1}(X,\rS^n(F))}}
\end{equation}
où $\cup$ est le cup-produit de la $A(X)$-algèbre $\oplus_{i\geq 0}\rH^i(X,A)$, 
les morphismes horizontaux sont induits par le cup-produit et la flèche verticale
de droite est le morphisme \eqref{higgs3-vck8e}, est commutatif.

\subsection{}\label{higgs3-vck6}
Soient $f\colon (X,A)\rightarrow (Y,B)$ un morphisme de topos annelés, 
\begin{equation}\label{higgs3-vck6a}
0\rightarrow A\rightarrow E\rightarrow F\rightarrow 0
\end{equation}
une suite exacte de $A$-modules localement projectifs de type fini.  On désigne par 
\begin{equation}\label{higgs3-vck6b}
u\colon f_*(F)\rightarrow \rR^1f_*(A)
\end{equation}
le morphisme bord de la suite exacte
longue de cohomologie déduite de la suite exacte courte \eqref{higgs3-vck6a}.
D'après \ref{higgs3-vck2}, la suite exacte \eqref{higgs3-vck6a} induit pour tout entier $n\geq 0$,  une suite exacte 
\begin{equation}\label{higgs3-vck6e}
0\rightarrow \rS^n(F)\rightarrow \rS^{n+1}(E)/\rS^{n-1}(E)\rightarrow \rS^{n+1}(F)\rightarrow 0.
\end{equation}

\begin{prop}\label{higgs3-vck7}
Sous les hypothèses de \eqref{higgs3-vck6}, pour tous entiers $n,q\geq 0$, on a un diagramme commutatif
\begin{equation}\label{higgs3-vck7a}
\xymatrix{
{\rS^{n+1}(f_*(F))\otimes_B\rR^qf_*(A)} \ar[d]_{\alpha\otimes \id}\ar[r]&{\rR^qf_*(\rS^{n+1}(F))}\ar[dd]^{\partial}\\
{\rS^n(f_*(F))\otimes_B\rR^1f_*(A)\otimes_B\rR^qf_*(A)}\ar[d]_{\id\otimes \cup}&\\
{\rS^n(f_*(F))\otimes_B\rR^{q+1}f_*(A)}\ar[r]&{\rR^{q+1}f_*(\rS^n(F))}}
\end{equation}
où $\partial$ est le morphisme bord de la suite exacte longue de cohomologie déduite de la suite exacte courte \eqref{higgs3-vck6a},
$\alpha$ est la restriction à $\rS^{n+1}(f_*(F))$ de la $B$-dérivation $d_u$ de 
$\rS(f_*(F))\otimes_B\wedge(\rR^1f_*(A))$ définie dans \eqref{higgs3-vck1a} 
relativement au morphisme $u$ \eqref{higgs3-vck6b}, $\cup$ est le cup-produit de la $B$-algèbre $\oplus_{i\geq 0}\rR^if_*(A)$
et les morphismes horizontaux sont induits par le cup-produit.
\end{prop}

En effet, $\rR^qf_*(\rS^{n+1}(F))$ est le faisceau sur $X$ (pour la topologie canonique) associé au préfaisceau
qui à tout $V\in \ob(Y)$ associe $\rH^q(f^*(V),\rS^{n+1}(F))$, et $d$ est induit par le morphisme 
\begin{equation}
\rH^q(f^*(V),\rS^{n+1}(F))\rightarrow \rH^{q+1}(f^*(V),\rS^{n}(F))
\end{equation}
bord de la suite exacte longue de cohomologie déduite de la suite exacte courte \eqref{higgs3-vck6e}.
La proposition résulte alors de \eqref{higgs3-vck8g}. 

\section{Catégories additives à isogénie près}\label{higgs3-caip}

\begin{defi}\label{higgs3-caip1}\index{Isogenie@Isogénie}
Soit $\bC$ une catégorie additive. 
\begin{itemize}
\item[(i)] Un morphisme $u\colon M\rightarrow N$ de $\bC$ est appelé {\em isogénie} s'il existe
un entier $n\not=0$ et un morphisme $v\colon N\rightarrow M$ de $\bC$ tels que $v\circ u=n\cdot \id_M$ 
et $u\circ v= n\cdot \id_N$. 
\item[(ii)] Un objet $M$ de $\bC$ est dit {\em d'exposant fini} s'il existe un entier $n\not=0$ tel que $n\cdot \id_M=0$. 
\end{itemize}
\end{defi}
On peut compléter la terminologie et faire les remarques suivantes~:

\subsubsection{}\label{higgs3-caip1a}\index{Categorie additive a isogenie pres@Catégorie additive à isogénie près}
\index{1000602@$\bC_\mQ$, $M_\mQ$}
\addtocounter{equation}{1}
La famille des isogénies de $\bC$ permet un calcul de fractions bilatéral (\cite{illusie1} I 1.4.2).
On appelle {\em catégorie des objets de $\bC$ à isogénie près}, et l'on note $\bC_\mQ$,  
la catégorie localisée de $\bC$ par rapport aux isogénies. On désigne par
\begin{equation}\label{higgs3-caip1aa}
F\colon \bC\rightarrow \bC_\mQ, \ \ \ M\mapsto M_\mQ
\end{equation}
le foncteur de localisation. On vérifie aisément que pour tous $M,N\in \ob(\bC)$, on a  
\begin{equation}\label{higgs3-caip1ab}
\Hom_{\bC_\mQ}(M_\mQ,N_\mQ)=\Hom_{\bC}(M,N)\otimes_{\mZ}\mQ.
\end{equation}
En particulier, la catégorie $\bC_\mQ$ est additive et le foncteur de localisation est additif. 
Pour qu'un objet $M$ de $\bC$ soit d'exposant fini, il faut et il suffit que $M_\mQ$ soit nul. 

\addtocounter{subsubsection}{2}

\subsubsection{}\label{higgs3-caip1b}
\addtocounter{equation}{1}
Si $\bC$ est une catégorie abélienne, la catégorie $\bC_\mQ$ est abélienne et le foncteur de localisation 
$F\colon \bC\rightarrow \bC_\mQ$ est exact. En fait, $\bC_\mQ$ s'identifie canoniquement
à la catégorie quotient de $\bC$ par la sous-catégorie épaisse $\bE$ des objets d'exposant fini.
En effet, notons $\bC/\bE$ la catégorie quotient de $\bC$ par $\bE$ et $T\colon \bC\rightarrow \bC/\bE$ 
le foncteur canonique  (\cite{gabriel} III §1). Pour tout $M\in \ob(\bE)$, on a $F(M)=0$. 
Par suite, il existe un et un unique foncteur $F'\colon \bC/\bE\rightarrow \bC_\mQ$ tel que $F=F'\circ T$. 
Par ailleurs, pour tout $M\in \ob(\bC)$ et tout entier $n\not=0$, $T(n\cdot \id_M)$ est un isomorphisme. 
Il existe donc un et un unique foncteur $T'\colon \bC_\mQ\rightarrow \bC/\bE$ tel que $T=T'\circ F$. 
On voit aussitôt que $T'$ et $F'$ sont des équivalences de catégories quasi-inverses l'une de l'autre.

\subsubsection{}\label{higgs3-caip1c}
\addtocounter{equation}{1}
Tout foncteur additif (resp. exact)  entre catégories additives (resp. abéliennes) $\bC\rightarrow \bC'$ 
s'étend de manière unique en un foncteur additif (resp. exact) $\bC_\mQ\rightarrow \bC'_\mQ$, compatible aux foncteurs 
de localisation.  

\subsubsection{}\label{higgs3-caip1d}
\addtocounter{equation}{1}
Si $\bC$ est une catégorie abélienne, le foncteur de localisation $\bC\rightarrow \bC_\mQ$ transforme 
les objets injectifs en des objets injectifs (\cite{gabriel} III cor.~1 à prop.~1). 
En particulier, si $\bC$ possède suffisamment d'injectifs, il en est de même de $\bC_\mQ$. 

\subsection{}\label{higgs3-caip4}\index{1000610@$M\otimes_{A_\mQ}N$, $\cHom_{A_\mQ}(M,N)$}
Soit $(X,A)$ un topos annelé. On désigne par $\bMod(A)$ la catégorie des $A$-modules de $X$
et par $\bMod_{\mQ}(A)$, au lieu de $\bMod(A)_\mQ$, la catégorie des $A$-modules de $X$ à isogénie 
près \eqref{higgs3-caip1a}. Le produit tensoriel des $A$-modules induit un bifoncteur 
\begin{equation}\label{higgs3-caip4a}
\bMod_{\mQ}(A)\times \bMod_{\mQ}(A)\rightarrow \bMod_{\mQ}(A),\ \ \ (M,N)\mapsto M\otimes_{A_\mQ}N,
\end{equation} 
faisant de $\bMod_{\mQ}(A)$ une catégorie monoïdale symétrique, ayant $A_\mQ$ pour objet unité. 
Les objets de $\bMod_{\mQ}(A)$ seront aussi appelés des {\em $A_\mQ$-modules}. Cette terminologie
se justifie en considérant $A_\mQ$ comme un monoïde de $\bMod_{\mQ}(A)$. 
Si $M$ et $N$ sont deux $A_\mQ$-modules, on note $\Hom_{A_\mQ}(M,N)$ 
le groupe des morphismes de $M$ dans $N$ dans $\bMod_{\mQ}(A)$. Le bifoncteur ``faisceau des morphismes''
de la catégorie des $A$-modules de $X$ induit un bifoncteur 
\begin{equation}\label{higgs3-caip4b}
\bMod_{\mQ}(A)\times \bMod_{\mQ}(A)\rightarrow \bMod_{\mQ}(A),\ \ \ (M,N)\mapsto \cHom_{A_\mQ}(M,N).
\end{equation} 
Les bifoncteurs \eqref{higgs3-caip4a} et \eqref{higgs3-caip4b} héritent des mêmes propriétés d'exactitudes que les bifoncteurs 
sur la catégorie $\bMod(A)$ qui leurs ont donné naissance.

\subsection{}\label{higgs3-caip2}
Pour tout morphisme de topos annelés $f\colon (Y,B)\rightarrow (X,A)$, on désigne encore par 
\begin{eqnarray}
f^*\colon \bMod_\mQ(A)&\rightarrow& \bMod_{\mQ}(B),\label{higgs3-caip2a}\\
f_*\colon \bMod_\mQ(B)&\rightarrow &\bMod_{\mQ}(A),\label{higgs3-caip2aa}
\end{eqnarray}
les foncteurs induits par les foncteurs image inverse et image directe par $f$,
de sorte que le premier est un adjoint à gauche du second. Le premier foncteur est exact et le second 
foncteur est exact à gauche. On note 
\begin{eqnarray}
\rR f_*\colon \bD^+(\bMod_\mQ(B))&\rightarrow& \bD^+(\bMod_{\mQ}(A)),\label{higgs3-caip2b}\\
\rR^q f_*\colon \bMod_\mQ(B)&\rightarrow& \bMod_{\mQ}(A), \ \ \ (q\in \mN),\label{higgs3-caip2bb}
\end{eqnarray}
les foncteurs dérivés droits de $f_*$ \eqref{higgs3-caip2aa}. 
Ces notations n'induisent aucune confusion avec celles des foncteurs dérivés 
droits du foncteur $f_*\colon \bMod(B)\rightarrow \bMod(A)$, puisque le foncteur de localisation
$\bMod(B)\rightarrow \bMod_\mQ(B)$ est exact et transforme les objets injectifs en des objets injectifs.

\begin{defi}\label{higgs3-caip5}
Soit $(X,A)$ un topos annelé. On dit qu'un $A_\mQ$-module $M$ est {\em plat} (ou {\em $A_\mQ$-plat}) 
si le foncteur $N\mapsto M\otimes_{A_\mQ}N$ de la catégorie $\bMod_{\mQ}(A)$
dans elle-même est exact.
\end{defi}

On peut faire les remarques suivantes~:

\subsubsection{}\label{higgs3-caip5a}
Soit $M$ un $A$-module. Pour que $M_\mQ$ soit $A_\mQ$-plat, il faut et il suffit que 
pour tout morphisme injectif de $A$-modules $u\colon N\rightarrow N'$, 
le noyau de $u\otimes\id_M$ soit d'exposant fini \eqref{higgs3-caip1b}. 

\subsubsection{}\label{higgs3-caip5b}
Si $M$ est un $A$-module plat, alors $M_\mQ$ est $A_\mQ$-plat. 

\subsubsection{}\label{higgs3-caip5c}
Soient $B$ une $A$-algèbre, $M$ un $A_\mQ$-module plat. Alors $M\otimes _{A_\mQ}B_\mQ$ est $B_\mQ$-plat.

\subsubsection{}\label{higgs3-caip5d}
Soit $B$ une $A$-algèbre telle que le foncteur 
\[
\bMod(A)\rightarrow \bMod(B), \ \ \ N\mapsto N\otimes_AB
\]
soit exact et fidèle. Pour qu'un $A_\mQ$-module $M$ soit plat, 
il faut et il suffit que le $B_\mQ$-module $M\otimes_{A_\mQ}B_\mQ$ soit plat.

\subsection{}\label{higgs3-caip3}
Soient $(X,A)$ un topos annelé, $U$ un objet de $X$. On désigne par 
$j_U\colon X_{/U}\rightarrow X$ le morphisme de localisation de $X$ en $U$.
Pour tout $F\in \ob(X)$, le faisceau $j_U^*(F)$ sera aussi noté $F|U$. Le topos $X_{/U}$ sera annelé par $A|U$. 
Le foncteur prolongement par zéro $j_{U!}\colon \bMod(A|U)\rightarrow \bMod(A)$ étant exact et fidèle (\cite{sga4} IV 11.3.1), 
il induit un foncteur exact et fidèle que l'on note encore
\begin{equation}\label{higgs3-caip3a}
j_{U!}\colon \bMod_\mQ(A|U)\rightarrow \bMod_\mQ(A), \ \ \ P\mapsto j_{U!}(P),
\end{equation}
et que l'on appelle encore le {\em prolongement par zéro}. 
C'est un adjoint à gauche du foncteur  \eqref{higgs3-caip2a}
\begin{equation}\label{higgs3-caip3c}
j_U^*\colon \bMod_\mQ(A)\rightarrow \bMod_\mQ(A|U).
\end{equation}
Pour tout $A_\mQ$-module $M$, le $(A|U)_\mQ$-module $j_U^*(M)$ sera aussi noté $M|U$. 
Pour tout $A$-module $N$, on a par définition $(N|U)_\mQ=N_\mQ|U$.

\begin{lem}\label{higgs3-caip6}
Soient $(X,A)$ un topos annelé, $U$ un objet de $X$.  Alors~:
\begin{itemize}
\item[{\rm (i)}] Pour tout $(A_\mQ|U)$-module plat $P$, $j_{U!}(P)$ est $A_\mQ$-plat. 
\item[{\rm (ii)}] Pour tout $A_\mQ$-module plat $M$, $j_U^*(M)$ est $(A_\mQ|U)$-plat. 
\end{itemize}
\end{lem}

En effet, il résulte aussitôt de (\cite{sga4} IV 12.11) que pour tout $A_\mQ$-module $M$
et tout $(A_\mQ|U)$-module $P$, on a un isomorphisme canonique fonctoriel
\begin{equation}\label{higgs3-caip6a}
j_{U!}(P\otimes_{(A_\mQ|U)}j_U^*(M))\stackrel{\sim}{\rightarrow} j_{U!}(P)\otimes_{A_\mQ}M.
\end{equation}

(i) Cela résulte de \eqref{higgs3-caip6a} et du fait que les foncteurs $j_U^*$ et $j_{U!}$ sont exacts. 

(ii) Il résulte de \eqref{higgs3-caip6a} et du fait que le foncteur $j_{U!}$ est exact que le foncteur 
\begin{equation}
\bMod_\mQ(A|U)\rightarrow \bMod_\mQ(A), \ \ \ P\mapsto j_{U!}(P\otimes_{(A_\mQ|U)}j_U^*(M))
\end{equation}
est exact. Comme le foncteur $j_{U!}$ est de plus fidèle \eqref{higgs3-caip3}, on en déduit que
le foncteur $P\mapsto P\otimes_{(A_\mQ|U)}j_U^*(M)$ sur la catégorie $\bMod_\mQ(A|U)$ est exact~; 
d'où la proposition. 

\begin{lem}\label{higgs3-caip7}
Soient $(X,A)$ un topos annelé, $(U_i)_{1\leq i\leq n}$ un recouvrement 
fini de l'objet final de $X$, $M, N$ deux $A_\mQ$-modules. Pour tous $1\leq i,j\leq n$, 
on pose $U_{ij}=U_i\times U_j$. Alors~:
\begin{itemize}
\item[{\rm (i)}] Le diagramme d'applications d'ensembles
\begin{equation}\label{higgs3-caip7a}
\Hom_{A_\mQ}(M,N)\rightarrow \prod_{1\leq i\leq n}\Hom_{(A_\mQ|U_i)}(M|U_i,N|U_i)\rightrightarrows 
\prod_{1\leq i,j\leq n}\Hom_{(A_\mQ|U_{i,j})}(M|U_{ij},N|U_{ij})
\end{equation}
est exact.
\item[{\rm (ii)}] Pour que $M$ soit nul, il faut et il suffit que pour tout $1\leq i\leq n$, $M|U_i$ soit nul.
\item[{\rm (iii)}] Pour que $M$ soit $A_\mQ$-plat, il faut et il suffit que pour tout $1\leq i\leq n$, $M|U_i$ soit $(A_\mQ|U_i)$-plat.
\end{itemize}
\end{lem}

(i) Soient $M^\circ,N^\circ$ deux $A$-modules tels que $M=M^\circ_\mQ$
et $N=N^\circ_\mQ$, $u,v\colon M^\circ\rightarrow N^\circ$ deux morphismes $A$-linéaires. 
Supposons que pour tout $1\leq i\leq n$, on ait $u_\mQ|U_i=v_\mQ|U_i$. Il existe alors un entier $m\not=0$
tel que $m\cdot u|U_i=m\cdot v|U_i$. On en déduit que $m\cdot u=m\cdot v$, d'où l'exactitude à gauche de \eqref{higgs3-caip7a}.
Par ailleurs, soient, pour tout $1\leq i\leq n$, $u_i\colon M^\circ|U_i\rightarrow N^\circ|U_i$ un morphisme $(A|U_i)$-linéaire 
tels que $(u_{i,\mQ})_{1\leq i\leq n}$ soit dans le noyau de
la double flèche de \eqref{higgs3-caip7a}. Il existe alors un entier $m'\not=0$ tel que 
pour tout $1\leq i,j\leq n$, on ait $m'\cdot u_i|U_{ij}=m'\cdot u_j|U_{ij}$. Par suite, les morphismes
$(m'\cdot u_i)_{1\leq i\leq n}$ se recollent en un morphisme $A$-linéaire $w\colon M^\circ\rightarrow N^\circ$. 
Il est clair que $(u_{i,\mQ})_{1\leq i\leq n}$ est l'image canonique de $m'^{-1}w_\mQ$, d'où  
l'exactitude au centre de \eqref{higgs3-caip7a}.

(ii) En effet, $M$ est nul si et seulement si $\id_M=0$. L'assertion résulte donc de (i). 

(iii) En effet, la condition est nécessaire en vertu de \ref{higgs3-caip6}(ii), et elle est suffisante compte tenu de
(ii) et (\cite{sga4} IV 12.11).

\subsection{}\label{higgs3-imh1}\index{Isogenie de Higgs@Isogénie de Higgs}\index{1000630@$\bIH(A,E)$, $\bIH_\mQ(A,E)$}
Soient $(X,A)$ un topos annelé, $E$ un $A$-module. On appelle {\em $A$-isogénie de Higgs à coefficients dans $E$}
la donnée d'un quadruplet
\begin{equation}\label{higgs3-imh1a}
(M,N,u\colon M\rightarrow N,\theta\colon M\rightarrow N\otimes_AE)
\end{equation}
formé de deux $A$-modules $M$ et $N$ et de deux morphismes $A$-linéaires $u$ et $\theta$  
vérifiant la propriété suivante~: il existe un entier $n\not=0$ et un morphisme $A$-linéaire $v\colon N\rightarrow M$ tels que 
$v\circ u=n\cdot \id_M$, $u\circ v=n\cdot \id_N$ et que $(M,(v\otimes \id_E)\circ \theta)$ et $(N,\theta\circ v)$ 
soient des $A$-modules de Higgs à coefficients dans $E$ (\cite{ag1} 2.8). 
On notera que $u$ induit une isogénie de modules de Higgs de $(M,(v\otimes \id_E)\circ \theta)$
dans $(N,\theta\circ v)$ \eqref{higgs3-caip1}, d'où la terminologie.  
Soient $(M,N,u,\theta)$, $(M',N',u',\theta')$ deux $A$-isogénies de Higgs à coefficients dans $E$. 
Un morphisme de $(M,N,u,\theta)$ dans $(M',N',u',\theta')$ est la donnée de deux morphismes $A$-linéaires
$\alpha\colon M\rightarrow M'$ et $\beta\colon N\rightarrow N'$ tels que $\beta\circ u=u'\circ \alpha$
et $(\beta\otimes \id_E)\circ \theta=\theta'\circ \alpha$. 
On désigne par $\bIH(A,E)$ la catégorie des $A$-isogénies de Higgs à coefficients dans $E$. C'est une catégorie 
additive. On note $\bIH_\mQ(A,E)$ la catégorie des objets de $\bIH(A,E)$ à isogénie près. 

\subsection{}\label{higgs3-imh2}\index{Complexe de Dolbeault d'une isogenie de Higgs@Complexe de Dolbeault d'une isogénie de Higgs}
Soient $(X,A)$ un topos annelé, $E$ un $A$-module, 
$(M,N,u,\theta)$ une $A$-isogénie de Higgs à coefficients dans $E$. Pour tout $i\geq 1$, on désigne par
\begin{equation}\label{higgs3-imh2a}
\theta_i\colon M\otimes_A \wedge^iE \rightarrow N\otimes_A \wedge^{i+1}E
\end{equation}
le morphisme $A$-linéaire défini pour toutes sections locales 
$m$ de $M$ et $\omega$ de $\wedge^iE$ par $\theta_i(m\otimes \omega)=\theta(m)\wedge \omega$.
On note 
\begin{equation}\label{higgs3-imh2b}
\otheta_i\colon M_\mQ\otimes_{A_\mQ} (\wedge^iE)_\mQ \rightarrow M_\mQ\otimes_{A_\mQ} (\wedge^{i+1}E)_\mQ
\end{equation}
le morphisme de $\bMod_\mQ(A)$ composé de l'image de $\theta_i$ et de l'inverse de l'image  
de $u\otimes \id_{ \wedge^{i+1}E}$. Soient $v\colon N\rightarrow M$ un morphisme $A$-linéaire, 
$n$ un entier non nul tels que $v\circ u=n\cdot \id_M$ et que $(M,(v\otimes \id_E)\circ \theta)$ 
soit un $A$-module de Higgs à coefficients dans $E$. Notons
\begin{equation}\label{higgs3-imh2c}
\vartheta_i\colon M\otimes_A \wedge^iE \rightarrow M\otimes_A \wedge^{i+1}E
\end{equation}
le morphisme $A$-linéaire induit par $(v\otimes \id_E)\circ \theta$ (\cite{ag1} (2.8.3)).
L'image canonique de $\vartheta_i$ dans $\bMod_\mQ(A)$ est alors égale à $n\cdot\otheta_i$. 
On en déduit que $\otheta_{i+1}\circ \otheta_i=0$ (cf. \cite{ag1} (2.8.2)). 
On appelle complexe de {\em Dolbeault} de $(M,N,u,\theta)$
et l'on note $\mK^\bullet(M,N,u,\theta)$ le complexe de cochaînes de $\bMod_\mQ(A)$ 
\begin{equation}\label{higgs3-imh2d}
M_\mQ\stackrel{\otheta_0}{\longrightarrow}M_\mQ\otimes_{A_\mQ}E_\mQ\stackrel{\otheta_1}{\longrightarrow} 
M_\mQ\otimes_{A_\mQ}(\wedge^2E)_\mQ\rightarrow \dots,
\end{equation}
où $M_\mQ$ est placé en degré $0$ et les différentielles sont de degré $1$. On obtient ainsi un foncteur
de la catégorie $\bIH(A,E)$ dans la catégorie des complexes de $\bMod_\mQ(A)$.
Toute isogénie de $\bIH(A,E)$ 
induit un isomorphisme des complexes de Dolbeault associés. Le foncteur ``complexe de Dolbeault''
induit donc un foncteur de $\bIH_\mQ(A,E)$ dans la catégorie des complexes de $\bMod_\mQ(A)$.

\subsection{}\label{higgs3-isoco1}\index{isoconnexion (integrable)@$\lambda$-isoconnexion (intégrable)}
Soient $(X,A)$ un topos annelé, $B$ une $A$-algèbre, $\lambda\in \Gamma(X,A)$. 
On appelle {\em $\lambda$-isoconnexion relativement à l'extension $B/A$}
(ou simplement {\em $\lambda$-isoconnexion} lorsqu'il n'y a aucun risque de confusion) la donnée d'un quadruplet   
\begin{equation}\label{higgs3-isoco1a}
(M,N,u\colon M\rightarrow N,\nabla\colon M\rightarrow \Omega^1_{B/A}\otimes_BN)
\end{equation}
où $M$ et $N$ sont des $B$-modules, $u$ est une isogénie de $B$-modules \eqref{higgs3-caip1}
et $\nabla$ est un morphisme $A$-linéaire tel que pour toutes sections locales 
$x$ de $B$ et $t$ de $M$, on ait 
\begin{equation}\label{higgs3-isoco1b}
\nabla(xt)=\lambda d(x) \otimes u(t)+x\nabla(t).
\end{equation} 
Pour tout morphisme $B$-linéaire $v\colon N\rightarrow M$ pour lequel il existe un entier $n$ 
tel que  $u\circ v=n\cdot \id_N$ et $v\circ u=n\cdot \id_M$, 
les couples $(M,(\id\otimes v)\circ \nabla)$ et $(N,\nabla\circ v)$ sont des modules à
$(n\lambda)$-connexions (\cite{ag1} 2.10),
et $u$ est un morphisme de $(M,(\id\otimes v)\circ \nabla)$ dans $(N,\nabla\circ v)$. On dit que la $\lambda$-isoconnexion 
$(M,N,u,\nabla)$ est {\em intégrable} s'il existe un morphisme $B$-linéaire $v\colon N\rightarrow M$ et un entier $n\not= 0$ 
tels que  $u\circ v=n\cdot \id_N$, $v\circ u=n\cdot \id_M$ et que les $(n\lambda)$-connexions
$(\id\otimes v)\circ \nabla$ sur $M$ et $\nabla\circ v$ sur $N$ soient intégrables. 

Soient $(M,N,u,\nabla)$, $(M',N',u',\nabla')$ deux $\lambda$-isoconnexions. 
Un morphisme de $(M,N,u,\nabla)$ dans $(M',N',u',\nabla')$ est la donnée de 
deux morphismes $B$-linéaires $\alpha\colon M\rightarrow M'$ et $\beta\colon N\rightarrow N'$
tels que  $\beta\circ u=u'\circ \alpha$ et $(\id \otimes \beta)\circ \nabla=\nabla'\circ \alpha$. 

\subsection{}\label{higgs3-isoco2}
Soient $f\colon (X',A')\rightarrow (X,A)$ un morphisme de topos annelés, $B$ une $A$-algèbre, 
$B'$ une $A'$-algèbre, $\alpha\colon f^*(B)\rightarrow B'$ un homomorphisme de $A'$-algèbres,
$\lambda\in \Gamma(X,A)$, $(M,N,u,\nabla)$ une $\lambda$-isoconnexion relativement à l'extension $B/A$. 
Notons $\lambda'$ l'image canonique de $\lambda$ dans $\Gamma(X',A')$, 
$d'\colon B'\rightarrow \Omega^1_{B'/A'}$ la $A'$-dérivation universelle de $B'$ et 
\begin{equation}
\gamma\colon f^*(\Omega^1_{B/A}) \rightarrow \Omega^1_{B'/A'}
\end{equation}
le morphisme $\alpha$-linéaire canonique. On voit aussitôt que 
$(f^*(M),f^*(N),f^*(u),f^*(\nabla))$ est une $\lambda'$-isoconnexion relativement à l'extension $f^*(B)/A'$,
qui est intégrable si $(M,N,u,\nabla)$ l'est.

Il existe un unique morphisme $A'$-linéaire 
\begin{equation}
\nabla'\colon B'\otimes_{f^*(B)}f^*(M)\rightarrow  \Omega^1_{B'/A'}\otimes_{f^*(B)}f^*(N)
\end{equation}
tel que pour toutes sections locales $x'$ de $B'$ et $t$ de $f^*(M)$, on ait 
\begin{equation}
\nabla'(x'\otimes t)=\lambda' d'(x')\otimes f^*(u)(t)+x'(\gamma\otimes\id_{f^*(N)})(f^*(\nabla)(t)).
\end{equation}
Le quadruplet $(B'\otimes_{f^*(B)}f^*(M),B'\otimes_{f^*(B)}f^*(N),\id_{B'}\otimes_{f^*(B)} f^*(u),\nabla')$ 
est une $\lambda'$-isoconnexion relativement à l'extension 
$B'/A'$, qui est intégrable si $(M,N,u,\nabla)$ l'est. 

\subsection{}\label{higgs3-isoco3}
Soient $(X,A)$ un topos annelé, $B$ une $A$-algèbre, $\lambda\in \Gamma(X,A)$, $(M,N,u,\nabla)$ une 
$\lambda$-isocon\-nexion intégrable relativement à l'extension $B/A$. Supposons qu'il existe un $A$-module $E$
et un $B$-isomorphisme $\gamma\colon E\otimes_AB\stackrel{\sim}{\rightarrow}\Omega^1_{B/A}$ 
tels que pour toute section locale $\omega$ de $E$, on ait $d(\gamma(\omega\otimes 1))=0$.
Notons $\theta\colon M\rightarrow E\otimes_AN$  le morphisme induit par $\nabla$ et $\gamma$.
Alors $(M,N,u,\theta)$ est une $A$-isogénie de Higgs à coefficients dans $E$ (cf. \cite{ag1} 2.12).

Soit $(M',N',u',\theta')$ une $A$-isogénie de Higgs à coefficients dans $E$.  
Il existe un unique morphisme $A$-linéaire 
\begin{equation}
\nabla'\colon M\otimes_AM'\rightarrow  \Omega^1_{B/A}\otimes_BN\otimes_AN'
\end{equation}
tel que pour toutes sections locales $t$ de $M$ et $t'$ de $M'$, on ait 
\begin{equation}
\nabla'(t\otimes t')=\nabla(t)\otimes_Au'(t')+ (\gamma\otimes_B\id_{N\otimes_A N'})(u(t)\otimes_A \theta'(t')).
\end{equation}
Le quadruplet $(M\otimes_AM',N\otimes_AN',u\otimes u',\nabla')$ est une $\lambda$-isoconnexion intégrable.

\subsection{}\label{higgs3-isoco4}\index{isoconnexion adique (integrable)@$\lambda$-isoconnexion adique (intégrable)}
Soient $A$ un anneau adique, $I$ un idéal de définition de $A$, $\lambda\in A$, $B$ une $A$-algèbre adique, 
{\em i.e.}, $B$ est une $A$-algèbre complète et séparée pour la topologie $(IB)$-adique.  On rappelle que la topologie 
canonique du $B$-module  $\Omega^1_{B/A}$ est déduite de celle de $B$ (\cite{ega4} 0.20.4.5).
On désigne par $\hOmega^1_{B/A}$ son séparé complété et on note encore 
\begin{equation}
d\colon B\rightarrow \hOmega^1_{B/A}
\end{equation}
la $A$-dérivation continue universelle de $B$. On appelle {\em $\lambda$-isoconnexion adique} 
(ou {\em $\lambda$-isoconnexion $I$-adique}) relativement à l'extension $B/A$ la donnée d'un quadruplet   
\begin{equation}
(M,N,u\colon M\rightarrow N,\nabla\colon M\rightarrow \hOmega^1_{B/A}\hotimes_BN)
\end{equation}
où $M$ et $N$ sont des $B$-modules complets et séparés pour les topologies $(IB)$-adiques, 
$u$ est une isogénie de $B$-modules \eqref{higgs3-caip1} et $\nabla$ est un morphisme $A$-linéaire tel que pour tous 
$x\in B$ et $t\in M$, on ait 
\begin{equation}
\nabla(xt)=\lambda d(x) \hotimes u(t)+x\nabla(t).
\end{equation} 
Pour tout morphisme $B$-linéaire $v\colon N\rightarrow M$ pour lequel il existe un entier $n$ 
tel que  $u\circ v=n\cdot \id_N$ et $v\circ u=n\cdot \id_M$, 
les couples $(M,(\id\hotimes v)\circ \nabla)$ et $(N,\nabla\circ v)$ sont des modules à
$(n\lambda)$-connexions adiques (\cite{ag1} 2.14),
et $u$ est un morphisme de $(M,(\id\hotimes v)\circ \nabla)$ dans $(N,\nabla\circ v)$. On dit que la $\lambda$-isoconnexion adique
$(M,N,u,\nabla)$ est {\em intégrable} s'il existe un morphisme $B$-linéaire $v\colon N\rightarrow M$ et un entier $n\not= 0$ 
tels que  $u\circ v=n\cdot \id_N$, $v\circ u=n\cdot \id_M$ et que les $(n\lambda)$-connexions adiques 
$(\id\hotimes v)\circ \nabla$ sur $M$ et $\nabla\circ v$ sur $N$ soient intégrables (cf. \cite{ag1} 2.14).

Soient $(M,N,u,\nabla)$, $(M',N',u',\nabla')$ deux $\lambda$-isoconnexions adiques relativement à l'extension $B/A$.
Un morphisme de $(M,N,u,\nabla)$ dans $(M',N',u',\nabla')$ est la donnée de 
deux morphismes $B$-linéaires $\alpha\colon M\rightarrow M'$ et $\beta\colon N\rightarrow N'$
tels que  $\beta\circ u=u'\circ \alpha$ et $(\id \hotimes \beta)\circ \nabla=\nabla'\circ \alpha$. 

\subsection{}\label{higgs3-isoco5}
Soient $A$ un anneau adique, $\lambda\in A$, $B$ une $A$-algèbre adique, 
$(M,N,u,\nabla)$ une $\lambda$-isoconnexion adique intégrable relativement à l'extension $B/A$. 
Supposons les conditions suivantes remplies~:
\begin{itemize}
\item[(i)] $A$ admet un idéal de définition de type fini $I$, et posant $A_1=A/I$ et $B_1=B\otimes_AA_1$, 
le $B_1$-module $\Omega^1_{B_1/A_1}$ est de type fini.
\item[(ii)] Il existe un $A$-module {\em libre de type fini} $E$
et un isomorphisme $B$-linéaire $\gamma\colon E\otimes_AB\stackrel{\sim}{\rightarrow}\hOmega^1_{B/A}$ 
tels que $\gamma(E)\subset d(B)$. 
\end{itemize}
On observera que $\hOmega^1_{B/A}\hotimes_BN=\hOmega^1_{B/A}\otimes_BN=E\otimes_AN$.
Notons $\theta\colon M\rightarrow E\otimes_AN$  le morphisme induit par 
$\nabla$ et $\gamma$. Il résulte alors de (\cite{ag1} 2.16) que $(M,N,u,\theta)$ est une $A$-isogénie 
de Higgs à coefficients dans $E$.

\subsection{}\label{higgs3-formel1}
On rappelle que $\cS$ désigne le schéma formel $\Spf(\co_C)$ \eqref{higgs3-not1}.
Soient $\fX$ un $\cS$-schéma formel localement de présentation finie (cf. \cite{egr1} 2.3.15), 
$\cJ$ un idéal de définition cohérent de $\fX$, $\cF$ un $\co_\fX$-module. 
Le schéma formel $\fX$ est donc idyllique (\cite{egr1} 2.6.13). 
Suivant (\cite{egr1} 2.10.1), on appelle {\em clôture rigide} de $\cF$ et l'on note $\crig(\cF)$, 
le $\co_\fX$-module 
\begin{equation}\label{higgs3-formel1a}
\crig(\cF)=\underset{\underset{n\geq 0}{\longrightarrow}}{\lim}\ \cHom_{\co_\fX}(\cJ^n,\cF).
\end{equation}
Cette notion ne dépend pas de l'idéal $\cJ$. Par ailleurs, on pose 
\begin{equation}\label{higgs3-formel1b}
\cF[\frac 1 p]=\cF_{\mQ_p}=\cF\otimes_{\mZ_p}\mQ_p.
\end{equation} 
Comme $p\co_\fX$ est un idéal de définition de $\fX$, le morphisme canonique 
$\cF_{\mQ_p}\rightarrow \crig(\cF)$ est un isomorphisme d'après (\cite{egr1} 2.10.5). 
On dit que $\cF$ est {\em rig-nul} si le morphisme canonique 
$\cF\rightarrow \cF_{\mQ_p}$ est nul  (cf. \cite{egr1} 2.10.1.4). 

Considérons les conditions suivantes~:
\begin{itemize}
\item[(i)] $\cF_{\mQ_p}=0$.
\item[(ii)] $\cF$ est rig-nul. 
\item[(iii)] Il existe un entier $n\geq 1$ tel que $p^n\cF=0$. 
\end{itemize}
D'après (\cite{egr1} 2.10.10), on a alors {\rm (iii)}$\Rightarrow${\rm (i)}$\Leftrightarrow${\rm (ii)}. 
De plus, si $\fX$ est quasi-compact et si $\cF$ est de type fini, les trois conditions sont équivalentes. 
On en déduit que si $\fX$ est quasi-compact et si $\cF$ est de type fini, 
pour tout $\co_\fX$-module $\cG$, l'homomorphisme canonique 
\begin{equation}\label{higgs3-formel1c}
\Hom_{\co_\fX}(\cF,\cG)\otimes_{\mZ_p}\mQ_p\rightarrow \Hom_{\co_\fX[\frac 1 p]}(\cF_{\mQ_p},\cG_{\mQ_p})
\end{equation}
est injectif.

\begin{lem}\label{higgs3-formel2}
Soit $\fX$ un $\cS$-schéma formel de présentation finie. On note
$\bMod^\coh(\co_{\fX})$ (resp. $\bMod^{\coh}(\co_{\fX}[\frac 1 p])$) la catégorie des $\co_{\fX}$-modules 
(resp. $\co_{\fX}[\frac 1 p]$-modules) cohérents et 
$\bMod^{\coh}_\mQ(\co_{\fX})$ la catégorie des $\co_\fX$-modules cohérents à isogénie près. 
Alors le foncteur canonique
\begin{equation}\label{higgs3-formel2a}
\bMod^{\coh}(\co_{\fX})\rightarrow \bMod^{\coh}(\co_{\fX}[\frac 1 p]),\ \ \ \cF\mapsto \cF_{\mQ_p},
\end{equation}
induit une équivalence de catégories abéliennes
\begin{equation}\label{higgs3-formel2b}
\bMod^{\coh}_\mQ(\co_{\fX})\stackrel{\sim}{\rightarrow} \bMod^{\coh}(\co_{\fX}[\frac 1 p]). 
\end{equation} 
\end{lem}
On notera d'abord que le foncteur \eqref{higgs3-formel2a} est bien défini en vertu de (\cite{egr1} 2.10.24(i))
et qu'il induit un foncteur exact 
\begin{equation}\label{higgs3-formel2c}
\bMod^{\coh}_\mQ(\co_{\fX})\rightarrow \bMod^{\coh}(\co_{\fX}[\frac 1 p]).
\end{equation}
Celui-ci est essentiellement surjectif en vertu de (\cite{egr1} 2.10.24(ii)). 
Montrons qu'il est pleinement fidèle. Soient $\cF$, $\cG$ deux $\co_\fX$-modules cohérents. 
L'homomorphisme canonique 
\begin{equation}\label{higgs3-formel2d}
\Hom_{\co_\fX}(\cF,\cG)\otimes_{\mZ_p}\mQ_p\rightarrow \Hom_{\co_\fX[\frac 1 p]}(\cF_{\mQ_p},\cG_{\mQ_p})
\end{equation}
est injectif d'après \eqref{higgs3-formel1c}. 
D'autre part, pour tout morphisme $\co_\fX$-linéaire $v\colon\cF\rightarrow \cG_{\mQ_p}$,
il existe un entier $n\geq 0$ tel que $v(p^n\cF)$ soit contenu dans l'image du morphisme 
canonique $c_\cG\colon \cG\rightarrow \cG_{\mQ_p}$. Comme $\cG_\tor=\ker(c_\cG)$ est cohérent (\cite{egr1} 2.10.14), 
il existe un entier $m\geq 0$ tel que $p^m\cG_\tor=0$. On en déduit qu'il existe un morphisme $\co_\fX$-linéaire
$w\colon \cF\rightarrow \cG$ tel que $c_\cG\circ w=p^{n+m}v$. 
L'homomorphisme \eqref{higgs3-formel2d} est donc surjectif~; d'où la proposition. 

\begin{lem}\label{higgs3-formel20}
Soient $\fX=\Spf(R)$ un $\cS$-schéma formel affine de présentation finie, $\cF$ un $\co_\fX[\frac 1 p]$-module cohérent.
Pour que $\cF$ soit un $\co_\fX[\frac 1 p]$-module localement projectif de type fini \eqref{higgs3-not3},
il faut et il suffit que $\Gamma(\fX,\cF)$ soit un $R[\frac 1 p]$-module projectif de type fini.
\end{lem}

On rappelle que l'anneau $R[\frac 1 p]$ est noethérien (\cite{egr1} 1.10.2(i)). 
D'après \ref{higgs3-formel2} et (\cite{egr1} 2.7.2), il existe un $R$-module cohérent $M$ tel que $\cF=(M^\Delta)_{\mQ_p}$.   
On a $\Gamma(\fX,\cF)=M_{\mQ_p}$ en vertu de (\cite{egr1}  (2.10.5.1)). 

Supposons d'abord que $\cF$ soit un $\co_\fX[\frac 1 p]$-module localement projectif de type fini et montrons 
que  $\Gamma(\fX,\cF)$ est un $R[\frac 1 p]$-module projectif de type fini. D'après 
(\cite{egr1}  2.7.4, (2.10.5.1) et 5.1.11), il existe une $R$-algèbre fidèlement plate et topologiquement
de présentation finie $R'$ telle que le $R'$-module $\Gamma(\fX,\cF)\otimes_RR'$ soit projectif de type fini.
On en déduit par descente fidèlement plate que $\Gamma(\fX,\cF)$ est un $R[\frac 1 p]$-module projectif de type fini.

Supposons ensuite que $\Gamma(\fX,\cF)$ soit un $R[\frac 1 p]$-module projectif de type fini et montrons que 
$\cF$ est un $\co_\fX[\frac 1 p]$-module localement projectif de type fini.
Compte tenu de (\cite{egr1} 1.10.2(iii) et (2.10.5.1)), on peut supposer le $R[\frac 1 p]$-module $\Gamma(\fX,\cF)$ libre de type fini.
Il existe alors un entier $n\geq 0$ et un morphisme $R$-linéaire $R^n\rightarrow M$ 
dont le noyau et le conoyau sont de torsion. Par suite, $\cF$ est un $\co_\fX[\frac 1 p]$-module libre de type fini.

\begin{lem}\label{higgs3-formel21}
Soient $A$ un anneau, $t\in A$, $M$ un $A$-module de type fini. On suppose que $A$ est complet 
et séparé pour la topologie $(tA)$-adique, que $t$ n'est pas un diviseur de zéro dans $A$ et 
que $M_t$ est un $A_t$-module projectif. On note $\hM$ le séparé complété de $M$ pour la topologie
 $(tA)$-adique. Alors le morphisme canonique $M\rightarrow \hM$ induit un isomorphisme $M_t\stackrel{\sim}{\rightarrow} \hM_t$.
 \end{lem}
 
En effet, notons $M'$ l'image du morphisme canonique $M\rightarrow M_t$. 
Il existe un entier $n\geq 1$ et un morphisme $A_t$-linéaire injectif $u\colon M_t\rightarrow A_t^n$. 
Comme $t$ n'est pas un diviseur de zéro dans $A$, il existe un entier $i\geq 0$
tel que $M'\subset u^{-1}(t^{-i}A)^n$. Par suite, $M'$ est contenu dans un $A$-module libre de type fini,
et est donc séparé  pour la topologie $(tA)$-adique. 
Par ailleurs, $M$ et $M'$ sont complets pour les topologies $(tA)$-adiques 
d'après (\cite{ac} chap.~III § 2.12 cor.~1 de prop.~16). On en déduit que le morphisme canonique 
$M\rightarrow \hM$ est surjectif et que $M'$ est complet et séparé pour la topologie $(tA)$-adique. 
Par suite, le morphisme surjectif canonique $M\rightarrow M'$ se factorise en $M\rightarrow \hM\rightarrow M'$. 
Comme $M_t\rightarrow M'_t$ est un isomorphisme, $M_t\rightarrow \hM_t$ est aussi un isomorphisme.

\subsection{}\label{higgs3-formel3}
Soient $\fX$ un $\cS$-schéma formel de présentation finie, $\cE$ un $\co_\fX$-module. 
On désigne par $\bMH(\co_\fX[\frac 1 p], \cE_{\mQ_p})$ la catégorie des $\co_\fX[\frac 1 p]$-modules 
de Higgs à coefficients dans $\cE_{\mQ_p}$ (\cite{ag1} 2.8), par $\bMH^\coh(\co_\fX[\frac 1 p], \cE_{\mQ_p})$
la sous-catégorie pleine formée des modules de Higgs dont le $\co_\fX[\frac 1 p]$-module sous-jacent est cohérent,
par $\bIH(\co_\fX,\cE)$ la catégorie des $\co_\fX$-isogénies de Higgs à coefficients dans $\cE$ \eqref{higgs3-imh1}
et par $\bIH^\coh(\co_\fX,\cE)$ la sous-catégorie pleine formée des quadruplets $(\cM,\cN,u,\theta)$ tels que 
les $\co_\fX$-modules $\cM$ et $\cN$ soient cohérents. Ce sont des catégories additives. 
On note $\bIH_\mQ(\co_\fX,\cE)$ (resp. $\bIH^\coh_\mQ(\co_\fX,\cE)$) la 
catégorie des objets de $\bIH(\co_\fX,\cE)$ (resp. $\bIH^\coh(\co_\fX,\cE)$) à isogénie près \eqref{higgs3-caip1a}.
On a un foncteur 
\begin{equation}\label{higgs3-formel3a}
\bIH(\co_\fX,\cE)\rightarrow \bMH(\co_\fX[\frac 1 p], \cE_{\mQ_p}), \ \ \ 
(\cM,\cN,u,\theta)\mapsto (\cM_{\mQ_p}, (u_{\mQ_p}^{-1}\otimes \id_{\cE_{\mQ_p}})\circ\theta_{\mQ_p}).
\end{equation}

\begin{lem}\label{higgs3-formel4}
Les hypothèses étant celles de \eqref{higgs3-formel3}, soient, de plus, $(\cM,\cN,u,\theta)$, 
$(\cM',\cN',u',\theta')$ deux objets de $\bIH(\co_\fX,\cE)$ tels que $\cM$ et $\cN$ soient des $\co_\fX$-modules 
de type fini, $(\cM_{\mQ_p},\ttheta)$ et $(\cM'_{\mQ_p},\ttheta')$ leurs images respectives par le foncteur \eqref{higgs3-formel3a}. 
Alors l'homomorphisme canonique 
\begin{eqnarray}\label{higgs3-formel4a}
\lefteqn{\Hom_{\bIH(\co_\fX,\cE)}((\cM,\cN,u,\theta),(\cM',\cN',u',\theta'))\otimes_{\mZ_p}\mQ_p\rightarrow}  \\
&&\ \ \ \ \ \ \ \ \ \ \ \ \ \ \ \ \ \ \ \ \ \ \ \ \ \ \ \ 
\Hom_{\bMH(\co_\fX[\frac 1 p], \cE_{\mQ_p})}((\cM_{\mQ_p},\ttheta),(\cM_{\mQ_p},\ttheta'))\nonumber 
\end{eqnarray}
est injectif.
\end{lem}
En effet, soient $\alpha\colon \cM\rightarrow \cM'$ et $\beta\colon \cN\rightarrow \cN'$ deux morphismes 
$\co_\fX$-linéaires définissant un morphisme de $(\cM,\cN,u,\theta)$ dans $(\cM',\cN',u',\theta')$ de 
$\bIH(\co_\fX,\cE)$, dont l'image $\alpha_{\mQ_p}$ par l'homomorphisme \eqref{higgs3-formel4a} est nulle. 
Comme $(\alpha(\cM))_{\mQ_p}=0$ et que $\cM$ est de type fini sur $\co_\fX$, 
il existe un entier $n\geq 0$ tel que $p^n\alpha=0$ \eqref{higgs3-formel1}. De même, comme $\beta_{\mQ_p}=0$,
il existe un entier $m\geq 0$ tel que $p^m\beta=0$. La proposition s'ensuit.

\begin{lem}\label{higgs3-formel5}
Les hypothèses étant celles de \eqref{higgs3-formel3}, supposons de plus $\cE$ cohérent.  
Alors le foncteur 
\begin{equation}\label{higgs3-formel5a}
\bIH^\coh_\mQ(\co_\fX,\cE)\rightarrow \bMH^\coh(\co_\fX[\frac 1 p], \cE_{\mQ_p})
\end{equation}
induit par \eqref{higgs3-formel3a} est une équivalence de catégories.
\end{lem}

Soient $N$ un $\co_\fX[\frac 1 p]$-module cohérent, $\theta$ un $\co_\fX[\frac 1 p]$-champ de Higgs 
sur $N$ à coefficients dans $\cE_{\mQ_p}$. D'après (\cite{egr1} 2.10.24(ii)), il existe 
un $\co_\fX$-module cohérent $\cN$ et un isomorphisme $u\colon \cN_{\mQ_p}\stackrel{\sim}{\rightarrow}N$. 
On peut supposer $\cN$ sans $p$-torsion (\cite{egr1} 2.10.14).
D'après la preuve de \ref{higgs3-formel2}, il existe un entier $n\geq 0$ et un morphisme 
$\co_\fX$-linéaire $\vartheta \colon \cN\rightarrow \cN\otimes_{\co_\fX}\cE$ tels que  le diagramme 
\begin{equation}\label{higgs3-formel5d}
\xymatrix{
\cN\ar[r]\ar[d]_\vartheta&N\ar[d]^{p^n\theta}\\
{\cN\otimes_{\co_\fX}\cE}\ar[r]&{N\otimes_{\co_\fX[\frac 1 p]}\cE_{\mQ_p}}}
\end{equation}
où les flèches horizontales sont induites par $u$, soit commutatif. Quitte à multiplier $\vartheta$
par une puissance de $p$, on peut supposer que $\vartheta\wedge \vartheta=0$ \eqref{higgs3-formel1c}
({\em i.e.}, que $\vartheta$ est un champ de Higgs). Comme $\cN$ est sans $p$-torsion, 
le morphisme $\vartheta$ se factorise en deux morphismes $\co_\fX$-linéaires
\begin{equation}
\xymatrix{
\cN\ar[r]^-(0.5){\nu_n}&{p^n\cN}\ar[r]^-(0.5){\vartheta_n}&{\cN\otimes_{\co_\fX}\cE}},
\end{equation}
où $\nu_n$ est l'isomorphisme induit par la multiplication par $p^n$ sur $\cN$. 
Le composé 
\begin{equation}
\xymatrix{
{p^n\cN}\ar[r]^-(0.5){\vartheta_n}&{\cN\otimes_{\co_\fX}\cE}\ar[rr]^{\nu_n\otimes \id_\cE}&&{(p^n\cN)\otimes_{\co_\fX}\cE}}
\end{equation}
est alors aussi un champ de Higgs. Notons $\iota_n\colon p^n\cN\rightarrow \cN$ l'injection canonique. 
On a $\vartheta\circ \iota_n =p^n\vartheta_n$, de sorte que le diagramme
\begin{equation}
\xymatrix{
{p^n\cN}\ar[r]^{\iota_n}\ar[d]_{\vartheta_n}\ar[r]&\cN\ar[r]&N\ar[d]^{\theta}\\
{\cN\otimes_{\co_\fX}\cE}\ar[rr]&&{N\otimes_{\co_\fX[\frac 1 p]}\cE_{\mQ_p}}}
\end{equation}
est commutatif. 
Par suite, $(p^n\cN,\cN,\iota_n,\vartheta_n)$ est  un objet de $\bIH^\coh(\co_\fX,\cE)$
dont l'image par le foncteur \eqref{higgs3-formel3a} est isomorphe à $(N,\theta)$. 
Le foncteur \eqref{higgs3-formel5a} est donc essentiellement surjectif. On sait d'après \ref{higgs3-formel4} qu'il est fidèle. 
Montrons qu'il est plein. Soient $(\cM,\cN,u,\theta)$ et $(\cM',\cN',u',\theta')$ deux objets de 
$\bIH^\coh(\co_\fX,\cE)$,
$(\cM_{\mQ_p},\ttheta)$ et $(\cM'_{\mQ_p},\ttheta')$ leurs images respectives par le foncteur \eqref{higgs3-formel5a}, 
$\lambda\colon \cM_{\mQ_p}\rightarrow \cM'_{\mQ_p}$ un morphisme $\co_\fX[\frac 1 p]$-linéaire tel que 
$(\lambda \otimes \id_{\cE_{\mQ_p}})\circ \ttheta=\ttheta'\circ \lambda$. D'après la preuve de \ref{higgs3-formel2}, 
il existe un entier $n\geq 0$ et un morphisme $\co_\fX$-linéaire $\alpha\colon \cM\rightarrow \cM'$ tels que  le diagramme 
\begin{equation}\label{higgs3-formel5c}
\xymatrix{
\cM\ar[r]\ar[d]_\alpha&{\cM_{\mQ_p}}\ar[d]^{p^n\lambda}\\
{\cM'}\ar[r]&{\cM'_{\mQ_p}}}
\end{equation}
où les flèches horizontales sont les morphismes canoniques, soit commutatif. 
D'après (\cite{egr1} 2.10.22(i) et 2.10.10), il existe un entier $m\geq 0$ tel que $p^m$
annule le noyau et le conoyau de $u$.  
Il existe donc un morphisme $\co_\fX$-linéaire $\beta\colon \cN\rightarrow \cN'$ tel que $\beta\circ u=u'\circ (p^{2m}\alpha)$. 
Notons 
\[
c\colon \cN'\otimes_{\co_\fX}\cE\rightarrow \cN'_{\mQ_p}\otimes_{\co_\fX[\frac 1 p]}\cE_{\mQ_p}
\] 
le morphisme canonique. 
Comme $(\cN'\otimes_{\co_\fX}\cE)_\tor=\ker(c)$ est cohérent (\cite{egr1} 2.10.14), 
il existe un entier $q\geq 0$ tel que $p^q \ker(c)=0$. 
La relation $(\lambda \otimes \id_{\cE_{\mQ_p}})\circ \ttheta=\ttheta'\circ \lambda$ implique alors que 
$( (p^q\beta) \otimes \id_\cE)\circ \theta=\theta'\circ (p^{q+2m}\alpha)$. 
Le foncteur \eqref{higgs3-formel5a} est donc plein. 

\begin{prop}\label{higgs3-formel6}
Soient $\fX$ un $\cS$-schéma formel de présentation finie, $f\colon \fX'\rightarrow \fX$ un morphisme de présentation finie
et fidèlement plat {\rm (\cite{egr1} 5.1.7)}, 
$\fX''=\fX'\times_\fX\fX'$, $\rp_1,\rp_2\colon \fX''\rightarrow \fX'$ les projections canoniques. 
\begin{itemize}
\item[{\rm (i)}] Soient $\cF$ et $\cG$ deux $\co_\fX[\frac 1 p]$-modules cohérents, $\cF'$ et $\cG'$ leurs images inverses sur 
$\fX'$, $\cF''$ et $\cG''$ leurs images inverses sur $\fX''$. Alors, le diagramme d'applications d'ensembles 
\begin{equation}
\Hom_{\co_\fX[\frac 1 p]}(\cF,\cG)\rightarrow \Hom_{\co_{\fX'}[\frac 1 p]}(\cF',\cG')\rightrightarrows
\Hom_{\co_{\fX''}[\frac 1 p]}(\cF'',\cG'')
\end{equation}
défini par les foncteurs de changement de base par $f$, $\rp_1$ et $\rp_2$ est exact. 
\item[{\rm (ii)}] Pour tout $\co_\fX[\frac 1 p]$-module cohérent $\cF'$, toute donnée de descente sur $\cF'$
relative à $f$ est effective.  
\end{itemize}
\end{prop}

Cela résulte de la preuve de (\cite{egr1} 5.11.12).

\begin{rema}
Soient $\fX$ un $\cS$-schéma formel de présentation finie, $\cF$ et $\cG$ deux $\co_\fX[\frac 1 p]$-modules cohérents.
Notons $\fX^\rig$ l'espace rigide associé à $\fX$ (\cite{egr1} 4.1.6),
$\fX^\rig_\ad$ le topos admissible de $\fX^\rig$ (\cite{egr1} 4.4.1) et 
\begin{equation}
\varrho_{\fX}\colon (\fX^\rig_\ad,\co_{\fX^\rig})\rightarrow (\fX_\zar,\co_\fX[\frac 1 p])
\end{equation}
le morphisme canonique de topos annelés  (\cite{egr1} (4.7.5.1)). 
D'après (\cite{egr1} 2.10.24(ii), 4.7.8 et 4.7.28), le morphisme d'adjonction $\cF\rightarrow \varrho_{\fX*}(\varrho_\fX^*(\cF))$
est un isomorphisme. Par suite, l'application 
\begin{equation}
\Hom_{\co_{\fX^\rig}}(\varrho_\fX^*(\cF),\varrho_\fX^*(\cG))\rightarrow\Hom_{\co_\fX[\frac 1 p]}(\cF,\cG),
\ \ \ u\mapsto \varrho_{\fX*}(u)
\end{equation}
est bijective~; c'est l'isomorphisme d'adjonction. 
La proposition \ref{higgs3-formel6}(i) résulte alors aussitôt de (\cite{egr1} 5.11.12(i)).
La proposition \ref{higgs3-formel6}(ii) ne résulte pas formellement de l'énoncé de (\cite{egr1} 5.11.12(ii)), 
mais plutôt de la même preuve. 
\end{rema}

\section{Systèmes projectifs d'un topos}\label{higgs3-spsa}

\subsection{}\label{higgs3-spsa1}\index{1000701@$X^{I^\circ}$ ($X$ topos, $I$ catégorie)}\index{1000702@$\alpha_i\colon X\rightarrow X^{I^\circ}$}
Dans cette section, $X$ désigne un $\mU$-topos et $I$ une $\mU$-petite catégorie \eqref{higgs3-not0}. 
On considère toujours $X$ comme muni de sa topologie canonique, qui en fait un $\mU$-site. 
On munit $X\times I$ de la topologie totale relative au site fibré constant $X\times I\rightarrow I$ 
de fibre $X$ (cf. \cite{sga4} VI 7.4.1 et \cite{ag2} 7.1), qui en fait un $\mU$-site. 
On rappelle (\cite{sga4} VI 7.4.7) que le topos des faisceaux 
de $\mU$-ensembles sur $X\times I$ est canoniquement équivalent à la catégorie 
$\bHom(I^\circ,X)$  des foncteurs de $I^\circ$ dans $X$, que l'on note encore $X^{I^\circ}$. 
En particulier, $X^{I^\circ}$ est un $\mU$-topos. Ce dernier fait peut se voir directement (\cite{sga4} IV 1.2).
On renvoie à (\cite{ag2} 7.4) pour la description des anneaux et des modules de $X^{I^\circ}$.

Pour tout $i\in \ob(I)$, on note 
\begin{equation}\label{higgs3-spsa1'a}
\alpha_{i!}\colon X\rightarrow X\times I
\end{equation} 
le foncteur qui à un objet $F$ de $X$ associe le couple $(F,i)$. 
Celui-ci étant cocontinu (\cite{sga4} VI 7.4.2), il définit un morphisme de topos (\cite{sga4} IV 4.7)
\begin{equation}\label{higgs3-spsa1a}
\alpha_i\colon X\rightarrow X^{I^\circ}.
\end{equation}
D'après (\cite{sga4} VI 7.4.7), pour tous $F\in \ob(X^{I^\circ})$ et $i\in \ob(I)$, on a
\begin{equation}\label{higgs3-spsa1b}
\alpha_i^*(F)=F(i). 
\end{equation}

On note encore 
\begin{equation}\label{higgs3-spsa1c}
\alpha_{i!}\colon X\rightarrow X^{I^\circ}
\end{equation}
le composé de $\alpha_{i!}$ \eqref{higgs3-spsa1'a} et du foncteur canonique 
$X\times I\rightarrow X^{I^\circ}$. Pour tous $j\in \ob(I)$ et $F\in \ob(X)$, on a 
\begin{equation}\label{higgs3-spsa1ca}
\alpha_{i!}(F)(j)=F\times (\Hom_I(j,i))_X, 
\end{equation}
où $(\Hom_I(j,i))_X$ est le faisceau constant sur $X$ de valeur $\Hom_I(j,i)$. 
D'après (\cite{sga4} VI 7.4.3(4)), le foncteur $\alpha_{i!}$ \eqref{higgs3-spsa1c} est un adjoint à gauche de $\alpha_i^*$. 

Pour tout morphisme $f\colon i\rightarrow j$ de $I$, on a un morphisme 
\begin{equation}\label{higgs3-spsa1bb}
\rho_f\colon \alpha_i\rightarrow \alpha_j,
\end{equation}
défini au niveau des images inverses, pour tout  $F\in \ob(X^{I^\circ})$, 
par le morphisme $F(j)\rightarrow F(i)$ induit par $f$ (\cite{sga4} VI (7.4.5.2)). Si $f$
et $g$ sont deux morphismes composables de $I$, on a 
$\rho_{gf}=\rho_g\rho_f$.

\begin{remas}\label{higgs3-spsa98}
(i)\  Soient $U\in \ob(X)$, $i\in \ob(I)$. Pour qu'une famille $((U_n,i_n)\rightarrow (U,i))_{n\in N}$ 
soit couvrante pour la topologie totale de $X\times I$, il faut et il suffit qu'elle soit raffinée 
par une famille $((V_m,i)\rightarrow (U,i))_{m\in M}$, où $(V_m\rightarrow U)_{m\in M}$ est une famille 
couvrante de $X$ (\cite{sga4} VI 7.4.2(1)). 

(ii)\  Il résulte de (i) que pour qu'un crible $\cR$ de $X^{I^\circ}$ soit épimorphique strict universel, 
i.e. couvre l'objet final de $X^{I^\circ}$ pour la topologie canonique, 
il faut et il suffit que pour tout $i\in \ob(I)$, il existe un raffinement $(U_{i,n})_{n\in N_i}$ de l'objet final de $X$ 
tel que pour tout $n\in N_i$, $\alpha_{i!}(U_{i,n})$ soit un objet de $\cR$ \eqref{higgs3-spsa1c}.

(iii)\ Supposons que les produits fibrés soient représentables dans $I$. 
La topologie totale sur $X\times I$ coïncide alors avec la topologie co-évanescente 
relative au site fibré constant $X\times I\rightarrow I$ de fibre $X$, 
lorsque l'on munit $I$ de la topologie chaotique (\cite{ag2} 5.4). 
Autrement dit, la topologie totale sur $X\times I$ est engendrée par 
la {\em prétopologie} formée des recouvrements verticaux (\cite{ag2} 5.3 et 5.7).  
\end{remas}

\begin{remas}\label{higgs3-spsa99}
(i)\ Les $\mU$-limites inductives (resp. les limites projectives finies) dans $X^{I^\circ}$ 
se calculent terme à terme, autrement dit, 
pour tout foncteur $\varphi\colon N\rightarrow X^{I^\circ}$ tel que la catégorie $N$ soit $\mU$-petite (resp. finie), 
si $F$ est un objet de $X^{I^\circ}$ qui représente la limite inductive (resp. projective) de $\varphi$, 
alors pour tout $i\in \ob(I)$, la limite inductive (resp. projective) de $\alpha_i^*\circ \varphi$ est représentable 
par $\alpha_i^*(F)$. 

(ii)\ Soit $\iota$ un objet final de $I$. D'après \eqref{higgs3-spsa1ca}, pour tout $j\in \ob(I)$, 
$\alpha_j^*\alpha_{\iota!}$ est le foncteur identique de $X$. Donc en vertu de (i), $\alpha_{\iota!}$ est exact à gauche,
et le couple $(\alpha_{\iota!},\alpha_\iota^*)$ forme un morphisme de topos 
\begin{equation}\label{higgs3-spsa99a}
\beta_\iota\colon X^{I^\circ}\rightarrow X.
\end{equation}
Le morphisme $\beta_\iota \alpha_\iota\colon X\rightarrow X$ est isomorphe au morphisme identique 
(cf. \cite{sga4} VI 7.4.12). 

(iii)\ Soit $i$ un objet de $I$ qui n'est pas final. Il existe alors $j \in \ob(I)$ tel que $\Hom_I(j,i)$
ne soit pas un singleton. En particulier, $\alpha_j^*\alpha_{i!}$ ne transforme pas 
l'objet final en l'objet final \eqref{higgs3-spsa1ca}. Par suite, $\alpha_{i!}$ n'est pas exact à gauche, et 
le couple $(\alpha_{i!},\alpha_{i}^*)$ ne forme pas un morphisme 
de topos, contrairement à ce qui a été affirmé dans (\cite{ekedahl2} ligne 20 page 59). 
Toutefois, l'isomorphisme (4.4) et la suite spectrale (4.5) de {\em loc. cit.} sont corrects en vertu de (\cite{ag2} 7.7 et 7.8). 
\end{remas}

\subsection{}\label{higgs3-spsa3}\index{1000715@$\lambda\colon X^{I^\circ}\rightarrow X$}
Le foncteur 
\begin{equation}\label{higgs3-spsa3a}
\lambda^*\colon X\rightarrow X^{I^\circ}
\end{equation}
qui à un objet $F$ de $X$ associe le foncteur constant $I^\circ\rightarrow X$ de valeur $F$,
est exact à gauche en vertu de \ref{higgs3-spsa99}(i). Il admet pour adjoint à droite le foncteur 
\begin{equation}\label{higgs3-spsa3b}
\lambda_*\colon X^{I^\circ}\rightarrow X
\end{equation}
qui a un foncteur $I^\circ\rightarrow X$ associe sa limite projective (\cite{sga4} II 4.1(3)). 
Le couple $(\lambda^*,\lambda_*)$ définit donc un morphisme de topos 
\begin{equation}\label{higgs3-spsa3c}
\lambda\colon X^{I^\circ}\rightarrow X.
\end{equation}

\subsection{}\label{higgs3-spsa35}\index{1000720@$f^{I^\circ}$ ($f$ morphisme de topos, $I$ catégorie)}
Tout morphisme de $\mU$-topos  $f\colon X\rightarrow Y$ induit un morphisme cartésien de topos fibrés
constants au-dessus de $I$
\begin{equation}\label{higgs3-spsa1d}
f\times \id_I \colon X\times I\rightarrow Y\times I. 
\end{equation}
D'après (\cite{sga4} VI 7.4.10), celui-ci induit un morphisme de topos 
\begin{equation}\label{higgs3-spsa1e}
f^{I^\circ}\colon X^{I^\circ}\rightarrow Y^{I^\circ}
\end{equation}
tel que pour tout $F\in X^{I^\circ}$, on ait (\cite{sga4} VI (7.4.9.2))
\begin{equation}\label{higgs3-spsa1f}
(f^{I^\circ})_*(F)= f_*\circ F.
\end{equation}
On voit aussitôt que pour tout  $G\in Y^{I^\circ}$, on a
\begin{equation}\label{higgs3-spsa1g}
(f^{I^\circ})^*(G)= f^*\circ G.
\end{equation}

Notons $\rR^q(f^{I^\circ})_*$ $(q\in \mN)$ les foncteurs dérivés droits du foncteur $(f^{I^\circ})_*$ 
pour les groupes abéliens. D'après (\cite{ag2} 7.7), pour tout groupe abélien 
$F$ de $X^{I^\circ}$ et tout $i\in \ob(I)$, on a un isomorphisme canonique fonctoriel
\begin{equation}\label{higgs3-spsa1h}
\rR^q(f^{I^\circ})_*(F)(i)\stackrel{\sim}{\rightarrow} \rR^qf_*(F(i)).
\end{equation}

\begin{prop}\label{higgs3-spsa100}
Supposons que les produits fibrés soient représentables dans $I$~; soient,
de plus, $U$ un objet de $X$, $j_{U!}\colon X_{/U}\rightarrow X$ le foncteur canonique,
$j_U\colon X_{/U}\rightarrow X$ le morphisme de localisation de $X$ en $U$, 
$j_{\lambda^*(U)}\colon (X^{I^\circ})_{/\lambda^*(U)}\rightarrow X^{I^\circ}$ 
le morphisme de localisation de $X^{I^\circ}$ en $\lambda^*(U)$. Alors~:
\begin{itemize}
\item[{\rm (i)}] La topologie totale sur $X_{/U}\times I$ est induite par la topologie totale 
sur $X\times I$ via le foncteur $j_{U!}\times \id_I$.
\item[{\rm (ii)}] Il existe une équivalence canonique de topos 
\begin{equation}\label{higgs3-spsa100a}
h\colon (X_{/U})^{I^\circ}\stackrel{\sim}{\rightarrow} (X^{I^\circ})_{/\lambda^*(U)},
\end{equation}
telle que $(j_U)^{I^\circ}=j_{\lambda^*(U)} \circ h$ \eqref{higgs3-spsa1e}. En particulier, pour tout $F\in \ob(X^{I^\circ})$, 
$F\times \lambda^*(U)$ s'identifie au foncteur $I^\circ\rightarrow X_{/U}, i\mapsto F(i)\times U$. 
\end{itemize}
\end{prop}

On notera d'abord que la topologie canonique de $X_{/U}$
est induite par la topologie canonique de $X$ via le foncteur $j_{U!}$, et que le foncteur ``extension
par le vide'' s'identifie dans ce cas au foncteur $j_{U!}$ (\cite{sga4} IV 1.2, III 3.5 et 5.4), d'où la notation.  

(i) La topologie totale sur $X\times I$ 
est engendrée par la prétopologie formée des recouvrements verticaux d'après \ref{higgs3-spsa98}(iii); 
et de même pour $X_{/U}\times I$.
La proposition résulte donc de (\cite{sga4} III 3.3 et II 1.4). 

(ii) Pour tous $V\in \ob(X)$ et $i\in \ob(I)$, on a un isomorphisme canonique \eqref{higgs3-spsa1b}
\begin{equation}\label{higgs3-spsa100b}
\lambda^*(U)(V\times i)=\Hom_{X^{I^\circ}}(\alpha_{i!}(V),\lambda^*(U))\stackrel{\sim}{\rightarrow}
\Hom_{X}(V,\alpha_{i}^*(\lambda^*(U)))=U(V).
\end{equation}
Par suite, le foncteur $j_{U!}\times \id_I\colon X_{/U}\times I\rightarrow X\times I$ se factorise canoniquement en 
\begin{equation}\label{higgs3-spsa100c}
\xymatrix{
{X_{/U} \times I}\ar[r]^-(0.4)e_-(0.4)\sim&{(X\times I)_{/\lambda^*(U)}}\ar[r]^-(0.4){j'_{\lambda^*(U)}}&{X\times I}},
\end{equation}
où $e$ est une équivalence de catégories et $j'_{\lambda^*(U)}$ est le foncteur canonique. 
D'après (i) et (\cite{sga4} III 5.4), $e$ induit une équivalence de topos  
\begin{equation}\label{higgs3-spsa100d}
h\colon (X_{/U})^{I^\circ}\stackrel{\sim}{\rightarrow} (X^{I^\circ})_{/\lambda^*(U)}.
\end{equation}
Pour tout $F\in \ob(X^{I^\circ})$, on a 
$j_{\lambda^*(U)}^*(F)=F\circ j'_{\lambda^*(U)}$ d'après (\cite{sga4} III 2.3 et 5.2(2)). 
Donc $h^*(j_{\lambda^*(U)}^*(F))=F\circ (j_{U!}\times \id_I)$.
Comme $F\circ (j_{U!}\times \id_I)=((j_U)^{I^\circ})^*(F)$ \eqref{higgs3-spsa1g},
on en déduit que $(j_U)^{I^\circ}=j_{\lambda^*(U)} \circ h$.

\subsection{}\label{higgs3-spsa2}\index{1000730@$X^{\mN^\circ}$, $X^{[n]^\circ}$}
Soit $J$ un ensemble ordonné $\mU$-petit. On note encore $J$ la catégorie définie par $J$, 
c'est-à-dire la catégorie ayant pour objets les éléments de $J$, avec au plus une flèche de source et de but donnés,
et pour tous $i,j\in J$, l'ensemble $\Hom_J(i,j)$ est non-vide si et seulement si $i\leq j$.   
Il est souvent commode d'utiliser pour les objets $F\colon J^\circ\rightarrow X$ de $X^{J^\circ}$
la notation indicielle $(F_j)_{j\in J}$ ou même $(F_j)$, où $F_j=F(j)$ pour tout $j\in J$.
On dit qu'un objet $(F_j)$ de $X^{J^\circ}$ est {\em strict} si pour tous éléments $i\geq j$ de $J$,
le morphisme de transition $F_i\rightarrow F_j$ est un épimorphisme.

Nous nous limitons dans cet article aux cas où $J$ est soit l'ensemble ordonné des entiers naturels $\mN$, 
soit l'un des sous-ensembles ordonnés $[n]=\{0,1,\dots,n\}$. 
On observera que dans chacun de ces cas, les produits fibrés sont représentables dans $J$.

\begin{lem}\label{higgs3-spsa6}
Soient $n$ un entier $\geq 0$, $\iota_n\colon X\times [n]\rightarrow X\times \mN$ le foncteur d'injection canonique. Alors~:
\begin{itemize}
\item[{\rm (i)}] La topologie totale sur $X\times [n]$ est induite par la topologie totale sur $X\times \mN$ via le foncteur~$\iota_n$. 
\item[{\rm (ii)}] Le foncteur $\iota_n$ est continu et cocontinu pour les topologies totales. Notons 
\begin{equation}\label{higgs3-spsa6a}
\varphi_n\colon X^{[n]^\circ}\rightarrow X^{\mN^\circ}
\end{equation}
le morphisme associé de topos. Pour tout $F=(F_i)_{i\in \mN}\in \ob(X^{\mN^\circ})$, on a 
\begin{equation}\label{higgs3-spsa6b}
\varphi_n^*(F)=(F_i)_{i\in [n]}.
\end{equation}
\item[{\rm (iii)}] Le foncteur $\varphi_n^*\colon \colon X^{\mN^\circ}\rightarrow X^{[n]^\circ}$ admet pour adjoint à gauche le foncteur
\begin{equation}\label{higgs3-spsa6c}
\varphi_{n!}\colon X^{[n]^\circ}\rightarrow X^{\mN^\circ},
\end{equation}
défini pour tout $F=(F_i)_{i\in [n]}\in \ob(X^{[n]^\circ})$ par
\begin{equation}\label{higgs3-spsa6d}
\varphi_{n!}(F)=(F_i)_{i\in \mN},
\end{equation}
où pour tout $i\geq n+1$, $F_i=\emptyset$ est l'objet initial de $X$. 
\end{itemize}
\end{lem}

(i) La topologie totale sur $X\times\mN$ 
est engendrée par la prétopologie formée des recouvrements verticaux d'après \ref{higgs3-spsa98}(iii)~; 
et de même pour $X\times [n]$. La proposition résulte donc de (\cite{sga4} III 3.3 et II 1.4). 

(ii) Notons $\hX$ (resp. $(X\times \mN)^\wedge$, resp. $(X\times [n])^\wedge$) 
la catégorie des préfaisceaux de $\mU$-ensembles sur $X$ (resp. $X\times \mN$, resp. $X\times [n]$). 
On a alors une équivalence de catégories \eqref{higgs3-spsa1'a}
\begin{equation}\label{higgs3-spsa6e}
(X\times \mN)^\wedge\stackrel{\sim}{\rightarrow} 
\hX^{\mN^\circ}=\bHom(\mN^\circ,\hX), \ \ \ F\mapsto (F\circ \alpha_{i!})_{i\in \mN};
\end{equation}
et de même pour $(X\times [n])^\wedge$.
Le foncteur $\iota_n$ induit par composition le foncteur 
\begin{equation}\label{higgs3-spsa6f}
\hiota_n^*\colon \hX^{\mN^\circ}\rightarrow \hX^{[n]^\circ}, \ \ \ (F_i)_{i\in \mN}\mapsto (F_i)_{i\in [n]}.
\end{equation}
Celui-ci admet pour adjoint à droite le foncteur
\begin{equation}\label{higgs3-spsa6g}
\hiota_{n*}\colon \hX^{[n]^\circ}\rightarrow \hX^{\mN^\circ}, \ \ \ (F_i)_{i\in [n]}\mapsto (F_i)_{i\in \mN},
\end{equation}
où pour tout $i\geq n+1$, $F_i=F_n$ et le morphisme $F_i\rightarrow F_{i-1}$ est l'identité de $F_n$. 
Le foncteur $\hiota_n^*$ admet pour adjoint à gauche le foncteur 
\begin{equation}\label{higgs3-spsa6h}
\hiota_{n!}\colon \hX^{[n]^\circ}\rightarrow \hX^{\mN^\circ}, \ \ \ (F_i)_{i\in [n]}\mapsto (F_i)_{i\in \mN},
\end{equation}
où pour tout $i\geq n+1$, $F_i=\emptyset$ est l'objet initial de $\hX$. 
Il est clair que $\hiota_n^*$ transforme les faisceaux sur $X\times \mN$ en des faisceaux sur $X\times [n]$
et que $\hiota_{n*}$ transforme les faisceaux sur $X\times [n]$ en des faisceaux sur $X\times \mN$.
Par suite, $i_n$ est continu et cocontinu. La continuité résulte aussi de (i). La formule \eqref{higgs3-spsa6b}
est une conséquence de \eqref{higgs3-spsa6f} et (\cite{sga4} III 2.3). 

(iii) Cela résulte de \eqref{higgs3-spsa6h} et (\cite{sga4} III 1.3).

\begin{lem}\label{higgs3-spsa7}
Soient $U$ un objet de $X$, $j_U\colon X_{/U}\rightarrow X$ le morphisme de localisation de $X$ en $U$, 
$n\in \mN$, $j_{(U,n)}\colon (X^{\mN^\circ})_{/\alpha_{n!}(U)}\rightarrow X^{\mN^\circ}$ le morphisme de localisation
de $X^{\mN^\circ}$ en  $\alpha_{n!}(U)$ \eqref{higgs3-spsa1c}. On a alors une équivalence canonique de topos 
\begin{equation}\label{higgs3-spsa7a}
h\colon (X_{/U})^{[n]^\circ}\stackrel{\sim}{\rightarrow}(X^{\mN^\circ})_{/\alpha_{n!}(U)},
\end{equation}
telle que $j_{(U,n)}\circ h$ soit le composé 
\begin{equation}\label{higgs3-spsa7b}
(X_{/U})^{[n]^\circ}\stackrel{\varphi_n}{\longrightarrow} (X_{/U})^{\mN^\circ}
\stackrel{(j_U)^{\mN^\circ}}{\longrightarrow}X^{\mN^\circ},
\end{equation}
où la première flèche est le morphisme \eqref{higgs3-spsa6a} et la seconde flèche est le morphisme \eqref{higgs3-spsa1e}. 
\end{lem}

Notons $u\colon \alpha_{n!}(U)\rightarrow \lambda^*(U)$ l'adjoint de l'isomorphisme canonique
$U\stackrel{\sim}{\rightarrow}\alpha_n^*(\lambda^*(U))$ et $\talpha_{n!}(U)$ l'image de $(U,n)$ par le foncteur canonique 
$X_{/U}\times \mN\rightarrow (X_{/U})^{\mN^\circ}$ \eqref{higgs3-spsa1c}. D'après \ref{higgs3-spsa100}(ii), 
on a une équivalence canonique de topos 
\begin{equation}\label{higgs3-spsa7c}
g\colon (X_{/U})^{\mN^\circ}\stackrel{\sim}{\rightarrow} (X^{\mN^\circ})_{/\lambda^*(U)}.
\end{equation}
Compte tenu de \eqref{higgs3-spsa100c} et (\cite{sga4} III 5.4), on a $g(\talpha_{n!}(U))=u$. 
On peut donc se borner au cas où $U$ est l'objet final $e_X$ de $X$ en vertu de \ref{higgs3-spsa100}(ii) et (\cite{sga4} IV 5.6). 
 
Le foncteur d'injection canonique $\iota_n\colon X\times [n]\rightarrow X\times \mN$ se factorise canoniquement en  
\begin{equation}\label{higgs3-spsa7d}
\xymatrix{
{X\times [n]}\ar[r]^-(0.4)\nu_-(0.4)\sim&{(X\times \mN)_{/(e_X,n)}}\ar[r]^-(0.4){j'_n}&{X\times \mN}},
\end{equation}
où $\nu$ est une équivalence de catégories et $j'_n$ est le foncteur canonique. 
D'après \ref{higgs3-spsa6}(i) et (\cite{sga4} III 5.4), $\nu$ induit une équivalence de topos  
\begin{equation}\label{higgs3-spsa7e}
h\colon X^{[n]^\circ}\stackrel{\sim}{\rightarrow} (X^{\mN^\circ})_{/\alpha_{n!}(e_X)}.
\end{equation}
Pour tout $F=(F_i)_{i\in \mN}\in \ob(X^{\mN^\circ})$, on a $j^*_{(e_X,n)}(F)=F\circ j'_n$
(\cite{sga4} III 2.3 et 5.2(2)). Donc $h^*(j^*_{(e_X,n)}(F))=F\circ \iota_n$.
Comme $F\circ \iota_n=(F_i)_{i\in [n]}=\varphi_n^*(F)$ \eqref{higgs3-spsa6b}, on en déduit que 
$j_{(e_X,n)}\circ h=\varphi_n$.

\begin{prop}[\cite{ag2} 7.9] \label{higgs3-spsa55}
Soient $n$ un entier $\geq 0$, $A=(A_i)_{i\in [n]}$ un anneau de $X^{[n]^\circ}$, 
$M=(M_i)_{i\in [n]}$ un $A$-module de $X^{[n]^\circ}$.
Pour tout entier $q\geq 0$, on a alors un isomorphisme canonique
\begin{equation}\label{higgs3-spsa55a}
\rH^q(X^{[n]^\circ},M)\stackrel{\sim}{\rightarrow } \rH^q(X,M_n).
\end{equation}
\end{prop}

\begin{prop}[\cite{ag2} 7.10] \label{higgs3-spsa5}
Soient $A=(A_n)_{n\in \mN}$ un anneau de $X^{\mN^\circ}$, 
$U$ un objet de $X$, $M=(M_n)_{n\in \mN}$ un $A$-module de $X^{\mN^\circ}$.
Pour tout entier $q\geq 0$, on a alors une suite exacte canonique et fonctorielle
\begin{equation}\label{higgs3-spsa5a}
0\rightarrow \rR^1\underset{\underset{n\in \mN^\circ}{\longleftarrow}}{\lim}\ \rH^{q-1}(U,M_n)\rightarrow 
\rH^q(\lambda^*(U),M)\rightarrow \underset{\underset{n\in \mN^\circ}{\longleftarrow}}{\lim}\ \rH^q(U,M_n)\rightarrow 0,
\end{equation}
où l'on a posé $\rH^{-1}(U,M_n)=0$ pour tout $n\in \mN$. 
\end{prop}

En effet, on peut se borner au cas où $U$ est l'objet final de $X$ \eqref{higgs3-spsa100},
auquel cas la proposition est un cas particulier de (\cite{ag2} 7.10).

\subsection{}\label{higgs3-spsa200}
Soient $A=(A_n)_{n\in \mN}$ un anneau de $X^{\mN^\circ}$, 
$M=(M_n)_{n\in \mN}$ et $N=(N_n)_{n\in \mN}$ deux $A$-modules de $X^{\mN^\circ}$.
On a alors un isomorphisme canonique bifonctoriel 
\begin{equation}\label{higgs3-spsa2a}
M\otimes_AN\stackrel{\sim}{\rightarrow} (M_n\otimes_{A_n}N_n)_{n\in \mN}.
\end{equation}
En effet, pour tout $n\in \mN$, on a un isomorphisme canonique (\cite{sga4} IV 13.4)
\begin{equation}\label{higgs3-spsa2b}
\alpha_n^*(M\otimes_AN)\stackrel{\sim}{\rightarrow} \alpha_n^*(M)\otimes_{\alpha_n^*(A)}\alpha_n^*(N).
\end{equation}
De même, pour tout entier $q\geq 0$, on a des isomorphismes canoniques fonctoriels \eqref{higgs3-not9}
\begin{eqnarray}
\rS^q_A(M)&\stackrel{\sim}{\rightarrow}&(\rS^q_{A_n}(M_n))_{n\in \mN},\label{higgs3-spsa2c}\\
\wedge^q_A(M)&\stackrel{\sim}{\rightarrow}&(\wedge^q_{A_n}(M_n))_{n\in \mN}.\label{higgs3-spsa2d}
\end{eqnarray}

\begin{lem}\label{higgs3-spad14}
Soient $A=(A_n)_{n\in \mN}$ un anneau de $X^{\mN^\circ}$, 
$M=(M_n)_{n\in \mN}$ un $A$-module de $X^{\mN^\circ}$. Pour que 
$M$ soit $A$-plat, il faut et il suffit que pour tout entier $n\geq 0$, $M_n$ soit $A_n$-plat.
\end{lem}

En effet, la condition est nécessaire en vertu de \eqref{higgs3-spsa1b} et (\cite{sga4} V 1.7.1),
et elle est suffisante d'après \ref{higgs3-spsa99}(i) et \eqref{higgs3-spsa2a}.

\begin{prop}\label{higgs3-spad15}
Soient $A=(A_n)_{n\in \mN}$ un anneau de $X^{\mN^\circ}$, 
$M=(M_n)_{n\in \mN}$ un $A$-module strict de $X^{\mN^\circ}$ \eqref{higgs3-spsa2}. Pour que 
$M$ soit de type fini sur $A$, il faut et il suffit que pour tout entier $n\geq 0$, $M_n$ soit de type fini sur $A_n$.
\end{prop}

Il n'y a évidemment que la suffisance de la condition à montrer \eqref{higgs3-spsa1b}. Supposons que pour tout entier $n\geq 0$, 
le $A_n$-module $M_n$ soit type fini. Soient $n$ un entier $\geq 0$, $U\in \ob(X)$
tels que $M_n|U$ soit engendré sur $A_n|U$ par un nombre fini de sections 
$s_1,\dots,s_{\ell}\in \Gamma(U,M_n)$. Comme $M$ est strict, pour tout entier 
$0\leq i\leq n$, $M_i|U$ est engendré sur $A_i|U$ par les images canoniques des sections 
$s_1,\dots,s_{\ell}$ dans $\Gamma(U,M_i)$. En vertu de \ref{higgs3-spsa7}, on a une équivalence canonique de catégories 
\begin{equation}
h\colon (X_{/U})^{[n]^\circ}\stackrel{\sim}{\rightarrow} (X^{\mN^\circ})_{/\alpha_{n!}(U)},
\end{equation}
telle que $h^*(A|\alpha_{n!}(U))$ soit isomorphe à l'anneau $(A_i|U)_{i\in [n]}$ de $(X_{/U})^{[n]^\circ}$
et que $h^*(M|\alpha_{n!}(U))$ soit isomorphe au module $(M_i|U)_{i\in [n]}$ de $(X_{/U})^{[n]^\circ}$.
D'après \ref{higgs3-spsa55}, on a un isomorphisme canonique
\begin{equation}
\Gamma((X_{/U})^{[n]^\circ}, (M_i|U)_{i\in [n]})\stackrel{\sim}{\rightarrow}\Gamma(U,M_n). 
\end{equation}
Compte tenu de \ref{higgs3-spsa99}(i), 
le module $(M_i|U)_{i\in [n]}$ est alors engendré sur $(A_i|U)_{i\in [n]}$
par les sections $s_1,\dots,s_{\ell}\in \Gamma(U,M_n)$.

Pour tout entier $m\geq 0$, soit $(U_{m,j})_{j\in J_m}$ un raffinement de l'objet final de $X$
tel que pour tout $j\in J_m$, $M_m|U_{m,j}$ soit engendré sur $A_m|U_{m,j}$ par un nombre fini de sections 
de $\Gamma(U_{m,j},M_m)$.  
D'après ce qui précède, pour tout $m\in \mN$ et tout $j\in J_m$, le module $M|\alpha_{m!}(U_{m,j})$ est engendré sur 
$A|\alpha_{m!}(U_{m,j})$ par un nombre fini de sections de $\Gamma(U_{m,j},M_m)$.
D'autre part, la famille $(\alpha_{m!}(U_{m,j}))_{m\in \mN,j\in J_m}$ est un raffinement de 
l'objet final de $X^{\mN^\circ}$ en vertu de \ref{higgs3-spsa98}(ii). 
Par suite, $M$ est de type fini sur $A$.

\begin{defi}\label{higgs3-spad300}
Soient $n$ un entier $\geq 0$, $A=(A_i)_{i\in [n]}$ un anneau de $X^{[n]^\circ}$. 
On dit qu'un $A$-module $(M_i)_{i\in [n]}$ de $X^{[n]^\circ}$ est {\em adique} 
si pour tous entiers $i$ et $j$ tels que $0\leq i\leq j\leq n$, le morphisme 
$M_{j}\otimes_{A_j}A_i\rightarrow M_i$ déduit du morphisme de transition $M_j\rightarrow M_i$ 
est un isomorphisme. 
\end{defi}

\begin{defi}\label{higgs3-spad3}
Soit $A=(A_i)_{i\in \mN}$ un anneau de $X^{\mN^\circ}$. 
On dit qu'un $A$-module $(M_i)_{i\in \mN}$ de $X^{\mN^\circ}$ est {\em adique} 
si pour tous entiers $i$ et $j$ tels que $0\leq i\leq j$, le morphisme 
$M_{j}\otimes_{A_j}A_i\rightarrow M_i$ déduit du morphisme de transition $M_j\rightarrow M_i$ 
est un isomorphisme. 
\end{defi}

Les notions de modules adiques définies ci-dessus sont des cas particuliers de la notion de module cocartésien 
introduite dans (\cite{ag2} 7.11).

\begin{lem}\label{higgs3-spad1}
Supposons que $X$ ait suffisamment de points. 
Soient de plus, $R$ un anneau de $X$, 
$J$ un idéal de $R$. Pour tout entier $n\geq 0$, on pose $R_n=R/J^{n+1}$. On note $\bvR$
l'anneau $(R_n)_{n\in \mN}$ de $X^{\mN^\circ}$.  
Alors pour qu'un $\bvR$-module adique $(M_n)_{n\in \mN}$ de $X^{\mN^\circ}$ soit de type fini, 
il faut et il suffit que le $R_0$-module $M_0$ soit de type fini.
\end{lem}
Il n'y a évidemment que la suffisance de la condition à montrer. Supposons le $R_0$-module $M_0$ 
de type fini. Soient $\varphi\colon X\rightarrow \Ens$ un foncteur fibre, $n$ un entier $\geq 1$. Comme on a 
$\varphi(R_n)=\varphi(R)/(\varphi(J))^{n+1}$ et $\varphi(M_0)=\varphi(M_n)/(\varphi(J)\varphi(M_n))$,
$\varphi(M_n)$ est de type fini sur $\varphi(R_n)$ d'après (\cite{egr1} 1.8.5).
Par suite, $M_n$ est de type fini sur $R_n$, d'où la proposition en vertu de \ref{higgs3-spad15}.

\subsection{}\label{higgs3-spad11}
Soient $f\colon Y\rightarrow X$ un morphisme de $\mU$-topos, 
$A=(A_n)_{n\in \mN}$ un anneau de $X^{\mN^\circ}$,
$B=(B_n)_{n\in \mN}$ un anneau de $Y^{\mN^\circ}$, 
$u\colon A\rightarrow (f^{\mN^\circ})_*(B)$
un homomorphisme d'anneaux. On considère $f^{\mN^\circ}\colon Y^{\mN^\circ}\rightarrow X^{\mN^\circ}$ 
comme un morphisme de topos annelés respectivement, par $A$ et $B$. 
Nous utilisons pour les modules la notation $(f^{\mN^\circ})^{-1}$ pour désigner l'image
inverse au sens des faisceaux abéliens et nous réservons la notation 
$(f^{\mN^\circ})^*$ pour l'image inverse au sens des modules.
La donnée de $u$ est équivalente à la donnée pour tout $n\in \mN$ d'un  
homomorphisme d'anneaux $u_n\colon A_n\rightarrow f_*(B_n)$ telle que ces homomorphismes soient compatibles 
aux morphismes de transition de $A$ et $B$ \eqref{higgs3-spsa1f}. L'homomorphisme 
$(f^{\mN^\circ})^{-1}(A)\rightarrow B$ adjoint de $u$ correspond au système d'homomorphismes 
$f^*(A_n)\rightarrow B_n$ adjoints des $u_n$ ($n\in \mN$) \eqref{higgs3-spsa1g}. Pour tout $n\in \mN$, 
on désigne par 
\begin{equation}
f_n\colon (Y,B_n)\rightarrow (X,A_n)
\end{equation}
le morphisme de topos annelés défini par $f$ et $u_n$. 
D'après \eqref{higgs3-spsa1g} et \eqref{higgs3-spsa2a}, pour tout $A$-module $M=(M_n)_{n\in \mN}$ de $X^{\mN^\circ}$, 
on a un isomorphisme canonique fonctoriel
\begin{equation}\label{higgs3-spad11b}
(f^{\mN^\circ})^*(M)\stackrel{\sim}{\rightarrow} (f_n^*(M_n)).
\end{equation}
Par suite, si $M$ est adique, il en est de même de $(f^{\mN^\circ})^*(M)$.

\begin{lem}\label{higgs3-spad160}
Soient $n$ un entier $\geq 0$, $A=(A_i)_{i\in [n]}$ un anneau de $X^{[n]^\circ}$, 
$(M_i)_{i\in [n]}$ un $A$-module de $X^{[n]^\circ}$. Pour que le $A$-module 
$M$ soit localement projectif de type fini \eqref{higgs3-not3}, il faut et il suffit que $M$ soit adique et que 
le $A_n$-module $M_n$ soit localement projectif de type fini.
\end{lem}

Supposons d'abord $M$ localement projectif de type fini sur $A$. 
Le $A_n$-module $M_n$ est alors localement projectif de type fini \eqref{higgs3-spsa1b}. 
Montrons que $M$ est adique. D'après \ref{higgs3-spsa98}(ii), il existe un raffinement $(U_j)_{j\in J}$ de l'objet final de $X$
tel que pour tout $j\in J$, $M|\alpha_{n!}(U_j)$ soit un facteur direct d'un $A|\alpha_{n!}(U_j)$-module libre de type fini.
On a $\lambda^*(U_j)=\alpha_{n!}(U_j)$ \eqref{higgs3-spsa1ca}. Compte tenu de \ref{higgs3-spsa100}(ii), 
on peut alors se borner au cas où $M$ est un facteur direct d'un $A$-module libre de type fini.
Il existe donc un entier $d\geq 1$, un $A$-module $N=(N_i)_{i\in [n]}$ et un isomorphisme $A$-linéaire 
$A^d\stackrel{\sim}{\rightarrow} M\oplus N$.
On en déduit que pour tous entiers $i$ et $j$ tels que $0\leq i\leq j\leq n$, les morphismes canoniques
$M_{j}\otimes_{A_j}A_i\rightarrow M_i$ et $N_{j}\otimes_{A_j}A_i\rightarrow N_i$ sont des isomorphismes~;
autrement dit, $M$ et $N$ sont adiques. 

Supposons ensuite que $M$ soit adique et que le $A_n$-module $M_n$ soit localement projectif de type fini.
Montrons que le $A$-module $M$ est localement projectif de type fini. 
Compte tenu de \ref{higgs3-spsa100}(ii), on peut se borner au cas où $M_n$ est un facteur direct d'un 
$A_n$-module libre de type fini. Il existe donc un entier $d\geq 1$, un $A_n$-module $N_n$ et 
un isomorphisme $A_n$-linéaire $A_n^d\stackrel{\sim}{\rightarrow} M_n\oplus N_n$.
Celui-ci induit un isomorphisme $A$-linéaire 
\begin{equation}
A^d\stackrel{\sim}{\rightarrow} M\oplus (N_n\otimes_{A_n}A_i)_{i\in [n]}.
\end{equation}

\begin{prop}\label{higgs3-spad16}
Soient $A=(A_n)_{n\in \mN}$ un anneau de $X^{\mN^\circ}$, 
$M=(M_n)_{n\in \mN}$ un $A$-module de $X^{\mN^\circ}$. Pour que le $A$-module 
$M$ soit localement projectif de type fini \eqref{higgs3-not3}, il faut et il suffit que $M$ soit adique 
et que pour tout entier $n\geq 0$, le $A_n$-module $M_n$ soit localement projectif de type fini.
\end{prop}

Supposons d'abord $M$ localement projectif de type fini sur $A$. Pour tout entier $n\geq 0$, 
le $(A_i)_{i\in [n]}$-module $(M_i)_{i\in [n]}$ est localement projectif de type fini d'après \ref{higgs3-spsa6}(ii). 
On en déduit que $M$ est adique et que pour tout entier $n\geq 0$, 
le $A_n$-module $M_n$ est localement projectif de type fini en vertu de \ref{higgs3-spad160}.
Supposons ensuite que $M$ soit adique et que pour tout entier $n\geq 0$, 
le $A_n$-module $M_n$ soit localement projectif de type fini. Soient $n$ un entier $\geq 0$, $U\in \ob(X)$. 
D'après \ref{higgs3-spsa7}, on a une équivalence canonique de catégories 
\begin{equation}
h\colon (X_{/U})^{[n]^\circ}\stackrel{\sim}{\rightarrow} (X^{\mN^\circ})_{/\alpha_{n!}(U)},
\end{equation}
telle que $h^*(A|\alpha_{n!}(U))$ soit isomorphe à l'anneau $(A_i|U)_{i\in [n]}$ de $(X_{/U})^{[n]^\circ}$
et que $h^*(M|\alpha_{n!}(U))$ soit isomorphe au module $(M_i|U)_{i\in [n]}$ de $(X_{/U})^{[n]^\circ}$.
Par suite, le $A|\alpha_{n!}(U)$-module $M|\alpha_{n!}(U)$ est localement projectif de type fini en vertu de \ref{higgs3-spad160}. 
On en déduit, compte tenu de \ref{higgs3-spsa98}(ii), que le $A$-module $M$ est localement projectif de type fini.

\subsection{}\label{higgs3-spad12}
Soient $Y$ un schéma connexe, $\oy$ un point géométrique de $Y$, $Y_\fet$ le topos fini étale de $Y$ \eqref{higgs3-not2}, 
$\bB_{\pi_1(Y,\oy)}$ le topos classifiant du groupe profini $\pi_1(Y,\oy)$, 
\begin{equation}\label{higgs3-spad12a}
\nu_\oy\colon Y_\fet\rightarrow \bB_{\pi_1(Y,\oy)}
\end{equation}
le foncteur fibre de $Y_\fet$ en $\oy$ \eqref{higgs3-not6}. Soit $R$ un anneau de $Y_\fet$. 
Posons $R_\oy=\nu_\oy(R)$, qui est un anneau muni de la topologie discrète et d'une action 
continue de $\pi_1(Y,\oy)$ par des homomorphismes d'anneaux. On désigne par $\hR_\oy$ le séparé complété 
$p$-adique de $R_\oy$ que l'on munit de la topologie $p$-adique et de l'action de $\pi_1(Y,\oy)$ 
induite par celle sur $R_\oy$. 
D'après (\cite{ac} III § 2.11 prop.~14 et cor.~1; cf. aussi \cite{egr1} 1.8.7), pour tout  entier $n\geq 1$, on a 
\begin{equation}\label{higgs3-spad12b}
\hR_\oy/p^n\hR_\oy\simeq R_\oy/p^nR_\oy.
\end{equation}
L'action de $\pi_1(Y,\oy)$ sur $\hR_\oy$ est donc continue.

Pour tout anneau topologique $A$
muni d'une action continue de $\pi_1(Y,\oy)$ par des homomorphismes d'anneaux, 
on désigne par $\bRep_A^{\cont}(\pi_1(Y,\oy))$  
la catégorie des $A$-représentations continues de $\pi_1(Y,\oy)$  (\cite{ag1} 3.1). 
Considérons les catégories suivantes~:
\begin{itemize}
\item[(a)] $\bRep_{R_\oy}^{\disc}(\pi_1(Y,\oy))$ la sous-catégorie pleine de  
$\bRep_{R_\oy}^{\cont}(\pi_1(Y,\oy))$ formée des $R_\oy$-représenta\-tions continues de $\pi_1(Y,\oy)$ 
pour lesquelles la topologie est discrète. 

\item[(b)] $\bpad(\bRep_{R_\oy}^{\disc}(\pi_1(Y,\oy)))$ la catégorie des 
systèmes projectifs $p$-adiques de $\bRep_{R_\oy}^{\disc}(\pi_1(Y,\oy))$~: 
un système projectif $(M_n)_{n\in \mN}$ de $\bRep_{R_\oy}^{\disc}(\pi_1(Y,\oy))$ est {\em $p$-adique} 
si pour tout entier $n\geq 0$, $p^{n+1}M_n=0$ et pour tous entiers $m\geq n\geq 0$, le morphisme 
$M_m/p^{n+1}M_m\rightarrow M_n$ induit par le morphisme de transition $M_m\rightarrow M_n$ 
est un isomorphisme (\cite{sga5} V 3.1.1).

\item[(c)] $\bpad_\tf(\bRep_{R_\oy}^{\disc}(\pi_1(Y,\oy)))$ la sous-catégorie pleine de 
$\bpad(\bRep_{R_\oy}^{\disc}(\pi_1(Y,\oy)))$ formée
des systèmes projectifs $p$-adiques de type fini~: on dit qu'un système projectif $p$-adique $(M_n)_{n\in \mN}$ de 
$\bRep_{R_\oy}^{\disc}(\pi_1(Y,\oy))$ est {\em de type fini} si pour tout entier $n\geq 0$, 
$M_n$ est un $R_\oy$-module de type fini, ou ce qui revient au même, si $M_0$ est un $R_\oy$-module de type fini  
(\cite{egr1} 1.8.5).

\item[(d)] $\bRep_{\hR_\oy}^{p\aatf}(\pi_1(Y,\oy))$ la sous-catégorie pleine de 
$\bRep_{\hR_\oy}^{\cont}(\pi_1(Y,\oy))$ formée des $\hR_\oy$-représenta\-tions 
$p$-adiques de type fini de $\pi_1(Y,\oy)$~:
une $\hR_\oy$-représentation de $\pi_1(Y,\oy)$ est {\em $p$-adique de type fini} 
si elle est continue pour la topologie 
$p$-adique et si le $\hR_\oy$-module sous-jacent est séparé de type fini. 
On notera que tout $\hR_\oy$-module de type fini est complet pour la topologie $p$-adique 
(\cite{ac} chap. III § 2.11 cor.~1 de prop.~16). 
\end{itemize}

La limite projective définit un foncteur 
\begin{equation}\label{higgs3-spad12c}
\bpad(\bRep_{R_\oy}^{\disc}(\pi_1(Y,\oy)))\rightarrow\bRep_{\hR_\oy}^\cont(\pi_1(Y,\oy)),
\end{equation}
qui induit une équivalence de catégories
\begin{equation}\label{higgs3-spad12d}
\bpad_\tf(\bRep_{R_\oy}^{\disc}(\pi_1(Y,\oy)))\stackrel{\sim}{\rightarrow} \bRep_{\hR_\oy}^{p\aatf}(\pi_1(Y,\oy)).
\end{equation}
En effet, pour tout objet $(M_n)$ 
de $\bpad_\tf(\bRep_{R_\oy}^{\disc}(\pi_1(Y,\oy)))$, le $\hR_\oy$-module 
\[
\hM=\underset{\underset{n\in \mN}{\longleftarrow}}{\lim}\ M_n 
\]
est de type fini et pour tout $n\in \mN$, on a $\hM/p^{n+1}\hM\simeq M_n$ d'après (\cite{ac} III § 2.11 prop.~14 et cor.~1). 

En restreignant le foncteur $\nu_\oy$ \eqref{higgs3-spad12a} aux $R$-modules, on obtient une équivalence 
de catégories que l'on note encore
\begin{equation}\label{higgs3-spad12e}
\nu_\oy\colon \bMod(R)\stackrel{\sim}{\rightarrow} \bRep_{R_\oy}^{\disc}(\pi_1(Y,\oy)). 
\end{equation}
On note $\bvR$ l'anneau $(R/p^{n+1}R)_{n\in \mN}$ de $Y_\fet^{\mN^\circ}$ et $\bMod^\ad(\bvR)$
(resp. $\bMod^\atf(\bvR)$) la catégorie des $\bvR$-modules adiques (resp. adiques de type fini) de $Y_\fet^{\mN^\circ}$. 
Pour qu'un $\bvR$-module $(M_n)_{n\in \mN}$ de $Y_\fet^{\mN^\circ}$ soit adique (resp. adique de type fini), 
il faut et il suffit que le système projectif $(\nu_\oy(M_n))$ de $\bRep_{R_\oy}^{\disc}(\pi_1(Y,\oy))$
soit $p$-adique (resp. $p$-adique de type fini en vertu de \ref{higgs3-spad15}).
Les foncteurs \eqref{higgs3-spad12e} et \eqref{higgs3-spad12c} induisent donc un foncteur 
\begin{equation}\label{higgs3-spad12f}
\bMod^\ad(\bvR)\rightarrow\bRep_{\hR_\oy}^\cont(\pi_1(Y,\oy))
\end{equation}
et une équivalence de catégories 
\begin{equation}\label{higgs3-spad12g}
\bMod^\atf(\bvR)\stackrel{\sim}{\rightarrow}\bRep_{\hR_\oy}^{p\aatf}(\pi_1(Y,\oy)).
\end{equation}

\section{Topos annelé de Faltings}\label{higgs3-tfa}

\subsection{}\label{higgs3-tfa1}\index{1000801@$X$, $\oX$, $Y$, $h$, $\hbar$, $j$}\index{1000802@$\oU$, $U_Y$ ($U$ un $X$-schéma)}
\index{1000805@$\pi\colon E\rightarrow \Et_{/X}$}\index{1000806@$(V\rightarrow U)$}\index{1000807@$\fF$, $\fF^\vee$, $\cP^\vee$}
\index{Site fibre de Faltings@Site fibré de Faltings}
Dans cette section, on se donne un diagramme commutatif de morphismes de schémas 
\begin{equation}\label{higgs3-tfa1aa}
\xymatrix{Y\ar[r]^j\ar[rd]_-(0.4)h&\oX\ar[d]^{\hbar}\\
&X}
\end{equation}
tel que $\oX$ soit normal et localement irréductible \eqref{higgs3-sli1} et que $j$ soit une immersion ouverte quasi-compacte. 
On notera que $\oX$ et par suite $Y$ sont étale-localement connexes d'après \ref{higgs3-sli2}(iii). 
Pour tout $X$-schéma $U$, on pose 
\begin{equation}\label{higgs3-tfa1ab}
\oU=U\times_X\oX \ \ \ {\rm et}\ \ \ U_Y=U\times_XY.
\end{equation} 

On désigne par 
\begin{equation}\label{higgs3-tfa1a}
\pi\colon E\rightarrow \Et_{/X}
\end{equation}
le $\mU$-site fibré de Faltings associé au morphisme $h$ \eqref{higgs3-not2} (\cite{ag2} 10.1). 
On rappelle que les objets de $E$ sont les morphismes de schémas $V\rightarrow U$ au-dessus de $h$ 
tels que le morphisme
$U\rightarrow X$ soit étale et que le morphisme $V\rightarrow U_Y$ soit étale fini. 
Soient $(V'\rightarrow U')$ et $(V\rightarrow U)$ deux objets de $E$. Un morphisme 
de $(V'\rightarrow U')$ dans $(V\rightarrow U)$ est la donnée d'un $X$-morphisme $U'\rightarrow U$ et 
d'un $Y$-morphisme $V'\rightarrow V$ tels que le diagramme
\begin{equation}\label{higgs3-tfa1b}
\xymatrix{
V'\ar[r]\ar[d]&U'\ar[d]\\
V\ar[r]&U}
\end{equation}
soit commutatif. Le foncteur $\pi$ est alors défini pour tout $(V\rightarrow U)\in \ob(E)$, par
\begin{equation}\label{higgs3-tfa1c}
\pi(V\rightarrow U)=U.
\end{equation}
Pour tout $U\in \ob(\Et_{/X})$, la fibre de $E$ au-dessus de $U$ s'identifie canoniquement au 
site fini étale de $U_Y$ \eqref{higgs3-not2}. On note 
\begin{equation}\label{higgs3-tfa1d}
\alpha_{U!}\colon \Et_{\rf/U_Y}\rightarrow E, \ \ \ V\mapsto (V\rightarrow U)
\end{equation}
le foncteur canonique (\cite{ag2} (5.1.2)). 

On désigne par 
\begin{equation}\label{higgs3-tfa1e}
\fF\rightarrow  \Et_{/X}
\end{equation}
le $\mU$-topos fibré associé à $\pi$. La catégorie fibre de $\fF$ au-dessus de tout $U\in \ob(\Et_{/X})$ est 
canoniquement équivalente au topos fini étale $(U_Y)_\fet$ de $U_Y$ \eqref{higgs3-not2}
et le foncteur image inverse 
pour tout morphisme $f\colon U'\rightarrow U$ de $\Et_{/X}$ s'identifie au foncteur 
$(f_Y)_{\fet}^*\colon (U_Y)_\fet\rightarrow (U'_Y)_\fet$ image inverse par le morphisme 
de topos $(f_Y)_\fet\colon (U'_Y)_\fet\rightarrow (U_Y)_\fet$ (\cite{ag2} 9.3). On désigne par 
\begin{equation}\label{higgs3-tfa1f}
\fF^\vee\rightarrow (\Et_{/X})^\circ
\end{equation}
la catégorie fibrée obtenue en associant à tout $U\in \ob(\Et_{/X})$ la catégorie $(U_Y)_\fet$, et à tout morphisme 
$f\colon U'\rightarrow U$ de $\Et_{/X}$ le foncteur 
$(f_Y)_{\fet*}\colon (U'_Y)_\fet\rightarrow (U_Y)_\fet$ 
image directe par le morphisme de topos $(f_Y)_\fet$. On désigne par
\begin{equation}\label{higgs3-tfa1g}
\cP^\vee\rightarrow (\Et_{/X})^\circ
\end{equation}
la catégorie fibrée obtenue en associant à tout $U\in \ob(\Et_{/X})$ la catégorie $(\Et_{\rf/U_Y})^\wedge$ 
des préfaisceaux de $\mU$-ensembles
sur $\Et_{\rf/U_Y}$, et à tout morphisme $f\colon U'\rightarrow U$ de $\Et_{/X}$ 
le foncteur 
\begin{equation}\label{higgs3-tfa1h}
(f_Y)_{\fet*}\colon (\Et_{\rf/U'_Y})^\wedge\rightarrow (\Et_{\rf/U_Y})^\wedge
\end{equation} 
obtenu en composant 
avec le foncteur image inverse $f^+_Y\colon \Et_{\rf/U_Y}\rightarrow \Et_{\rf/U'_Y}$.

\subsection{}\label{higgs3-tfa2}\index{1000810@$\hE$}\index{1000811@$\tE$}
On désigne par $\hE$ la catégorie des préfaisceaux de $\mU$-ensembles sur $E$. On a 
alors une équivalence de catégories (\cite{ag2} 5.2)
\begin{eqnarray}\label{higgs3-tfa2a}
\hE&\rightarrow& \bHom_{(\Et_{/X})^\circ}((\Et_{/X})^\circ,\cP^\vee)\\
F&\mapsto &\{U\mapsto F\circ \alpha_{U!}\}.\nonumber
\end{eqnarray}
On identifiera dans la suite $F$ à la section $\{U\mapsto F\circ \alpha_{U!}\}$ qui lui est associée par cette équivalence.

On munit $E$ de la topologie co-évanescente définie par $\pi$ (\cite{ag2} 5.3) 
et on note $\tE$ le topos des faisceaux de $\mU$-ensembles sur $E$. 
Les site et topos ainsi définis sont appelés {\em site et topos de Faltings}  
associés à $h$ (\cite{ag2} 10.1). Si $F$ est un préfaisceau sur $E$, on note $F^a$ le faisceau associé. 
D'après (\cite{ag2} 5.11), le foncteur \eqref{higgs3-tfa2a} induit un foncteur pleinement fidèle 
\begin{equation}\label{higgs3-tfa2b}
\tE\rightarrow \bHom_{(\Et_{/X})^\circ}((\Et_{/X})^\circ,\fF^\vee)
\end{equation}
d'image essentielle les sections $\{U\mapsto F_U\}$ vérifiant une condition de recollement.

\subsection{}\label{higgs3-tfa3}\index{1000815@$\sigma$}\index{1000816@$\beta$}\index{1000817@$\Psi$}
Le foncteur $\alpha_{X!}\colon \Et_{\rf/Y}\rightarrow E$ 
est continu et exact à gauche (\cite{ag2} 5.32). Il définit donc un morphisme de topos (\cite{ag2} (10.6.3))
\begin{equation}\label{higgs3-tfa3a}
\beta\colon \tE\rightarrow Y_\fet.
\end{equation}
Le foncteur 
\begin{equation}\label{higgs3-tfa3b}
\sigma^+\colon \Et_{/X}\rightarrow E, \ \ \ U\mapsto (U_Y\rightarrow U)
\end{equation}
est continu et exact à gauche (\cite{ag2} 5.32). Il définit donc un morphisme de topos (\cite{ag2} (10.6.4))
\begin{equation}\label{higgs3-tfa3c}
\sigma\colon \tE \rightarrow X_\et.
\end{equation}
Par ailleurs, le foncteur 
\begin{equation}\label{higgs3-tfa3d}
\Psi^+\colon E\rightarrow \Et_{/Y}, \ \ \ (V\rightarrow U)\mapsto V
\end{equation}
est continu et exact à gauche (\cite{ag2} 10.7). Il définit donc un morphisme de topos 
\begin{equation}\label{higgs3-tfa3e}
\Psi\colon Y_\et\rightarrow \tE.
\end{equation}
On a des morphismes canoniques (\cite{ag2} (10.8.3) et (10.8.4))
\begin{eqnarray}
\sigma^*&\rightarrow& \Psi_* h_\et^*, \label{higgs3-tfa3f}\\
\beta^*&\rightarrow&\Psi_* \rho^*_Y, \label{higgs3-tfa3g}
\end{eqnarray} 
où $\rho_Y\colon Y_\et\rightarrow Y_\fet$ est le morphisme canonique \eqref{higgs3-not2a}. 
Si $X$ est quasi-séparé et si $Y$ est cohérent, \eqref{higgs3-tfa3g} est un isomorphisme en vertu de (\cite{ag2} 10.9(iii)). 

\begin{rema}\label{higgs3-tfa31}
Il résulte aussitôt de (\cite{ag2} 5.10) que $\{U\mapsto U_Y\}$ est un faisceau sur $E$. 
C'est donc un objet final de $\tE$ et on a des isomorphismes canoniques 
\begin{equation}\label{higgs3-tfa31a}
\sigma^*(X)\stackrel{\sim}{\rightarrow} \{U\mapsto U_Y\}\stackrel{\sim}{\leftarrow} \beta^*(Y).
\end{equation}
\end{rema}

\subsection{}\label{higgs3-tfa4}
Considérons un diagramme commutatif 
\begin{equation}\label{higgs3-tfa4c}
\xymatrix{
Y'\ar[r]^{h'}\ar[d]&X'\ar[d]\\
Y\ar[r]^h&X}
\end{equation}
et notons $E'$ le topos de Faltings associé au morphisme $h'$ et $\tE'$ 
le topos des faisceaux de $\mU$-ensembles sur $E'$ (\cite{ag2} 10.1). 
On a alors un foncteur continu et exact à gauche (\cite{ag2} 10.12)
\begin{equation}\label{higgs3-tfa4a}
\Phi^+\colon E\rightarrow E', \ \ \ (V\rightarrow U)\mapsto (V\times_YY'\rightarrow U\times_XX').
\end{equation}
Il définit un morphisme de topos 
\begin{equation}\label{higgs3-tfa4b}
\Phi\colon \tE'\rightarrow \tE.
\end{equation}

\subsection{}\label{higgs3-tfa41}\index{1000820@$\rho\colon X_\et\gtimes_{X_\et}Y_\et\rightarrow \tE$}\index{1000822@$\rho(\oy\rightsquigarrow \ox)$}
On désigne par $D$ le site co-évanescent associé au foncteur $h^+\colon \Et_{/X}\rightarrow \Et_{/Y}$ 
induit par $h$ (\cite{ag2} 4.1). Le topos des 
faisceaux de $\mU$-ensembles sur $D$ est le topos co-évanescent $X_\et\gtimes_{X_\et}Y_\et$ 
du morphisme $h_\et\colon Y_\et\rightarrow X_\et$ induit par $h$ (\cite{ag2} 3.12 et 4.10). 
Tout objet de $E$ est naturellement un objet de $D$. 
On définit ainsi un foncteur pleinement fidèle et exact à gauche
\begin{equation}\label{higgs3-tfa41a}
\rho^+\colon E\rightarrow D.
\end{equation} 
Celui-ci est continu et exact à gauche (\cite{ag2} 10.15). 
Il définit donc un morphisme de topos
\begin{equation}\label{higgs3-tfa41b}
\rho\colon X_\et\gtimes_{X_\et}Y_\et\rightarrow \tE.
\end{equation}

La donnée d'un point de $X_\et\gtimes_{X_\et}Y_\et$ 
est équivalente à la donnée d'une paire de points géométriques $\ox$ de $X$ et $\oy$ de $Y$
et d'une flèche de spécialisation $u$ de $h(\oy)$ vers $\ox$, c'est-à-dire, d'un $X$-morphisme 
$u\colon \oy\rightarrow X_{(\ox)}$, où $X_{(\ox)}$ désigne le localisé strict de $X$ en $\ox$ (\cite{ag2} 10.18). 
Un tel point sera noté $(\oy\rightsquigarrow \ox)$ ou encore $(u\colon \oy\rightsquigarrow \ox)$. 
On désigne par $\rho(\oy\rightsquigarrow \ox)$ son image par  $\rho$, qui est donc un point de $\tE$.  

Si $X$ et $Y$ sont cohérents, lorsque $(\oy \rightsquigarrow \ox)$ décrit la famille des points de $X_\et\gtimes_{X_\et}Y_\et$, 
la famille des foncteurs fibres de $\tE$ associés aux points $\rho(\oy \rightsquigarrow \ox)$ est conservative en vertu de (\cite{ag2} 10.21).

\subsection{}\label{higgs3-tfa5}\index{1000825@$\theta\colon Y_\fet\rightarrow \tE$}
Supposons $X$ strictement local, de point fermé $x$. 
On désigne par $E_{\scoh}$ la sous-catégorie pleine de $E$ formée des objets $(V\rightarrow U)$ tels que 
le morphisme $U\rightarrow X$ soit séparé et cohérent, que l'on munit de la topologie induite par celle de $E$. 
Le foncteur d'injection canonique $E_\scoh\rightarrow E$ induit alors par restriction 
une équivalence de catégories entre $\tE$ et le topos des faisceaux de $\mU$-ensembles sur $E_\scoh$ (\cite{ag2} 10.4). 

Pour tout morphisme étale, séparé et cohérent $U\rightarrow X$, 
on désigne par $U^\rf$ la somme disjointe des localisés stricts de $U$ en les points de $U_x$;
c'est un sous-schéma ouvert et fermé de $U$, qui est fini sur $X$ (\cite{ega4} 18.5.11). 
D'après (\cite{ag2} 10.23), le foncteur
\begin{equation}\label{higgs3-tfa5a}
\theta^+\colon E_{\scoh}\rightarrow \Et_{\rf/Y}, \ \ \ (V\rightarrow U)\mapsto V\times_UU^\rf
\end{equation}
est continu et exact à gauche. Il définit donc  un morphisme de topos
\begin{equation}\label{higgs3-tfa5b}
\theta\colon Y_\fet\rightarrow \tE.
\end{equation}
On a un isomorphisme canonique (\cite{ag2} (10.24.3))
\begin{equation}\label{higgs3-tfa5c}
\beta\theta\stackrel{\sim}{\rightarrow} \id_{Y_\fet}.
\end{equation}
On en déduit un morphisme de changement de base (\cite{ag2} (10.24.4))
\begin{equation}\label{higgs3-tfa5d}
\beta_*\rightarrow \theta^*.
\end{equation}
Celui-ci est un isomorphisme en vertu de (\cite{ag2} 10.27)~; en particulier, le foncteur $\beta_*$ est exact.

\subsection{}\label{higgs3-tfa6}
Soient $\ox$ un point géométrique de $X$, $X'$ le localisé strict de $X$ en $\ox$, $Y'=Y\times_XX'$, 
$h'\colon Y'\rightarrow X'$ la projection canonique. On désigne par $E'$ le site de Faltings associé au morphisme $h'$, 
par $\tE'$ le topos des faisceaux de $\mU$-ensembles sur $E'$, par  
\begin{equation}\label{higgs3-tfa6a}
\beta'\colon \tE'\rightarrow Y'_\fet
\end{equation}
le morphisme canonique \eqref{higgs3-tfa3a}, par 
\begin{equation}\label{higgs3-tfa6b}
\Phi\colon \tE'\rightarrow \tE
\end{equation}
le morphisme de fonctorialité  \eqref{higgs3-tfa4b} et par 
\begin{equation}\label{higgs3-tfa6c}
\theta\colon Y'_\fet\rightarrow \tE'
\end{equation}
le morphisme \eqref{higgs3-tfa5b}. On note 
\begin{equation}\label{higgs3-tfa6d}
\varphi_\ox\colon \tE\rightarrow Y'_\fet
\end{equation}
le foncteur composé $\theta^*\circ\Phi^*$. 

Si $X$ et $Y$ sont cohérents, la famille des foncteurs $\varphi_\ox$,
lorsque $\ox$ décrit l'ensemble des points géométriques de $X$, est conservative en vertu de (\cite{ag2} 10.32).

On note $\fV_\ox$ la catégorie des $X$-schémas étales $\ox$-pointés (\cite{sga4} VIII 3.9), 
ou ce qui revient au même, la catégorie des voisinages du point de $X_\et$ associé à $\ox$ 
dans le site $\Et_{/X}$ (\cite{sga4} IV 6.8.2). 
Pour tout objet $(U,\fp\colon \ox\rightarrow U)$ de $\fV_\ox$, on désigne encore
par $\fp\colon X'\rightarrow U$ le morphisme déduit de $\fp$ (\cite{sga4} VIII 7.3) et par
\begin{equation}\label{higgs3-tfa6e}
\fp_Y\colon Y' \rightarrow U_Y
\end{equation}
son changement de base par $h$. D'après (\cite{ag2} 10.37), pour tout faisceau $F=\{U\mapsto F_U\}$ de $\tE$, on a un isomorphisme canonique fonctoriel 
\begin{equation}\label{higgs3-tfa6f}
\underset{\underset{(U,\fp)\in \fV_\ox^\circ}{\longrightarrow}}{\lim}\ (\fp_Y)^*_\fet(F_U)
\stackrel{\sim}{\rightarrow}\varphi_\ox(F).
\end{equation}

\subsection{}\label{higgs3-tfa7}\index{1000830@$\cB$}
On désigne par $\cB$ le préfaisceau sur $E$ défini pour tout $(V\rightarrow U)\in \ob(E)$, par 
\begin{equation}\label{higgs3-tfa7a}
\cB((V\rightarrow U))=\Gamma(\oU,\co_{\oU})
\end{equation} 
et par $\cB^a$ le faisceau associé. D'après (\cite{ag2} 5.34(ii)) et avec les conventions de notation de \ref{higgs3-not2}, 
on a un isomorphisme canonique 
\begin{equation}\label{higgs3-tfa67b}
\sigma^*(\hbar_*(\co_\oX))\stackrel{\sim}{\rightarrow}\cB^a.
\end{equation}

\subsection{}\label{higgs3-tfa8}\index{1000835@$\oU^V$}\index{1000836@$\ocB$, $\ocB_{U}$}
Pour tout $(V\rightarrow U)\in \ob(E)$, on note $\oU^V$ la fermeture intégrale de $\oU=U\times_X\oX$ dans $V$. 
Pour tout morphisme $(V'\rightarrow U')\rightarrow (V\rightarrow U)$ de $E$, on a un morphisme canonique 
$\oU'^{V'}\rightarrow \oU^V$ qui s'insère dans un diagramme commutatif 
\begin{equation}\label{higgs3-tfa8a}
\xymatrix{
V'\ar[r]\ar[d]&{\oU'^{V'}}\ar[d]\ar[r]&\oU'\ar[r]\ar[d]&U'\ar[d]\\
V\ar[r]&{\oU^V}\ar[r]&\oU\ar[r]&U}
\end{equation} 
On désigne par $\ocB$ le préfaisceau sur $E$ défini pour tout $(V\rightarrow U)\in \ob(E)$, par 
\begin{equation}\label{higgs3-tfa8b}
\ocB((V\rightarrow U))=\Gamma(\oU^V,\co_{\oU^V}).
\end{equation}
Pour tout $U\in \ob(\Et_{/X})$, on pose \eqref{higgs3-tfa1d}
\begin{equation}\label{higgs3-tfa8d}
\ocB_{U}=\ocB\circ \alpha_{U!}.
\end{equation}

\begin{remas}\label{higgs3-tfa9}
Soit $(V\rightarrow U)$ un objet de $E$. Alors~:
\begin{itemize}
\item[{\rm (i)}] Comme $V$ est entier sur $U_Y$, 
le morphisme canonique $V\rightarrow U_Y\times_{\oU}\oU^V$ est un isomorphisme.
En particulier, le morphisme canonique $V\rightarrow \oU^V$ est une immersion ouverte 
schématiquement dominante.
\item[{\rm (ii)}] Le schéma $\oU^V$ est normal et localement irréductible \eqref{higgs3-sli1}. 
En effet, $\oU$ et $V$ sont normaux et localement irréductibles d'après \ref{higgs3-sli3}.
Soit $U_0$ un ouvert de $\oU$ n'ayant qu'un nombre fini de composantes irréductibles. 
Alors $U_0\times_{\oU}\oU^V$ est la somme finie des fermetures intégrales de $U_0$ dans 
les points génériques de $V$ qui sont au-dessus de $U_0$, dont chacune est un schéma intègre et normal 
en vertu de (\cite{ega2} 6.3.7). 
\item[{\rm (iii)}] Pour tout $U$-schéma étale $U'$, posant $V'=V\times_UU'$, le morphisme canonique \eqref{higgs3-tfa8a}
\begin{equation}\label{higgs3-tfa9a}
\oU'^{V'}\rightarrow\oU^V\times_UU'
\end{equation} 
est un isomorphisme. En effet, $\oU^V\times_UU'$ est normal et localement irréductible d'après (ii) et \ref{higgs3-sli3},
et le morphisme canonique $V'\rightarrow \oU^V\times_UU'$ est une immersion ouverte 
schématiquement dominante en vertu de (i) et (\cite{ega4} 11.10.5). L'assertion s'ensuit car $\oU'^{V'}$
s'identifie à la fermeture intégrale de $\oU^V\times_UU'$ dans $V'$. 
\end{itemize}
\end{remas}

\begin{lem}\label{higgs3-tfa10}
Soient $(V\rightarrow U)$ un objet de $E$, $f\colon W\rightarrow V$ un torseur 
pour la topologie étale de $V$ sous un groupe constant fini $G$. On désigne par $\of\colon \oU^W\rightarrow \oU^V$
le morphisme induit par $f$. Alors l'action naturelle de $G$ sur $\oU^W$ est admissible {\rm (\cite{sga1} V Déf.~1.7)}
et $(\oU^V,\of)$ est un schéma quotient de $\oU^W$ par $G$, autrement dit, le morphisme canonique 
\begin{equation}\label{higgs3-tfa10a}
\co_{\oU^V}\rightarrow\of_*(\co_{\oU^W})^G
\end{equation}
est un isomorphisme. En particulier, le morphisme canonique 
\begin{equation}\label{higgs3-tfa10b}
\Gamma(\oU^V,\co_{\oU^V})\rightarrow\Gamma(\oU^W,\co_{\oU^W})^G
\end{equation}
est un isomorphisme. 
\end{lem}

En effet, l'action de $G$ sur $\oU^W$ est admissible en vertu de (\cite{sga1} V Cor.~1.8). 
Notons $Z$ le quotient de $\oU^W$ par $G$.  
D'après \ref{higgs3-tfa9}(i) et (\cite{sga1} Prop.~1.9), $f$ induit un isomorphisme 
$V\stackrel{\sim}{\rightarrow}U_Y\times_{\oU}Z$,
et le morphisme $V\rightarrow Z$ est une immersion ouverte schématiquement dominante. 
Comme $Z$ est entier sur $\oU$, on en déduit que le morphisme
$Z\rightarrow \oU^V$ induit par $\of$  est un isomorphisme, d'où l'isomorphisme \eqref{higgs3-tfa10a}.
L'isomorphisme \eqref{higgs3-tfa10b} s'en déduit puisque le foncteur $\Gamma(\oU^V,-)$ sur le topos de Zariski de $\oU^V$
est exact à gauche. 

\subsection{}\label{higgs3-tfa112}\index{1000840@$\oR^\oy_U$}
Soient $U$ un objet de $\Et_{/X}$, $\oy$ un point géométrique de $U_Y$. 
Le schéma $\oU$ étant localement irréductible \eqref{higgs3-sli3}, 
il est la somme des schémas induits sur ses composantes irréductibles. On note $\oU^\star$
la composante irréductible de $\oU$ (ou ce qui revient au même, sa composante connexe) contenant $\oy$. 
De même, $U_Y$ est la somme des schémas induits sur ses composantes irréductibles.
On pose $V=\oU^\star\times_{\oX}Y$ qui est la composante irréductible de $U_Y$ contenant $\oy$. 
On désigne par $\bB_{\pi_1(V,\oy)}$ le topos classifiant du groupe profini 
$\pi_1(V,\oy)$, par $(V_i)_{i\in I}$ le revêtement universel normalisé de $V$ en $\oy$ \eqref{higgs3-not6} et par
\begin{equation}\label{higgs3-tfa112a}
\nu_\oy\colon V_\fet \stackrel{\sim}{\rightarrow}\bB_{\pi_1(V,\oy)}, \ \ \ 
F\mapsto \underset{\underset{i\in I}{\longrightarrow}}{\lim}\ F(V_i)
\end{equation}
le foncteur fibre de $V_\fet$ en $\oy$ \eqref{higgs3-not6c}. 
Pour chaque $i\in I$, $(V_i\rightarrow U)$ est naturellement un objet de $E$. 
On peut donc considérer le système projectif filtrant des schémas $(\oU^{V_i})_{i\in I}$. 
On pose 
\begin{equation}\label{higgs3-tfa112c}
\oR^\oy_U=\underset{\underset{i\in I}{\longrightarrow}}{\lim}\ \Gamma(\oU^{V_i},\co_{\oU^{V_i}}),
\end{equation} 
qui est une représentation discrète continue de $\pi_1(V,\oy)$ en vertu de \ref{higgs3-tfa10},
autrement dit, c'est un objet de $\bB_{\pi_1(V,\oy)}$.

\begin{rema}\label{higgs3-tfa114}
Conservons les hypothèses de \eqref{higgs3-tfa112} et notons $\oU^\oy$
la limite projective dans la catégorie des $\oU$-schémas 
du système projectif filtrant $(\oU^{V_i})_{i\in I}$, qui existe en vertu de (\cite{ega4} 8.2.3). 
D'après (\cite{sga4} VI 5.3), si $\oU$ est cohérent, on a un isomorphisme canonique  
\begin{equation}\label{higgs3-tfa114a}
\Gamma(\oU^\oy,\co_{\oU^\oy}) \stackrel{\sim}{\rightarrow} \oR^\oy_U.
\end{equation}
\end{rema}

\begin{lem}\label{higgs3-tfa113}
Sous les hypothèses de \eqref{higgs3-tfa112}, $\ocB_U$ \eqref{higgs3-tfa8d} est un faisceau
pour la topologie étale de $\Et_{\rf/U_Y}$, et on a un isomorphisme canonique 
\begin{equation}\label{higgs3-tfa113a}
\nu_\oy(\ocB_U|V)\stackrel{\sim}{\rightarrow}\oR^\oy_U.
\end{equation}
\end{lem}

Il suffit de montrer que la restriction du préfaisceau $\ocB_U$ au site $\Et_{\rf/V}$ est un faisceau. 
L'isomorphisme \eqref{higgs3-tfa113a} résultera alors aussitôt des définitions. Notons 
\begin{equation}\label{higgs3-tfa113b}
\mu_\oy^+\colon \Et_{\rf/V}\rightarrow \bB_{\pi_1(V,\oy)}
\end{equation}
le foncteur fibre en $\oy$ \eqref{higgs3-not6b}. 
Pour tout $W\in \ob(\Et_{\rf/V})$, on a un isomorphisme fonctoriel 
\begin{equation}\label{higgs3-tfa113c}
\mu_\oy^+(W) \stackrel{\sim}{\rightarrow} \underset{\underset{i\in I}{\longrightarrow}}{\lim}\ \Hom_{V}(V_i,W).
\end{equation}
On en déduit un morphisme fonctoriel en $W$
\begin{equation}\label{higgs3-tfa113d}
\Gamma(\oU^W,\co_{\oU^W})\rightarrow  
\Hom_{\bB_{\pi_1(V,\oy)}}(\mu^+_\oy(W),\oR^\oy_U).
\end{equation}
Celui-ci est un isomorphisme en vertu de \ref{higgs3-tfa10}, d'où l'assertion recherchée.

\begin{prop}\label{higgs3-tfa11}
Le préfaisceau $\ocB$ sur $E$ est un faisceau pour la topologie co-évanescente.
\end{prop}

Soient $(V\rightarrow U)\in \ob(E)$, 
$(U_i\rightarrow U)_{i\in I}$ un recouvrement de $\Et_{/X}$. Pour tout $(i,j)\in I^2$, posons $V_i=V\times_UU_i$,
$U_{ij}=U_i\times_UU_j$ et $V_{ij}=U_{ij}\times_UV$. On a 
$\oU_i^{V_i}\simeq \oU^V\times_UU_i$ et $\oU_{ij}^{V_{ij}}\simeq \oU^V\times_UU_{ij}$ d'après \ref{higgs3-tfa9}(iii).
Considérant $\co_{\oU^V}$ comme un faisceau du topos étale de $\oU^V$, le
recouvrement étale $(\oU_i^{V_i}\rightarrow \oU^V)_{i\in I}$ induit alors une suite exacte d'applications d'ensembles
\begin{equation}
\ocB((V\rightarrow U))\rightarrow \prod_{i\in I}\ocB((V_i\rightarrow U_i))\rightrightarrows 
\prod_{(i,j)\in I^2}\ocB((V_{ij}\rightarrow U_{ij})).
\end{equation}
La proposition s'ensuit compte tenu de \ref{higgs3-tfa113} et (\cite{ag2} 5.10).

\subsection{}\label{higgs3-tfa111}
On dira dans la suite que $\ocB$ est l'anneau de $\tE$ associé à $\oX$.
D'après \ref{higgs3-tfa7}, l'homomorphisme canonique $\cB\rightarrow \ocB$ induit un homomorphisme 
$\sigma^*(\hbar_*(\co_\oX))\rightarrow \ocB$  (cf. \ref{higgs3-not2}). Sauf mention explicite du contraire, 
on considère $\sigma\colon \tE\rightarrow X_\et$ \eqref{higgs3-tfa3c} 
comme un morphisme de topos annelés, respectivement par $\ocB$ et $\hbar_*(\co_\oX)$. 
Nous utilisons pour les modules la notation $\sigma^{-1}$ pour désigner l'image
inverse au sens des faisceaux abéliens et nous réservons la notation 
$\sigma^*$ pour l'image inverse au sens des modules.

\begin{rema}\label{higgs3-tfa12}
Le préfaisceau $\ocB$ n'est pas en général un faisceau pour la topologie de $E$ définie originellement par Faltings dans 
(\cite{faltings2} page 214). En effet, supposons que $\hbar$ soit entier, que $j$ ne soit pas fermée
et qu'il existe une immersion ouverte affine $Z\rightarrow X$
telle que $Y=Z\times_X\oX$, $j\colon Y\rightarrow \oX$ étant la projection canonique. 
Considérons deux schémas affines $U$ et $U'$ de $\Et_{/X}$ 
et un morphisme surjectif $U'\rightarrow U$ tel que $U'_Y\rightarrow U_Y$ soit fini
et que $\oU'\rightarrow \oU$ ne soit pas entier. 
On pourrait par exemple prendre pour $U$ un sous-schéma ouvert affine de $X$ tel que l'immersion ouverte 
$j_U\colon U_Y\rightarrow \oU$ ne soit pas fermée et pour $U'$ la somme de $U$ et $U_Z$. 
Posons $V=U'_Y$. Il est clair que $(V\rightarrow U)$ est un objet de $E$
et que $(V\rightarrow U')\rightarrow (V\rightarrow U)$ est un recouvrement pour la topologie 
de $E$ définie par Faltings dans (\cite{faltings2} page 214). Mais la suite 
\begin{equation}\label{higgs3-tfa12a}
\ocB((V\rightarrow U))\rightarrow \ocB((V\rightarrow U')) \rightrightarrows \ocB((V\rightarrow U'\times_UU'))
\end{equation}
n'est cependant pas exacte. En effet, supposons qu'elle le soit. Le diagramme 
\begin{equation}\label{higgs3-tfa12b}
V\rightarrow U'\stackrel{\Delta}{\rightarrow}U'\times_UU'\rightrightarrows U',
\end{equation}
où $\Delta$ est le morphisme diagonal et la double flèche représente les deux projections canoniques,
est commutatif. Par suite, la double flèche de \eqref{higgs3-tfa12a} est faite de deux copies du même morphisme. 
Donc l'application canonique 
\begin{equation}\label{higgs3-tfa12c}
\ocB((V\rightarrow U))\rightarrow\ocB((V\rightarrow U'))
\end{equation}
est un isomorphisme. Comme $\Gamma(\oU',\co_{\oU'})\subset  \ocB((V\rightarrow U'))$, on en déduit que 
$\Gamma(\oU',\co_{\oU'})$ est entier sur $\Gamma(\oU,\co_\oU)$, et par suite que $\oU'\rightarrow \oU$ est entier
puisque $\oU$ et $\oU'$ sont affines, ce qui contredit les hypothèses.  
\end{rema}

\begin{prop}\label{higgs3-tfa18}
Supposons $X$ et $\oX$ strictement locaux et soit de plus 
$\oy$ un point géométrique de $Y$. Notons $x$ le point fermé de $X$
et $(\oy\rightsquigarrow x)$ le point de $X_\et\gtimes_{X_\et}Y_\et$
défini par l'unique flèche de spécialisation de $h(\oy)$ dans $x$ \eqref{higgs3-tfa41}. 
Alors~:
\begin{itemize}
\item[{\rm (i)}] La fibre $\ocB_{\rho(\oy \rightsquigarrow x)}$ de $\ocB$ en le point $\rho(\oy \rightsquigarrow x)$ 
de $\tE$ \eqref{higgs3-tfa41b} est un anneau normal et strictement local. 
\item[{\rm (ii)}] On a un isomorphisme canonique 
\begin{equation}\label{higgs3-tfa18a}
(\hbar_*(\co_\oX))_{x}\stackrel{\sim}{\rightarrow}\Gamma(\oX,\co_\oX).
\end{equation}
\item[{\rm (iii)}]  L'homomorphisme 
\begin{equation}\label{higgs3-tfa18b}
(\hbar_*(\co_{\oX}))_x\rightarrow \ocB_{\rho(\oy \rightsquigarrow x)}
\end{equation} 
induit par l'homomorphisme canonique $\sigma^{-1}(\hbar_*(\co_\oX))\rightarrow \ocB$  \eqref{higgs3-tfa111} est injectif et local.
\end{itemize}
\end{prop}

(i) On notera d'abord que $Y$ est intègre. Soit $(V_i)_{i\in I}$ le revêtement universel normalisé 
de $Y$ au point $\oy$ \eqref{higgs3-not6}. D'après (\cite{ag2} (10.36.2)), on a un isomorphisme canonique 
\begin{equation}\label{higgs3-tfa18c}
\underset{\underset{i\in I}{\longrightarrow}}{\lim}\ 
\ocB((V_i\rightarrow X)) \stackrel{\sim}{\rightarrow} \ocB_{\rho(\oy \rightsquigarrow x)}.
\end{equation}
Pour chaque $i\in I$, le schéma $\oX^{V_i}$ est normal, intègre et entier sur $\oX$ \eqref{higgs3-tfa9}. 
Il est donc strictement local en vertu de \ref{higgs3-sli5}. Par ailleurs, pour tous $(i,j)\in I^2$ avec $j\geq i$,
le morphisme de transition $\oX^{V_j}\rightarrow \oX^{V_i}$ est entier et dominant. 
En particulier, l'homomorphisme de transition $\ocB((V_i\rightarrow X))\rightarrow \ocB((V_j\rightarrow X))$ est local.  
On en déduit que l'anneau $\ocB_{\rho(\oy \rightsquigarrow x)}$ est local, normal et hensélien (\cite{ega1n} 0.6.5.12(ii)
et \cite{raynaud1} I § 3 prop.~1). 
Comme l'homomorphisme $\Gamma(\oX,\co_{\oX})\rightarrow \ocB_{\rho(\oy \rightsquigarrow x)}$ est entier et donc local, 
le corps résiduel de $\ocB_{\rho(\oy \rightsquigarrow x)}$ est une extension algébrique de celui de 
$\Gamma(\oX,\co_{\oX})$.
Il est donc séparablement clos. 

(ii) Cela résulte aussitôt du fait que $X$ est strictement local. 

(iii) Rappelons qu'on a un isomorphisme canonique (\cite{ag2} (10.18.1)) 
\begin{equation}\label{higgs3-tfa18d}
(\sigma^{-1}(\hbar_*(\co_\oX)))_{\rho(\oy \rightsquigarrow x)}\stackrel{\sim}{\rightarrow} \hbar_*(\co_{\oX})_x.
\end{equation}
Compte tenu de (ii) et \eqref{higgs3-tfa18c}, la fibre de l'homomorphisme canonique 
$\sigma^{-1}(\hbar_*(\co_{\oX}))\rightarrow \ocB$ en $\rho(\oy \rightsquigarrow x)$ s'identifie à l'homomorphisme canonique
\begin{equation}\label{higgs3-tfa18e}
\Gamma(\oX,\co_{\oX})\rightarrow \ocB_{\rho(\oy \rightsquigarrow x)}=\underset{\underset{i\in I}{\longrightarrow}}{\lim}\ 
\ocB((V_i\rightarrow X)), 
\end{equation} 
qui est clairement injectif et entier et donc local.

\subsection{}\label{higgs3-tfa13}
Considérons un diagramme commutatif 
\begin{equation}\label{higgs3-tfa13a}
\xymatrix{
Y'\ar[r]^{j'}\ar[d]_{g'}&{\oX'}\ar[r]^{\hbar'}\ar[d]^{\ogg}&X'\ar[d]^g\\
Y\ar[r]^j&{\oX}\ar[r]^{\hbar}&X}
\end{equation}
et posons $h'=\hbar'\circ j'$. Supposons que $\oX'$ soit normal et localement irréductible 
et que $j'$ soit une immersion ouverte quasi-compacte. On désigne par $E'$ (resp. $\tE'$) le site (resp. topos) de 
Faltings associé au morphisme $h'$ \eqref{higgs3-tfa2}, par $\ocB'$ l'anneau de $\tE'$ associé à $\oX'$ \eqref{higgs3-tfa111} et par 
\begin{equation}\label{higgs3-tfa13b}
\Phi\colon \tE'\rightarrow \tE
\end{equation}
le morphisme de fonctorialité \eqref{higgs3-tfa4b}. 
Pour tout $(V'\rightarrow U')\in \ob(E')$, on pose $\oU'=U'\times_{X'}\oX'$ et on note
$\oU'^{V'}$ la fermeture intégrale de $\oU'$ dans $V'$, de sorte que 
\begin{equation}\label{higgs3-tfa13c}
\ocB'((V'\rightarrow U'))=\Gamma(\oU'^{V'},\co_{\oU'^{V'}}).
\end{equation}
Pour tout $(V\rightarrow U)\in \ob(E)$, posons $V'=V\times_YY'$ et $U'=U\times_XX'$, 
de sorte $(V'\rightarrow U')$ est un objet de $E'$ et que l'on a un diagramme commutatif 
\begin{equation}\label{higgs3-tfa13d}
\xymatrix{Y'\ar[d]&V'\ar@{}[ld]|{\Box}\ar[l]\ar[r]\ar[d]&\oU'\ar[d]\ar[r]\ar@{}[rd]|{\Box}&\oX'\ar[d]\\
Y&V\ar[l]\ar[r]&\oU\ar[r]&\oX}
\end{equation}
On en déduit un morphisme 
\begin{equation}\label{higgs3-tfa13e}
\oU'^{V'}\rightarrow \oU^V,
\end{equation}
et par suite un homomorphisme d'anneaux de $\tE$
\begin{equation}\label{higgs3-tfa13f}
\ocB\rightarrow \Phi_*(\ocB').
\end{equation}
Nous considérons dans la suite $\Phi$ comme un morphisme de topos annelés (respectivement par $\ocB'$ et $\ocB$). 
Nous utilisons pour les modules la notation $\Phi^{-1}$ pour désigner l'image
inverse au sens des faisceaux abéliens et nous réservons la notation 
$\Phi^*$ pour l'image inverse au sens des modules.

\begin{lem}\label{higgs3-tfa14}
Les hypothèses étant celles de \eqref{higgs3-tfa13}, supposons de plus que $g$ soit étale 
et que les deux carrés du diagramme \eqref{higgs3-tfa13a} soient cartésiens, 
de sorte que $(Y'\rightarrow X')$ est un objet de $E$. Alors~:
\begin{itemize} 
\item[{\rm (i)}] Le morphisme 
\begin{equation}\label{higgs3-tfa14a}
\Phi^{-1}(\ocB)\rightarrow \ocB'
\end{equation}
adjoint du morphisme \eqref{higgs3-tfa13f} est un isomorphisme. 
\item[{\rm (ii)}] Le morphisme de topos annelés $\Phi\colon (\tE',\ocB')\rightarrow (\tE,\ocB)$ 
s'identifie au morphisme de localisation de $(\tE,\ocB)$ en $\sigma^*(X')$.
\end{itemize}
\end{lem}

En effet, le morphisme de topos $\Phi\colon \tE'\rightarrow \tE$ s'identifie au morphisme de localisation 
de $\tE$ en $\sigma^*(X')=(Y'\rightarrow X')^a$ en vertu de (\cite{ag2} 10.14). 
Il suffit donc de montrer la première proposition. Le foncteur $\Phi^*\colon \tE\rightarrow \tE'$ s'identifie au 
foncteur de restriction par le foncteur canonique $E'\rightarrow E$. 
Pour tout $(V\rightarrow U)\in \ob(E')$, on a un isomorphisme canonique 
\begin{equation}\label{higgs3-tfa14b}
U\times_{X'}\oX'\stackrel{\sim}{\rightarrow}U\times_X\oX=\oU. 
\end{equation}
On en déduit un isomorphisme (fonctoriel en $(V\rightarrow U)$)
\begin{equation}\label{higgs3-tfa14c}
\Phi^{-1}(\ocB)((V\rightarrow U))\stackrel{\sim}{\rightarrow}\ocB'((V\rightarrow U)).
\end{equation}
Il reste à montrer que celui-ci est adjoint du morphisme \eqref{higgs3-tfa13f}.
Posons $U'=U\times_XX'$ et $V'=V\times_YY'$. Les morphismes structuraux 
$U\rightarrow X'$ et $V\rightarrow Y'$ induisent des sections $U\rightarrow U'$ et $V\rightarrow V'$ des projections
canoniques $U'\rightarrow U$ et $V'\rightarrow V$.  
On en déduit un diagramme commutatif 
\begin{equation}\label{higgs3-tfa14d}
\xymatrix{
V\ar[r]\ar[d]\ar@/_1pc/[dd]_{\id_V}&\oU\ar[d]\ar@/^1pc/[dd]^{\id_\oU}\\
V'\ar[r]\ar[d]&\oU'\ar[d]\\
V\ar[r]&\oU}
\end{equation}
et par suite des morphismes $\oU^V\rightarrow \oU'^{V'}\rightarrow \oU^V$ dont le composé est l'identité de $\oU^V$. 
Le composé 
\begin{equation}\label{higgs3-tfa14e}
\Phi^{-1}(\ocB)((V\rightarrow U))\rightarrow \Phi^{-1}(\Phi_*(\ocB'))((V\rightarrow U))
= \ocB'((V'\rightarrow U'))\rightarrow \ocB'((V\rightarrow U))
\end{equation}
où la première flèche est induite par \eqref{higgs3-tfa13f} et la dernière flèche est le morphisme d'adjonction, est donc
l'isomorphisme \eqref{higgs3-tfa14c}; d'où la première proposition.

\begin{lem}\label{higgs3-tfa17}
Les hypothèses étant celles de \eqref{higgs3-tfa13}, supposons de plus que  
$X'$ soit le localisé strict de $X$ en un point géométrique $\ox$ et 
que les deux carrés du diagramme \eqref{higgs3-tfa13a} soient cartésiens.
Soient $(V\rightarrow U)$ un objet de $E$, $U'=U\times_XX'$, $V'=V\times_YY'$. 
Alors le morphisme canonique $\oU'^{V'}\rightarrow \oU^V\times_XX'$ \eqref{higgs3-tfa13e}
est un isomorphisme. 
\end{lem}

En effet, on a un diagramme commutatif à carrés cartésiens 
\begin{equation}
\xymatrix{
V'\ar[r]\ar[d]&\oU'\ar[r]\ar[d]&U'\ar[r]\ar[d]&X'\ar[d]\\
V\ar[r]&\oU\ar[r]&U\ar[r]&X}
\end{equation}
On peut donc identifier $\oU'^{V'}$ à la fermeture intégrale de $\oU^V\times_XX'$ dans $V'$. 
Il suffit alors de montrer que le morphisme canonique $V'\rightarrow \oU^V\times_XX'$ est  une immersion ouverte 
quasi-compacte et schématiquement dominante et que $\oU^V\times_XX'$ est normal et localement irréductible.
La première assertion résulte de \ref{higgs3-tfa9}(i) et (\cite{ega4} 11.10.5). La seconde assertion étant locale, pour la prouver,
on peut se borner au cas où $X$ est affine. Considérons $X'$ comme une limite projective cofiltrante 
de voisinages étales affines $(X_i)_{\in I}$ de $\ox$ dans $X$ (cf. \cite{sga4} VIII 4.5). 
Alors $\oU^V\times_XX'$ est canoniquement isomorphe à la limite projective des schémas 
$(\oU^V\times_XX_i)_{i\in I}$ (\cite{ega4} 8.2.5). Comme chaque $\oU^V\times_XX_i$ est normal d'après 
\ref{higgs3-tfa9}(ii) et (\cite{raynaud1} VII prop.~2), $\oU^V\times_XX'$ est normal d'après (\cite{ega1n} 0.6.5.12(ii)).
D'autre part, comme $\oX'$ est normal et localement irréductible par hypothèse, il en est de même de $V'$
d'après \ref{higgs3-sli3}. On en déduit que $\oU^V\times_XX'$ est localement irréductible en vertu de \ref{higgs3-sli4}~; d'où la 
proposition.

\begin{prop}\label{higgs3-tfa16}
Les hypothèses étant celles de \eqref{higgs3-tfa13}, supposons de plus que $\hbar$ soit cohérent, que 
$X'$ soit le localisé strict de $X$ en un point géométrique $\ox$ et 
que les deux carrés du diagramme \eqref{higgs3-tfa13a} soient cartésiens.
Alors le morphisme 
\begin{equation}\label{higgs3-tfa16a}
\Phi^{-1}(\ocB)\rightarrow \ocB'
\end{equation}
adjoint du morphisme \eqref{higgs3-tfa13f} est un isomorphisme. 
\end{prop}
Choisissons un voisinage ouvert affine $X_0$ de l'image de $\ox$ dans $X$ et 
notons $I$ la catégorie des $X_0$-schémas étales $\ox$-pointés
qui sont affines au-dessus de $X_0$ (cf. \cite{sga4} VIII 3.9 et 4.5). On désigne par 
\begin{equation}\label{higgs3-tfa16b}
\cE\rightarrow I
\end{equation}
le $\mU$-site fibré défini dans (\cite{ag2} (11.2.2)): la catégorie fibre de $\cE$   
au-dessus d'un objet $U$ de $I$ est le site de Faltings associé à la projection canonique 
$h_U\colon U_Y\rightarrow U$, et le foncteur image inverse associé à un morphisme $f\colon U'\rightarrow U$ de $I$
est le foncteur $\Phi_f^+\colon \cE_U\rightarrow \cE_{U'}$ défini dans \eqref{higgs3-tfa4a}. On note
\begin{equation}\label{higgs3-tfa16c}
\cF\rightarrow I
\end{equation}
le $\mU$-topos fibré associé à $\cE/I$. La catégorie fibre de $\cF$ au-dessus d'un objet $U$ de $I$
est le topos $\tcE_U$ des faisceaux de $\mU$-ensembles sur le site co-évanescent $\cE_U$, et le foncteur image 
inverse relatif à un morphisme $f\colon U'\rightarrow U$ de $I$
est le foncteur image inverse par le morphisme de topos $\Phi_f\colon \tcE_{U'}\rightarrow \tcE_U$ défini dans \eqref{higgs3-tfa4b}.  
On note 
\begin{equation}\label{higgs3-tfa16d}
\cF^\vee\rightarrow I^\circ
\end{equation}
la catégorie fibrée obtenue en associant à tout objet $U$ de $I$ la catégorie $\cF_U=\tcE_U$ et 
à tout morphisme $f\colon U'\rightarrow U$ de $I$ le foncteur image directe par le morphisme de topos $\Phi_f$.  
On rappelle (\cite{ag2} 10.14) que pour tout $U\in \ob(I)$, $\tcE_U$ est canoniquement équivalent au topos 
$\tE_{/(U_Y\rightarrow U)^a}$, où $(U_Y\rightarrow U)^a$ désigne le faisceau 
associé à $(U_Y\rightarrow U)$. 

En vertu de (\cite{ag2} 11.3), le topos $\tE'$ est canoniquement équivalent à la limite projective du topos fibré $\cF/I$. 
Munissons $\cE$ de la topologie totale (\cite{sga4} VI 7.4.1) 
et notons $\Top(\cE)$ le topos des faisceaux de $\mU$-ensembles sur $\cE$. 
On a alors un morphisme canonique (\cite{ag2} (11.4.2))
\begin{equation}\label{higgs3-tfa16e}
\varpi\colon \tE'\rightarrow \Top(\cE)
\end{equation}
et un diagramme commutatif canonique de foncteurs (\cite{ag2} (11.4.3))
\begin{equation}\label{higgs3-tfa16f}
\xymatrix{
{\tE'}\ar[r]^-(0.5)\sim\ar[d]_{\varpi_*}&{\bHom_{\cart/I^\circ}(I^\circ, \cF^\vee)}\ar@{^(->}[d]\\
{\Top(\cE)}\ar[r]^-(0.5)\sim&{\bHom_{I^\circ}(I^\circ, \cF^\vee)}}
\end{equation}
où les flèches horizontales sont des équivalences de catégories et 
la flèche verticale de droite est l'injection canonique (cf. \cite{sga4} VI 8.2.9). 

Pour tout objet $U$ de $I$, on note $\cB_{/U}$ l'anneau de $\tcE_U$ associé à $\oU$ \eqref{higgs3-tfa111}.  
Pour tout morphisme $f\colon U'\rightarrow U$ de $I$, on a un homomorphisme 
canonique $\ocB_{/U}\rightarrow \Phi_{f*}(\ocB_{/U'})$ \eqref{higgs3-tfa13f}. Ces homomorphismes
vérifient une relation de compatibilité pour la composition des morphismes de $I$ du type (\cite{egr1} (1.1.2.2)).
Ils définissent donc un anneau de $\Top(\cE)$ que l'on note 
$\{U\mapsto \ocB_{/U}\}$ (à ne pas confondre avec l'anneau $\ocB=\{U\mapsto \ocB_U\}$ de $\tE$ \eqref{higgs3-tfa8d}). 
Pour tout $U\in \ob(I)$, on a un morphisme canonique $g_U\colon X'\rightarrow U$ qui induit un morphisme de topos
annelés \eqref{higgs3-tfa13f}
\begin{equation}
\Phi_U\colon (\tE',\ocB')\rightarrow (\tcE_U,\ocB_{/U}).
\end{equation} 
La collection des homomorphismes $\ocB_{/U}\rightarrow \Phi_{U*}(\ocB')$ 
définit un homomorphisme d'anneaux de $\Top(\cE)$
\begin{equation}\label{higgs3-tfa16g}
\{U\mapsto \ocB_{/U}\} \rightarrow \varpi_*(\ocB'). 
\end{equation}
Montrons que l'homomorphisme adjoint 
\begin{equation}\label{higgs3-tfa16h}
\varpi^*(\{U\mapsto \ocB_{/U}\}) \rightarrow \ocB'
\end{equation}
est un isomorphisme. Comme le foncteur $\varpi_*$ est pleinement fidèle \eqref{higgs3-tfa16f}, il suffit de montrer 
que l'homomorphisme induit 
\begin{equation}\label{higgs3-tfa16ha}
\varpi_*(\varpi^*(\{U\mapsto \ocB_{/U}\})) \rightarrow \varpi_*(\ocB') 
\end{equation}
est un isomorphisme. Compte tenu de (\cite{sga4} VI 8.5.3), cela revient à montrer que pour tout $U\in \ob(I)$, 
l'homomorphisme canonique de $\tcE_U$
\begin{equation}\label{higgs3-tfa16hb}
\underset{\underset{f\colon U'\rightarrow U}{\longrightarrow}}{\lim}\ \Phi_{f*}(\ocB_{/U'})\rightarrow \Phi_{U*}(\ocB'),
\end{equation}
où la limite est prise sur les morphismes $f\colon U'\rightarrow U$ de $I$, est un isomorphisme. 

Soit $(V_1\rightarrow U_1)$ un objet de $E_{/(U_Y\rightarrow U)}$ 
tel que le morphisme $U_1\rightarrow U$ soit cohérent. Le morphisme canonique
\begin{equation}\label{higgs3-tfa16hd}
\oU_1^{V_1}\times_UX'\rightarrow 
\underset{\underset{U'\in \ob(I_{/U})}{\longleftarrow}}{\lim}\ \oU_1^{V_1}\times_UU'
\end{equation}
est un isomorphisme d'après  (\cite{ega4} 8.2.5).  Comme $\hbar$ est cohérent, 
le schéma $\oU_1^{V_1}$ est cohérent. 
On en déduit par (\cite{sga4} VI 5.2), \ref{higgs3-tfa9}(iii) et \ref{higgs3-tfa17} que l'homomorphisme canonique 
\begin{equation}\label{higgs3-tfa16hc}
\underset{\underset{U'\in \ob(I_{/U})}{\longrightarrow}}{\lim}\ \ocB_{/U'}((V_1\times_UU'\rightarrow U_1\times_UU'))\rightarrow 
\ocB'((V_1\times_UX'\rightarrow U_1\times_UX'))
\end{equation}
est un isomorphisme. Par ailleurs, $h$ étant cohérent, le faisceau $(V_1\rightarrow U_1)^a$ de $\tcE_U$ 
associé à $(V_1\rightarrow U_1)$ est cohérent d'après (\cite{ag2} 10.5(i)).   
Donc en vertu de (\cite{ag2} 10.4 et 10.5(ii)) et (\cite{sga4} VI 5.3), l'isomorphisme \eqref{higgs3-tfa16hc}
montre que \eqref{higgs3-tfa16hb} est un isomorphisme. 

On déduit de l'isomorphisme \eqref{higgs3-tfa16h} et de (\cite{ag2} (11.4.4)) que l'homomorphisme canonique 
\begin{equation}\label{higgs3-tfa16i}
\underset{\underset{i\in I^\circ}{\longrightarrow}}{\lim}\ \Phi_U^{-1}(\ocB_{/U}) \rightarrow \ocB' 
\end{equation}
est un isomorphisme. 
Pour tout $U\in \ob(I)$, le morphisme canonique $U\rightarrow X$ induit un morphisme de topos annelés
\begin{equation}
\rho_U\colon (\tcE_U,\ocB_{/U}) \rightarrow  (\tE,\ocB).
\end{equation} 
Comme l'homomorphisme $\rho_U^{-1}(\ocB)\rightarrow \ocB_{/U}$ est un isomorphisme en vertu de \ref{higgs3-tfa14}(i), 
la proposition résulte de l'isomorphisme \eqref{higgs3-tfa16i}.

\section{Topos de Faltings au-dessus d'un trait}\label{higgs3-TFT}

\subsection{}\label{higgs3-TFT1}\index{1000902@$X$, $X^\circ$}\index{1000904@$U^\circ=U\times_XX^\circ$}\index{1000906@$\ocB_n$, $\ocB_{U,n}$}
Dans cette section, $S$ désigne le trait fixé dans \eqref{higgs3-not1}. 
On se donne un $S$-schéma cohérent $X$ et un sous-schéma ouvert $X^\circ$ de $X_\eta$. 
Pour tout $X$-schéma $U$, on pose
\begin{equation}\label{higgs3-TFT1a}
U^\circ=U\times_XX^\circ.
\end{equation} 
On rappelle que pour tout $S$-schéma $Y$, on a posé
$\oY=Y\times_S\oS$ et pour tout entier $n\geq 1$, $Y_n=Y\times_SS_n$ \eqref{higgs3-not1c}.
On note $j\colon X^\circ\rightarrow X$ et $\hbar\colon \oX\rightarrow X$
les morphismes canoniques. On suppose $j$ quasi-compact et $\oX$ normal et localement irréductible \eqref{higgs3-sli1}. 
On observera que $\oX$ et $\oX^\circ$ sont cohérents et étale-localement connexes d'après \ref{higgs3-sli2}(iii).
On se propose d'appliquer les constructions de § \ref{higgs3-tfa}  aux morphismes de la ligne supérieure 
du diagramme commutatif suivant~:
\begin{equation}\label{higgs3-TFT1b}
\xymatrix{
{\oX^\circ}\ar[r]^{j_\oX}\ar[d]&{\oX}\ar[r]^{\hbar}\ar[d]\ar@{}[rd]|{\Box}&X\ar[d]\\
{\oeta}\ar[r]&\oS\ar[r]&S}
\end{equation}
On désigne par $E$ (resp. $\tE$) le site (resp. topos) de Faltings associé au morphisme 
$\hbar\circ j_\oX\colon \oX^\circ\rightarrow X$ \eqref{higgs3-tfa2}, et par $\ocB$ l'anneau 
de $\tE$ associé à $\oX$ \eqref{higgs3-tfa111}, qui est alors une $\co_\oK$-algèbre. On note
\begin{eqnarray}
\sigma\colon \tE &\rightarrow& X_\et,\label{higgs3-TFT1e}\\
\rho\colon X_\et\gtimes_{X_\et}\oX^\circ_\et&\rightarrow &\tE,\label{higgs3-TFT1c}
\end{eqnarray}
les morphismes canoniques \eqref{higgs3-tfa3c} et \eqref{higgs3-tfa41b}, respectivement.

Pour tout entier $n\geq 0$, on pose 
\begin{equation}\label{higgs3-TFT1d}
\ocB_n=\ocB/p^n\ocB.
\end{equation}
Pour tout $U\in \ob(\Et_{/X})$, on pose $\ocB_U=\ocB\circ \alpha_{U!}$ \eqref{higgs3-tfa1d} et
\begin{equation}\label{higgs3-TFT1dd}
\ocB_{U,n}=\ocB_U/p^n\ocB_U.
\end{equation}
On notera que l'homomorphisme canonique $\ocB_{U,n}\rightarrow \ocB_n\circ \alpha_{U!}$ 
n'est pas en général un isomorphisme~; c'est pourquoi nous n'utiliserons pas la notation $\ocB_{n,U}$.
Toutefois, les correspondances $\{U\mapsto p^n\ocB_U\}$ et $\{U\mapsto \ocB_{U,n}\}$
forment naturellement des préfaisceaux sur $E$ \eqref{higgs3-tfa2a}, et les morphismes canoniques 
\begin{eqnarray}
\{U\mapsto p^n\ocB_U\}^a&\rightarrow&p^n\ocB,\label{higgs3-TFT1f}\\
\{U\mapsto \ocB_{U,n}\}^a&\rightarrow&\ocB_n,\label{higgs3-TFT1g}
\end{eqnarray}
où  les termes de gauche désignent les faisceaux associés dans $\tE$, sont des isomorphismes
d'après (\cite{ag2} 8.2 et 8.9).

\begin{lem}\label{higgs3-TFT2}
L'anneau $\ocB$ est $\co_\oK$-plat. 
\end{lem}

Pour tout $(V\rightarrow U)\in \ob(E)$, comme $V$ est un $\oeta$-schéma, 
$\oU^V$ est $\oS$-plat en vertu de \ref{higgs3-tfa9}(i). Par suite, $\ocB$ n'a pas 
de $\co_\oK$-torsion et est donc $\co_\oK$-plat (\cite{ac} Chap.~VI §3.6 lem.~1).

\subsection{}\label{higgs3-TFT3}\index{1000911@$\tE_s$, $\delta$}\index{1000910@$\gamma$}
Comme $X_\eta$ est un ouvert de $X_\et$, {\em i.e.}, un sous-objet de l'objet final $X$ (\cite{sga4} IV 8.3),
$\sigma^*(X_\eta)=(\oX^\circ\rightarrow X_\eta)^a$ est un ouvert de $\tE$ \eqref{higgs3-tfa3c}. 
On rappelle que le topos $\tE_{/\sigma^*(X_\eta)}$
est canoniquement équivalent au topos de Faltings associé au morphisme $\oX^\circ\rightarrow X_\eta$ (\cite{ag2} 10.14). 
On note 
\begin{equation}\label{higgs3-TFT3a}
\gamma\colon \tE_{/\sigma^*(X_\eta)}\rightarrow \tE
\end{equation}
le morphisme de localisation de $\tE$ en $\sigma^*(X_\eta)$, que l'on identifie au morphisme de fonctorialité 
induit par l'injection canonique $X_\eta\rightarrow X$ \eqref{higgs3-tfa4}. On a alors une suite de trois foncteurs 
adjoints 
\begin{equation}\label{higgs3-TFT3b}
\gamma_!\colon \tE_{/\sigma^*(X_\eta)}\rightarrow \tE, \ \ \ 
\gamma^* \colon \tE\rightarrow \tE_{/\sigma^*(X_\eta)}, \ \ \ 
\gamma_*\colon \tE_{/\sigma^*(X_\eta)}\rightarrow \tE, 
\end{equation}
dans le sens que pour deux foncteurs consécutifs de la suite, celui de droite est adjoint à droite
de l'autre. Les foncteurs $\gamma_!$ et $\gamma_*$ sont pleinement fidèles (\cite{sga4} IV 9.2.4). 

On désigne par $\tE_s$ le sous-topos fermé de $\tE$ complémentaire de l'ouvert $\sigma^*(X_\eta)$, 
c'est-à-dire la sous-catégorie pleine de $\tE$ formée des faisceaux $F$ tels que $\gamma^*(F)$
soit un objet final de $\tE_{/\sigma^*(X_\eta)}$  (\cite{sga4} IV 9.3.5), et par 
\begin{equation}\label{higgs3-TFT3c}
\delta\colon \tE_s\rightarrow \tE
\end{equation} 
le plongement canonique, c'est-à-dire le morphisme de topos tel que  
$\delta_*\colon \tE_s\rightarrow \tE$ soit le foncteur d'injection canonique. 
Pour tout $F\in \ob(\tE)$, on pose $F_s=\delta^*(F)$. 

On désigne par $\Pt(\tE)$, $\Pt(\tE_{/\sigma^*(X_\eta)})$ et $\Pt(\tE_s)$ les catégories des points de $\tE$,
$\tE_{/\sigma^*(X_\eta)}$ et $\tE_s$, respectivement, et par 
\begin{equation}\label{higgs3-TFT3d}
u\colon \Pt(\tE_{/\sigma^*(X_\eta)})\rightarrow \Pt(\tE) \ \ \ {\rm et}\ \ \ v\colon \Pt(\tE_s)\rightarrow \Pt(\tE)
\end{equation}
les foncteurs induits par $\gamma$ et $\delta$, respectivement. Ces foncteurs sont pleinement fidèles,
et tout point de $\tE$ appartient à l'image essentielle de l'un ou l'autre de ces foncteurs exclusivement 
(\cite{sga4} IV 9.7.2). 

\begin{rema}\label{higgs3-TFT4}
Par définition, pour tout $F\in \ob(\tE_s)$, la projection canonique 
\begin{equation}
\sigma^*(X_\eta)\times \delta_*(F)\rightarrow \sigma^*(X_\eta) 
\end{equation}
est un isomorphisme. Il existe donc un unique morphisme $\sigma^*(X_\eta)\rightarrow \delta_*(F)$ de $\tE$, 
autrement dit, $\delta^*(\sigma^*(X_\eta))$ est un objet initial de $\tE_s$.  
\end{rema}

\begin{lem}\label{higgs3-TFT5}
{\rm (i)}\ Soit $(\oy\rightsquigarrow \ox)$ un point de $X_\et\gtimes_{X_\et}\oX^\circ_\et$ \eqref{higgs3-tfa41}. 
Pour que $\rho(\oy\rightsquigarrow \ox)$ appartienne à l'image essentielle de $u$ (resp. $v$) \eqref{higgs3-TFT3d}, 
il faut et il suffit que $\ox$ soit au-dessus de $\eta$ (resp.~$s$). 

{\rm (ii)}\ La famille des points de $\tE_{/\sigma^*(X_\eta)}$ (resp. $\tE_s$) définie par la famille des points 
$\rho(\oy\rightsquigarrow \ox)$ de $\tE$ tels que $\ox$ soit au-dessus de $\eta$ (resp. $s$)
est conservative. 
\end{lem}

(i) En effet, pour que $\rho(\oy\rightsquigarrow \ox)$ appartienne à l'image essentielle de $u$ (resp. $v$), 
il faut et il suffit que $(\sigma^*(X_\eta))_{\rho(\oy\rightsquigarrow \ox)}$ soit un singleton (resp. vide). 
Par ailleurs, on a un isomorphisme
canonique (\cite{ag2} (10.18.1))
\begin{equation}\label{higgs3-TFT5a}
(\sigma^*(X_\eta))_{\rho(\oy\rightsquigarrow \ox)}\stackrel{\sim}{\rightarrow} (X_\eta)_{\ox},
\end{equation}
d'où la proposition.  

(ii) Cela résulte de (i), (\cite{ag2} 10.21) et (\cite{sga4} IV 9.7.3).

\begin{lem}\label{higgs3-TFT6}
Pour tout faisceau $F=\{U\mapsto F_U\}$ de $\tE$, les propriétés suivantes sont équivalentes~:
\begin{itemize}
\item[{\rm (i)}] $F$ est un objet de $\tE_s$. 
\item[{\rm (ii)}] Pour tout $U\in \Et_{/X_\eta}$, $F_U$ est un objet final de $\oU^\circ_\fet$, 
i.e., est représentable par $\oU^\circ$. 
\item[{\rm (iii)}] Pour tout point $(\oy\rightsquigarrow \ox)$ de $X_\et\gtimes_{X_\et}\oX^\circ_\et$ \eqref{higgs3-tfa41}
tel que $\ox$ soit au-dessus de $\eta$, la fibre $F_{\rho(\oy\rightsquigarrow \ox)}$ de $F$ en $\rho(\oy\rightsquigarrow \ox)$
est un singleton. 
\end{itemize}
\end{lem}

En effet, d'après (\cite{ag2} 5.38), on a un isomorphisme canonique 
\begin{equation}\label{higgs3-TFT6a}
\gamma^*(F)\stackrel{\sim}{\rightarrow} \{U'\mapsto F_{U'}\}, \ \ \ (U'\in \ob(\Et_{/X_\eta})). 
\end{equation}
Comme $\{U'\mapsto \oU'^\circ\}$, pour $U'\in \ob(\Et_{/X_\eta})$, est un objet final de $\tE_{/\sigma^*(X_\eta)}$ \eqref{higgs3-tfa31},
les conditions (i) et (ii) sont équivalentes. 
Par ailleurs, les conditions (i) et (iii) sont équivalentes en vertu de \ref{higgs3-TFT5}(ii).

\begin{lem}\label{higgs3-TFT7}
Pour tout $n\geq 0$, l'anneau $\ocB_n$ \eqref{higgs3-TFT1d} est un objet de $\tE_s$.
\end{lem}

En effet, pour tout point $(\oy\rightsquigarrow \ox)$ de $X_\et\gtimes_{X_\et}\oX^\circ_\et$ \eqref{higgs3-tfa41}
tel que $\ox$ soit au-dessus de $\eta$, l'homomorphisme canonique $\sigma^{-1}(\co_X)\rightarrow \ocB$ \eqref{higgs3-tfa111}
induit un homomorphisme $\co_{X,\ox}\rightarrow \ocB_{\rho(\oy\rightsquigarrow \ox)}$ (\cite{ag2} (10.18.1)).
Par suite, $p$ est inversible dans $\ocB_{\rho(\oy\rightsquigarrow \ox)}$, d'où la proposition en vertu de \ref{higgs3-TFT6}.

\subsection{}\label{higgs3-TFT8}\index{1000915@$\sigma_s$, $\sigma_\eta$}
On note $a\colon X_s\rightarrow X$ et $b\colon X_\eta\rightarrow X$ les injections canoniques. 
Le topos $\tE_{/\sigma^*(X_\eta)}$ étant canoniquement équivalent au topos de Faltings associé au morphisme 
$\oX^\circ\rightarrow X_\eta$, notons 
\begin{equation}\label{higgs3-TFT8a}
\sigma_\eta\colon \tE_{/\sigma^*(X_\eta)}\rightarrow X_{\eta,\et}
\end{equation}
le morphisme canonique \eqref{higgs3-tfa3c}. D'après (\cite{ag2} (10.12.6)), le diagramme 
\begin{equation}\label{higgs3-TFT8b}
\xymatrix{
{\tE_{/\sigma^*(X_\eta)}}\ar[r]^{\sigma_\eta}\ar[d]_{\gamma}&{X_{\eta,\et}}\ar[d]^{b}\\
{\tE}\ar[r]^\sigma&{X_\et}}
\end{equation}
est commutatif à isomorphisme canonique près. En vertu de (\cite{sga4} IV 9.4.3), il existe un morphisme 
\begin{equation}\label{higgs3-TFT8c}
\sigma_s\colon \tE_s\rightarrow X_{s,\et}
\end{equation}
unique à isomorphisme près tel que le diagramme 
\begin{equation}\label{higgs3-TFT8d}
\xymatrix{
{\tE_s}\ar[r]^{\sigma_s}\ar[d]_{\delta}&{X_{s,\et}}\ar[d]^{a}\\
{\tE}\ar[r]^\sigma&{X_\et}}
\end{equation}
soit commutatif à isomorphisme près. Par définition, on a un isomorphisme canonique 
\begin{equation}\label{higgs3-TFT8f}
\sigma^*\circ a_*\stackrel{\sim}{\rightarrow}\delta_*\circ \sigma_s^*.
\end{equation}
Les foncteurs $a_*$ et $\delta_*$ étant exacts, 
pour tout groupe abélien $F$ de $\tE_s$ et tout entier $i\geq 0$, on a un isomorphisme canonique 
\begin{equation}\label{higgs3-TFT8e}
a_*(\rR^i\sigma_{s*}(F))\stackrel{\sim}{\rightarrow}\rR^i\sigma_*(\delta_*F). 
\end{equation}

\subsection{}\label{higgs3-TFT9}\index{1000920@$a$, $a_n$, $\oa_n$, $\iota_n$, $\oiota_n$}\index{1000922@$\sigma_n$}
Pour tout entier $\geq 1$, on note $a\colon X_s\rightarrow X$, $a_n\colon X_s\rightarrow X_n$, $\iota_n\colon X_n\rightarrow X$
et $\oiota_n\colon \oX_n\rightarrow \oX$ les injections canoniques. 
Le corps résiduel de $\co_K$ étant algébriquement clos, il existe un unique $S$-morphisme $s\rightarrow \oS$. 
Celui-ci induit des immersions fermées $\oa\colon X_s\rightarrow \oX$ et $\oa_n\colon X_s\rightarrow \oX_n$
qui relèvent $a$ et $a_n$, respectivement. 
\begin{equation}\label{higgs3-TFT9a}
\xymatrix{
{X_s}\ar[r]_{\oa_n}\ar@{=}[d]\ar@/^1pc/[rr]^{\oa}&{\oX_n}\ar[r]_{\oiota_n}\ar[d]^{\hbar_n}&\oX\ar[d]^\hbar\\
{X_s}\ar[r]^{a_n}\ar@/_1pc/[rr]_{a}&{X_n}\ar[r]^{\iota_n}&X}
\end{equation}
L'homomorphisme canonique $\sigma^{-1}(\hbar_*(\co_\oX))\rightarrow \ocB$ \eqref{higgs3-tfa111} induit 
un homomorphisme \eqref{higgs3-TFT1d}
\begin{equation}\label{higgs3-TFT9b}
\sigma^{-1}(\iota_{n*}(\hbar_{n*}(\co_{\oX_n})))\rightarrow \ocB_n.
\end{equation}
Comme $\oa_n$ est un homéomorphisme universel, on peut considérer $\co_{\oX_n}$
comme un faisceau de $X_{s,\et}$ (cf. \ref{higgs3-not2}). 
On peut alors identifier les anneaux $\iota_{n*}(\hbar_{n*}(\co_{\oX_n}))$ et $a_*(\co_{\oX_n})$ 
de $X_\et$ et par suite les anneaux $\sigma^{-1}(\iota_{n*}(\hbar_{n*}(\co_{\oX_n})))$ et 
$\delta_*(\sigma^*_s(\co_{\oX_n}))$ de $\tE_s$ \eqref{higgs3-TFT8f}. 
Comme $\ocB_n$ est un objet de $\tE_s$ \eqref{higgs3-TFT7}, on peut considérer 
\eqref{higgs3-TFT9b} comme un homomorphisme de $\tE_s$ 
\begin{equation}\label{higgs3-TFT9c}
\sigma_s^*(\co_{\oX_n})\rightarrow \ocB_n.
\end{equation}
Le morphisme $\sigma_s$ \eqref{higgs3-TFT8c} est donc sous-jacent à un morphisme de topos annelés, que l'on note
\begin{equation}\label{higgs3-TFT9d}
\sigma_n\colon (\tE_s,\ocB_n)\rightarrow (X_{s,\et},\co_{\oX_n}).
\end{equation}

\subsection{}\label{higgs3-TFT14}\index{1000925@$\tE_s^{\mN^\circ}$} \index{1000926@$\bvocB$}
\index{1000927@$\bvsigma$,$\bvu$}\index{1000929@$\top\colon (\tE_s^{\mN^\circ},\bvocB)\rightarrow (X_{s,\zar},\co_{\fX})$}\index{1000928@$\fX$}
On désigne par $\bvocB$ l'anneau $(\ocB_{n+1})_{n\in \mN}$ de $\tE_s^{\mN^\circ}$ \eqref{higgs3-spsa2} et par 
$\co_{\bvoX}$ l'anneau $(\co_{\oX_{n+1}})_{n\in \mN}$ de $X_{s,\zar}^{\mN^\circ}$ ou de $X_{s,\et}^{\mN^\circ}$,
selon le contexte. Cet abus de notation n'induit aucune confusion (cf. \ref{higgs3-not2}). 
D'après \ref{higgs3-spsa35}, les morphismes $(\sigma_{n+1})_{n\in \mN}$ \eqref{higgs3-TFT9d} induisent un morphisme de topos annelés 
\begin{equation}\label{higgs3-TFT14a}
\bvsigma\colon (\tE_s^{\mN^\circ},\bvocB)\rightarrow(X_{s,\et}^{\mN^\circ},\co_{\bvoX}).
\end{equation}
Pour tout entier $n\geq 1$, on désigne par
\begin{equation}\label{higgs3-TFT9e}
u_n\colon (X_{s,\et},\co_{\oX_n})\rightarrow (X_{s,\zar},\co_{\oX_n})
\end{equation}
le morphisme canonique de topos annelés (cf. \ref{higgs3-not2} et \ref{higgs3-TFT9}).  
Les morphismes $(u_{n+1})_{n\in \mN}$ induisent un morphisme de topos annelés
\begin{equation}\label{higgs3-TFT14b}
\bvu\colon (X_{s,\et}^{\mN^\circ},\co_{\bvoX})\rightarrow (X_{s,\zar}^{\mN^\circ},\co_{\bvoX}).
\end{equation}

On désigne par $\fX$ le schéma formel complété $p$-adique de $\oX$. 
L'espace topologique sous-jacent à $\fX$ est canoniquement isomorphe à $X_s$.
Il est annelé par le faisceau d'anneaux topologiques $\co_\fX$ limite 
projective des faisceaux d'anneaux pseudo-discrets $(\co_{\oX_{n+1}})_{n\in \mN}$ (\cite{ega1n} 0.3.9.1). 
On désigne par
\begin{equation}\label{higgs3-TFT14d}
\lambda\colon X_{s,\zar}^{\mN^\circ}\rightarrow X_{s,\zar}
\end{equation}
le morphisme défini dans \eqref{higgs3-spsa3c}. Le faisceau d'anneaux $\lambda_*(\co_{\bvoX})$ est canoniquement 
isomorphe au faisceau d'anneaux (sans topologies) sous-jacent à $\co_\fX$ (\cite{ega1n} 0.3.9.1 et 0.3.2.6). 
On considère alors $\lambda$ comme un morphisme de topos annelés, respectivement par $\co_{\bvoX}$ et $\co_\fX$. 
On note 
\begin{equation}\label{higgs3-TFT14e}
\top\colon (\tE_s^{\mN^\circ},\bvocB)\rightarrow (X_{s,\zar},\co_{\fX})
\end{equation}
le morphisme composé de topos annelés  $\lambda\circ \bvu\circ \bvsigma$.

\subsection{}\label{higgs3-TFT12}
Soient $g\colon X'\rightarrow X$ un morphisme cohérent, $X'^\rhd$ un sous-schéma ouvert de $X'^\circ=X^\circ\times_XX'$. 
Pour tout $X'$-schéma $U'$, on pose
\begin{equation}
U'^\rhd=U'\times_{X'}X'^\rhd.
\end{equation} 
On note $j'\colon X'^\rhd\rightarrow X'$ l'injection canonique 
et $\hbar'\colon \oX'\rightarrow X'$ le morphisme canonique \eqref{higgs3-TFT1a}.  
Supposons $\oX'$ normal et localement irréductible. 
On désigne par $E'$ (resp. $\tE'$) le site (resp. topos) de 
Faltings associé au morphisme  $\hbar'\circ j'_{\oX'}\colon \oX'^\rhd\rightarrow X'$ \eqref{higgs3-tfa2},
par $\ocB'$ l'anneau de $\tE'$ associé à $\oX'$ \eqref{higgs3-tfa111} et par 
\begin{equation}\label{higgs3-TFT12a}
\sigma'\colon (\tE',\ocB') \rightarrow (X'_\et,\hbar'_*(\co_{\oX'}))
\end{equation}
le morphisme canonique de topos annelés \eqref{higgs3-tfa111}. On note
$\tE'_s$ le sous-topos fermé de $\tE'$ complémentaire de l'ouvert $\sigma'^*(X'_\eta)$, 
\begin{equation}\label{higgs3-TFT12d}
\delta'\colon \tE'_s\rightarrow \tE'
\end{equation}
le plongement canonique et 
\begin{equation}\label{higgs3-TFT12cd}
\sigma'_s\colon \tE'_s\rightarrow X'_{s,\et}
\end{equation}
le morphisme canonique de topos \eqref{higgs3-TFT8c}. 
Pour tout entier $n\geq 1$, on pose $\ocB'_n=\ocB'/p^n\ocB'$, et on désigne par
\begin{equation}
\sigma'_n\colon (\tE'_s,\ocB'_n)\rightarrow (X'_{s,\et},\co_{\oX'_n})
\end{equation}
le morphisme de topos annelés induit par $\sigma'$ \eqref{higgs3-TFT9d}.

Soit
\begin{equation}\label{higgs3-TFT12c}
\Phi\colon (\tE',\ocB')\rightarrow (\tE,\ocB)
\end{equation}
le morphisme de topos annelés déduit de $g$ par fonctorialité \eqref{higgs3-tfa13}. D'après (\cite{ag2} (10.12.6)) et 
les définitions \ref{higgs3-tfa111} et \ref{higgs3-tfa13}, 
le diagramme de morphismes de topos annelés 
\begin{equation}\label{higgs3-TFT12cc}
\xymatrix{
{(\tE',\ocB')}\ar[r]^{\Phi}\ar[d]_{\sigma'}&{(\tE,\ocB)}\ar[d]^{\sigma}\\
{(X'_{\et},\hbar'_*(\co_{\oX'}))}\ar[r]^{\ogg}&{(X_{\et},\hbar_*(\co_{\oX}))}}
\end{equation}
où $\ogg$ est le morphisme induit par $g$, 
est commutatif à isomorphisme canonique près, dans le sens de (\cite{egr1} 1.2.3).
On a un isomorphisme canonique $\Phi^*(\sigma^*(X_\eta))\simeq \sigma'^*(X'_\eta)$ \eqref{higgs3-tfa4a}.
En vertu de (\cite{sga4} IV 9.4.3), il existe donc un morphisme de topos
\begin{equation}\label{higgs3-TFT12h}
\Phi_s\colon \tE'_s\rightarrow \tE_s
\end{equation}
unique à isomorphisme près tel que le diagramme 
\begin{equation}\label{higgs3-TFT12i}
\xymatrix{
{\tE'_s}\ar[r]^{\Phi_s}\ar[d]_{\delta'}&{\tE_s}\ar[d]^{\delta}\\
{\tE'}\ar[r]^\Phi&{\tE}}
\end{equation}
soit commutatif à isomorphisme près, et même $2$-cartésien. Il résulte de (\cite{ag2} (10.12.6)) et (\cite{sga4} IV 9.4.3) 
que le diagramme de morphisme de topos 
\begin{equation}\label{higgs3-TFT12ij}
\xymatrix{
{\tE'_s}\ar[r]^{\Phi_s}\ar[d]_{\sigma'_s}&{\tE_s}\ar[d]^{\sigma_s}\\
{X'_{s,\et}}\ar[r]^{g_s}&{X_{s,\et}}}
\end{equation}
est commutatif à isomorphisme canonique près.

L'homomorphisme canonique $\Phi^{-1}(\ocB)\rightarrow \ocB'$
induit un homomorphisme $\Phi_s^*(\ocB_n)\rightarrow \ocB'_n$. 
Le morphisme $\Phi_s$ est donc sous-jacent à un morphisme de topos annelés, que l'on note 
\begin{equation}\label{higgs3-TFT12j}
\Phi_n\colon (\tE'_s,\ocB'_n)\rightarrow (\tE_s,\ocB_n).
\end{equation}
Il résulte de \eqref{higgs3-TFT12cc} et \eqref{higgs3-TFT12ij} que le diagramme de morphismes de topos annelés 
\begin{equation}\label{higgs3-TFT12jj}
\xymatrix{
{(\tE'_s,\ocB'_n)}\ar[r]^{\Phi_n}\ar[d]_{\sigma'_n}&{(\tE_s,\ocB_n)}\ar[d]^{\sigma_n}\\
{(X'_{s,\et},\co_{\oX'_n})}\ar[r]^{\ogg_n}&{(X_{s,\et},\co_{\oX_n})}}
\end{equation}
où $\ogg_n$ est le morphisme induit par $g$, est commutatif à isomorphisme canonique près, 
dans le sens de (\cite{egr1} 1.2.3). 

Posons $\bvocB'=(\ocB'_{n+1})_{n\in \mN}$, qui est un anneau de $\tE'^{\mN^\circ}_s$.
D'après \ref{higgs3-spsa35}, les morphismes $(\Phi_{n+1})_{n\in \mN}$ définissent un morphisme de topos annelés
\begin{equation}\label{higgs3-TFT12ja}
\bvPhi\colon (\tE'^{\mN^\circ}_s,\bvocB')\rightarrow (\tE^{\mN^\circ}_s,\bvocB).
\end{equation} 
On désigne par $\fX'$ le schéma formel complété $p$-adique de $\oX'$ et par 
\begin{equation}\label{higgs3-TFT12jb}
\top'\colon (\tE'^{\mN^\circ}_s,\bvocB')\rightarrow (\fX'_{\zar},\co_{\fX'})
\end{equation}
le morphisme de topos annelés défini dans \eqref{higgs3-TFT14e} relativement à $(X',X'^\rhd)$.
Il résulte aussitôt de \eqref{higgs3-TFT12jj} et du caractère fonctoriel des morphismes 
$\bvu$ \eqref{higgs3-TFT14b} et $\lambda$ \eqref{higgs3-TFT14d} que le diagramme de morphismes de topos
annelés 
\begin{equation}\label{higgs3-TFT12jc}
\xymatrix{
{(\tE'^{\mN^\circ}_s,\bvocB')}\ar[r]^{\bvPhi}\ar[d]_{\top'}&{(\tE^{\mN^\circ}_s,\bvocB)}\ar[d]^\top\\
{(X'_{s,\zar},\co_{\fX'})}\ar[r]^{\fgg}&{(X_{s,\zar},\co_{\fX})}}
\end{equation}
où $\fgg$ est le morphisme induit par $g$, 
est commutatif à isomorphisme canonique près, 
dans le sens de (\cite{egr1} 1.2.3). On en déduit, pour tout $\bvocB$-module $\cF$ et tout 
entier $q\geq 0$, un morphisme de changement de base 
\begin{equation}\label{higgs3-TFT12jd}
\fgg^*(\rR^q\top_*(\cF))\rightarrow \rR^q\top'_*(\bvPhi^*(\cF)).
\end{equation}

\begin{lem}\label{higgs3-TFT16}
Les hypothèses étant celles de \eqref{higgs3-TFT12}, supposons de plus $g$ étale et $X'^\rhd=X'^\circ$. 
Alors $\Phi_s$ \eqref{higgs3-TFT12h} est canoniquement isomorphe 
au morphisme de localisation de $\tE_s$ en $\sigma^*_s(X'_s)$.
\end{lem}

On observera d'abord que $(\oX'^\rhd\rightarrow X')$ est un objet de $E$ et que le faisceau associé n'est autre
que $\sigma^*(X')$. Par ailleurs, on a un isomorphisme canonique $\delta^*(\sigma^*(X'))\stackrel{\sim}{\rightarrow}\sigma^*_s(X'_s)$ \eqref{higgs3-TFT8d}. Notons $j\colon \tE_{/\sigma^*(X')}\rightarrow \tE$ 
(resp. $j_s\colon (\tE_s)_{/\sigma^*_s(X'_s)}\rightarrow \tE_s$) 
le morphisme de localisation de $\tE$ en $\sigma^*(X')$ (resp. de $\tE_s$ en $\sigma^*_s(X'_s)$). 
D'après (\cite{sga4} IV 5.10), le morphisme $\delta$ induit un morphisme 
\begin{equation}
\delta_{/\sigma^*(X')}\colon (\tE_s)_{/\sigma^*_s(X'_s)}\rightarrow \tE_{/\sigma^*(X')}
\end{equation}
qui s'insère dans un diagramme commutatif à isomorphisme canonique près
\begin{equation}
\xymatrix{
{(\tE_s)_{/\sigma^*_s(X'_s)}}\ar[r]^-(0.4){j_s}\ar[d]_{\delta_{/\sigma^*(X')}}&{\tE_s}\ar[d]^{\delta}\\
{\tE_{/\sigma^*(X')}}\ar[r]^-(0.5)j&{\tE}}
\end{equation}
Ce diagramme est en fait 2-cartésien d'après (\cite{sga4} IV 5.11). 
D'autre part, les topos $\tE'$ et $\tE_{/\sigma^*(X')}$ sont canoniquement équivalents 
et le morphisme $\Phi$ s'identifie à $j$ en vertu de (\cite{ag2} 10.14). 
Comme le diagramme \eqref{higgs3-TFT12i} est aussi 2-cartésien, la proposition s'ensuit. 

\begin{lem}\label{higgs3-TFT13}
Les hypothèses étant celles de \eqref{higgs3-TFT12}, supposons de plus que $X'^\rhd=X'^\circ$
et que l'une des deux conditions suivantes soit remplie~:
\begin{itemize}
\item[{\rm (i)}] $g$ est étale.
\item[{\rm (ii)}]  $X'$ est le localisé strict de $X$ en un point géométrique $\ox$.
\end{itemize}
Alors pour tout entier $n\geq 1$, l'homomorphisme $\Phi_s^*(\ocB_n)\rightarrow \ocB'_n$ est un isomorphisme. 
\end{lem}

Cela résulte de \ref{higgs3-tfa14}(i) et \ref{higgs3-tfa16}. 

\begin{prop}\label{higgs3-TFT17}
Les hypothèses étant celles de \eqref{higgs3-TFT12}, supposons de plus $g$ étale et $X'^\rhd=X'^\circ$,
et notons encore $\lambda\colon \tE_s^{\mN^\circ}\rightarrow \tE_s$ le morphisme de topos défini dans \eqref{higgs3-spsa3c}.
Alors $\bvPhi$ \eqref{higgs3-TFT12ja} est canoniquement isomorphe 
au morphisme de localisation du topos annelé $(\tE_s^{\mN^\circ},\bvocB)$ en $\lambda^*(\sigma^*_s(X'_s))$.
\end{prop}

En effet, le morphisme de topos $\Phi_s^{\mN^\circ}\colon \tE'^{\mN^\circ}_s
\rightarrow   \tE^{\mN^\circ}_s$ s'identifie au morphisme de localisation de $\tE^{\mN^\circ}_s$ en 
$\lambda^*(\sigma^*_s(X'_s))$, d'après \ref{higgs3-TFT16} et \ref{higgs3-spsa100}(ii). D'autre part, l'homomorphisme canonique 
$(\Phi_s^{\mN^\circ})^*(\bvocB)\rightarrow \bvocB$ est un isomorphisme en vertu de \ref{higgs3-TFT13} 
et \eqref{higgs3-spsa1g}; d'où la proposition. 

\begin{cor}\label{higgs3-TFT15}
Les hypothèses étant celles de \eqref{higgs3-TFT12}, supposons de plus que 
$g$ soit une immersion ouverte et que $X'^\rhd=X'^\circ$. 
Alors pour tout $\bvocB$-module $\cF$ et tout 
entier $q\geq 0$, le morphisme de changement de base \eqref{higgs3-TFT12jd} 
\begin{equation}\label{higgs3-TFT15a}
\fgg^*(\rR^q\top_*(\cF))\rightarrow \rR^q\top'_*(\bvPhi^*(\cF))
\end{equation}
est un isomorphisme.  
\end{cor}

En effet, les carrés du diagramme de morphismes de topos 
\begin{equation}
\xymatrix{
{\tE_s^{\mN^\circ}}\ar[r]^{\lambda}\ar[d]_{\sigma_s^{\mN^\circ}}&{\tE_s}\ar[d]^{\sigma_s}\\
{X_{s,\et}^{\mN^\circ}}\ar[r]^{\lambda}\ar[d]_{u_s^{\mN^\circ}}&{X_{s,\et}}\ar[d]^{u_s}\\
{X_{s,\zar}^{\mN^\circ}}\ar[r]^{\lambda}&{X_{s,\zar}}}
\end{equation}
où $u_s$ est le morphisme canonique \eqref{higgs3-not2b} et $\lambda$ désigne (abusivement) les morphismes 
définis dans \eqref{higgs3-spsa3c}, sont commutatifs à isomorphismes canoniques près \eqref{higgs3-spsa1g}. 
Notant $\ttT= \lambda\circ u_s^{\mN^\circ}\circ \sigma_s^{\mN^\circ} $, 
qui est le morphisme de topos sous-jacent à $\top$ \eqref{higgs3-TFT14e}, on en déduit un isomorphisme
\begin{equation}
\ttT^*(X'_s)\stackrel{\sim}{\rightarrow} \lambda^*(\sigma_s^*(X'_s)).
\end{equation} 
Il résulte alors de \eqref{higgs3-TFT12jc} et \ref{higgs3-TFT17}
que $\top'$ s'identifie au morphisme $\top_{/X'_s}$ (cf. \cite{sga4} IV 5.10);
d'où la proposition.

\section{Les algèbres de Higgs-Tate}\label{higgs3-RGG}

\subsection{}\label{higgs3-RGG1}\index{10001001@$(X,\cM_X)$, $X^\circ$}\index{1000904@$U^\circ=U\times_XX^\circ$}
\index{10001003@$\tOmega^1_{X/S}$, $\tOmega^1_{\oX_n/\oS_n}$}\index{10001005@$(\oX,\cM_{\oX})$, $(\coX,\cM_{\coX})$}
\index{10001007@$(\tX,\cM_\tX)$}\index{deformation@$(\cA_2(\oS),\cM_{\cA_2(\oS)})$-déformation}
Dans cette section, $(S,\cM_S)$ désigne le trait logarithmique fixé dans \eqref{higgs3-not1} et 
\begin{equation}
f\colon (X,\cM_X)\rightarrow (S,\cM_S)
\end{equation} 
un morphisme adéquat de schémas logarithmiques \eqref{higgs3-slad6}. 
On note $X^\circ$ le sous-schéma ouvert maximal de $X$
où la structure logarithmique $\cM_X$ est triviale~; c'est un sous-schéma ouvert de $X_\eta$. 
Pour tout $X$-schéma usuel $U$, on pose  
\begin{equation}\label{higgs3-RGG1a}
U^\circ=U\times_XX^\circ.
\end{equation} 
On rappelle que pour tout $S$-schéma $Y$, on a posé
$\oY=Y\times_S\oS$, $\coY=Y\times_S\coS$ et pour tout entier $n\geq 1$, $Y_n=Y\times_SS_n$ \eqref{higgs3-not1c}.
On note $j\colon X^\circ\rightarrow X$ et $\hbar\colon \oX\rightarrow X$
les morphismes canoniques. Pour alléger les notations, on pose 
\begin{equation}\label{higgs3-RGG1aa}
\tOmega^1_{X/S}=\Omega^1_{(X,\cM_X)/(S,\cM_S)},
\end{equation}
que l'on considère comme un faisceau de $X_\zar$ ou $X_\et$, selon le contexte (cf. \ref{higgs3-not2}).
Pour tout entier $n\geq 1$,  on pose 
\begin{equation}\label{higgs3-RGG1ab}
\tOmega^1_{\oX_n/\oS_n}=\tOmega^1_{X/S}\otimes_{\co_X}\co_{\oX_n},
\end{equation}
que l'on considère comme un faisceau de $X_{s,\zar}$ ou $X_{s,\et}$, selon le contexte (cf. \ref{higgs3-not2} et \ref{higgs3-TFT9}). 
On munit $\oX$ et $\coX$ des structures logarithmiques $\cM_{\oX}$ et $\cM_{\coX}$ 
images inverses de $\cM_X$. On a alors des isomorphismes canoniques  \eqref{higgs3-not1}
\begin{eqnarray}
(\oX,\cM_{\oX})\stackrel{\sim}{\rightarrow}(X,\cM_X)\times_{(S,\cM_S)}(\oS,\cM_\oS),\label{higgs3-RGG1bb}\\
(\coX,\cM_{\coX})\stackrel{\sim}{\rightarrow}(X,\cM_X)\times_{(S,\cM_S)}(\coS,\cM_\coS),\label{higgs3-RGG1b}
\end{eqnarray} 
le produit étant indifféremment pris dans la catégorie des schémas logarithmiques ou 
dans celle des schémas logarithmiques fins. 

Une {\em $(\cA_2(\oS),\cM_{\cA_2(\oS)})$-déformation lisse} de $(\coX,\cM_{\coX})$ \eqref{higgs3-not4} 
est la donnée d'un morphisme lisse de schémas logarithmiques fins 
\begin{equation}\label{higgs3-RGG1c}
\ttf\colon (\tX,\cM_\tX)\rightarrow (\cA_2(\oS),\cM_{\cA_2(\oS)})
\end{equation}
et d'un $(\coS,\cM_\coS)$-isomorphisme
\begin{equation}\label{higgs3-RGG1d}
(\coX,\cM_{\coX})\stackrel{\sim}{\rightarrow}(\tX,\cM_\tX)\times_{(\cA_2(\oS),\cM_{\cA_2(\oS)})}(\coS,\cM_\coS),
\end{equation}
le produit étant indifféremment pris dans la catégorie des schémas logarithmiques ou  
dans celle des schémas logarithmiques fins (cf. \cite{kato1} 3.14).
Dans la suite de cette section, {\em on suppose qu'il existe une telle déformation $(\tX,\cM_\tX)$ que l'on fixe}.

\subsection{}\label{higgs3-RGG2}\index{1000805@$\pi\colon E\rightarrow \Et_{/X}$|textbf}\index{Site fibre de Faltings@Site fibré de Faltings}
\index{1000815@$\sigma$|textbf}\index{1000820@$\rho\colon X_\et\gtimes_{X_\et}Y_\et\rightarrow \tE$|textbf}\index{1000911@$\tE_s$, $\delta$|textbf}
\index{1000906@$\ocB_n$, $\ocB_{U,n}$|textbf}\index{1000922@$\sigma_n$}\index{1000811@$\tE$|textbf}
D'après \ref{higgs3-slad4}(iii), les schémas $X$ et $\oX$ sont normaux et localement irréductibles. 
Par ailleurs, $X$ étant noethérien, $j$ est quasi-compact. On peut donc 
appliquer les constructions de § \ref{higgs3-TFT}  aux morphismes de la ligne supérieure 
du diagramme commutatif suivant~:
\begin{equation}\label{higgs3-RGG2a}
\xymatrix{
{\oX^\circ}\ar[r]^{j_\oX}\ar[d]&{\oX}\ar[r]^{\hbar}\ar[d]\ar@{}[rd]|{\Box}&X\ar[d]\\
{\oeta}\ar[r]&\oS\ar[r]&S}
\end{equation}
On note
\begin{equation}\label{higgs3-RGG2d}
\pi\colon E\rightarrow \Et_{/X}
\end{equation}
le $\mU$-site fibré de Faltings associé au morphisme $h=\hbar\circ j_\oX\colon \oX^\circ\rightarrow X$ \eqref{higgs3-tfa1}.
On munit $E$ de la topologie co-évanescente et on désigne par $\tE$ le topos des faisceaux de $\mU$-ensembles sur 
$E$ \eqref{higgs3-tfa2}, par $\ocB$ l'anneau de $\tE$ associé à $\oX$ \eqref{higgs3-tfa111} et par
\begin{eqnarray}
\sigma\colon \tE \rightarrow X_\et,\label{higgs3-RGG2b}\\
\rho\colon X_\et\gtimes_{X_\et}\oX^\circ_\et\rightarrow \tE,\label{higgs3-RGG2c}
\end{eqnarray}
les morphismes canoniques \eqref{higgs3-tfa3c} et \eqref{higgs3-tfa41b}. 
On note $\tE_s$ le sous-topos fermé de $\tE$ complémentaire de l'ouvert $\sigma^*(X_\eta)$ \eqref{higgs3-TFT3}, 
\begin{equation}\label{higgs3-RGG2bc}
\delta\colon \tE_s\rightarrow \tE
\end{equation} 
le plongement canonique et 
\begin{equation}\label{higgs3-RGG2e}
\sigma_s\colon \tE_s\rightarrow X_{s,\et}
\end{equation} 
le morphisme de topos induit par $\sigma$ \eqref{higgs3-TFT8c}.
Pour tout entier $n\geq 0$, on pose 
\begin{equation}\label{higgs3-RGG2f}
\ocB_n=\ocB/p^n\ocB.
\end{equation}
Pour tout $U\in \ob(\Et_{/X})$, on pose $\ocB_U=\ocB\circ \alpha_{U!}$ \eqref{higgs3-tfa1d} et
\begin{equation}\label{higgs3-RGG2ff}
\ocB_{U,n}=\ocB_U/p^n\ocB_U.
\end{equation}
On notera que l'homomorphisme canonique $\ocB_{U,n}\rightarrow \ocB_n\circ \alpha_{U!}$ 
n'est pas en général un isomorphisme (cf. \eqref{higgs3-TFT1g}). 
On rappelle \eqref{higgs3-TFT7} que $\ocB_n$ est un anneau de $\tE_s$. Si $n\geq 1$, on note
\begin{equation}\label{higgs3-RGG2g}
\sigma_n\colon (\tE_s,\ocB_n)\rightarrow (X_{s,\et},\co_{\oX_n})
\end{equation}
le morphisme canonique de topos annelés \eqref{higgs3-TFT9d}.

\subsection{}\label{higgs3-RGG8}\index{10001012@$\bP$}
On désigne par $\bP$ la sous-catégorie pleine de $\Et_{/X}$ formée des schémas affines $U$ tels que
le morphisme $(U,\cM_X|U)\rightarrow (S,\cM_S)$ induit par $f$ admette une carte adéquate \eqref{higgs3-slad3}. 
On munit $\bP$ de la topologie induite par celle de $\Et_{/X}$. Comme $X$ est noethérien et donc quasi-séparé, 
tout objet de $\bP$ est cohérent sur $X$. Par suite, $\bP$ est 
une famille $\mU$-petite, topologiquement génératrice du site $\Et_{/X}$ et est stable par produits fibrés. 
On désigne par 
\begin{equation}\label{higgs3-RGG8a}
\pi_\bP\colon E_\bP\rightarrow \bP
\end{equation} 
le site fibré déduit de $\pi$ \eqref{higgs3-RGG2d}
par changement de base par le foncteur d'injection canonique $\bP\rightarrow \Et_{/X}$. 
On munit $E_\bP$ de 
la topologie co-évanescente définie par $\pi_\bP$ et on note $\tE_\bP$ le topos des faisceaux de $\mU$-ensembles 
sur $E_\bP$. D'après (\cite{ag2} 5.21 et 5.22), 
la topologie de $E_\bP$ est induite par celle de $E$ au moyen du foncteur de projection canonique 
$E_\bP\rightarrow E$, et celui-ci induit par restriction une équivalence de catégories 
\begin{equation}\label{higgs3-RGG8b}
\tE\stackrel{\sim}{\rightarrow}\tE_\bP.
\end{equation}

\begin{rema}\label{higgs3-RGG65}
Soient $U$ un objet de $\bP$, $\oy$ un point géométrique générique de $\oU^\circ$, $\oR^\oy_U$
l'anneau défini dans \eqref{higgs3-tfa112c}. Les schémas $U$ et $\oU$
étant localement irréductibles \eqref{higgs3-sli3}, ils sont les sommes des schémas induits sur 
leurs composantes irréductibles. Notons $U^c$ (resp. $\oU^\star$) la composante irréductible de $U$ (resp. $\oU$)
contenant $\oy$. De même, $\oU^\circ$
est la somme des schémas induits sur ses composantes irréductibles et $\oU^{\star \circ}=\oU^\star\times_XX^\circ$ 
est la composante irréductible de $\oU^\circ$ contenant $\oy$.  
Alors $U^c$ est naturellement un objet de $\bP$ au-dessus de $U$, et l'homomorphisme canonique 
$\oR^\oy_{U}\rightarrow \oR^\oy_{U^c}$ est un isomorphisme $\pi_1(\oU^{\star \circ},\oy)$-équivariant. 
Si $U^c_s=\emptyset$, $\oR^\oy_{U}$ est une $\oK$-algèbre.
Supposons $U^c_s\not=\emptyset$, de sorte que
le morphisme $(U^c,\cM_X|U^c)\rightarrow (S,\cM_S)$ induit par $f$ vérifie les hypothèses de (\cite{ag1} 6.2). 
L'algèbre $\oR^\oy_U$ munie de l'action de $\pi_1(\oU^{\star \circ},\oy)$ correspond alors à l'algèbre $\oR$ 
munie de l'action de $\Delta$ introduite dans (\cite{ag1} 6.7 et 6.10); d'où la notation. 
L'anneau du schéma affine $\oU^\star$ correspond à l'algèbre $R_1$ dans {\em loc. cit.},
d'après (\cite{ag1} 6.8(i)). 
\end{rema}

\subsection{}\label{higgs3-RGG85}\index{10001012@$\bQ$, $E_\bQ$}
On désigne par $\bQ$ la sous-catégorie pleine de $\bP$ \eqref{higgs3-RGG8} formée des schémas affines connexes $U$ tels qu'il 
existe une carte fine et saturée $M\rightarrow \Gamma(U,\cM_X)$ pour $(U,\cM_X|U)$
induisant un isomorphisme 
\begin{equation}\label{higgs3-RGG85a}
M\stackrel{\sim}{\rightarrow} \Gamma(U,\cM_X)/\Gamma(U,\co^\times_X).
\end{equation}
Cette carte est a priori indépendante de la carte adéquate requise dans la définition des objets de $\bP$. 
On munit $\bQ$ de la topologie induite par celle de $\Et_{/X}$. 
Il résulte de (\cite{ag1} 5.17) que $\bQ$ est une sous-catégorie 
topologiquement génératrice de $\Et_{/X}$. On désigne par 
\begin{equation}\label{higgs3-RGG85b}
\pi_\bQ\colon E_\bQ\rightarrow \bQ
\end{equation}
le site fibré déduit de $\pi$ \eqref{higgs3-RGG2d}
par changement de base par le foncteur d'injection canonique $\bQ\rightarrow \Et_{/X}$. 
Le foncteur de projection canonique $E_\bQ\rightarrow E$ est pleinement fidèle et 
la catégorie $E_\bQ$ est $\mU$-petite et topologiquement génératrice du site $E$. On munit $E_\bQ$ de la 
topologie induite par celle de $E$. Par restriction, le topos $\tE$ est alors équivalent à la catégorie 
des faisceaux de $\mU$-ensembles sur $E_\bQ$ (\cite{sga4} III 4.1). 
On prendra garde qu'en général, $\bQ$ n'étant pas stable par produits fibrés, 
on ne peut pas parler de topologie co-évanescente sur $E_\bQ$ associée à $\pi_\bQ$, 
et encore moins appliquer  (\cite{ag2} 5.21 et 5.22). 

\subsection{}\label{higgs3-RGG86}
On désigne par $\hE_\bQ$ la catégorie des préfaisceaux de $\mU$-ensembles sur $E_\bQ$
et par 
\begin{equation}\label{higgs3-RGG86a}
\cP^\vee_\bQ\rightarrow \bQ^\circ
\end{equation}
la catégorie fibrée obtenue en associant à tout $U\in \ob(\bQ)$ la catégorie $(\Et_{\rf/\oU^\circ})^\wedge$ 
des préfaisceaux de $\mU$-ensembles
sur $\Et_{\rf/\oU^\circ}$, et à tout morphisme $f\colon U'\rightarrow U$ de $\bQ$ 
le foncteur 
\begin{equation}\label{higgs3-RGG86b}
\of^\circ_{\fet*}\colon (\Et_{\rf/\oU'^\circ})^\wedge\rightarrow (\Et_{\rf/\oU^\circ})^\wedge
\end{equation} 
obtenu en composant avec le foncteur image inverse $\Et_{\rf/\oU^\circ} \rightarrow \Et_{\rf/\oU'^\circ}$
par le morphisme $\of^\circ\colon \oU'^\circ\rightarrow \oU^\circ$ déduit de $f$;
autrement dit, $\cP^\vee_\bQ$ est la catégorie fibrée sur $\bQ^\circ$ déduite de la catégorie fibrée \eqref{higgs3-tfa1g} 
par changement de base par le foncteur d'injection canonique $\bQ\rightarrow \Et_{/X}$.
Pour tout $U\in \ob(\bQ)$, on note $\alpha_{U!}\colon \Et_{\rf/\oU^\circ}\rightarrow E_\bQ$ 
le foncteur canonique \eqref{higgs3-tfa1d}. D'après (\cite{sga1} VI 12; cf. aussi \cite{egr1} 1.1.2),
on a une équivalence de catégories 
\begin{eqnarray}\label{higgs3-RGG86c}
\hE_\bQ&\stackrel{\sim}{\rightarrow}& \bHom_{\bQ^\circ}(\bQ^\circ,\cP^\vee_\bQ)\\
F&\mapsto &\{U\mapsto F\circ \alpha_{U!}\}.\nonumber
\end{eqnarray}
On identifiera dans la suite $F$ à la section $\{U\mapsto F\circ \alpha_{U!}\}$ qui lui est associée par cette équivalence.

Comme $E_\bQ$ est une sous-catégorie topologiquement génératrice de $E$, le foncteur ``faisceau associé'' sur $E_\bQ$
induit un foncteur que l'on note aussi 
\begin{equation}\label{higgs3-RGG86d}
\hE_\bQ\rightarrow \tE, \ \ \ F\mapsto F^a.
\end{equation}

Soient $F=\{W\mapsto G_W\}$ $(W\in \ob(\Et_{/X}))$ un objet de $\hE$ \eqref{higgs3-tfa2a}, $F_\bQ=\{U\mapsto G_U\}$ 
$(U\in \ob(\bQ))$ l'objet de $\hE_\bQ$ obtenu en restreignant $F$ à $E_\bQ$. 
Il résulte aussitôt de (\cite{sga4} II 3.0.4) et de la définition du foncteur ``faisceau associé'' (\cite{sga4} II 3.4)
qu'on a un isomorphisme canonique de $\tE$
\begin{equation}\label{higgs3-RGG86e}
(F_\bQ)^a\stackrel{\sim}{\rightarrow} F^a.
\end{equation}

\begin{rema}\label{higgs3-RGG865}
Soit $F=\{U\mapsto F_U\}$ un préfaisceau sur $E_\bQ$. Pour chaque $U\in \ob(\bQ)$,
notons $F_U^a$ le faisceau de $\oU^\circ_\fet$ associé à $F_U$. Alors $\{U\mapsto F_U^a\}$ est un 
préfaisceau sur $E_\bQ$ et on a un morphisme canonique 
$\{U\mapsto F_U\}\rightarrow \{U\mapsto F_U^a\}$ de $\hE_\bQ$, 
induisant un isomorphisme entre les faisceaux associés. 
Cette assertion ne résulte pas directement de (\cite{ag2} 5.17)
puisque $\bQ$ n'est pas stable par produits fibrés. La preuve est toutefois similaire. 
Le seul point à vérifier est l'isomorphisme (\cite{ag2} (5.17.3)). 
Soient $G=\{W\mapsto G_W\}$ $(W\in \ob(\Et_{/X}))$ un objet de $\tE$, $G_\bQ=\{U\mapsto G_U\}$ 
$(U\in \ob(\bQ))$ l'objet de $\hE_\bQ$ obtenu en restreignant $G$ à $E_\bQ$. Pour tout $U\in \ob(\bQ)$,
$G_U$ est un faisceau de $\oU^\circ_\fet$ \eqref{higgs3-tfa2b}. Par suite, l'application 
\begin{equation}\label{higgs3-RGG865a}
\Hom_{\hE_\bQ}(\{U\mapsto F^a_U\},\{U\mapsto G_U\})\rightarrow 
\Hom_{\hE_\bQ}(\{U\mapsto F_U\},\{U\mapsto G_U\})
\end{equation}
induite par le morphisme canonique $\{U\mapsto F_U\}\rightarrow \{U\mapsto F^a_U\}$
est un isomorphisme. L'assertion s'ensuit. 
\end{rema}

\subsection{}\label{higgs3-RGG3}\index{10001015@$\fV_\ox$, $\fV_\ox(\bP)$, $\fV_\ox(\bQ)$}\index{10001016@$(U,\fp)$}
\index{10001017@$\varphi_\ox\colon \tE\rightarrow \oX'^\circ_\fet$}
Soient $(\oy\rightsquigarrow \ox)$ un point de $X_\et\gtimes_{X_\et}\oX^\circ_\et$ \eqref{higgs3-tfa41} tel que $\ox$ soit
au-dessus de $s$, $X'$ le localisé strict de $X$ en $\ox$. 
On rappelle que la donnée d'un voisinage du point de $X_\et$ associé à $\ox$ 
dans le site $\Et_{/X}$ (resp. $\bP$ \eqref{higgs3-RGG8}, resp. $\bQ$ \eqref{higgs3-RGG85})
est équivalente à la donnée d'un $X$-schéma étale $\ox$-pointé (resp. de $\bP$, resp. de $\bQ$) (\cite{sga4} IV 6.8.2). 
Ces objets forment naturellement une catégorie cofiltrante, que l'on note $\fV_\ox$ (resp. $\fV_\ox(\bP)$,
resp. $\fV_\ox(\bQ)$). Les catégories $\fV_\ox(\bP)$ et $\fV_\ox(\bQ)$ sont $\mU$-petites, et 
les foncteurs d'injection canoniques $\bQ\rightarrow \bP\rightarrow \Et_{/X}$ induisent des foncteurs pleinement 
fidèles et cofinaux $\fV_\ox(\bQ)\rightarrow \fV_\ox(\bP)\rightarrow \fV_\ox$. 
Pour tout objet $(U,\fp\colon \ox\rightarrow U)$ de $\fV_\ox$, 
on note encore $\fp\colon X'\rightarrow U$ le $X$-morphisme déduit de $\fp$ (\cite{sga4} VIII 7.3) et on pose
\begin{equation}\label{higgs3-RGG3a} 
\ofp^\circ=\fp\times_X\oX^\circ\colon \oX'^\circ \rightarrow \oU^\circ.
\end{equation}

En vertu de \ref{higgs3-sli6}, $\oX'$ est normal et strictement local (et en particulier intègre). 
Le $X$-morphisme $u\colon \oy\rightarrow X'$
définissant $(\oy\rightsquigarrow \ox)$ se relève en un $\oX^\circ$-morphisme $v\colon \oy\rightarrow \oX'^\circ$ et 
induit donc un point géométrique de $\oX'^\circ$ que l'on note aussi (abusivement) $\oy$.
Pour tout $(U,\fp)\in \ob(\fV_\ox)$, on note aussi (abusivement) $\oy$ le point géométrique $\ofp^\circ(v(\oy))$ de $\oU^\circ$. 
Comme $\oU$ est localement irréductible \eqref{higgs3-sli3}, il est la somme des schémas induits sur ses composantes irréductibles. 
Notons $\oU^\star$  la composante irréductible de $\oU$ contenant $\oy$. De même, $\oU^\circ$
est la somme des schémas induits sur ses composantes irréductibles et $\oU^{\star\circ}=\oU^\star\times_XX^\circ$ 
est la composante irréductible de $\oU^\circ$ contenant $\oy$. 
Le morphisme $\ofp^\circ\colon \oX'^\circ \rightarrow \oU^\circ$ se factorise donc à travers $\oU^{\star\circ}$.

On désigne par
\begin{equation}\label{higgs3-RGG3b} 
\varphi_\ox\colon \tE\rightarrow \oX'^\circ_\fet
\end{equation}
le foncteur canonique \eqref{higgs3-tfa6d} et par 
\begin{equation}\label{higgs3-RGG3c} 
\nu_{\oy}\colon \oX'^\circ_\fet \stackrel{\sim}{\rightarrow}\bB_{\pi_1(\oX'^\circ,\oy)}
\end{equation}
le foncteur fibre de $\oX'^\circ_\fet$ en $\oy$  \eqref{higgs3-not6c}. 
D'après \eqref{higgs3-tfa6f}, on a un isomorphisme canonique 
\begin{equation}\label{higgs3-RGG3d}
\varphi_\ox(\ocB)\stackrel{\sim}{\rightarrow}
\underset{\underset{(U,\fp)\in \fV_\ox^\circ}{\longrightarrow}}{\lim}\ (\ofp^{\circ})^*_\fet(\ocB_U),
\end{equation}
où $\ocB_U$ est le faisceau de $\oU^\circ_\fet$ défini dans  \eqref{higgs3-tfa8d}. 
Compte tenu de \eqref{higgs3-tfa113a}, on en déduit un isomorphisme de $\co_\oK$-algèbres de $\bB_{\pi_1(\oX'^\circ,\oy)}$
\begin{equation}\label{higgs3-RGG3e}
\nu_{\oy}(\varphi_\ox(\ocB))\stackrel{\sim}{\rightarrow} 
\underset{\underset{(U,\fp)\in \fV_\ox^\circ}{\longrightarrow}}{\lim}\ \oR^{\oy}_U,
\end{equation} 
où $\oR^{\oy}_U$ la $\co_\oK$-algèbre de $\bB_{\pi_1(\oU^{\star \circ},\oy)}$ définie dans \eqref{higgs3-tfa112c}. 
D'après (\cite{ag2} 10.31 et 9.9), l'anneau sous-jacent à $\nu_{\oy}(\varphi_\ox(\ocB))$ est canoniquement isomorphe à 
la fibre $\ocB_{\rho(\oy \rightsquigarrow \ox)}$.

\begin{rema}\label{higgs3-RGG30} 
Pour tout point géométrique $\ox$ de $X$ au-dessus de $s$, 
il existe un point $(\oy\rightsquigarrow \ox)$ de $X_\et\gtimes_{X_\et}\oX^\circ_\et$ \eqref{higgs3-tfa41}. 
En effet, notons $X'$ le localisé strict de $X$ en $\ox$. 
D'après \ref{higgs3-sli6}, $\oX'$ est normal et strictement local (et en particulier intègre). 
Comme $X^\circ$ est schématiquement dense dans $X$ d'après \ref{higgs3-slad4}(iv), $\oX'^\circ$ est intègre et non-vide 
(\cite{ega4} 11.10.5). Soit $v\colon \oy\rightarrow \oX'^\circ$ un point géométrique de $\oX'^\circ$.
On note encore $\oy$ le point géométrique de $\oX^\circ$ et 
$u\colon \oy\rightarrow X'$ le $X$-morphisme induits par $v$. 
On obtient ainsi un point $(\oy\rightsquigarrow \ox)$ de $X_\et\gtimes_{X_\et}\oX^\circ_\et$. 
\end{rema}

\begin{prop}\label{higgs3-RGG4}
Soient $(\oy\rightsquigarrow \ox)$ un point de $X_\et\gtimes_{X_\et}\oX^\circ_\et$ \eqref{higgs3-tfa41} tel que $\ox$ soit
au-dessus de $s$, $X'$  le localisé strict de $X$ en $\ox$. Alors~:
\begin{itemize}
\item[{\rm (i)}] La fibre $\ocB_{\rho(\oy \rightsquigarrow \ox)}$ de $\ocB$ en $\rho(\oy \rightsquigarrow \ox)$
est un anneau normal et strictement local. 
\item[{\rm (ii)}] La fibre $\hbar_*(\co_{\oX})_{\ox}$ de $\hbar_*(\co_{\oX})$ en $\ox$ 
est un anneau normal et strictement local.
\item[{\rm (iii)}]  L'homomorphisme 
\begin{equation}\label{higgs3-RGG4a}
\hbar_*(\co_{\oX})_{\ox}\rightarrow \ocB_{\rho(\oy \rightsquigarrow \ox)}
\end{equation}
induit par l'homomorphisme canonique $\sigma^{-1}(\hbar_*(\co_\oX))\rightarrow \ocB$ \eqref{higgs3-tfa111}
est injectif et local.
\end{itemize}
\end{prop}

Notons $j'\colon X'^\circ\rightarrow X'$ l'injection canonique 
et $g$, $\ogg$ et $\hbar'$ les flèches canoniques du diagramme cartésien suivant~:
\begin{equation}\label{higgs3-RGG4b}
\xymatrix{
{\oX'}\ar[r]^{\ogg}\ar[d]_{\hbar'}\ar@{}[rd]|{\Box}&{\oX}\ar[d]^{\hbar}\\
X'\ar[r]^g&X}
\end{equation}
D'après \ref{higgs3-sli6}, $\oX'$ est normal et strictement local (et en particulier intègre).
On désigne par $E'$ (resp. $\tE'$) le site (resp. topos) de 
Faltings associé au morphisme  $h'=\hbar'\circ j'_{\oX'}\colon \oX'^\circ\rightarrow X'$ \eqref{higgs3-tfa2} et 
par $\ocB'$ l'anneau de $\tE'$ associé à $\oX'$ \eqref{higgs3-tfa111}. On note 
\begin{eqnarray}
\sigma'\colon \tE' &\rightarrow& X'_\et,\label{higgs3-RGG4c}\\
\rho'\colon X'_\et\gtimes_{X'_\et}\oX'^\circ_\et&\rightarrow &\tE',\label{higgs3-RGG4d}
\end{eqnarray}
les morphismes canoniques \eqref{higgs3-tfa3c} et \eqref{higgs3-tfa41b}, respectivement, et
\begin{equation}\label{higgs3-RGG4e}
\Phi\colon (\tE',\ocB')\rightarrow (\tE,\ocB)
\end{equation}
le morphisme de topos annelés déduit de $g$ par fonctorialité \eqref{higgs3-tfa13}.
Le $X$-morphisme $u\colon \oy\rightarrow X'$
définissant $(\oy\rightsquigarrow \ox)$ induit un $X'$-morphisme $v\colon \oy\rightarrow \oX'^\circ$.
On note encore (abusivement) $\ox$ le point fermé de $X'$, $\oy$ le point géométrique de $\oX'^\circ$ défini par $v$,
et $(\oy\rightsquigarrow \ox)$ le point de $X'_\et\gtimes_{X'_\et}\oX'^\circ_\fet$ défini par $u$.
Les points $\rho(\oy \rightsquigarrow \ox)$ et $\Phi(\rho'(\oy\rightsquigarrow \ox))$ de $\tE$
sont alors canoniquement isomorphes (\cite{ag2} 10.17). 

(i) L'homomorphisme canonique $\Phi^{-1}(\ocB)\rightarrow \ocB'$ étant un isomorphisme d'après \ref{higgs3-tfa16}, 
il induit un isomorphisme 
\begin{equation}\label{higgs3-RGG4f}
\ocB_{\rho(\oy \rightsquigarrow \ox)}\stackrel{\sim}{\rightarrow}
\ocB'_{\rho'(\oy \rightsquigarrow \ox)}.
\end{equation}
La proposition résulte alors de \ref{higgs3-tfa18}(i).

(ii) D'après \ref{higgs3-sli7}, le morphisme canonique $\ogg^{-1}(\co_{\oX})\rightarrow \co_{\oX'}$  
est un isomorphisme de $\oX'_\et$. 
On en déduit par (\cite{sga4} VIII 5.2) un isomorphisme canonique 
\begin{equation}\label{higgs3-RGG4g}
\hbar_*(\co_{\oX})_{\ox}\stackrel{\sim}{\rightarrow}\Gamma(\oX',\co_{\oX'}).
\end{equation}
La proposition s'ensuit en vertu de \ref{higgs3-sli6}. 

(iii) Le diagramme de morphismes de topos 
\begin{equation}\label{higgs3-RGG4h}
\xymatrix{
{\tE'}\ar[r]^\Phi\ar[d]_{\sigma'}&{\tE}\ar[d]^{\sigma}\\
{X'_\et}\ar[r]^g&{X_\et}}
\end{equation}
est commutatif à isomorphisme canonique près (\cite{ag2} (10.12.6)). De plus, le diagramme 
\begin{equation}\label{higgs3-RGG4i}
\xymatrix{
{\sigma'^{-1}(g^{-1}(\hbar_*(\co_\oX)))}\ar@{=}[r]\ar[d]_c&{\Phi^{-1}(\sigma^{-1}(\hbar_*\co_\oX))}\ar[r]&
{\Phi^{-1}(\ocB)}\ar[d]^b\\
{\sigma'^{-1}(\hbar'_*(\ogg^{-1}(\co_{\oX})))}\ar[r]^-(0.5)a&{\sigma'^{-1}(\hbar'_*(\co_{\oX'}))}\ar[r]&{\ocB'}}
\end{equation}
où $c$ est le morphisme de changement de base relativement au diagramme \eqref{higgs3-RGG4b}
et les autres flèches sont les morphismes canoniques, est commutatif. Comme $\hbar$ est entier, $c$ est un 
isomorphisme (\cite{sga4} VIII 5.6). Par ailleurs, $a$ est un isomorphisme \eqref{higgs3-sli7}
et $b$ est un isomorphisme \eqref{higgs3-tfa16}. 
La proposition résulte alors de \ref{higgs3-tfa18}(iii).

\begin{cor}\label{higgs3-RGG5}
\begin{itemize}
\item[{\rm (i)}] Le topos $\tE_s$ est localement annelé par $\ocB_s=\ocB|\tE_s$. 
\item[{\rm (ii)}] Pour tout entier $n\geq 1$, $\sigma_n\colon (\tE_s,\ocB_n)\rightarrow (X_{s,\et},\co_{\oX_n})$ 
\eqref{higgs3-RGG2g} est un morphisme de topos localement annelés.
\end{itemize} 
\end{cor} 

Cela résulte de \ref{higgs3-TFT5}, \ref{higgs3-RGG4} et  (\cite{sga4} IV 13.9).

\begin{prop}\label{higgs3-RGG7}
L'endomorphisme de Frobenius absolu de $\ocB_1$ est surjectif.
\end{prop}
Soient $(\oy\rightsquigarrow \ox)$ un point de $X_\et\gtimes_{X_\et}\oX^\circ_\et$ \eqref{higgs3-tfa41} tel que $\ox$ soit
au-dessus de $s$ et que $\oy$ soit un point géométrique générique de $\oX^\circ$, $X'$ le localisé strict de $X$ en $\ox$. 
Avec les notations de \ref{higgs3-RGG3}, on a un isomorphisme canonique \eqref{higgs3-RGG3e}
\begin{equation}\label{higgs3-RGG7a}
\nu_\oy(\varphi_\ox(\ocB_1))\stackrel{\sim}{\rightarrow}\underset{\underset{(U,\fp)\in \fV_\ox(\bP)^\circ}{\longrightarrow}}{\lim}\ 
\oR^{\oy}_{U}/p\oR^{\oy}_{U}.
\end{equation}
Par fonctorialité de l'isomorphisme \eqref{higgs3-tfa6f}, 
celui-ci est compatible aux endomorphismes de Frobenius absolus de $\ocB_1$ et 
$\oR^{\oy}_U/p\oR^{\oy}_U$. 
Pour tout $(U,\fp)\in \ob(\fV_\ox(\bP))$, l'endomorphisme de Frobenius absolu de 
$\oR^{\oy}_U/p\oR^{\oy}_U$ est surjectif en vertu de \ref{higgs3-RGG65} et (\cite{ag1} 9.10).
La proposition s'ensuit en vertu de  \ref{higgs3-TFT5} et \ref{higgs3-TFT7}.

\subsection{}\label{higgs3-RGG10}\index{10001030@$\oY^\oy$, $\hoY^\oy$}\index{10001032@$(\cA_2(\oY^\oy),\cM_{\cA_2(\oY^\oy)})$}
\index{10001035@$\cL^\oy_Y$, $\cF^\oy_Y$, $\cC^\oy_Y$}
Soient $Y$ un objet de $\bQ$ \eqref{higgs3-RGG85} tel que $Y_s\not=\emptyset$, 
$\oy$ un point géométrique de $\oY^\circ$. 
Comme $\oY$ est localement irréductible \eqref{higgs3-sli3}, 
il est la somme des schémas induits sur ses composantes irréductibles. On note  
$\oY^\star$ la composante irréductible de $\oY$ contenant $\oy$. De même, $\oY^\circ$
est la somme des schémas induits sur ses composantes irréductibles, et $\oY^{\star \circ}=\oY^\star\times_XX^\circ$ 
est la composante irréductible de $\oY^\circ$ contenant $\oy$.
On désigne par $\oR^\oy_Y$ l'anneau défini dans 
\eqref{higgs3-tfa112c} et par $\hoR^\oy_Y$ son séparé complété $p$-adique. On pose \eqref{higgs3-not5a}
\begin{equation}\label{higgs3-RGG10a}
\cR_{\oR^\oy_Y}=\underset{\underset{x\mapsto x^p}{\longleftarrow}}{\lim}\oR^\oy_Y/p\oR^\oy_Y,
\end{equation} 
et on note $\theta_Y\colon \rW(\cR_{\oR^\oy_Y})\rightarrow \hoR^\oy_Y$
l'homomorphisme de Fontaine défini dans (\cite{ag1} (9.3.4)). On pose 
\begin{equation}\label{higgs3-RGG10c}
\cA_2(\oR^\oy_Y)=\rW(\cR_{\oR^\oy_Y})/\ker(\theta_Y)^2,
\end{equation}
et on note encore $\theta_Y\colon \cA_2(\oR^\oy_Y)\rightarrow \hoR^\oy_Y$ l'homomorphisme induit par 
$\theta_Y$ (\cite{ag1} (9.3.5)). On pose enfin
\begin{eqnarray}
\oY^\oy&=&\Spec(\oR^\oy_Y),\label{higgs3-RGG10cd}\\
\hoY^\oy&=&\Spec(\hoR^\oy_Y),\label{higgs3-RGG10d}\\
\cA_2(\oY^\oy)&=&\Spec(\cA_2(\oR^\oy_Y)).\label{higgs3-RGG10e}
\end{eqnarray}
On observera que $\oY$ étant affine, $\oY^\oy$ n'est autre que le schéma défini dans \ref{higgs3-tfa114}. 
On munit $\oY^\oy$ (resp. $\hoY^\oy$) de la structure logarithmique $\cM_{\oY^\oy}$ (resp. $\cM_{\hoY^\oy}$) 
image inverse de $\cM_X$ (\cite{ag1} 5.10) et $\cA_2(\oY^\oy)$ de la structure logarithmique $\cM_{\cA_2(\oY^\oy)}$
définie comme suit. Soient $Q_Y$ le monoïde et $q_Y\colon Q_Y\rightarrow \rW(\cR_{\oR^\oy_Y})$
l'homomorphisme définis dans (\cite{ag1} 9.6) (notés $Q$ et $q$ dans {\em loc. cit.}) 
en prenant pour $u$ l'homomorphisme canonique 
$\Gamma(Y,\cM_X)\rightarrow \Gamma(\oY^\oy,\cM_{\oY^\oy})$. 
On désigne par $\cM_{\cA_2(\oY^\oy)}$ la structure logarithmique 
sur $\cA_2(\oY^\oy)$ associée à la structure pré-logarithmique définie par l'homomorphisme 
$Q_Y\rightarrow \cA_2(\oR^\oy_Y)$ induit par $q_Y$. 
L'homomorphisme $\theta_Y$ induit alors un morphisme (\cite{ag1} (9.6.4)) 
\begin{equation}\label{higgs3-RGG10f}
i_Y\colon (\hoY^\oy,\cM_{\hoY^\oy})\rightarrow (\cA_2(\oY^\oy),\cM_{\cA_2(\oY^\oy)}).
\end{equation}
Le schéma logarithmique $(\cA_2(\oY^\oy),\cM_{\cA_2(\oY^\oy)})$ 
est fin et saturé et $i_Y$ est une immersion fermée exacte. En effet, tous les foncteurs fibre de $\oY^{\star \circ}_\fet$
étant isomorphes, il suffit de montrer cette assertion dans le cas où $\oy$ est localisé en un point générique de $\oY$. 
Compte tenu de \ref{higgs3-RGG65}, les notations ci-dessus correspondent alors
à celles introduites dans (\cite{ag1} 9.11), à l'exception de $\cM_{\cA_2(\oY^\oy)}$ qui 
correspond plutôt à la structure logarithmique $\cM'_{\cA_2(\oY^\oy)}$ dans {\em loc. cit.}
Mais comme $Y$ est un objet de $\bQ$, cette dernière est canoniquement 
isomorphe à la structure logarithmique $\cM_{\cA_2(\oY^\oy)}$ introduite dans (\cite{ag1} 9.12) en vertu de (\cite{ag1} 9.13); 
d'où l'assertion (et la notation). 

On pose 
\begin{equation}\label{higgs3-RGG10h}
\rT^\oy_Y=\Hom_{\hoR^\oy_Y}(\tOmega^1_{X/S}(Y)\otimes_{\co_X(Y)}\hoR^\oy_Y,\xi\hoR^\oy_Y)
\end{equation} 
et on identifie le $\hoR^\oy_Y$-module dual à $\xi^{-1}\tOmega^1_{X/S}(Y)\otimes_{\co_X(Y)}\hoR^\oy_Y$ (cf. \ref{higgs3-not4}).
On désigne par $\hoY^\oy_\zar$ le topos de Zariski de $\hoY^\oy$, par $\trT^\oy_Y$ le $\co_{\hoY^\oy}$-module
associé à $\rT^\oy_Y$ et par $\bT^\oy_Y$ le $\hoY^\oy$-fibré vectoriel associé à son dual, autrement dit,  \eqref{higgs3-not9}
\begin{equation}
\bT^\oy_Y=\Spec(\rS_{\hoR^\oy_Y}(\xi^{-1}\tOmega^1_{X/S}(Y)\otimes_{\co_X(Y)}\hoR^\oy_Y)).
\end{equation}

Soient $U$ un ouvert de Zariski de $\hoY^\oy$, $\tU$ l'ouvert de $\cA_2(\oY^\oy)$ défini par $U$. 
On note $\cL_Y^\oy(U)$ l'ensemble des morphismes représentés par des flèches pointillées qui complètent  le diagramme 
\begin{equation}\label{higgs3-RGG10g}
\xymatrix{
{(U,\cM_{\hoY^\oy}|U)}\ar[r]^-(0.5){i_Y|U}\ar[d]&{(\tU,\cM_{\cA_2(\oY^\oy)}|\tU)}\ar@{.>}[d]\ar@/^2pc/[dd]\\
{(\coX,\cM_{\coX})}\ar[r]\ar[d]&{(\tX,\cM_\tX)}\ar[d]\\
{(\coS,\cM_\coS)}\ar[r]^-(0.5){i_\oS}&{(\cA_2(\oS),\cM_{\cA_2(\oS)})}}
\end{equation}
de façon à le laisser commutatif. D'après (\cite{ag1} 5.23), le foncteur $U\mapsto \cL_Y^\oy(U)$ 
est un $\trT^\oy_Y$-torseur de $\hoY^\oy_\zar$. 
On désigne par $\cF^\oy_Y$ le $\hoR^\oy_Y$-module des fonctions affines sur $\cL^\oy_Y$ (cf. \cite{ag1} 4.9). 
Celui-ci s'insère dans une suite exacte canonique (\cite{ag1} (4.9.1))
\begin{equation}\label{higgs3-RGG10m}
0\rightarrow \hoR^\oy_Y\rightarrow \cF^\oy_Y\rightarrow \xi^{-1}\tOmega^1_{X/S}(Y)\otimes_{\co_X(Y)}
\hoR^\oy_Y \rightarrow 0. 
\end{equation}
Cette suite induit pour tout entier $m\geq 1$, une suite exacte \eqref{higgs3-vck2b}
\begin{equation}
0\rightarrow \rS^{m-1}_{\hoR^\oy_Y}(\cF^\oy_Y)\rightarrow \rS^m_{\hoR^\oy_Y}(\cF^\oy_Y)\rightarrow 
\rS^m_{\hoR^\oy_Y}(\xi^{-1}\tOmega^1_{X/S}(Y)\otimes_{\co_X(Y)} \hoR^\oy_Y)\rightarrow 0.
\end{equation}
Les $\hoR^\oy_Y$-modules $(\rS^m_{\hoR^\oy_Y}(\cF^\oy_Y))_{m\in \mN}$ 
forment donc un système inductif filtrant, dont la limite inductive 
\begin{equation}\label{higgs3-RGG10n}
\cC^\oy_Y=\underset{\underset{m\geq 0}{\longrightarrow}}\lim\ \rS^m_{\hoR^\oy_Y}(\cF^\oy_Y)
\end{equation}
est naturellement munie d'une structure de $\hoR^\oy_Y$-algèbre. D'après (\cite{ag1} 4.10), le $\hoY^\oy$-schéma 
\begin{equation}
\bL^\oy_Y=\Spec(\cC^\oy_Y)
\end{equation}
est naturellement un $\bT^\oy_Y$-fibré principal homogène sur $\hoY^\oy$ qui représente canoniquement $\cL^\oy_Y$. 
On notera que $\cL^\oy_Y$, $\cF^\oy_Y$, $\cC^\oy_Y$ et $\bL^\oy_Y$ dépendent du choix de la déformation 
$(\tX,\cM_{\tX})$ fixée dans \ref{higgs3-RGG1}.

Le groupe $\pi_1(\oY^{\star \circ},\oy)$ agit naturellement à gauche sur les schémas
logarithmiques $(\hoY^\oy,\cM_{\hoY^\oy})$ et $(\cA_2(\oY^\oy),\cM_{\cA_2(\oY^\oy)})$, et le morphisme $i_Y$ 
est $\pi_1(\oY^{\star \circ},\oy)$-équivariant (cf. \cite{ag1} 9.11).
Procédant comme dans (\cite{ag1} 10.4), on munit $\trT^\oy_Y$ d'une structure canonique de $\co_{\hoY^\oy}$-module 
$\pi_1(\oY^{\star \circ},\oy)$-équivariant et $\cL^\oy_Y$ d'une structure canonique de $\trT^\oy_Y$-torseur 
$\pi_1(\oY^{\star \circ},\oy)$-équivariant (cf. \cite{ag1} 4.18). D'après (\cite{ag1} 4.21), 
ces deux structures induisent une structure $\pi_1(\oY^{\star \circ},\oy)$-équivariante 
sur le $\co_{\hoY^\oy}$-module associé à $\cF^\oy_Y$, 
ou, ce qui revient au même, une action $\hoR^\oy_Y$-semi-linéaire de $\pi_1(\oY^{\star \circ},\oy)$ sur $\cF^\oy_Y$, 
telle que les morphismes de la suite \eqref{higgs3-RGG10m} soient $\pi_1(\oY^{\star \circ},\oy)$-équivariants. 
On en déduit une action de $\pi_1(\oY^{\star \circ},\oy)$ sur $ \cC^\oy_Y$ 
par des automorphismes d'anneaux, compatible avec son action sur $\hoR^\oy_Y$.

\begin{lem}\label{higgs3-RGG105}
Sous les hypothèses de \eqref{higgs3-RGG10}, les actions de $\pi_1(\oY^{\star \circ},\oy)$ sur $\cF^\oy_Y$
et sur $\cC^\oy_Y$ sont continues pour les topologies $p$-adiques.
\end{lem}
En effet, compte tenu de \ref{higgs3-tfa113} et du fait que tous les foncteurs fibres de $\oY^{\star \circ}_\fet$
sont isomorphes \eqref{higgs3-not6c}, on peut supposer que $\oy$
est localisé en un point générique de $\oY$. 
Notons $\tY\rightarrow \tX$ l'unique morphisme étale 
qui relève $\coY\rightarrow \coX$ (\cite{ega4} 18.1.2)
et $\cM_{\coY}$ (resp. $\cM_\tY$) la structure logarithmique sur $\coY$
(resp. $\tY$) image inverse de $\cM_\coX$ (resp. $\cM_\tX$). 
Soient $U$ un ouvert de Zariski de $\hoY^\oy$, $\tU$ l'ouvert de $\cA_2(\oY^\oy)$ défini par $U$. 
L'ensemble $\cL_Y^\oy(U)$ est alors canoniquement isomorphe à 
l'ensemble des morphismes représentés par des flèches pointillées qui complètent  le diagramme 
\begin{equation}\label{higgs3-RGG105a}
\xymatrix{
{(U,\cM_{\hoY^\oy}|U)}\ar[r]^-(0.5){i_Y|U}\ar[d]&{(\tU,\cM_{\cA_2(\oY^\oy)}|\tU)}\ar@{.>}[d]\ar@/^2pc/[dd]\\
{(\coY,\cM_{\coY})}\ar[r]\ar[d]&{(\tY,\cM_\tY)}\ar[d]\\
{(\coS,\cM_\coS)}\ar[r]^-(0.5){i_\oS}&{(\cA_2(\oS),\cM_{\cA_2(\oS)})}}
\end{equation}
de façon à le laisser commutatif. La $\hoR^\oy_Y$-algèbre $ \cC^\oy_Y$ munie de l'action de 
$\pi_1(\oY^{\star \circ},\oy)$ s'identifie donc à l'algèbre de Higgs-Tate associée à $(Y,\cM_Y,\tY,\cM_\tY)$ 
définie dans (\cite{ag1} 10.5) (cf. \ref{higgs3-RGG65}). La proposition s'ensuit en vertu de (\cite{ag1} 12.4).

\subsection{}\label{higgs3-RGG11}\index{10001040@$\cF_{Y,n}$, $\cC_{Y,n}$ ($Y\in \ob(\bQ)$, $n\in \mN$)}
Soient $Y$ un objet de $\bQ$ \eqref{higgs3-RGG8} tel que $Y_s\not=\emptyset$, $n$ un entier $\geq 0$. 
Si $A$ est un anneau et $M$ un $A$-module, on note encore $A$ (resp. $M$) le faisceau constant 
de valeur $A$ (resp. $M$) de $\oY^\circ_\fet$. 
Le schéma $\oY^\circ$ étant localement irréductible, 
il est la somme des schémas induits sur ses composantes irréductibles.
Soient $W$ une composante irréductible de $\oY^\circ$, $\Pi(W)$ 
son groupoïde fondamental (\cite{ag2} 9.10). 
Compte tenu de \ref{higgs3-tfa113} et (\cite{ag2} 9.11), le faisceau $\ocB_Y|W$ de $W_\fet$ définit un foncteur 
\begin{equation}\label{higgs3-RGG11a}
\Pi(W)\rightarrow  \Ens, \ \ \ \oy\mapsto \oR^\oy_Y.
\end{equation}
On en déduit un foncteur 
\begin{equation}\label{higgs3-RGG11c}
\Pi(W)\rightarrow \Ens, \ \ \ \oy\mapsto \cF^\oy_Y/p^n\cF^\oy_Y.
\end{equation}
D'après \ref{higgs3-RGG105}, pour tout point géométrique $\oy$ de $W$, 
$\cF^\oy_Y/p^n\cF^\oy_Y$ est une représentation discrète et continue de $\pi_1(W,\oy)$.
Par suite, en vertu de (\cite{ag2} 9.11), le foncteur \eqref{higgs3-RGG11c} définit un $(\ocB_{Y,n}|W)$-module 
$\cF_{W,n}$ de $W_\fet$, unique à isomorphisme canonique près, où $\ocB_{Y,n}=\ocB_Y/p^n\ocB_Y$ \eqref{higgs3-RGG2ff}. 
Par descente (\cite{giraud2} II 3.4.4), il existe un $\ocB_{Y,n}$-module $\cF_{Y,n}$ de $\oY^\circ_\fet$, unique 
à isomorphisme canonique près, tel que 
pour toute composante irréductible $W$ de $\oY^\circ$, on ait $\cF_{Y,n}|W=\cF_{W,n}$.

La suite exacte \eqref{higgs3-RGG10m} induit une suite exacte de 
$\ocB_{Y,n}$-modules 
\begin{equation}\label{higgs3-RGG14a}
0\rightarrow \ocB_{Y,n}\rightarrow \cF_{Y,n}\rightarrow 
\xi^{-1}\tOmega^1_{X/S}(Y)\otimes_{\co_X(Y)}\ocB_{Y,n} \rightarrow 0. 
\end{equation}
Celle-ci induit pour tout entier $m\geq 1$, une suite exacte \eqref{higgs3-vck2b}
\[
0\rightarrow \rS^{m-1}_{\ocB_{Y,n}}(\cF_{Y,n})\rightarrow \rS^{m}_{\ocB_{Y,n}}(\cF_{Y,n})\rightarrow 
\rS^m_{\ocB_{Y,n}}(\xi^{-1}\tOmega^1_{X/S}(Y)\otimes_{\co_X(Y)} \ocB_{Y,n})\rightarrow 0.
\]
Les $\ocB_{Y,n}$-modules $(\rS^{m}_{\ocB_{Y,n}}(\cF_{Y,n}))_{m\in \mN}$ forment donc un système inductif 
dont la limite inductive 
\begin{equation}\label{higgs3-RGG14b}
 \cC_{Y,n}=\underset{\underset{m\geq 0}{\longrightarrow}}\lim\ \rS^m_{\ocB_{Y,n}}(\cF_{Y,n})
\end{equation}
est naturellement munie d'une structure de $\ocB_{Y,n}$-algèbre de $\oY^\circ_\fet$. 
On notera que $\cF_{Y,n}$ et $\cC_{Y,n}$ dépendent du choix de la déformation 
$(\tX,\cM_{\tX})$ fixée dans \ref{higgs3-RGG1}. 

\subsection{}\label{higgs3-RGG12}
Soient $g\colon Y\rightarrow Z$ un morphisme de $\bQ$ tel que $Y_s\not=\emptyset$,
$\oy$ un point géométrique de $\oY^\circ$, $\oz=\ogg(\oy)$. 
On rappelle que $\oY$ et $\oZ$ sont sommes des schémas induits sur leurs composantes irréductibles \eqref{higgs3-sli3}. 
On note $\oY^\star$ la composante irréductible de $\oY$ contenant $\oy$ et 
$\oZ^\star$ la composante irréductible de $\oZ$ contenant $\oz$, de sorte que $\ogg(\oY^\star)\subset \oZ^\star$.
Reprenons les notations de \ref{higgs3-RGG10} pour $Y$ et pour $Z$.  
Le morphisme $\ogg^\circ\colon \oY^\circ\rightarrow \oZ^\circ$ 
induit un homomorphisme de groupes $\pi_1(\oY^{\star \circ},\oy) \rightarrow\pi_1(\oZ^{\star \circ},\oz)$. 
Le morphisme canonique $(\ogg^\circ)^*_\fet(\ocB_Z)\rightarrow \ocB_Y$ 
induit un homomorphisme d'anneaux $\pi_1(\oY^{\star \circ},\oy)$-équivariant \eqref{higgs3-tfa113a}
\begin{equation}\label{higgs3-RGG12a}
\oR^\oz_Z\rightarrow \oR^{\oy}_Y, 
\end{equation} 
et par suite un morphisme $\pi_1(\oY^{\star \circ},\oy)$-équivariant de schémas 
$h\colon \hoY^{\oy}\rightarrow \hoZ^\oz$. Comme $g$ est étale, on a un morphisme canonique 
$\co_{\hoZ^\oz}$-linéaire et $\pi_1(\oY^{\star \circ},\oy)$-équivariant $u\colon \trT^\oz_Z\rightarrow h_*(\trT^\oy_Y)$
tel que le morphisme adjoint $u^\sharp\colon h^*(\trT^\oz_Z)\rightarrow \trT^\oy_Y$ soit un isomorphisme.
Il résulte aussitôt des définitions \eqref{higgs3-RGG10g}
qu'on a un morphisme canonique $u$-équivariant et $\pi_1(\oY^{\star \circ},\oy)$-équivariant
\begin{equation}\label{higgs3-RGG12b}
v\colon \cL_Z^\oz\rightarrow h_*(\cL_Y^\oy).
\end{equation}
D'après (\cite{ag1} 4.22), le couple $(u,v)$ induit un isomorphisme $\hoR^\oz_Y$-linéaire et 
$\pi_1(\oY^{\star \circ},\oy)$-équivariant 
\begin{equation}\label{higgs3-RGG12bb}
\cF^\oy_Y\stackrel{\sim}{\rightarrow} \cF^\oz_Z\otimes_{\hoR^\oz_Z}\hoR^\oy_Y,
\end{equation}
et par suite, un morphisme $\hoR^\oz_Z$-linéaire et $\pi_1(\oY^{\star \circ},\oy)$-équivariant 
\begin{equation}\label{higgs3-RGG12d}
\cF^\oz_Z\rightarrow \cF^\oy_Y
\end{equation}
qui s'insère dans un diagramme commutatif 
\begin{equation}\label{higgs3-RGG12e}
\xymatrix{
0\ar[r]&{\hoR^\oz_Z}\ar[r]\ar[d]&{\cF^\oz_Z}\ar[r]\ar[d]&
{\xi^{-1}\tOmega^1_{X/S}(Z)\otimes_{\co_X(Z)}\hoR^\oz_Z}\ar[r]\ar[d]&0\\
0\ar[r]&{\hoR^\oy_Y}\ar[r]&{\cF^\oy_Y}\ar[r]&
{\xi^{-1}\tOmega^1_{X/S}(Y)\otimes_{\co_X(Y)}\hoR^\oy_Y}\ar[r]&0}
\end{equation}
On en déduit un homomorphisme $\pi_1(\oY^{\star \circ},\oy)$-équivariant de $\hoR^\oz_Z$-algèbres
\begin{equation}\label{higgs3-RGG12c}
\cC^{\oz}_Z\rightarrow\cC^\oy_Y.
\end{equation}

On désigne par $\Pi(\oY^{\star \circ})$ et $\Pi(\oZ^{\star \circ})$ les groupoïdes fondamentaux de $\oY^{\star \circ}$ et 
$\oZ^{\star \circ}$ et par
\begin{equation}\label{higgs3-RGG12f}
\gamma\colon \Pi(\oY^{\star \circ})\rightarrow \Pi(\oZ^{\star \circ})
\end{equation}
le foncteur induit par le foncteur image inverse $\Et_{\rf/\oZ^{\star \circ}}\rightarrow \Et_{\rf/\oY^{\star \circ}}$.
Pour tout entier $n\geq 0$, on note $F_{Y,n}\colon \Pi(\oY^{\star \circ})\rightarrow \Ens$ et 
$F_{Z,n}\colon \Pi(\oZ^{\star \circ})\rightarrow \Ens$ les foncteurs associés par 
(\cite{ag2} 9.11) aux objets $\cF_{Y,n}|\oY^{\star \circ}$ de $\oY^{\star \circ}_\fet$ et 
$\cF_{Z,n}|\oZ^{\star \circ}$ de $\oZ^{\star \circ}_\fet$, respectivement. 
Le morphisme \eqref{higgs3-RGG12d} induit clairement un morphisme de foncteurs
\begin{equation}\label{higgs3-RGG12g}
F_{Z,n}\circ \gamma \rightarrow F_{Y,n}. 
\end{equation} 
On en déduit par (\cite{ag2} 9.11) un morphisme $(\ogg^{\circ})_\fet^*(\ocB_{Z,n})$-linéaire 
\begin{equation}\label{higgs3-RGG12ii}
(\ogg^{\circ})^*_\fet(\cF_{Z,n}) \rightarrow\cF_{Y,n},
\end{equation} 
et donc par adjonction, un morphisme $\ocB_{Z,n}$-linéaire 
\begin{equation}\label{higgs3-RGG12i}
\cF_{Z,n} \rightarrow \ogg^{\circ}_{\fet *}(\cF_{Y,n}). 
\end{equation} 
Il résulte de \eqref{higgs3-RGG12e} que le diagramme 
\begin{equation}\label{higgs3-RGG12ih}
\xymatrix{
0\ar[r]&{(\ogg^{\circ})^*_\fet(\ocB_{Z,n})}\ar[r]\ar[d]&{(\ogg^{\circ})^*_\fet(\cF_{Z,n})}\ar[r]\ar[d]&
{\xi^{-1}\tOmega^1_{X/S}(Z)\times_{\co_X(Z)}(\ogg^{\circ})^*_\fet(\ocB_{Z,n})}\ar[r]\ar[d]&0\\
0\ar[r]&{\ocB_{Y,n}}\ar[r]&{\cF_{Y,n}}\ar[r]&
{\xi^{-1}\tOmega^1_{X/S}(Y)\times_{\co_X(Y)}\ocB_{Y,n}}\ar[r]&0}
\end{equation}
est commutatif. 
On en déduit un homomorphisme de $(\ogg^{\circ})_\fet^*(\ocB_{Z,n})$-algèbres 
\begin{equation}\label{higgs3-RGG12hh}
(\ogg^{\circ})^*_\fet(\cC_{Z,n}) \rightarrow\cC_{Y,n},
\end{equation} 
et donc par adjonction un homomorphisme de $\ocB_{Z,n}$-algèbres
\begin{equation}\label{higgs3-RGG12h}
\cC_{Z,n} \rightarrow\ogg^{\circ}_{\fet *}(\cC_{Y,n}). 
\end{equation}

\subsection{}
Pour tout entier $n\geq 0$ et tout objet $Y$ de $\bQ$ tel que $Y_s=\emptyset$, on pose $\cC_{Y,n}=\cF_{Y,n}=0$.  
La suite exacte \eqref{higgs3-RGG14a} vaut encore dans ce cas, puisque $\ocB_Y$ est une $\oK$-algèbre.
Les morphismes \eqref{higgs3-RGG12i} et \eqref{higgs3-RGG12h} sont alors définis pour tout morphisme 
de $\bQ$, et ils vérifient des relations de cocycles du type (\cite{egr1} (1.1.2.2)).

\subsection{}\label{higgs3-RGG22}\index{10001045@$\cF^{(r)}_{Y,n}$, $\cC^{(r)}_{Y,n}$ ($Y\in \ob(\bQ)$, $r\in \mQ_{\geq 0}$, $n\in \mN$)}
\index{10001047@$\tta_{Y,n}^{r,r'}$, $\alpha_{Y,n}^{r,r'}$}
Soient $r$ un nombre rationnel $\geq 0$, $n$ un entier $\geq 0$, $Y$ un objet de $\bQ$.   
On désigne par $\cF^{(r)}_{Y,n}$ l'extension de $\ocB_{Y,n}$-modules de $\oY^\circ_\fet$ déduite de 
$\cF_{Y,n}$ \eqref{higgs3-RGG14a} par image inverse 
par le morphisme de multiplication par $p^r$ sur $\xi^{-1}\tOmega^1_{X/S}(Y)\otimes_{\co_X(Y)}\ocB_{Y,n}$,
de sorte qu'on a une suite exacte canonique de $\ocB_{Y,n}$-modules
\begin{equation}\label{higgs3-RGG22a}
0\rightarrow \ocB_{Y,n}\rightarrow \cF^{(r)}_{Y,n}\rightarrow 
\xi^{-1}\tOmega^1_{X/S}(Y)\otimes_{\co_X(Y)}\ocB_{Y,n} \rightarrow 0. 
\end{equation}
Celle-ci induit pour tout entier $m\geq 1$, une suite exacte de $\ocB_{Y,n}$-modules \eqref{higgs3-vck2b}
\[
0\rightarrow \rS^{m-1}_{\ocB_{Y,n}}(\cF^{(r)}_{Y,n})\rightarrow \rS^{m}_{\ocB_{Y,n}}(\cF^{(r)}_{Y,n})\rightarrow 
\rS^m_{\ocB_{Y,n}}(\xi^{-1}\tOmega^1_{X/S}(Y)\otimes_{\co_X(Y)} \ocB_{Y,n})\rightarrow 0.
\]
Les $\ocB_{Y,n}$-modules $(\rS^{m}_{\ocB_{Y,n}}(\cF^{(r)}_{Y,n}))_{m\in \mN}$ forment donc un système inductif 
dont la limite inductive 
\begin{equation}\label{higgs3-RGG22b}
 \cC^{(r)}_{Y,n}=\underset{\underset{m\geq 0}{\longrightarrow}}\lim\ \rS^m_{\ocB_{Y,n}}(\cF^{(r)}_{Y,n})
\end{equation}
est naturellement munie d'une structure de $\ocB_{Y,n}$-algèbre de $\oY^\circ_\fet$. 

Pour tous nombres rationnels $r\geq r'\geq 0$, on a un morphisme $\ocB_{Y,n}$-linéaire canonique 
\begin{equation}\label{higgs3-RGG22g}
\tta_{Y,n}^{r,r'}\colon\cF_{Y,n}^{(r)}\rightarrow \cF_{Y,n}^{(r')}
\end{equation} 
qui relève la multiplication par $p^{r-r'}$ sur  
$\xi^{-1}\tOmega^1_{X/S}(Y)\otimes_{\co_X(Y)}\ocB_{Y,n}$ et qui étend l'identité sur $\ocB_{Y,n}$ 
\eqref{higgs3-RGG22a}. Il induit un homomorphisme de $\ocB_{Y,n}$-algèbres 
\begin{equation}\label{higgs3-RGG22h}
\alpha_{Y,n}^{r,r'}\colon \cC_{Y,n}^{(r)}\rightarrow \cC_{Y,n}^{(r')}.
\end{equation}

\subsection{}\label{higgs3-RGG221}
Soient $r$ un nombre rationnel $\geq 0$, $n$ un entier $\geq 0$,
$g\colon Y\rightarrow Z$ un morphisme de $\bQ$.
Le diagramme \eqref{higgs3-RGG12ih} induit un morphisme $(\ogg^{\circ})^*_\fet(\ocB_{Z,n})$-linéaire
\begin{equation}\label{higgs3-RGG22c}
(\ogg^{\circ})^*_\fet(\cF^{(r)}_{Z,n}) \rightarrow\cF^{(r)}_{Y,n}
\end{equation} 
qui s'insère dans un diagramme commutatif
\begin{equation}\label{higgs3-RGG22cd}
\xymatrix{
0\ar[r]&{(\ogg^{\circ})^*_\fet(\ocB_{Z,n})}\ar[r]\ar[d]&{(\ogg^{\circ})^*_\fet(\cF^{(r)}_{Z,n})}\ar[r]\ar[d]&
{\xi^{-1}\tOmega^1_{X/S}(Z)\times_{\co_X(Z)}(\ogg^{\circ})^*_\fet(\ocB_{Z,n})}\ar[r]\ar[d]&0\\
0\ar[r]&{\ocB_{Y,n}}\ar[r]&{\cF^{(r)}_{Y,n}}\ar[r]&
{\xi^{-1}\tOmega^1_{X/S}(Y)\times_{\co_X(Y)}\ocB_{Y,n}}\ar[r]&0}
\end{equation}
On en déduit par adjonction un morphisme $\ocB_{Z,n}$-linéaire
\begin{equation}\label{higgs3-RGG22d}
\cF^{(r)}_{Z,n} \rightarrow(\ogg^{\circ})_{\fet*}(\cF^{(r)}_{Y,n}).
\end{equation} 
On en déduit aussi un morphisme de $(\ogg^{\circ})^*_\fet(\ocB_{Z,n})$-algèbres
\begin{equation}\label{higgs3-RGG22e}
(\ogg^{\circ})^*_\fet(\cC^{(r)}_{Z,n}) \rightarrow\cC^{(r)}_{Y,n},
\end{equation} 
et donc par adjonction un morphisme de $\ocB_{Z,n}$-algèbres
\begin{equation}\label{higgs3-RGG22f}
\cC^{(r)}_{Z,n} \rightarrow(\ogg^{\circ})_{\fet*}(\cC^{(r)}_{Y,n}). 
\end{equation} 
Les morphismes \eqref{higgs3-RGG22d} et \eqref{higgs3-RGG22f} vérifient des relations de cocycles 
du type (\cite{egr1} (1.1.2.2)).

\begin{lem}\label{higgs3-RGG23}
Pour tout nombre rationnel $r\geq 0$, tout entier $n\geq 0$ et 
tout morphisme $g\colon Y\rightarrow Z$ de $\bQ$, 
les morphismes \eqref{higgs3-RGG22c} et \eqref{higgs3-RGG22e} induisent des isomorphismes
\begin{eqnarray}
(\ogg^{\circ})^*_\fet(\cF^{(r)}_{Z,n})\otimes_{(\ogg^{\circ})^*_\fet(\ocB_{Z,n})}\ocB_{Y,n} &
\stackrel{\sim}{\rightarrow}&\cF^{(r)}_{Y,n},\label{higgs3-RGG23a}\\
(\ogg^{\circ})^*_\fet(\cC^{(r)}_{Z,n})\otimes_{(\ogg^{\circ})^*_\fet(\ocB_{Z,n})}\ocB_{Y,n} &
\stackrel{\sim}{\rightarrow}&\cC^{(r)}_{Y,n}.\label{higgs3-RGG23b}
\end{eqnarray} 
\end{lem}

En effet, le premier isomorphisme résulte du diagramme \eqref{higgs3-RGG22cd} et du fait que le morphisme canonique 
\begin{equation}
\tOmega^1_{X/S}(Z)\otimes_{\co_X(Z)}\co_X(Y)\rightarrow \tOmega^1_{X/S}(Y)
\end{equation} 
est un isomorphisme. Le second isomorphisme se déduit du premier. 

\subsection{}\label{higgs3-RGG24}\index{10001050@$\cF^{(r)}_n$, $\cC^{(r)}_n$ ($r\in \mQ_{\geq 0}$, $n\in \mN$)}
\index{10001051@$\cF_n$, $\cC_n$ ($n\in \mN$)}\index{10001052@$\tta_n^{r,r'}$, $\alpha_n^{r,r'}$}
\index{Extension de Higgs-Tate d'epaisseur1@$\ocB_n$-extension de Higgs-Tate d'épaisseur $r$ ($\cF^{(r)}_n$)}
\index{Algebre de Higgs-Tate d'epaisseur1@$\ocB_n$-algèbre de Higgs-Tate d'épaisseur $r$ ($\cC^{(r)}_n$)}
Soient $r$ un nombre rationnel $\geq 0$, $n$ un entier $\geq 0$. D'après \ref{higgs3-RGG221}, les correspondances 
\begin{equation}
\{Y\mapsto \cF^{(r)}_{Y,n}\} \ \ \ {\rm et}\ \ \ \{Y\mapsto \cC^{(r)}_{Y,n}\}, \ \ \ (Y\in \ob(\bQ)),
\end{equation} 
définissent des préfaisceaux  sur $E_\bQ$ \eqref{higgs3-RGG85b} de modules et d'algèbres, respectivement, 
relativement à l'anneau $\{Y\mapsto \ocB_{Y,n}\}$. On pose 
\begin{eqnarray}
\cF^{(r)}_n&=&\{Y\mapsto \cF^{(r)}_{Y,n}\}^a,\label{higgs3-RGG24a}\\
\cC^{(r)}_n&=&\{Y\mapsto \cC^{(r)}_{Y,n}\}^a, \label{higgs3-RGG24b}
\end{eqnarray}
les faisceaux associés dans $\tE$ \eqref{higgs3-RGG86d}. D'après \eqref{higgs3-TFT1g} et \eqref{higgs3-RGG86e},  
$\cF^{(r)}_n$ est un $\ocB_n$-module~; on l'appelle la {\em $\ocB_n$-extension de Higgs-Tate d'épaisseur $r$} 
associée à $(f,\tX,\cM_\tX)$.
De même, $\cC^{(r)}_n$ est une $\ocB_n$-algèbre~; on l'appelle la {\em $\ocB_n$-algèbre de Higgs-Tate d'épaisseur $r$} 
associée à $(f,\tX,\cM_\tX)$. On pose $\cF_n=\cF_n^{(0)}$ et $\cC_n=\cC_n^{(0)}$, et  
on les appelle la {\em $\ocB_n$-extension de Higgs-Tate} et  la {\em $\ocB_n$-algèbre de Higgs-Tate}, respectivement,
associées à $(f,\tX,\cM_\tX)$.

Pour tous nombres rationnels $r\geq r'\geq 0$, les morphismes \eqref{higgs3-RGG22g} induisent un morphisme $\ocB_n$-linéaire 
\begin{equation}\label{higgs3-RGG24c}
\tta_n^{r,r'}\colon \cF^{(r)}_n\rightarrow \cF_n^{(r')}.
\end{equation}
Les homomorphismes \eqref{higgs3-RGG22h} induisent un homomorphisme de $\ocB_n$-algèbres 
\begin{equation}\label{higgs3-RGG24d}
\alpha_n^{r,r'}\colon \cC_n^{(r)}\rightarrow \cC_n^{(r')}.
\end{equation}
Pour tous nombres rationnels $r\geq r'\geq r''\geq 0$, on a
\begin{equation}\label{higgs3-RGG18h}
\tta_n^{r,r''}=\tta_n^{r',r''} \circ \tta_n^{r,r'} \ \ \ {\rm et}\ \ \ \alpha_n^{r,r''}=\alpha_n^{r',r''} \circ \alpha_n^{r,r'}.
\end{equation}

\begin{prop}\label{higgs3-RGG18}
Soient $r$ un nombre rationnel $\geq 0$, $n$ un entier $\geq 1$. Alors~:
\begin{itemize}
\item[{\rm (i)}] Les faisceaux $\cF^{(r)}_n$ et $\cC^{(r)}_n$ sont des objets de $\tE_s$.  
\item[{\rm (ii)}] On a une suite exacte localement scindée canonique de $\ocB_n$-modules 
(\eqref{higgs3-RGG1ab} et \eqref{higgs3-RGG2g})
\begin{equation}\label{higgs3-RGG18a}
0\rightarrow \ocB_n\rightarrow \cF^{(r)}_n\rightarrow 
\sigma_n^*(\xi^{-1}\tOmega^1_{\oX_n/\oS_n})\rightarrow 0.
\end{equation}
Elle induit pour tout entier $m\geq 1$, une suite exacte de $\ocB_n$-modules \eqref{higgs3-not9}
\begin{equation}\label{higgs3-RGG18c}
0\rightarrow \rS^{m-1}_{\ocB_n}(\cF^{(r)}_n)\rightarrow \rS^m_{\ocB_n}(\cF^{(r)}_n)\rightarrow 
\sigma_n^*(\rS^m_{\co_{\oX_n}}(\xi^{-1}\tOmega^1_{\oX_n/\oS_n}))\rightarrow 0.
\end{equation}
En particulier, les $\ocB_n$-modules $(\rS^m_{\ocB_n}(\cF^{(r)}_n))_{m\in \mN}$ forment un système inductif filtrant. 
\item[{\rm (iii)}] On a un isomorphisme canonique de $\ocB_n$-algèbres  
\begin{equation}\label{higgs3-RGG18b}
\cC^{(r)}_n \stackrel{\sim}{\rightarrow}\underset{\underset{m\geq 0}{\longrightarrow}}\lim\ \rS^m_{\ocB_n}(\cF^{(r)}_n).
\end{equation}
\item[{\rm (iv)}] Pour tous nombres rationnels $r\geq r'\geq 0$, le diagramme 
\begin{equation}\label{higgs3-RGG18dd}
\xymatrix{
0\ar[r]&{\ocB_n}\ar[r]\ar@{=}[d]&
{\cF^{(r)}_n}\ar[r]\ar[d]^{\tta_n^{r,r'}}&{\sigma_n^*(\xi^{-1}\tOmega^1_{\oX_n/\oS_n})}\ar[r]\ar[d]^{\cdot p^{r-r'}}& 0\\
0\ar[r]&{\ocB_n}\ar[r]&{\cF^{(r')}_n}\ar[r]&{\sigma_n^*(\xi^{-1}\tOmega^1_{\oX_n/\oS_n})}\ar[r]& 0}
\end{equation}
où les lignes horizontales sont les suites exactes \eqref{higgs3-RGG18a}
et la flèche verticale de droite désigne la multiplication par $p^{r-r'}$, est commutatif. 
De plus, les morphismes $\tta_n^{r,r'}$ et $\alpha_n^{r,r'}$ sont compatibles avec les isomorphismes 
\eqref{higgs3-RGG18b} pour $r$ et $r'$. 
\end{itemize}
\end{prop}

(i) En effet, comme $\cF^{(r)}_{Y,n}=\cC^{(r)}_{Y,n}=0$ pour tout $Y\in \ob(\bQ)$ tel que $Y_s=\emptyset$, 
on a $\cF^{(r)}_n|\sigma^*(X_\eta)=\cC^{(r)}_n|\sigma^*(X_\eta)=0$ 
en vertu de (\cite{sga4} III 5.5).  

(ii) On a un isomorphisme canonique \eqref{higgs3-TFT8d}
\begin{equation}\label{higgs3-RGG18ee}
\sigma^*_n(\xi^{-1}\tOmega^1_{\oX_n/\oS_n})\stackrel{\sim}{\rightarrow}
\sigma^{-1}(\xi^{-1}\tOmega^1_{X/S})\otimes_{\sigma^{-1}(\co_X)}\ocB_n.
\end{equation}
Donc en vertu de (\cite{ag2} 5.34(ii), 8.9 et 5.17) et \eqref{higgs3-RGG86e}, $\sigma^*_n(\xi^{-1}\tOmega^1_{\oX_n/\oS_n})$ est 
le faisceau de $\tE$ associé au préfaisceau sur $E_\bQ$ défini par la correspondance
\begin{equation}
\{Y\mapsto \xi^{-1}\tOmega^1_{X/S}(Y)\otimes_{\co_X(Y)}\ocB_{Y,n}\}, \ \ \ (Y\in \ob(\bQ)).
\end{equation}
La suite exacte \eqref{higgs3-RGG18a} résulte alors de \eqref{higgs3-RGG22a} et \eqref{higgs3-RGG22cd} puisque le foncteur 
``faisceau associé'' est exact (\cite{sga4} II 4.1). Elle est localement scindée car le $\ocB_n$-module 
$\sigma_n^*(\xi^{-1}\tOmega^1_{\oX_n/\oS_n})$ est localement libre de type fini. La suite exacte \eqref{higgs3-RGG18c}
s'en déduit par \eqref{higgs3-vck2b}. 

(iii) On déduit facilement de \ref{higgs3-RGG865} et (\cite{sga4} IV 12.10) que pour tout entier $m\geq 0$, 
$\rS^m_{\ocB_n}(\cF^{(r)}_n)$ est le faisceau associé au préfaisceau 
\begin{equation}\label{higgs3-RGG18f}
\{Y\mapsto \rS^m_{\ocB_{Y,n}}(\cF^{(r)}_{Y,n})\}, \ \ \ (Y\in \ob(\bQ)).
\end{equation}
La proposition s'ensuit compte tenu de \ref{higgs3-RGG865} et du fait que 
le foncteur ``faisceau associé'' commute aux limites inductives (\cite{sga4} II 4.1). 

(iv) Cela résulte aussitôt des preuves de (ii) et (iii).

\subsection{}\index{10001055@$d_n^{(r)}$}
Soient $r$ un nombre rationnel $\geq 0$, $n$ un entier $\geq 1$. 
D'après \ref{higgs3-RGG18}, on a un isomorphisme canonique $\cC_n^{(r)}$-linéaire
\begin{equation}\label{higgs3-RGG18d}
\Omega^1_{\cC_n^{(r)}/\ocB_n}\stackrel{\sim}{\rightarrow} 
\sigma_n^*(\xi^{-1}\tOmega^1_{\oX_n/\oS_n})\otimes_{\ocB_n}\cC_n^{(r)}.
\end{equation}
La $\ocB_n$-dérivation universelle de $\cC_n^{(r)}$ correspond via cet isomorphisme à l'unique $\ocB_n$-dérivation 
\begin{equation}\label{higgs3-RGG18e}
d_n^{(r)}\colon \cC_n^{(r)}\rightarrow \sigma_n^*(\xi^{-1}\tOmega^1_{\oX_n/\oS_n})\otimes_{\ocB_n}\cC_n^{(r)}
\end{equation}
qui prolonge le morphisme canonique $\cF_n^{(r)}\rightarrow \sigma_n^*(\xi^{-1}\tOmega^1_{\oX_n/\oS_n})$ \eqref{higgs3-RGG18a}. 
Il résulte de \ref{higgs3-RGG18}(iv) que pour tous nombres rationnels $r\geq r'\geq 0$, on a 
\begin{equation}\label{higgs3-RGG18i}
p^{r-r'}(\id \otimes \alpha^{r,r'}_n) \circ d^{(r)}_n=d^{(r')}_n\circ \alpha^{r,r'}_n.
\end{equation}

\subsection{}\label{higgs3-RGG19}\index{10001060@$\cF_Y^{\oy,(r)}$, $\cC_Y^{\oy,(r)}$}
Soient $Y$ un objet de $\bQ$ \eqref{higgs3-RGG85} tel que $Y_s\not=\emptyset$, 
$\oy$ un point géométrique de $\oY^\circ$. Reprenons les notations de \eqref{higgs3-RGG10}.
Pour tout nombre rationnel $r\geq 0$, on désigne par $\cF_Y^{\oy,(r)}$ 
l'extension de $\hoR_Y^\oy$-modules déduite de $\cF_Y^\oy$ 
\eqref{higgs3-RGG10m} par image inverse par le morphisme de multiplication par $p^r$ sur 
$\xi^{-1}\tOmega^1_{X/S}(Y)\otimes_{\co_X(Y)} \hoR^\oy_Y$,
de sorte qu'on a une suite exacte de $\hoR^\oy_Y$-modules 
\begin{equation}\label{higgs3-RGG19a}
0\rightarrow \hoR^\oy_Y\rightarrow \cF_Y^{\oy,(r)}\rightarrow \xi^{-1}\tOmega^1_{X/S}(Y)\otimes_{\co_X(Y)} \hoR^\oy_Y
\rightarrow 0.
\end{equation}
Celle-ci induit pour tout entier $m\geq 1$, une suite exacte de $\hoR^\oy_Y$-modules \eqref{higgs3-vck2b}
\begin{equation}\label{higgs3-RGG19aa}
0\rightarrow \rS^{m-1}_{\hoR^\oy_Y}(\cF^{\oy,(r)}_Y)\rightarrow \rS^m_{\hoR^\oy_Y}(\cF^{\oy,(r)}_Y)\rightarrow 
\rS^m_{\hoR^\oy_Y}(\xi^{-1}\tOmega^1_{X/S}(Y)\otimes_{\co_X(Y)} \hoR^\oy_Y)\rightarrow 0.
\end{equation}
En particulier, les $\hoR^\oy_Y$-modules $(\rS^m_{\hoR^\oy_Y}(\cF^{\oy,(r)}_Y))_{m\in \mN}$ 
forment un système inductif filtrant dont la limite inductive 
\begin{equation}\label{higgs3-RGG19b}
\cC_Y^{\oy,(r)}= \underset{\underset{m\geq 0}{\longrightarrow}}\lim\ \rS^m_{\hoR^\oy_Y}(\cF_Y^{\oy,(r)}),
\end{equation}
est naturellement munie d'une structure de $\hoR^\oy_Y$-algèbre.

Il résulte de \ref{higgs3-RGG12} que les formations de $\cF_Y^{\oy,(r)}$ et $\cC_Y^{\oy,(r)}$ sont fonctorielles 
en la paire $(Y,\oy)$. Plus précisément, soient $g\colon Y\rightarrow Z$ un morphisme de $\bQ$, $\oz$ l'image
de $\oy$ par le morphisme $\ogg^\circ\colon \oY^\circ\rightarrow \oZ^\circ$.
Il résulte aussitôt de \ref{higgs3-RGG12} que le diagramme canonique
\begin{equation}\label{higgs3-RGG19i}
\xymatrix{
0\ar[r]&{\hoR^\oz_Z}\ar[r]\ar[d]&{\cF^{\oz,(r)}_Z}\ar[r]\ar[d]&
{\xi^{-1}\tOmega^1_{X/S}(Z)\otimes_{\co_X(Z)}\hoR^\oz_Z}\ar[r]\ar[d]&0\\
0\ar[r]&{\hoR^\oy_Y}\ar[r]&{\cF^{\oy,(r)}_Y}\ar[r]&
{\xi^{-1}\tOmega^1_{X/S}(Y)\otimes_{\co_X(Y)}\hoR^\oy_Y}\ar[r]&0}
\end{equation}
est commutatif. Le morphisme canonique 
\begin{equation}\label{higgs3-RGG19j}
\tOmega^1_{X/S}(Z)\otimes_{\co_X(Z)}\co_X(Y)\rightarrow \tOmega^1_{X/S}(Y)
\end{equation} 
étant un isomorphisme, on en déduit que les morphismes canoniques
\begin{eqnarray}
\cF^{\oz,(r)}_Z\otimes_{\hoR_Z^\oz}\hoR_Y^\oy&\rightarrow&\cF^{\oy,(r)}_Y,\label{higgs3-RGG19k}\\
\cC^{\oz,(r)}_Z\otimes_{\hoR_Z^\oz}\hoR_Y^\oy&\rightarrow&\cC^{\oy,(r)}_Y,\label{higgs3-RGG19l}
\end{eqnarray} 
sont des isomorphismes. 

Comme $\oY$ est localement irréductible \eqref{higgs3-sli3}, 
il est la somme des schémas induits sur ses composantes irréductibles. On note  
$\oY^\star$ la composante irréductible de $\oY$ contenant $\oy$. De même, $\oY^\circ$
est la somme des schémas induits sur ses composantes irréductibles, et $\oY^{\star \circ}=\oY^\star\times_XX^\circ$ 
est la composante irréductible de $\oY^\circ$ contenant $\oy$. On note $\bB_{\pi_1(\oY^{\star \circ},\oy)}$ 
le topos classifiant du groupe profini $\pi_1(\oY^{\star \circ},\oy)$ et 
\begin{equation}\label{higgs3-RGG19e}
\nu_\oy\colon \oY^{\star \circ}_\fet \stackrel{\sim}{\rightarrow} \bB_{\pi_1(\oY^{\star \circ},\oy)}
\end{equation}
le foncteur fibre en $\oy$ \eqref{higgs3-not6c}. On a alors un isomorphisme canonique d'anneaux \eqref{higgs3-tfa113a} 
\begin{equation}\label{higgs3-RGG19f}
\nu_\oy(\ocB_Y|\oY^{\star \circ})\stackrel{\sim}{\rightarrow} \oR^\oy_Y.
\end{equation}
Comme $\nu_\oy$ est exact et qu'il commute aux limites inductives, 
pour tout entier $n\geq 0$, on a des isomorphismes canoniques 
de $\oR^\oy_Y$-modules et de $\oR^\oy_Y$-algèbres, respectivement,
\begin{eqnarray}
\nu_\oy(\cF^{(r)}_{Y,n}|\oY^{\star \circ})&\stackrel{\sim}{\rightarrow}& \cF^{\oy,(r)}_Y/p^n \cF^{\oy,(r)}_Y,\label{higgs3-RGG19g}\\
\nu_\oy(\cC^{(r)}_{Y,n}|\oY^{\star \circ})&\stackrel{\sim}{\rightarrow}& \cC^{\oy,(r)}_Y/p^n \cC^{\oy,(r)}_Y.\label{higgs3-RGG19h}
\end{eqnarray}

\subsection{}\label{higgs3-RGG16}\index{10001065@$\oR^\oy_{X'}$, $\cF^{\oy,(r)}_{X'}$,  $\cC^{\oy,(r)}_{X'}$}
Soient $(\oy\rightsquigarrow \ox)$ un point de $X_\et\gtimes_{X_\et}\oX^\circ_\et$ \eqref{higgs3-tfa41} tel que $\ox$ soit
au-dessus de $s$, $X'$ le localisé strict de $X$ en $\ox$. Reprenons les notations de \ref{higgs3-RGG3}. 
On pose
\begin{equation}\label{higgs3-RGG16d}
\oR^\oy_{X'}=\underset{\underset{(U,\fp)\in \fV_\ox^\circ}{\longrightarrow}}{\lim}\ \oR^{\oy}_U,
\end{equation}
où $\oR^{\oy}_U$ est l'anneau défini dans \eqref{higgs3-tfa112c}.  
On note $\hoR^\oy_{X'}$ le complété séparé $p$-adique de $\oR^\oy_{X'}$. 
On a un isomorphisme canonique \eqref{higgs3-RGG3e}
\begin{equation}\label{higgs3-RGG16m}
\nu_{\oy}(\varphi_\ox(\ocB)) \stackrel{\sim}{\rightarrow} \oR^{\oy}_{X'}. 
\end{equation} 

Pour tout nombre rationnel $r\geq 0$, on pose 
\begin{eqnarray}
\cF^{\oy,(r)}_{X'}&=&\underset{\underset{(U,\fp)\in \fV_\ox(\bQ)^\circ}{\longrightarrow}}{\lim}\ 
\cF^{\oy,(r)}_U\otimes_{\hoR_U^\oy}\hoR_{X'}^\oy,\label{higgs3-RGG16e}\\
\cC^{\oy,(r)}_{X'}&=&\underset{\underset{(U,\fp)\in \fV_\ox(\bQ)^\circ}{\longrightarrow}}{\lim}\ 
\cC^{\oy,(r)}_U\otimes_{\hoR_U^\oy}\hoR_{X'}^\oy,\label{higgs3-RGG16f}
\end{eqnarray}
où $\cF^{\oy,(r)}_U$ est le $\hoR^{\oy}_U$-module 
défini dans \eqref{higgs3-RGG19a} et $\cC^{\oy,(r)}_U$ est la $\hoR^{\oy}_U$-algèbre définie dans \eqref{higgs3-RGG19b}. 
On note $\hcC^{\oy,(r)}_{X'}$ le complété séparé $p$-adique de $\cC^{\oy,(r)}_{X'}$. 
D'après \eqref{higgs3-RGG19k} et \eqref{higgs3-RGG19l}, pour tout objet $(U,\fp)$ de $\fV_\ox(\bQ)$, les homomorphismes canoniques
\begin{eqnarray}
\cF^{\oy,(r)}_U\otimes_{\hoR_U^\oy}\hoR_{X'}^\oy&\rightarrow& \cF^{\oy,(r)}_{X'},\label{higgs3-RGG16k}\\
\cC^{\oy,(r)}_U\otimes_{\hoR_U^\oy}\hoR_{X'}^\oy&\rightarrow& \cC^{\oy,(r)}_{X'},\label{higgs3-RGG16l}
\end{eqnarray}
sont des isomorphismes.  

\begin{rema}\label{higgs3-RGG172}
Sous les hypothèses de \ref{higgs3-RGG16}, pour tout entier $n\geq 1$, les morphismes canoniques 
\begin{eqnarray}
\underset{\underset{(U,\fp)\in \fV_\ox(\bQ)^\circ}{\longrightarrow}}{\lim}\ 
\cF^{\oy,(r)}_U/p^n\cF^{\oy,(r)}_U&\rightarrow&\cF^{\oy,(r)}_{X'}/p^n\cF^{\oy,(r)}_{X'},\label{higgs3-RGG172b}\\
\underset{\underset{(U,\fp)\in \fV_\ox(\bQ)^\circ}{\longrightarrow}}{\lim}\ 
\cC^{\oy,(r)}_U/p^n\cC^{\oy,(r)}_U&\rightarrow&\cC^{\oy,(r)}_{X'}/p^n\cC^{\oy,(r)}_{X'},\label{higgs3-RGG172a}
\end{eqnarray}
sont des isomorphismes.
\end{rema}

\begin{lem}\label{higgs3-RGG160}
Soient $(\oy\rightsquigarrow \ox)$ un point de $X_\et\gtimes_{X_\et}\oX^\circ_\et$ \eqref{higgs3-tfa41} tel que $\ox$ soit
au-dessus de $s$, $X'$ le localisé strict de $X$ en $\ox$, $r$ un nombre rationnel $\geq 0$. 
Reprenons les notations de \eqref{higgs3-RGG3} et \eqref{higgs3-RGG16}; posons, de plus, $R'_1=\Gamma(\oX',\co_{\oX'})$ et 
notons $\hR'_1$ son séparé complété $p$-adique. Alors, 
\begin{itemize}
\item[{\rm (i)}] Les anneaux $R'_1$ et $\oR^{\oy}_{X'}$ sont intègres, normaux et  $\co_\oK$-plats,
et l'homomorphisme canonique $R'_1\rightarrow \oR^{\oy}_{X'}$ est injectif et entier.
\item[{\rm (ii)}] Les anneaux $\hR'_1$, $\hoR^\oy_{X'}$, $\cC^{\oy,(r)}_{X'}$ et $\hcC^{\oy,(r)}_{X'}$ sont $\co_C$-plats.
\item[{\rm (iii)}] Pour tout entier $n\geq 1$, l'homomorphisme canonique $R'_1/p^nR'_1\rightarrow \oR^\oy_{X'}/p^n\oR^\oy_{X'}$ 
est injectif. 
\item[{\rm (iv)}] L'homomorphisme canonique $\hR'_1\rightarrow \hoR^\oy_{X'}$ est injectif,
et la topologie $p$-adique sur $\hR'_1[\frac 1 p]$ est induite par la topologie $p$-adique sur 
$\hoR^\oy_{X'}[\frac 1 p]$ {\rm (\cite{ag1} 2.2)}.
\end{itemize}
\end{lem}

(i) L'anneau $R'_1$ est intègre et normal en vertu de \ref{higgs3-sli6}, et il est clairement $\co_\oK$-plat. 
Pour tout objet $(U,\fp)$ de $\fV_\ox(\bP)$, l'anneau 
$\oR^\oy_U$ \eqref{higgs3-tfa112c} est intègre, normal et $\co_\oK$-plat d'après \ref{higgs3-tfa9} et (\cite{ega1n} 0.6.5.12(ii)).
De plus, notant $\oU^\star$  la composante irréductible de $\oU$ contenant $\oy$, l'homomorphisme canonique
$\Gamma(\oU^\star,\co_{\oU^\star})\rightarrow \oR^\oy_U$ est injectif et entier. On vérifie aussitôt que pour tout morphisme
$(U',\fp')\rightarrow (U,\fp)$ de $\fV_\ox(\bP)$, l'homomorphisme canonique 
$\oR^\oy_U\rightarrow \oR^\oy_{U'}$ est injectif. Par suite, $\oR^{\oy}_{X'}$ est intègre, normal et $\co_\oK$-plat. 
Comme $\oX'$ est intègre d'après \ref{higgs3-sli6}, on a 
\begin{equation}
R'_1\simeq \underset{\underset{(U,\fp)\in \fV_\ox(\bP)^\circ}{\longrightarrow}}{\lim}\ \Gamma(\oU^\star,\co_{\oU^\star}).
\end{equation}
L'homomorphisme canonique $R'_1\rightarrow \oR^{\oy}_{X'}$ est donc injectif et entier.

(ii) D'après (\cite{ac} chap.~III §2.11 prop.~14 et cor.~1), pour tout $n\geq 1$, on a 
\begin{equation}\label{higgs3-RGG160a}
\hR'_1/p^n\hR'_1\simeq R'_1/p^nR'_1.
\end{equation}
Soient $\alpha\in \hR'_1$ tel que $p\alpha=0$, $\oalpha$ sa classe dans $R'_1/p^nR'_1$ $(n\geq 1)$. 
Comme $R'_1$ est $\co_\oK$-plat, il résulte de \eqref{higgs3-RGG160a} que $\oalpha\in p^{n-1}\hR'_1/p^n\hR'_1$.
On en déduit que $\alpha\in \cap_{n\geq 0}p^n\hR'_1=\{0\}$. 
Par suite, $p$ n'est pas diviseur de zéro dans $\hR'_1$, et donc $\hR'_1$ est $\co_C$-plat (\cite{ac} Chap.~VI §3.6 lem.~1). 
On démontre de même que $\hoR^\oy_{X'}$ est $\co_C$-plat. 
Par suite, $\cC^{\oy,(r)}_{X'}$ est $\co_C$-plat \eqref{higgs3-RGG16l}. 
On en déduit comme plus haut que $\hcC^{\oy,(r)}_{X'}$ est $\co_C$-plat.

(iii) Cela résulte aussitôt de (i). 

(iv) La première assertion résulte aussitôt de (iii). D'après (iii) et \eqref{higgs3-RGG160a}, pour tout $n\geq 1$, on a 
$\hR'_1\cap p^n\hoR^{\oy}_{X'}=p^n\hR'_1$.  Comme $\hR'_1$ est $\co_C$-plat, 
on en déduit que $\hR'_1[\frac 1 p]\cap p^n\hoR^{\oy}_{X'}=p^n\hR'_1$; d'où la seconde assertion. 

\begin{prop}\label{higgs3-RGG15}
Soient $\ox$ un point de $X$ au-dessus de $s$, $X'$ le localisé strict de $X$ en $\ox$, 
$r$ un nombre rationnel $\geq 0$, $n$ un entier $\geq 0$.
Alors, avec les notations de \eqref{higgs3-RGG3}, on a des isomorphismes canoniques de $\oX'^\circ_\fet$
\begin{eqnarray}
\underset{\underset{(U,\fp)\in \fV_\ox(\bQ)^\circ}{\longrightarrow}}{\lim}\ 
(\ofp^\circ)^*_\fet(\cC^{(r)}_{U,n})&\stackrel{\sim}{\rightarrow}&\varphi_\ox(\cC^{(r)}_n),\label{higgs3-RGG15a}\\
\underset{\underset{(U,\fp)\in \fV_\ox(\bQ)^\circ}{\longrightarrow}}{\lim}\ 
(\ofp^\circ)^*_\fet(\cF^{(r)}_{U,n})&\stackrel{\sim}{\rightarrow}&\varphi_\ox(\cF^{(r)}_n).\label{higgs3-RGG15b}
\end{eqnarray}
\end{prop}
On notera que la proposition ne résulte pas directement de (\cite{ag2} 10.37) puisque $\bQ$ n'est pas
stable par produits fibrés. Posons \eqref{higgs3-tfa2a}
\begin{equation}\label{higgs3-RGG15c}
\cC^{(r)}_n=\{U\mapsto C_U\}, \ \ \ U\in \ob(\Et_{/X}).
\end{equation}
En vertu de (\cite{ag2} 10.37),  on a un isomorphisme canonique fonctoriel
\begin{equation}\label{higgs3-RGG15d}
\varphi_\ox(\cC^{(r)}_n)\stackrel{\sim}{\rightarrow} 
\underset{\underset{(U,\fp)\in \fV_\ox^\circ}{\longrightarrow}}{\lim}\ (\ofp^\circ)_\fet^*(C_U). 
\end{equation} 
On peut évidemment remplacer $\fV_\ox$ par la catégorie $\fV_\ox(\bQ)$ \eqref{higgs3-RGG3}. 
Pour démontrer \eqref{higgs3-RGG15a}, il suffit donc de montrer que le morphisme canonique
\begin{equation}\label{higgs3-RGG15e}
\gamma\colon \underset{\underset{(U,\fp)\in \fV_\ox(\bQ)^\circ}{\longrightarrow}}{\lim}\ (\ofp^\circ)_\fet^*(\cC^{(r)}_{U,n}) \rightarrow 
\underset{\underset{(U,\fp)\in \fV_\ox(\bQ)^\circ}{\longrightarrow}}{\lim}\ (\ofp^\circ)_\fet^*(C_U) 
\end{equation}
est un isomorphisme. 
Soient $\oy$ un point géométrique de $\oX'^\circ$, $\phi_{\oy}\colon \oX'^\circ_\fet\rightarrow \Ens$
le foncteur fibre associé au point $\rho_{\oX'^\circ}(\oy)$ de $\oX'^\circ_\fet$ \eqref{higgs3-not2a}. 
Il suffit encore de montrer que $\phi_\oy(\gamma)$ est un isomorphisme (\cite{ag2} 9.6). 

On note aussi $\oy$ le point géométrique de $\oX^\circ$ induit par $\oy$, $u\colon \oy\rightarrow X'$
le $X$-morphisme canonique et 
$(\oy \rightsquigarrow \ox)$ le point de $X_\et\gtimes_{X_\et}\oX^\circ_\et$ défini par $u$. 
Soit $\cQ_{\rho(\oy \rightsquigarrow \ox)}$ la catégorie des objets $\rho(\oy \rightsquigarrow \ox)$-pointés de $E_\bQ$ 
\eqref{higgs3-RGG85}, autrement dit la catégorie des triplets $((V\rightarrow U), \fp,\fq)$ 
formés d'un objet $(V\rightarrow U)$ de $E_\bQ$, d'un $X$-morphisme $\fp\colon \ox\rightarrow U$ et 
d'un $\oX^\circ$-morphisme $\fq\colon \oy\rightarrow V$ tels que notant encore 
$\fp\colon X'\rightarrow U$ le $X$-morphisme induit par $\fp$, le diagramme
\begin{equation}\label{higgs3-RGG15f}
\xymatrix{
\oy\ar[r]^u\ar[d]_{\fq}&{X'}\ar[d]^{\fp}\\
V\ar[r]&U}
\end{equation}
soit commutatif. D'après  (\cite{ag2} 10.20), $\cQ_{\rho(\oy \rightsquigarrow \ox)}$ 
est canoniquement équivalente à la catégorie des voisinages de 
$\rho(\oy \rightsquigarrow \ox)$ dans $E_\bQ$ (\cite{sga4} IV 6.8.2). Elle est donc cofiltrante et pour tout préfaisceau 
$F=\{U\mapsto F_U\}$ sur $E_\bQ$, on a un isomorphisme canonique fonctoriel (\cite{sga4} IV (6.8.4))
\begin{equation}\label{higgs3-RGG15g}
(F^a)_{\rho(\oy \rightsquigarrow \ox)}\stackrel{\sim}{\rightarrow} \underset{\underset{((V\rightarrow U), \fp,\fq)\in 
\cQ^\circ_{\rho(\oy \rightsquigarrow \ox)}}{\longrightarrow}}\lim\ F_U(V).
\end{equation} 
On a un foncteur 
\begin{equation}\label{higgs3-RGG15h}
\alpha\colon \cQ_{\rho(\oy \rightsquigarrow \ox)}\rightarrow \fV_\ox(\bQ), \ \ \ ((V\rightarrow U),\fp,\fq)\mapsto (U,\fp).
\end{equation}
Pour tout $(U,\fp)\in \ob(\fV_\ox(\bQ))$, la fibre de $\alpha$ au-dessus de $(U,\fp)$ est canoniquement équivalente 
à la catégorie $\cD^\oy_{(U,\fp)}$ des $\oU^\circ$-schémas finis étales, $\ofp^\circ(\oy)$-pointés \eqref{higgs3-RGG3a}. 

Compte tenu de (\cite{ag2} (9.3.4)) et (\cite{sga4} IV (6.8.4)), $\phi_\oy(\gamma)$ s'identifie à l'application canonique 
\begin{equation}
\phi_\oy(\gamma)\colon \underset{\underset{(U,\fp)\in \fV_\ox(\bQ)^\circ}{\longrightarrow}}{\lim}\ 
\underset{\underset{(V,\fq)\in (\cD^\oy_{(U,\fp)})^\circ}{\longrightarrow}}{\lim}\ \cC^{(r)}_{U,n}(V)
\rightarrow \underset{\underset{(U,\fp)\in \fV_\ox(\bQ)^\circ}{\longrightarrow}}{\lim}\ 
\underset{\underset{(V,\fq)\in (\cD^\oy_{(U,\fp)})^\circ}{\longrightarrow}}{\lim}\ C_U(V).
\end{equation}
On peut clairement remplacer chacune des doubles limites inductives ci-dessus par une limite inductive 
sur la catégorie $\cQ^\circ_{\rho(\oy \rightsquigarrow \ox)}$. Il résulte alors de \eqref{higgs3-RGG15g} que $\phi_\oy(\gamma)$
est bijective~; d'où l'isomorphisme \eqref{higgs3-RGG15a}. La preuve de \eqref{higgs3-RGG15b}
est similaire.

\begin{cor}\label{higgs3-RGG17}
Soient $(\oy\rightsquigarrow \ox)$ un point de $X_\et\gtimes_{X_\et}\oX^\circ_\et$ tel que $\ox$ soit
au-dessus de $s$, $X'$ le localisé strict de $X$ en $\ox$, $n$ un entier $\geq 0$, $r$ un nombre rationnel $\geq 0$.
Alors, avec les notations de \eqref{higgs3-RGG3} et \eqref{higgs3-RGG19}, on a des isomorphismes canoniques de $\bB_{\pi_1(\oX'^\circ,\oy)}$
\begin{eqnarray}
\nu_\oy(\varphi_\ox(\cC^{(r)}_n))&\stackrel{\sim}{\rightarrow}&\underset{\underset{(U,\fp)\in \fV_\ox(\bQ)^\circ}{\longrightarrow}}{\lim}\ 
\cC^{\oy,(r)}_U/p^n\cC^{\oy,(r)}_U,\label{higgs3-RGG17a}\\
\nu_\oy(\varphi_\ox(\cF^{(r)}_n))&\stackrel{\sim}{\rightarrow}&\underset{\underset{(U,\fp)\in \fV_\ox(\bQ)^\circ}{\longrightarrow}}{\lim}\ 
\cF^{\oy,(r)}_U/p^n\cF^{\oy,(r)}_U.\label{higgs3-RGG17b}
\end{eqnarray}
\end{cor}

Cela résulte de \ref{higgs3-RGG15}, \eqref{higgs3-RGG19g} et \eqref{higgs3-RGG19h}.

\begin{remas}\label{higgs3-RGG171}
Conservons les hypothèses de \eqref{higgs3-RGG17}~; supposons, de plus, $n\geq 1$. 
\begin{itemize}
\item[(i)] Les isomorphismes \eqref{higgs3-RGG16m}, \eqref{higgs3-RGG17a} et \eqref{higgs3-RGG17b} sont compatibles entre eux.  
\item[(ii)] On a un isomorphisme canonique de $\bB_{\pi_1(\oX'^\circ,\oy)}$
\begin{equation}\label{higgs3-RGG171a}
\nu_\oy(\varphi_\ox(\sigma_n^*(\xi^{-1}\tOmega^1_{\oX_n/\oS_n})))\stackrel{\sim}{\rightarrow}
\underset{\underset{(U,\fp)\in \fV_\ox^\circ}{\longrightarrow}}{\lim}\ 
\xi^{-1}\tOmega^1_{X/S}(U)\otimes_{\co_X(U)}(\oR^{\oy}_U/p^n\oR^{\oy}_U).
\end{equation}
En effet, d'après \eqref{higgs3-RGG18ee} et (\cite{ag2} 5.34(ii), 8.9 et 5.17), $\sigma^*_n(\xi^{-1}\tOmega^1_{\oX_n/\oS_n})$ est 
le faisceau de $\tE$ associé au préfaisceau sur $E_\bQ$ défini par la correspondance 
\begin{equation}\label{higgs3-RGG171b}
\{U\mapsto \xi^{-1}\tOmega^1_{X/S}(U)\otimes_{\co_X(U)}\ocB_{U,n}\}, \ \ \ (U\in \ob(\Et_{/X})).
\end{equation}
L'assertion résulte donc de (\cite{ag2} 10.37) et \eqref{higgs3-tfa113a}.
\item[(iii)] L'image de la suite exacte \eqref{higgs3-RGG18a} par le foncteur composé $\nu_\oy\circ \varphi_\ox$
s'identifie à la limite projective sur la catégorie $\fV_\ox(\bQ)^\circ$ des suites exactes
\begin{equation}\label{higgs3-RGG171c}
0\rightarrow\oR^{\oy}_U/p^n\oR^{\oy}_U \rightarrow 
\cF^{{\oy},(r)}_U/p^n\cF^{{\oy},(r)}_U \rightarrow 
\xi^{-1}\tOmega^1_{X/S}(U)\otimes_{\co_X(U)}(\oR^{\oy}_U/p^n\oR^{\oy}_U) \rightarrow 0
\end{equation}
déduites de \eqref{higgs3-RGG19a}. 
\item[(iv)] L'image de l'isomorphisme \eqref{higgs3-RGG18b} par le foncteur composé $\nu_\oy\circ \varphi_\ox$
s'identifie à l'isomorphisme
\begin{equation}\label{higgs3-RGG171d}
\underset{\underset{(U,\fp)\in \fV_\ox(\bQ)^\circ}{\longrightarrow}}{\lim}\ \cC^{{\oy},(r)}_U/p^n\cC^{{\oy},(r)}_U
\stackrel{\sim}{\rightarrow} \underset{\underset{(U,\fp)\in \fV_\ox(\bQ)^\circ}{\longrightarrow}}{\lim}\ 
\underset{\underset{m\geq 0}{\longrightarrow}}\lim\ 
\rS^m_{\oR_{U}^\oy}(\cF_{U}^{{\oy},(r)}/p^n\cF_{U}^{{\oy},(r)})
\end{equation}
déduit de \eqref{higgs3-RGG19b}. La preuve est similaire à celle \ref{higgs3-RGG15} et est laissée au lecteur. 
\end{itemize}
\end{remas}

\subsection{}\label{higgs3-RGG21}\index{10001070@$\bvsigma\colon (\tE_s^{\mN^\circ},\bvocB)\rightarrow(X_{s,\et}^{\mN^\circ},\co_{\bvoX})$}
\index{10001071@$\xi^{-1}\tOmega^1_{\bvoX/\bvoS}$}
\index{10001072@$\bvcF^{(r)}$, $\bvcC^{(r)}$ ($r\in \mQ_{\geq 0}$)}\index{10001073@$\bvcF$, $\bvcC$}\index{10001074@$\bvtta^{r,r'}$, $\bvalpha^{r,r'}$}
\index{10001075@$\bvd^{(r)}$}
\index{Algebre de Higgs-Tate d'epaisseur2@$\bvocB$-algèbre de Higgs-Tate d'épaisseur $r$ ($\bvcC^{(r)}$)}
\index{Extension de Higgs-Tate d'epaisseur2@$\bvocB$-extension de Higgs-Tate d'épaisseur $r$ ($\bvcF^{(r)}$)}
On désigne par $\bvocB$ l'anneau $(\ocB_{n+1})_{n\in \mN}$ de $\tE_s^{\mN^\circ}$  \eqref{higgs3-spsa2}, 
par $\co_{\bvoX}$ l'anneau $(\co_{\oX_{n+1}})_{n\in \mN}$ de $X_{s,\et}^{\mN^\circ}$
(à ne pas confondre avec $\co_{\coX}$ \eqref{higgs3-RGG1a}) et par $\xi^{-1}\tOmega^1_{\bvoX/\bvoS}$
le $\co_{\bvoX}$-module $(\xi^{-1}\tOmega^1_{\oX_{n+1}/\oS_{n+1}})_{n\in \mN}$ de $X_{s,\et}^{\mN^\circ}$ 
\eqref{higgs3-RGG1ab}. On note 
\begin{equation}\label{higgs3-RGG21c}
\bvsigma\colon (\tE_s^{\mN^\circ},\bvocB)\rightarrow(X_{s,\et}^{\mN^\circ},\co_{\bvoX})
\end{equation}
le morphisme de topos annelés induit par les $(\sigma_{n+1})_{n\in \mN}$ \eqref{higgs3-RGG2g}. 

Soit $r$ un nombre rationnel $\geq 0$.
Pour tous entiers $m\geq n\geq 1$, on a un morphisme $\ocB_m$-linéaire canonique
$\cF^{(r)}_m\rightarrow \cF^{(r)}_n$ et un homomorphisme canonique de $\ocB_m$-algèbres 
$\cC^{(r)}_m\rightarrow \cC^{(r)}_n$, 
compatibles avec la suite exacte \eqref{higgs3-RGG18a} et l'isomorphisme \eqref{higgs3-RGG18b} et tels que 
les morphismes induits 
\begin{equation}\label{higgs3-RGG21a}
\cF^{(r)}_m\otimes_{\ocB_m}\ocB_n\rightarrow \cF^{(r)}_n\ \ \ {\rm et}\ \ \ 
\cC^{(r)}_m\otimes_{\ocB_m}\ocB_n\rightarrow \cC^{(r)}_n
\end{equation}
soient des isomorphismes. Ces morphismes forment des systèmes compatibles lorsque $m$ et $n$ varient, 
de sorte que $(\cF^{(r)}_{n+1})_{n\in \mN}$ et $(\cC^{(r)}_{n+1})_{n\in \mN}$ sont des systèmes projectifs. 
On appelle {\em $\bvocB$-extension de Higgs-Tate d'épaisseur $r$} associée à $(f,\tX,\cM_\tX)$, 
et l'on note $\bvcF^{(r)}$, le $\bvocB$-module $(\cF^{(r)}_{n+1})_{n\in \mN}$ de $\tE_s^{\mN^\circ}$. 
On appelle {\em $\bvocB$-algèbre de Higgs-Tate d'épaisseur $r$} associée à $(f,\tX,\cM_\tX)$,
et l'on note $\bvcC^{(r)}$, la $\bvocB$-algèbre $(\cC^{(r)}_{n+1})_{n\in \mN}$ de $\tE_s^{\mN^\circ}$.
Ce sont des $\bvocB$-modules adiques \eqref{higgs3-spad3}. 
D'après \ref{higgs3-spsa99}(i), \eqref{higgs3-spsa1g} et \eqref{higgs3-spsa2a}, la suite exacte \eqref{higgs3-RGG18a}
induit une suite exacte de $\bvocB$-modules 
\begin{equation}\label{higgs3-RGG21b}
0\rightarrow \bvocB\rightarrow \bvcF^{(r)}\rightarrow 
\bvsigma^*(\xi^{-1}\tOmega^1_{\bvoX/\bvoS})\rightarrow 0.
\end{equation}
Comme le $\co_X$-module $\tOmega^1_{X/S}$ est localement libre de type fini, 
le $\bvocB$-module $\bvsigma^*(\xi^{-1}\tOmega^1_{\bvoX/\bvoS})$ est localement libre de type fini
et la suite \eqref{higgs3-RGG21b} est localement scindée. 
D'après \eqref{higgs3-vck2b}, elle induit pour tout entier $m\geq 1$, une suite exacte de $\bvocB$-modules 
\begin{equation}\label{higgs3-RGG21bc}
0\rightarrow \rS^{m-1}_{\bvocB}(\bvcF^{(r)})\rightarrow \rS^m_{\bvocB}(\bvcF^{(r)})\rightarrow 
\bvsigma^*(\rS^m_{\co_{\bvoX}}(\xi^{-1}\tOmega^1_{\bvoX/\bvoS}))\rightarrow 0.
\end{equation}
En particulier, les $\bvocB$-modules $(\rS^m_{\bvocB}(\bvcF^{(r)}))_{m\in \mN}$ forment un système inductif filtrant. 
D'après \ref{higgs3-spsa99}(i), \eqref{higgs3-spsa2c} et \eqref{higgs3-RGG18b}, on a un isomorphisme canonique de $\bvocB$-algèbres
\begin{equation}\label{higgs3-RGG21d}
\bvcC^{(r)}\stackrel{\sim}{\rightarrow}\underset{\underset{m\geq 0}{\longrightarrow}}\lim\ \rS^m_{\bvocB}(\bvcF^{(r)}). 
\end{equation}

On pose $\bvcF=\bvcF^{(0)}$ et $\bvcC=\bvcC^{(0)}$, et  
on les appelle la {\em $\bvocB$-extension de Higgs-Tate} et  la {\em $\bvocB$-algèbre de Higgs-Tate}, respectivement,
associées à $(f,\tX,\cM_\tX)$.
Pour tous nombres rationnels $r\geq r'\geq 0$, les morphismes $(\tta_n^{r,r'})_{n\in \mN}$ \eqref{higgs3-RGG24c}
induisent un morphisme $\bvocB$-linéaire 
\begin{equation}\label{higgs3-RGG180cc}
\bvtta^{r,r'}\colon \bvcF^{(r)}\rightarrow \bvcF^{(r')}.
\end{equation}
Les homomorphismes $(\alpha_n^{r,r'})_{n\in \mN}$ \eqref{higgs3-RGG24d}
induisent un homomorphisme de $\bvocB$-algèbres 
\begin{equation}\label{higgs3-RGG180c}
\bvalpha^{r,r'}\colon \bvcC^{(r)}\rightarrow \bvcC^{(r')}.
\end{equation}
Pour tous nombres rationnels $r\geq r'\geq r''\geq 0$, on a
\begin{equation}\label{higgs3-RGG180e}
\bvtta^{r,r''}=\bvtta^{r',r''} \circ \bvtta^{r,r'} \ \ \ {\rm et}\ \ \ \bvalpha^{r,r''}=\bvalpha^{r',r''} \circ \bvalpha^{r,r'}.
\end{equation}

Les dérivations $(d_{n+1}^{(r)})_{n\in \mN}$ \eqref{higgs3-RGG18e} définissent un morphisme
\begin{equation}\label{higgs3-RGG180b}
\bvd^{(r)}\colon \bvcC^{(r)}\rightarrow \bvsigma^*(\xi^{-1}\tOmega^1_{\bvoX/\bvoS})\otimes_{\bvocB}\bvcC^{(r)},
\end{equation}
qui n'est autre que la $\bvocB$-dérivation universelle de $\bvcC^{(r)}$. Elle prolonge le morphisme canonique 
$\bvcF^{(r)}\rightarrow \bvsigma^*(\xi^{-1}\tOmega^1_{\bvoX/\bvoS})$.
Pour tous nombres rationnels $r\geq r'\geq 0$, on a 
\begin{equation}\label{higgs3-RGG180d}
p^{r-r'}(\id \otimes \bvalpha^{r,r'}) \circ \bvd^{(r)}=\bvd^{(r')}\circ \bvalpha^{r,r'}.
\end{equation}

\begin{remas}\label{higgs3-RGG181}
Soient $r$ un nombre rationnel $\geq 0$, $n$ un  entier $\geq 1$. 
\begin{itemize}
\item[{\rm (i)}] Pour tout entier $m\geq 0$, les morphismes canoniques $\rS^m_{\ocB_n}(\cF^{(r)}_n)\rightarrow \cC_n^{(r)}$
et $\rS^m_{\bvocB}(\bvcF^{(r)})\rightarrow \bvcC^{(r)}$ sont injectifs. En effet, pour tout entier $m'\geq m$, le morphisme canonique $\rS^m_{\ocB_n}(\cF^{(r)}_n)\rightarrow  \rS^{m'}_{\ocB_n}(\cF^{(r)}_n)$
est injectif \eqref{higgs3-RGG18c}. Comme les limites injectives filtrantes commutent 
aux limites projectives finies dans $\tE_s^{\mN^\circ}$,
$\rS^m_{\ocB_n}(\cF^{(r)}_n)\rightarrow \cC_n^{(r)}$ est injectif. 
La seconde assertion se déduit de la première par \ref{higgs3-spsa99}(i). 
\item[{\rm (ii)}]  On a $\sigma_n^*(\xi^{-1}\tOmega^1_{\oX_n/\oS_n})=d_n^{(r)}(\cF_n^{(r)}) \subset d_n^{(r)}(\cC_n^{(r)})$
\eqref{higgs3-RGG18e}. Par suite, la dérivation $d_n^{(r)}$ est un $\ocB_n$-champ de Higgs à coefficients dans 
$\sigma_n^*(\xi^{-1}\tOmega^1_{\oX_n/\oS_n})$ d'après (\cite{ag1} 2.12).
\item[{\rm (iii)}] On a $\bvsigma^*(\xi^{-1}\tOmega^1_{\bvoX/\bvoS})=\bvd^{(r)}(\bvcF^{(r)}) \subset \bvd^{(r)}(\bvcC^{(r)})$
\eqref{higgs3-RGG180b}. Par suite, la dérivation $\bvd^{(r)}$ est un $\bvocB$-champ de Higgs à coefficients dans 
$\bvsigma^*(\xi^{-1}\tOmega^1_{\bvoX/\bvoS})$.  
\end{itemize}
\end{remas}

\begin{prop}\label{higgs3-plat2}
Pour tout nombre rationnel $r\geq 0$, le foncteur 
\begin{equation}\label{higgs3-plat2a}
\bMod(\bvocB)\rightarrow \bMod(\bvcC^{(r)}), \ \ \ M\mapsto M\otimes_{\bvocB} \bvcC^{(r)}
\end{equation}
est exact et fidèle~; en particulier, $\bvcC^{(r)}$ est $\bvocB$-plat. 
\end{prop}
Comme le $\co_X$-module $\tOmega^1_{X/S}$ est localement libre de type fini, 
la suite exacte \eqref{higgs3-RGG21b} est localement scindée. Un scindage local de cette suite 
induit, pour tout entier $m\geq 0$, un scindage local de la suite exacte \eqref{higgs3-RGG21bc}.
On en déduit que $\rS^m_{\bvocB}(\bvcF^{(r)})$ est $\bvocB$-plat et que l'homomorphisme 
canonique $\bvocB\rightarrow \bvcC^{(r)}$ admet localement des sections. 
La proposition s'ensuit compte tenu de \eqref{higgs3-RGG21d}.

\section{Calculs cohomologiques}\label{higgs3-coh}

\subsection{}\label{higgs3-coh0}\index{10001102@$\xi^{-i}\tOmega^i_{\bvoX/\bvoS}$}
\index{10001070@$\bvsigma\colon (\tE_s^{\mN^\circ},\bvocB)\rightarrow(X_{s,\et}^{\mN^\circ},\co_{\bvoX})$|textbf} 
\index{1000929@$\top\colon (\tE_s^{\mN^\circ},\bvocB)\rightarrow (X_{s,\zar},\co_{\fX})$|textbf}\index{1000928@$\fX$|textbf}
\index{1000927@$\bvsigma$,$\bvu$|textbf}\index{10001103@$\bvtau$}\index{1000922@$\sigma_n$}
Les hypothèses et notations générales de § \ref{higgs3-RGG} sont en vigueur dans cette section. 
On désigne, de plus, par $d=\dim(X/S)$ la dimension relative de $X$ sur $S$,
par $\bvocB$ l'anneau $(\ocB_{n+1})_{n\in \mN}$ de $\tE_s^{\mN^\circ}$  \eqref{higgs3-spsa2}
et par $\co_{\bvoX}$ l'anneau $(\co_{\oX_{n+1}})_{n\in \mN}$ 
de $X_{s,\zar}^{\mN^\circ}$ ou de $X_{s,\et}^{\mN^\circ}$, selon le contexte (cf. \ref{higgs3-not2} et \ref{higgs3-TFT9}). 
Pour tous entiers $i,n\geq 1$,  on pose
\begin{equation}\label{higgs3-coh0h}
\tOmega^i_{X/S}=\wedge^i(\tOmega^1_{X/S}) \ \ \ {\rm et} \ \ \  \tOmega^i_{\oX_n/\oS_n}=
\wedge^i(\tOmega^1_{\oX_n/\oS_n}). 
\end{equation}
On note $\xi^{-i}\tOmega^i_{\bvoX/\bvoS}$ le $\co_\bvoX$-module 
$(\xi^{-i}\tOmega^i_{\oX_{n+1}/\oS_{n+1}})_{n\in \mN}$. On a un isomorphisme canonique \eqref{higgs3-spsa2d} 
\begin{equation}
\xi^{-i}\tOmega^i_{\bvoX/\bvoS}\stackrel{\sim}{\rightarrow} \wedge^i(\xi^{-1}\tOmega^1_{\bvoX/\bvoS}).
\end{equation}

Pour tout entier $n\geq 1$, on désigne par
\begin{equation}\label{higgs3-coh0ab}
\sigma_n\colon (\tE_s,\ocB_n)\rightarrow (X_{s,\et},\co_{\oX_n})
\end{equation}
le morphisme canonique de topos annelés \eqref{higgs3-RGG2g}, par 
\begin{equation}\label{higgs3-coh0b}
u_n\colon (X_{s,\et},\co_{\oX_n})\rightarrow (X_{s,\zar},\co_{\oX_n})
\end{equation}
le morphisme canonique de topos annelés \eqref{higgs3-not2} et par  
\begin{equation}\label{higgs3-coh0a}
\tau_n\colon (\tE_s,\ocB_n)\rightarrow (X_{s,\zar},\co_{\oX_n})
\end{equation}
le morphisme composé $u_n\circ \sigma_n$. On note
\begin{eqnarray}
\bvsigma\colon (\tE_s^{\mN^\circ},\bvocB)&\rightarrow&(X_{s,\et}^{\mN^\circ},\co_{\bvoX}),\label{higgs3-coh0i}\\
\bvu\colon (X_{s,\et}^{\mN^\circ},\co_{\bvoX})&\rightarrow& (X_{s,\zar}^{\mN^\circ},\co_{\bvoX}),\label{higgs3-coh0f}\\
\bvtau\colon (\tE_s^{\mN^\circ},\bvocB)&\rightarrow& (X_{s,\zar}^{\mN^\circ},\co_{\bvoX}),\label{higgs3-coh0c}
\end{eqnarray}
les morphismes de topos annelés induit par les $(\sigma_{n+1})_{n\in \mN}$,  $(u_{n+1})_{n\in \mN}$ et 
$(\tau_{n+1})_{n\in \mN}$, respectivement, de sorte que $\bvtau=\bvu\circ \bvsigma$. 

On rappelle qu'on a posé $\cS=\Spf(\co_C)$ \eqref{higgs3-not1}. 
On désigne par $\fX$ le schéma formel complété $p$-adique de $\oX$. 
C'est un $\cS$-schéma formel de présentation finie (\cite{egr1} 2.3.15). Il est donc idyllique (\cite{egr1} 2.6.13). 
On désigne par $\xi^{-i}\tOmega^i_{\fX/\cS}$ le complété $p$-adique du $\co_\oX$-module 
$\xi^{-i}\tOmega^i_{\oX/\oS}=\xi^{-i}\tOmega^i_{X/S}\otimes_{\co_X}\co_{\oX}$ (\cite{egr1} 2.5.1).
On a un isomorphisme canonique (\cite{egr1} 2.5.5(ii))
\begin{equation}
\xi^{-i}\tOmega^i_{\fX/\cS}\stackrel{\sim}{\rightarrow}\wedge^i(\xi^{-1}\tOmega^1_{\fX/\cS}).
\end{equation}

On désigne par
\begin{equation}\label{higgs3-coh0d}
\lambda\colon (X_{s,\zar}^{\mN^\circ},\co_{\bvoX})\rightarrow (X_{s,\zar}, \co_\fX)
\end{equation}
le morphisme défini dans \eqref{higgs3-TFT14e} et par
\begin{equation}\label{higgs3-coh0e}
\top\colon (\tE_s^{\mN^\circ},\bvocB)\rightarrow (X_{s,\zar},\co_{\fX})
\end{equation}
le morphisme composé $\lambda\circ \bvtau$. 
Nous utilisons pour les modules la notation $\top^{-1}$ pour désigner l'image
inverse au sens des faisceaux abéliens et nous réservons la notation 
$\top^*$ pour l'image inverse au sens des modules~; et de même pour 
$\bvsigma$ et $\bvtau$. Pour tout $\co_\fX$-module $\cN$ de $X_{s,\zar}$, on a un isomorphisme canonique 
\begin{equation}\label{higgs3-coh0g}
\top^*(\cN)\stackrel{\sim}{\rightarrow}\bvtau^*((\cN/p^{n+1}\cN)_{n\in \mN}). 
\end{equation}
En particulier, $\top^*(\cN)$ est adique \eqref{higgs3-spad11}.

\subsection{}\label{higgs3-coh1}
Pour tout entier $n\geq 1$, le morphisme d'adjonction de $X_{s,\et}$
\begin{equation}\label{higgs3-coh1a}
\xi^{-1}\tOmega^1_{\oX_n/\oS_n}\rightarrow \sigma_{n*}(\sigma_n^*(\xi^{-1}\tOmega^1_{\oX_n/\oS_n}))
\end{equation}
induit un morphisme $\co_{\oX_n}$-linéaire 
\begin{equation}\label{higgs3-coh1aa}
\xi^{-1}\tOmega^1_{\oX_n/\oS_n}\rightarrow \tau_{n*}(\sigma_n^*(\xi^{-1}\tOmega^1_{\oX_n/\oS_n})).
\end{equation}
De même, compte tenu de \eqref{higgs3-spsa1f}, le morphisme d'adjonction 
\begin{equation}\label{higgs3-coh1b}
\xi^{-1}\tOmega^1_{\bvoX/\bvoS}\rightarrow \bvsigma_*(\bvsigma^*(\xi^{-1}\tOmega^1_{\bvoX/\bvoS}))
\end{equation}
induit un morphisme $\co_{\bvoX}$-linéaire 
\begin{equation}\label{higgs3-coh1bb}
\xi^{-1}\tOmega^1_{\bvoX/\bvoS}\rightarrow \bvtau_*(\bvsigma^*(\xi^{-1}\tOmega^1_{\bvoX/\bvoS})).
\end{equation}
Ce dernier induit un morphisme $\co_\fX$-linéaire 
\begin{equation}\label{higgs3-coh1c}
\xi^{-1}\tOmega^1_{\fX/\cS}\rightarrow \top_*(\bvsigma^*(\xi^{-1}\tOmega^1_{\bvoX/\bvoS})).
\end{equation}
On notera que le morphisme adjoint 
\begin{equation}\label{higgs3-coh1cc}
\top^*(\xi^{-1}\tOmega^1_{\fX/\cS})\rightarrow \bvsigma^*(\xi^{-1}\tOmega^1_{\bvoX/\bvoS})
\end{equation}
est un isomorphisme d'après \eqref{higgs3-coh0g} et la remarque suivant \eqref{higgs3-not2e}. 

On désigne par  
\begin{equation}\label{higgs3-coh1d}
\partial_n\colon \xi^{-1}\tOmega^1_{\oX_n/\oS_n}\rightarrow \rR^1\sigma_{n*}(\ocB_n)
\end{equation}
le morphisme $\co_{\oX_n}$-linéaire de $X_{s,\et}$ composé du 
morphisme \eqref{higgs3-coh1a} et du morphisme bord de la suite exacte longue de cohomologie déduite
de la suite exacte canonique \eqref{higgs3-RGG18a}
\begin{equation}\label{higgs3-coh1e}
0\rightarrow \ocB_n\rightarrow \cF_n\rightarrow 
\sigma_n^*(\xi^{-1}\tOmega^1_{\oX_n/\oS_n})\rightarrow 0.
\end{equation}
On note 
\begin{equation}\label{higgs3-coh1dd}
\bvpartial\colon \xi^{-1}\tOmega^1_{\bvoX/\bvoS}\rightarrow \rR^1\bvsigma_*(\bvocB)
\end{equation}
le morphisme $\co_{\bvoX}$-linéaire de $X_{s,\et}^{\mN^\circ}$ composé du 
morphisme \eqref{higgs3-coh1b} et du morphisme bord de la suite exacte longue de cohomologie déduite
de la suite exacte canonique \eqref{higgs3-RGG21b} 
\begin{equation}\label{higgs3-coh1ee}
0\rightarrow \bvocB\rightarrow \bvcF\rightarrow 
\bvsigma^*(\xi^{-1}\tOmega^1_{\bvoX/\bvoS})\rightarrow 0.
\end{equation}
Compte tenu de \eqref{higgs3-spsa1g}, \eqref{higgs3-spsa1h} et \eqref{higgs3-spsa2a}, 
on peut identifier $\bvpartial$ au morphisme $(\partial_n)_{n\geq 1}$. 

On désigne par 
\begin{equation}\label{higgs3-coh1f}
\delta\colon \xi^{-1}\tOmega^1_{\fX/\cS}\rightarrow \rR^1\top_*(\bvocB)
\end{equation}
le morphisme $\co_{\fX}$-linéaire de $X_{s,\zar}$ composé du 
morphisme \eqref{higgs3-coh1c} et du morphisme bord de la suite exacte longue de cohomologie déduite
de la suite exacte canonique \eqref{higgs3-RGG21b} 
\begin{equation}\label{higgs3-coh1g}
0\rightarrow \bvocB\rightarrow \bvcF\rightarrow 
\bvsigma^*(\xi^{-1}\tOmega^1_{\bvoX/\bvoS})\rightarrow 0.
\end{equation}

\begin{prop}\label{higgs3-coh2}
Soit $n$ un entier $\geq 1$. 
\begin{itemize}
\item[{\rm (i)}] Il existe un et un unique homomorphisme de $\co_{\oX_n}$-algèbres graduées de $\oX_{n,\et}$
\begin{equation}\label{higgs3-coh2a}
\wedge (\xi^{-1}\tOmega^1_{\oX_n/\oS_n})\rightarrow \oplus_{i\geq 0}\rR^i\sigma_{n*}(\ocB_n)
\end{equation}
dont la composante en degré un est le morphisme $\partial_n$ \eqref{higgs3-coh1d}. 
De plus, son noyau est annulé par $p^{\frac{2d}{p-1}}\fm_\oK$ 
et son conoyau est annulé par $p^{\frac{2d+1}{p-1}}\fm_\oK$ \eqref{higgs3-coh0}.
\item[{\rm (ii)}] Pour tout entier $i\geq d+1$, $\rR^i\sigma_{n*}(\ocB_n)$ est presque-nul. 
\end{itemize}
\end{prop}

Soient $\ox$ un point géométrique de $X$ au-dessus de $s$, $X'$ le localisé strict de $X$ en $\ox$, 
\begin{equation}\label{higgs3-coh2b}
\varphi_\ox\colon \tE\rightarrow \oX'^\circ_\fet
\end{equation}
le foncteur \eqref{higgs3-tfa6d}. D'après (\cite{ag2} 10.30) et \eqref{higgs3-TFT8e}, pour tout entier $i\geq 0$, on a un isomorphisme canonique 
\begin{equation}\label{higgs3-coh2d}
(\rR^i\sigma_{n*}(\ocB_n))_\ox\stackrel{\sim}{\rightarrow} \rH^i(\oX'^\circ_\fet, \varphi_\ox(\ocB_n)).
\end{equation}
En vertu de \ref{higgs3-sli6}, $\oX'$ est normal et strictement local (et en particulier intègre). 
Comme $\oX'$ est $\oS$-plat, $\oX'^\circ$ est intègre et non-vide. Soient 
$v\colon \oy\rightarrow \oX'^\circ$ un point géométrique générique de $\oX'^\circ$, 
\begin{equation}
\nu_{\oy}\colon \oX'^\circ_\fet \stackrel{\sim}{\rightarrow}\bB_{\pi_1(\oX'^\circ,\oy)}
\end{equation}
le foncteur fibre associé \eqref{higgs3-not6c}. 
On note encore $\oy$ le point géométrique de $\oX^\circ$ et 
$u\colon \oy\rightarrow X'$ le $X$-morphisme induits par $v$. 
On obtient ainsi un point $(\oy\rightsquigarrow \ox)$ de $X_\et\gtimes_{X_\et}\oX^\circ_\et$. 

On désigne par $\fV_\ox$ (resp. $\fV_\ox(\bQ)$) la catégorie des voisinages 
du point de $X_\et$ associé à $\ox$ dans le site $\Et_{/X}$ (resp. $\bQ$ \eqref{higgs3-RGG85}).
Pour tout $(U,\fp\colon \ox\rightarrow U)\in \ob(\fV_\ox)$, on note encore $\fp\colon X'\rightarrow U$ 
le morphisme déduit de $\fp$ et on pose
\begin{equation}\label{higgs3-coh2c}
\ofp^\circ=\fp\times_X\oX^\circ \colon \oX'^\circ \rightarrow \oU^\circ.
\end{equation}
On note encore (abusivement) $\oy$ le point géométrique $\ofp^\circ(v(\oy))$ de $\oU^\circ$.
On observera que $\oy$ est localisé en un point générique de $\oU^\circ$ car $\ofp^\circ$ est plat.
Comme $\oU$ est localement irréductible \eqref{higgs3-sli3}, 
il est la somme des schémas induits sur ses composantes irréductibles. 
Notons $\oU^\star$  la composante irréductible de $\oU$ contenant $\oy$. De même, $\oU^\circ$
est la somme des schémas induits sur ses composantes irréductibles, et $\oU^{\star \circ}=\oU^\star\times_XX^\circ$ 
est la composante irréductible de $\oU^\circ$ contenant $\oy$. 
Le morphisme $\ofp^\circ$ se factorise à travers $\oU^{\star \circ}$.
On a un isomorphisme canonique \eqref{higgs3-RGG3d}
\begin{equation}\label{higgs3-coh2e}
\varphi_\ox(\ocB)\stackrel{\sim}{\rightarrow}
\underset{\underset{(U,\fp)\in \fV_\ox^\circ}{\longrightarrow}}{\lim}\ (\ofp^\circ)^*_\fet(\ocB_U).
\end{equation}
Il induit pour tout entier $i\geq 0$, un isomorphisme 
\begin{equation}\label{higgs3-coh2f}
\rH^i(\oX'^\circ_\fet, \varphi_\ox(\ocB_n))
\stackrel{\sim}{\rightarrow} \underset{\underset{(U,\fp)\in \fV_\ox^\circ}{\longrightarrow}}{\lim}\
\rH^i((\oU^\circ_\fet)_{/\oU^{\star \circ}}, \ocB_{U,n}). 
\end{equation} 
En effet, si $(Z,\fq)$ est un objet de $\fV_\ox$ tel que $Z$ soit affine, il suffit d'appliquer (\cite{ag2} 11.10) 
à la sous-catégorie pleine de $(\fV_\ox)_{/(Z,\fq)}$ formée des objets $(U,\fp)$ tels que $U$ soit affine. 
Compte tenu de \eqref{higgs3-tfa113a} et (\cite{ag2} (9.8.6)), on déduit de \eqref{higgs3-coh2d} et \eqref{higgs3-coh2f} un isomorphisme 
\begin{equation}\label{higgs3-coh2g}
(\rR^i\sigma_{n*}(\ocB_n))_\ox\stackrel{\sim}{\rightarrow} \underset{\underset{(U,\fp)\in \fV_\ox^\circ}{\longrightarrow}}{\lim}\
\rH^i(\pi_1(\oU^{\star \circ},\oy),\oR^{\oy}_U/p^n\oR^{\oy}_U).
\end{equation}
Les isomorphismes \eqref{higgs3-coh2d} et \eqref{higgs3-coh2f} sont clairement compatibles aux cup-produits. 
Il en est donc de même de \eqref{higgs3-coh2g}.

D'autre part, $\oX'$ étant strictement local \eqref{higgs3-sli6}, il s'identifie au localisé strict de $\oX$ en $\ox$.
On a donc un isomorphisme canonique
\begin{equation}\label{higgs3-coh2gb}
\tOmega^1_{\oX_n/\oS_n,\ox}
\stackrel{\sim}{\rightarrow} \underset{\underset{(U,\fp)\in \fV^\circ_\ox}{\longrightarrow}}{\lim}\
\tOmega^1_{\oX_n/\oS_n}(\oU^\star), 
\end{equation} 
où l'on considère $\tOmega^1_{\oX_n/\oS_n}$ comme un faisceau de $\oX_\et$.

(i) Soit $(U,\fp)$ un objet de $\fV_\ox(\bQ)$. 
On notera que les schémas $U$, $\oU$ et $\oU^\star$ sont affines. 
La suite exacte de $\oR^{\oy}_U$-représentations de $\pi_1(\oU^{\star \circ},\oy)$
\begin{equation}\label{higgs3-coh2h}
0\rightarrow \oR^{\oy}_U/p^n\oR^{\oy}_U\rightarrow \cF^{\oy}_U/p^n\cF^{\oy}_U
\rightarrow \xi^{-1}\tOmega^1_{X/S}(U)\otimes_{\co_X(U)}(\oR^{\oy}_U/p^n\oR^{\oy}_U)\rightarrow 0
\end{equation}
déduite de \eqref{higgs3-RGG10m}, induit un morphisme $\co_{\oX_n}(\oU^\star)$-linéaire
\begin{equation}\label{higgs3-coh2i}
\alpha_{(U,\fp)}\colon\xi^{-1}\tOmega^1_{\oX_n/\oS_n}(\oU^\star)\rightarrow \rH^1(\pi_1(\oU^{\star \circ},\oy),
\oR^{\oy}_U/p^n\oR^{\oy}_U).
\end{equation}
On peut explicitement décrire ce morphisme en tenant compte de \ref{higgs3-RGG65}.
En effet, d'après (\cite{ag1} 10.16), on a un diagramme commutatif
\begin{equation}\label{higgs3-coh2k}
\xymatrix{
{\xi^{-1}\tOmega^1_{\oX_n/\oS_n}(\oU^\star)}\ar[d]_a\ar[r]^-(0.5){\alpha_{(U,\fp)}}&
{\rH^1(\pi_1(\oU^{\star \circ},\oy),\oR^{\oy}_U/p^n\oR^{\oy}_U)}\\
{\rH^1(\pi_1(\oU^{\star \circ},\oy),\xi^{-1}\oR^{\oy}_U(1)/p^n\xi^{-1}\oR^{\oy}_U(1))}\ar[r]^{-b}_\sim&
{\rH^1(\pi_1(\oU^{\star \circ},\oy),p^{\frac{1}{p-1}}\oR^{\oy}_U/p^{n+\frac{1}{p-1}}\oR^{\oy}_U)}\ar[u]_c}
\end{equation}
où $a$ est induit par le morphisme défini dans (\cite{ag1} (8.13.2)), 
$b$ est induit par l'isomorphisme 
\begin{equation}
\hoR^{\oy}_U(1)\stackrel{\sim}{\rightarrow} p^{\frac{1}{p-1}}\xi \hoR^{\oy}_U
\end{equation}
défini dans (\cite{ag1} 9.18) et $c$ est induit par l'injection canonique $p^{\frac{1}{p-1}}\oR^{\oy}_U \rightarrow 
\oR^{\oy}_U$.

En vertu de (\cite{ag1} 8.17(i)), il existe un et un unique homomorphisme de 
$\co_{\oX_n}(\oU^\star)$-algèbres graduées
\begin{equation}
\wedge (\xi^{-1}\tOmega^1_{\oX_n/\oS_n}(\oU^\star))\rightarrow 
\oplus_{i\geq 0}\rH^i(\pi_1(\oU^{\star \circ},\oy),\xi^{-i}\oR^{\oy}_U(i)/p^n\xi^{-i}\oR^{\oy}_U(i))
\end{equation}
dont la composante en degré un est le morphisme $a$. Celui-ci est presque-injectif et son conoyau est annulé par 
$p^{\frac{1}{p-1}}\fm_\oK$. On en déduit qu'il existe un et un unique homomorphisme de 
$\co_{\oX_n}(\oU^\star)$-algèbres graduées
\begin{equation}\label{higgs3-coh2j}
\wedge (\xi^{-1}\tOmega^1_{\oX_n/\oS_n}(\oU^\star))\rightarrow 
\oplus_{i\geq 0}\rH^i(\pi_1(\oU^{\star \circ},\oy),\oR^{\oy}_U/p^n\oR^{\oy}_U)
\end{equation}
dont la composante en degré $1$ est $\alpha_{(U,\fp)}$. Une chasse au diagramme \eqref{higgs3-coh2k} 
montre que le noyau de \eqref{higgs3-coh2j} est annulé par $p^{\frac{2d}{p-1}}\fm_\oK$.
Comme $\rH^i(\pi_1(\oU^{\star \circ},\oy),\oR^{\oy}_U/p^n\oR^{\oy}_U)$ 
est presque-nul pour tout $i\geq d+1$  en vertu de (\cite{ag1} 8.17(ii)), 
le conoyau de \eqref{higgs3-coh2j} est annulé par $p^{\frac{2d+1}{p-1}}\fm_\oK$.

Par ailleurs, d'après \ref{higgs3-RGG171}(iii), l'image de la suite exacte \eqref{higgs3-coh1e} par le foncteur $\nu_\oy\circ \varphi_\ox$ 
s'identifie à la limite inductive sur la catégorie $\fV_\ox(\bQ)^\circ$ des suites exactes \eqref{higgs3-coh2h}.  
Par suite, d'après (\cite{ag2} 10.30(iii)), la fibre du morphisme $\partial_n$ \eqref{higgs3-coh1d} en $\ox$ s'identifie à la limite
inductive sur la catégorie $\fV_\ox(\bQ)^\circ$ des morphismes $\alpha_{(U,\fp)}$.
On en déduit l'existence (et l'unicité) de l'homomorphisme \eqref{higgs3-coh2a} compte tenu de \ref{higgs3-TFT5}(ii). 
La fibre de ce dernier en $\ox$ s'identifie à la limite
inductive sur la catégorie $\fV_\ox(\bQ)^\circ$ des homomorphismes \eqref{higgs3-coh2j}. 
La proposition s'ensuit puisque les limites inductives filtrantes sont exactes.  

(ii) Cela résulte de \eqref{higgs3-coh2g}, \ref{higgs3-TFT5}(ii) et (\cite{ag1} 8.17(ii)). 

\begin{cor}\label{higgs3-coh200}
{\rm (i)}\ Il existe un et un unique homomorphisme de $\co_{\bvoX}$-algèbres graduées de $X_{s,\et}^{\mN^\circ}$
\begin{equation}\label{higgs3-coh200a}
\wedge (\xi^{-1}\tOmega^1_{\bvoX/\bvoS})\rightarrow \oplus_{i\geq 0}\rR^i\bvsigma_*(\bvocB)
\end{equation}
dont la composante en degré un est induite par le morphisme $\bvpartial$ \eqref{higgs3-coh1dd}. 
De plus, son noyau est annulé par $p^{\frac{2d}{p-1}}\fm_\oK$ 
et son conoyau est annulé par $p^{\frac{2d+1}{p-1}}\fm_\oK$. 

{\rm (ii)}\ Pour tout entier $i\geq d+1$, $\rR^i\bvsigma_*(\bvocB)$ est presque-nul. 
\end{cor}

Cela résulte de \ref{higgs3-coh2}, \ref{higgs3-spsa99}(i) et \eqref{higgs3-spsa1h}. 

\begin{prop}\label{higgs3-coh4}
{\rm (i)}\ Pour tous entiers $i\geq 1$ et $j\geq 1$, 
le $\co_{\bvoX}$-module $\rR^i\bvu_*(\rR^j\bvsigma_*(\bvocB))$ 
est annulé par $p^{\frac{4d+1}{p-1}}\fm_\oK$.

{\rm (ii)}\ Pour tout entier $q\geq 0$, le noyau du morphisme 
\begin{equation}\label{higgs3-coh4a}
\rR^q\bvtau_*(\bvocB)\rightarrow \bvu_*(\rR^q\bvsigma_*(\bvocB))
\end{equation}
induit par la suite spectrale de Cartan-Leray, est annulé par $p^{\frac{(d+1)(4d+1)}{p-1}}\fm_\oK$
et son conoyau est annulé par $p^{\frac{q(4d+1)}{p-1}}\fm_\oK$.
\end{prop}

(i) En effet, il résulte de \ref{higgs3-coh2}(i) et (\cite{sga4} VII 4.3)
que pour tout entier $n\geq 1$, les $\co_{\oX_n}$-modules $\rR^iu_{n*}(\rR^j\sigma_{n*}(\ocB_n))$ 
sont annulés par $p^{\frac{4d+1}{p-1}}\fm_\oK$. La proposition s'ensuit compte tenu de \eqref{higgs3-spsa1h}.

(ii) On peut clairement se borner au cas où $q\geq 1$.
Considérons la suite spectrale de Cartan-Leray (\cite{sga4} V 5.4)
\begin{equation}\label{higgs3-coh4b}
\rE_2^{i,j}=\rR^i\bvu_*(\rR^j\bvsigma_*(\bvocB))\Rightarrow \rR^{i+j}\bvtau_*(\bvocB), 
\end{equation}
et notons $(\rE^q_i)_{0\leq i\leq q}$ la filtration aboutissement sur $\rR^q\bvtau_*(\bvocB)$.
On observera que $\rE^q_1$ est le noyau du morphisme \eqref{higgs3-coh4a}.
On sait que $\rE_2^{i,j}$ est annulé par $p^{\frac{4d+1}{p-1}}\fm_\oK$ pour tous $i\geq 1$ et $j\geq 0$ d'après (i), 
et qu'il est presque-nul pour tous $i\geq 0$ et $j\geq d+1$ en vertu de \ref{higgs3-coh200}(ii). On en déduit que  
\begin{equation}\label{higgs3-coh4c}
\rE_i^q/\rE_{i+1}^q=\rE_\infty^{i,q-i}
\end{equation}
est annulé par $p^{\frac{4d+1}{p-1}}\fm_\oK$ pour tout $i\geq 1$ et qu'il est presque-nul pour tout $i\leq q-d-1$.
Par suite, $\rE^q_1$ est annulé par $p^{\frac{(4d+1)(d+1)}{p-1}}\fm_\oK$; d'où la première assertion. 
Par ailleurs, on a $\rE_\infty^{0,q}=\rE_{q+2}^{0,q}$, et le conoyau du morphisme \eqref{higgs3-coh4a} s'identifie au 
conoyau du composé des injections canoniques
\begin{equation}\label{higgs3-coh4d}
\rE_{q+2}^{0,q}\rightarrow \rE_{q+1}^{0,q}\rightarrow \dots \rightarrow \rE_2^{0,q}.
\end{equation}
Le conoyau de chacune de ces injections est annulé par $p^{\frac{4d+1}{p-1}}\fm_\oK$ d'après (i); d'où la 
seconde assertion.

\begin{lem}\label{higgs3-coh5}
Soit $\bvcM=(\cM_{n+1})_{n\in \mN}$ un $\co_{\bvoX}$-module de $X_{s,\zar}^{\mN^\circ}$
tel que pour tous entiers $m\geq n\geq 1$, $\cM_n$ soit un $\co_{\oX_n}$-module quasi-cohérent sur $\oX_n$
et que le morphisme canonique $\cM_m\rightarrow \cM_n$ soit surjectif. Alors $\rR^i\lambda_*(\bvcM)$ 
\eqref{higgs3-coh0d} est nul pour tout entier $i\geq 1$. 
\end{lem}

En effet, d'après (\cite{sga4} V 5.1),
$\rR^i\lambda_*(\bvcM)$ est le faisceau associé au préfaisceau sur $X_{s,\zar}$
qui à tout ouvert de Zariski $U$ de $X_s$, associe le $\co_{\bvoX}(U)$-module 
$\rH^i(\lambda^*(U),\bvcM)$. En vertu de \ref{higgs3-spsa5}, on a une suite exacte canonique
\begin{equation}\label{higgs3-coh5a}
0\rightarrow \rR^1\underset{\underset{n\geq 1}{\longleftarrow}}{\lim}\ \rH^{i-1}(U,\cM_n)
\rightarrow \rH^i(\lambda^*(U),\bvcM)
\rightarrow \underset{\underset{n\geq 1}{\longleftarrow}}{\lim}\ \rH^{i}(U,\cM_n)\rightarrow 0.
\end{equation}
On identifie dans la suite les sites de Zariski de $X_s$ et $\oX_1$ \eqref{higgs3-TFT9}. 
Soit $U$ un ouvert affine de $\oX_1$. 
Pour tout entier $n\geq 1$, $U$ détermine un ouvert affine de $\oX_n$ (\cite{ega1n} 2.3.5). 
Par suite, $\rH^{i}(U,\cM_n)=0$ et pour tous entiers $m\geq n\geq 1$,  le morphisme canonique 
\begin{equation}\label{higgs3-coh5b}
\rH^0(U,\cM_m)\rightarrow \rH^0(U,\cM_n)
\end{equation}
est surjectif, de sorte que le système projectif de groupes abéliens 
$(\rH^0(U,\cM_{n+1}))_{n\in \mN}$ vérifie la condition de Mittag-Leffler. 
On en déduit que $\rH^i(\lambda^*(U),\bvcM)=0$ (\cite{jannsen} 1.15 et \cite{roos} 3.1).
La proposition s'ensuit puisque les ouverts affines de $\oX_1$ forment une famille topologiquement génératrice 
du site de Zariski de $\oX_1$

\begin{prop}\label{higgs3-coh6}
Pour tout entier $j\geq 0$, posons $\alpha_j=\frac{(d+1+j)(4d+1)+6d+1}{p-1}$. Alors~: 
\begin{itemize}
\item[{\rm (i)}] Pour tous entiers $i\geq 1$ et $j\geq 0$,
le $\co_{\bvoX}$-module $\rR^i\lambda_*(\rR^j\bvtau_*(\bvocB))$ 
est annulé par $p^{\alpha_j}\fm_\oK$.\item[{\rm (ii)}] Pour tout entier $q\geq 0$, le noyau et le conoyau du morphisme 
\begin{equation}\label{higgs3-coh7a}
\rR^q\top_*(\bvocB)\rightarrow \lambda_*(\rR^q\bvtau_*(\bvocB))
\end{equation}
induit par la suite spectrale de Cartan-Leray \eqref{higgs3-coh0e}, sont annulés par $p^{\sum_{j=0}^{q-1}\alpha_j}\fm_\oK$. 
\end{itemize}
\end{prop}

(i) En effet,  en vertu de \ref{higgs3-coh200}(i), 
on a un morphisme $\co_\bvoX$-linéaire de $X_{s,\zar}^{\mN^\circ}$
\begin{equation}
\xi^{-j}\tOmega^j_{\bvoX/\bvoS}\rightarrow \bvu_*(\rR^j\bvsigma_*(\bvocB)),
\end{equation}
dont le noyau est annulé par $p^{\frac{2d}{p-1}}\fm_\oK$ 
et le conoyau est annulé par $p^{\frac{4d+1}{p-1}}\fm_\oK$.
On en déduit, compte tenu de \ref{higgs3-coh5}, que $\rR^i\lambda_*(\bvu_*(\rR^j\bvsigma_*(\bvocB)))$ est annulé par 
$p^{\frac{6d+1}{p-1}}\fm_\oK$.
Par ailleurs, d'après \ref{higgs3-coh4}(ii), on a un morphisme $\co_\bvoX$-linéaire 
\begin{equation}
\rR^j\bvtau(\bvocB)\rightarrow \bvu_*(\rR^j\bvsigma_*(\bvocB))
\end{equation}
dont le noyau est annulé par $p^{\frac{(d+1)(4d+1)}{p-1}}\fm_\oK$
et le conoyau est annulé par $p^{\frac{j(4d+1)}{p-1}}\fm_\oK$. La proposition s'ensuit. 

(ii) La preuve est similaire à celle de \ref{higgs3-coh4}(ii). On peut clairement supposer $q\geq 1$. 
Considérons la suite spectrale de Cartan-Leray 
\begin{equation}\label{higgs3-coh7b}
\rE_2^{i,j}=\rR^i\lambda_*(\rR^j\bvtau_*(\bvocB))\Rightarrow \rR^{i+j}\top_*(\bvocB),
\end{equation}
et notons $(\rE^q_i)_{0\leq i\leq q}$ la filtration aboutissement sur $\rR^q\top_*(\bvocB)$.
On observera que $\rE^q_1$ est le noyau du morphisme \eqref{higgs3-coh7a}. 
D'après (i), $\rE_2^{i,j}$ est annulé par $p^{\alpha_j}\fm_\oK$ 
pour tous $i\geq 1$ et $j\geq 0$. Il en est alors de même de 
\begin{equation}\label{higgs3-coh7c}
\rE_i^q/\rE_{i+1}^q=\rE_\infty^{i,q-i}
\end{equation}
pour tout $i\geq 1$. Donc $\rE^q_1$ est annulé par $p^{\sum_{j=0}^{q-1}\alpha_j}\fm_\oK$. 
Par ailleurs, on a $\rE_\infty^{0,q}=\rE_{q+2}^{0,q}$, et le conoyau du morphisme \eqref{higgs3-coh7a}  
s'identifie au conoyau du composé des injections canoniques
\begin{equation}\label{higgs3-coh7d}
\rE_{q+2}^{0,q}\rightarrow \rE_{q+1}^{0,q}\rightarrow \dots \rightarrow \rE_2^{0,q}.
\end{equation}
Pour tout entier $m$ tel que $2\leq m\leq q+1$,  
le conoyau de l'injection canonique $\rE_{m+1}^{0,q}\rightarrow  \rE_{m}^{0,q}$
est annulé par $p^{\alpha_{q-m+1}}\fm_\oK$ d'après (i); d'où la proposition.

\begin{cor}\label{higgs3-coh8}
Il existe un et un unique isomorphisme de $\co_\fX[\frac 1 p]$-algèbres graduées 
\eqref{higgs3-formel1b}
\begin{equation}\label{higgs3-coh8a}
\wedge(\xi^{-1}\tOmega^1_{\fX/\cS}[\frac 1 p])\stackrel{\sim}{\rightarrow} \oplus_{i\geq 0}\rR^i\top_*(\bvocB)[\frac 1 p]
\end{equation}
dont la composante en degré un est le morphisme $\delta\otimes_{\mZ_p}\mQ_p$ \eqref{higgs3-coh1f}.
\end{cor} 
En effet, l'homomorphisme \eqref{higgs3-coh200a} induit un homomorphisme de $\co_\fX$-algèbres graduées 
\begin{equation}\label{higgs3-coh8b}
\wedge(\xi^{-1}\tOmega^1_{\fX/\cS})\rightarrow \oplus_{i\geq 0}\lambda_*(\bvu_*(\rR^i\bvsigma_*(\bvocB)))
\end{equation}
dont le noyau et le conoyau sont rig-nuls d'après \ref{higgs3-coh200}(i) (cf. \ref{higgs3-formel1}). Par ailleurs, les morphismes 
\begin{equation}\label{higgs3-coh8c}
\oplus_{i\geq 0}\rR^i\top_*(\bvocB)\rightarrow \oplus_{i\geq 0} \lambda_*(\rR^i\bvtau_*(\bvocB))\rightarrow
\oplus_{i\geq 0} \lambda_*(\bvu_*(\rR^i\bvsigma_*(\bvocB)))
\end{equation}
induits par les suites spectrales de Cartan-Leray, sont des homomorphismes de $\co_\fX$-algèbres graduées 
(cf. \cite{ega3} 0.12.2.6), dont les noyaux et conoyaux sont rig-nuls en vertu de \ref{higgs3-coh4}(ii) et \ref{higgs3-coh6}(ii). 
On obtient l'isomorphisme \eqref{higgs3-coh8a} recherché en appliquant le foncteur $-\otimes_{\mZ_p}\mQ_p$.

\begin{lem}\label{higgs3-coh9}
Soient $\cM$ un $\co_\fX$-module localement libre de type fini, $q$ un entier $\geq 0$.
Pour tout entier $n\geq 1$, posons $\cM_n=\cM\otimes_{\co_\fX}\co_{\oX_n}$
et $\bvcM=(\cM_{n+1})_{n\in \mN}$ que l'on considère aussi comme un faisceau de $X_{s,\et}^{\mN^\circ}$ 
{\rm (\cite{sga4} VII 2 c))}. On a alors un isomorphisme canonique 
\begin{equation}\label{higgs3-coh9a}
\cM\otimes_{\co_\fX} \rR^q\top_*(\bvocB)\stackrel{\sim}{\rightarrow}\rR^q\top_*(\bvsigma^*(\bvcM)).
\end{equation}
\end{lem}
En effet, on a un isomorphisme canonique $\bvu_*(\bvcM)\stackrel{\sim}{\rightarrow}\bvcM$.
Par ailleurs, en vertu de (\cite{egr1} 2.8.5), on a un isomorphisme canonique 
$\cM\stackrel{\sim}{\rightarrow}\lambda_*(\bvcM)$. 
Le morphisme d'adjonction $\bvcM\rightarrow \bvsigma_*(\bvsigma^*(\bvcM))$ induit alors un morphisme $\co_\fX$-linéaire
\begin{equation}\label{higgs3-coh9c}
\cM\rightarrow\top_*(\bvsigma^*(\bvcM)).
\end{equation}
On en déduit par cup-produit un morphisme $\co_\fX$-linéaire
\begin{equation}\label{higgs3-coh9d}
\cM\otimes_{\co_\fX} \rR^q\top_*(\bvocB)\rightarrow\rR^q\top_*(\bvsigma^*(\bvcM)).
\end{equation}
Pour voir que c'est un isomorphisme, on peut se ramener au cas où 
$\cM=\co_{\fX}$ (la question étant locale pour la topologie de Zariski de $\fX$), 
auquel cas l'assertion est immédiate.

\subsection{}\label{higgs3-coh10}
D'après (\cite{illusie1} I 4.3.1.7), la suite exacte canonique \eqref{higgs3-RGG21b} 
\begin{equation}\label{higgs3-coh10a}
0\rightarrow \bvocB\rightarrow \bvcF\rightarrow 
\bvsigma^*(\xi^{-1}\tOmega^1_{\bvoX/\bvoS})\rightarrow 0
\end{equation}
induit, pour tout entier $m\geq 1$, une suite exacte \eqref{higgs3-not9}
\begin{equation}\label{higgs3-coh10b}
0\rightarrow \rS^{m-1}_{\bvocB}(\bvcF)\rightarrow \rS^m_{\bvocB}(\bvcF)\rightarrow 
\bvsigma^*(\rS^m_{\co_{\bvoX}}(\xi^{-1}\tOmega^1_{\bvoX/\bvoS}))\rightarrow 0.
\end{equation}
On munit $\rS^m_{\bvocB}(\bvcF)$ de la filtration décroissante 
exhaustive $(\rS^{m-i}_{\bvocB}(\bvcF))_{i\in \mN}$. On a alors une suite exacte canonique
\begin{equation}\label{higgs3-coh10c}
0\rightarrow \bvsigma^*(\rS^{m-1}_{\co_{\bvoX}}(\xi^{-1}\tOmega^1_{\bvoX/\bvoS}))\rightarrow 
\rS^m_{\bvocB}(\bvcF)/\rS^{m-2}_{\bvocB}(\bvcF)\rightarrow 
\bvsigma^*(\rS^m_{\co_{\bvoX}}(\xi^{-1}\tOmega^1_{\bvoX/\bvoS}))\rightarrow 0.
\end{equation}

Pour tous entiers $i$ et $j$, posons
\begin{equation}\label{higgs3-coh10d}
\rE_1^{i,j}=\rR^{i+j}\top_*(\bvsigma^*(\rS^{-i}_{\co_{\bvoX}}(\xi^{-1}\tOmega^1_{\bvoX/\bvoS}))),
\end{equation}
et notons 
\begin{equation}\label{higgs3-coh10f}
d_1^{i,j}\colon \rE_1^{i,j}\rightarrow \rE_1^{i+1,j}
\end{equation} 
le morphisme induit par la suite exacte \eqref{higgs3-coh10c}.
D'après \ref{higgs3-coh9}, on a un isomorphisme $\co_{\fX}$-linéaire canonique
\begin{equation}\label{higgs3-coh10e}
\rS^{-i}(\xi^{-1}\tOmega^1_{\fX/\cS})\otimes_{\co_{\fX}}\rR^{i+j}\top_*(\bvocB)
\stackrel{\sim}{\rightarrow}\rE_1^{i,j}.
\end{equation}
On en déduit, compte tenu de \ref{higgs3-coh8}, un isomorphisme 
\begin{equation}\label{higgs3-coh10h}
\rS^{-i}(\xi^{-1}\tOmega^1_{\fX/\cS}[\frac 1 p])\otimes_{\co_{\fX}[\frac 1 p]}
\wedge^{i+j}(\xi^{-1}\tOmega^1_{\fX/\cS}[\frac 1 p])
\stackrel{\sim}{\rightarrow}\rE_1^{i,j}[\frac 1 p].
\end{equation}

\begin{prop}\label{higgs3-coh12}
On a un diagramme commutatif
\begin{equation}
\xymatrix{
{\rS^{-i}(\xi^{-1}\tOmega^1_{\fX/\cS}[\frac 1 p])\otimes_{\co_\fX[\frac 1 p]}\wedge^{i+j}(\xi^{-1}\tOmega^1_{\fX/\cS}[\frac 1 p])} 
\ar[d]_{\phi^{i,j}}\ar[r]^-(0.5)\sim&{\rE_1^{i,j}[\frac 1 p]}\ar[d]^{d_1^{i,j}\otimes_{\mZ_p}\mQ_p}\\
{\rS^{-i-1}(\xi^{-1}\tOmega^1_{\fX/\cS}[\frac 1 p])\otimes_{\co_\fX[\frac 1 p]}
\wedge^{i+j+1}(\xi^{-1}\tOmega^1_{\fX/\cS}[\frac 1 p])}
\ar[r]^-(0.5)\sim&{\rE_1^{i+1,j}[\frac 1 p]}}
\end{equation}
où $\phi^{i,j}$ est la restriction de la $\co_\fX[\frac 1 p]$-dérivation de 
\begin{equation}
\rS(\xi^{-1}\tOmega^1_{\fX/\cS}[\frac 1 p])\otimes_{\co_\fX[\frac 1 p]}\wedge(\xi^{-1}\tOmega^1_{\fX/\cS}[\frac 1 p])
\end{equation} 
définie dans \eqref{higgs3-vck1a} relativement au morphisme identique de $\xi^{-1}\tOmega^1_{\fX/\cS}[\frac 1 p]$,
et les flèches horizontales sont les isomorphismes \eqref{higgs3-coh10h}.
\end{prop}
En effet, en vertu de \ref{higgs3-vck7} et compte tenu du morphisme \eqref{higgs3-coh1c}, on a un diagramme commutatif
\begin{equation}
\xymatrix{
{\rS^{-i}(\xi^{-1}\tOmega^1_{\fX/\cS})\otimes_{\co_\fX}\rR^{i+j}\top_*(\bvocB)} 
\ar[d]_{\alpha\otimes \id}\ar[r]^-(0.5)\sim&{\rE_1^{i,j}}\ar[dd]^{d_1^{i,j}}\\
{\rS^{-i-1}(\xi^{-1}\tOmega^1_{\fX/\cS})\otimes_{\co_\fX}\rR^1\top_*(\bvocB)\otimes_{\co_\fX}\rR^{i+j}\top_*(\bvocB)}
\ar[d]_{\id\otimes \cup}&\\
{\rS^{-i-1}(\xi^{-1}\tOmega^1_{\fX/\cS})\otimes_{\co_\fX}\rR^{i+j+1}\top_*(\bvocB)}\ar[r]^-(0.5)\sim&{\rE_1^{i+1,j}}}
\end{equation}
où $\alpha$ est la restriction à $\rS^{-i}(\xi^{-1}\tOmega^1_{\fX/\cS})$ de la $\co_\fX$-dérivation $d_{\delta}$ de 
$\rS(\xi^{-1}\tOmega^1_{\fX/\cS})\otimes_{\co_\fX}\wedge(\rR^1\top_*(\bvocB))$ définie dans \eqref{higgs3-vck1a} 
relativement au morphisme $\delta$ \eqref{higgs3-coh1f}, 
$\cup$ est le cup-produit de la $\co_\fX$-algèbre $\oplus_{i\geq 0} \rR^i\top_*(\bvocB)$
et les flèches horizontales sont les isomorphismes \eqref{higgs3-coh10e}. 
La proposition s'ensuit.

\begin{prop}\label{higgs3-coh11}
Soit $m$ un entier $\geq 1$.  
Alors~:
\begin{itemize}
\item[{\rm (i)}] Le morphisme  
\begin{equation}\label{higgs3-coh11a}
\top_*(\rS^{m-1}_{\bvocB}(\bvcF))[\frac 1 p]\rightarrow \top_*(\rS^m_{\bvocB}(\bvcF))[\frac 1 p]
\end{equation}
induit par \eqref{higgs3-coh10b} est un isomorphisme.
\item[{\rm (ii)}] Pour tout entier $q\geq 1$, le morphisme 
\begin{equation}\label{higgs3-coh11b}
\rR^q\top_*(\rS^{m-1}_{\bvocB}(\bvcF))[\frac 1 p]\rightarrow \rR^q\top_*(\rS^m_{\bvocB}(\bvcF))[\frac 1 p]
\end{equation}
induit par \eqref{higgs3-coh10b} est nul.
\end{itemize}
\end{prop}

Pour tous entiers $i$ et $j$, on pose \eqref{higgs3-coh10d}
\begin{equation}\label{higgs3-coh11c}
{_m\rE}_1^{i,j}=\left\{
\begin{array}{clcr}
\rE_1^{i-m,j+m} &{\rm si} \ i\geq 0,\\
0 \ \ \ \ \ \ \ \ \ &{\rm si} \ i<0.
\end{array}
\right.
\end{equation}
On désigne par 
\begin{equation}\label{higgs3-coh11d}
{_m\rE}_1^{i,j}\Rightarrow \rR^{i+j}\top_*(\rS^{m}_{\bvocB}(\bvcF))
\end{equation}
la suite spectrale d'hypercohomologie du $\bvocB$-module filtré $\rS^m_{\bvocB}(\bvcF)$ \eqref{higgs3-coh10b},
dont les différentielles ${_md}_1^{i,j}$ sont données par \eqref{higgs3-coh10f}
\begin{equation}\label{higgs3-coh11dc}
{_md}_1^{i,j}=\left\{
\begin{array}{clcr}
d_1^{i-m,j+m}&{\rm si} \ i\geq 0,\\
0 \ \ \ \ \ \ \ \ &{\rm si} \ i<0.
\end{array}
\right.
\end{equation}
Pour tout entier $q\geq 0$, 
notons $({_m\rE}^q_i)_{i\in \mZ}$ la filtration aboutissement sur 
$\rR^q\top_*(\rS^m_{\bvocB}(\bvcF))$, de sorte que l'on a 
\begin{equation}
{_m\rE}_\infty^{i,q-i}={_m\rE}^q_i/{_m\rE}^q_{i+1}.
\end{equation} 
On a alors
\begin{equation}
{_m\rE}^q_i=\left\{
\begin{array}{clcr}
\rR^q\top_*(\rS^m_{\bvocB}(\bvcF))&{\rm si} \ i\leq  0,\\
0 & {\rm si}\  i\geq m+1,
\end{array}
\right.
\end{equation} 
et ${_m\rE}_\infty^{0,q}\subset {_m\rE}_1^{0,q}$. On déduit que l'image du morphisme canonique 
\begin{equation}\label{higgs3-coh11e}
\rR^{q}\top_*(\rS^{m-1}_{\bvocB}(\bvcF))\rightarrow \rR^{q}\top_*(\rS^m_{\bvocB}(\bvcF))
\end{equation}
est ${_m\rE}^q_1$, et son conoyau est ${_m\rE}_\infty^{0,q}$.

Par ailleurs, il résulte de \ref{higgs3-coh12} et \ref{higgs3-vck1} que pour tous entiers $i$ et $q$ 
vérifiant l'une des deux conditions suivantes~: 
\begin{itemize}
\item[(i)] $q=0$ et $i<m$,
\item[(ii)] $q\geq 1$ et $i\geq 1$, 
\end{itemize}
on a  
\begin{equation}
{_m\rE}_\infty^{i,q-i}[\frac 1 p]={_m\rE}_2^{i,q-i}[\frac 1 p]=0.
\end{equation}
La proposition s'ensuit.

\begin{prop}\label{higgs3-coh13}
Soient $r,r'$ deux nombres rationnels tels que $r>r'>0$. Alors~:
\begin{itemize}
\item[{\rm (i)}] Pour tout entier $n\geq 1$, l'homomorphisme canonique \eqref{higgs3-RGG24b}
\begin{equation}\label{higgs3-coh13a}
\co_{\oX_n}\rightarrow \sigma_{n*}(\cC_n^{(r)})
\end{equation}
est presque-injectif. Notons $\cH^{(r)}_n$ son conoyau. 
\item[{\rm (ii)}] Il existe un nombre rationnel $\alpha>0$ tel que pour tout entier $n\geq 1$, 
le morphisme 
\begin{equation}\label{higgs3-coh13b}
\cH^{(r)}_n\rightarrow \cH^{(r')}_n
\end{equation} 
induit par l'homomorphisme 
$\alpha_n^{r,r'}\colon \cC_n^{(r)}\rightarrow \cC_n^{(r')}$ \eqref{higgs3-RGG24d} soit annulé par $p^{\alpha}$. 
\item[{\rm (iii)}] Il existe un nombre rationnel $\beta>0$ tel que pour tous entiers $n,q\geq 1$, 
le morphisme canonique
\begin{equation}\label{higgs3-coh13c}
\rR^q\sigma_{n*}(\cC_n^{(r)})\rightarrow \rR^q\sigma_{n*}(\cC_n^{(r')})
\end{equation}
soit annulé par $p^{\beta}$.
\end{itemize}
\end{prop}

Soient $n$ et $q$ deux entiers tels que $n\geq 1$ et $q\geq 0$, $\ox$ un point géométrique de $X$ au-dessus de $s$,
$X'$ le localisé strict de $X$ en $\ox$, 
\begin{equation}\label{higgs3-coh13d}
\varphi_\ox\colon \tE\rightarrow \oX'^\circ_\fet
\end{equation}
le foncteur \eqref{higgs3-tfa6d}. D'après (\cite{ag2} 10.30) et \eqref{higgs3-TFT8e}, on a un isomorphisme canonique 
\begin{equation}\label{higgs3-coh13f}
(\rR^q\sigma_{n*}(\cC^{(r)}_n))_\ox\stackrel{\sim}{\rightarrow} \rH^q(\oX'^\circ_\fet, \varphi_\ox(\cC^{(r)}_n)).
\end{equation}
En vertu de \ref{higgs3-sli6}, $\oX'$ est normal et strictement local (et en particulier intègre). 
Comme $\oX'$ est $\oS$-plat, $\oX'^\circ$ est intègre et non-vide. 
Soient $v\colon \oy\rightarrow \oX'^\circ$ un point géométrique générique de $\oX'^\circ$,
\begin{equation}\label{higgs3-coh13g}
\nu_{\oy}\colon \oX'^\circ_\fet \stackrel{\sim}{\rightarrow}\bB_{\pi_1(\oX'^\circ,\oy)}
\end{equation}
le foncteur fibre associé  \eqref{higgs3-not6c}.  On note encore $\oy$ le point géométrique de $\oX^\circ$ et 
$u\colon \oy\rightarrow X'$ le morphisme induits par $v$. 
On obtient ainsi un point $(\oy\rightsquigarrow \ox)$ de $X_\et\gtimes_{X_\et}\oX^\circ_\et$. 

On désigne par $\fV_\ox$ (resp. $\fV_\ox(\bQ)$) la catégorie des voisinages 
du point de $X_\et$ associé à $\ox$ dans le site $\Et_{/X}$ (resp. $\bQ$ \eqref{higgs3-RGG85}). 
Pour tout $(U,\fp\colon \ox\rightarrow U)\in \ob(\fV_\ox)$, on note encore $\fp\colon X'\rightarrow U$ le morphisme déduit de 
$\fp$, et on pose
\begin{equation}\label{higgs3-coh13e}
\ofp^\circ=\fp\times_X\oX^\circ\colon \oX'^\circ \rightarrow \oU^\circ.
\end{equation}
On note aussi (abusivement) $\oy$ le point géométrique $\ofp^\circ(v(\oy))$ de $\oU^\circ$.
On observera que $\oy$ est localisé en un point générique de $\oU^\circ$ car $\ofp^\circ$ est plat. 
Comme $\oU$ est localement irréductible \eqref{higgs3-sli3}, 
il est la somme des schémas induits sur ses composantes irréductibles. 
Notons $\oU^\star$  la composante irréductible de $\oU$ contenant $\oy$. De même, $\oU^\circ$
est la somme des schémas induits sur ses composantes irréductibles, et $\oU^{\star \circ}=\oU^\star\times_XX^\circ$ 
est la composante irréductible de $\oU^\circ$ contenant $\oy$. 
Le morphisme $\ofp^\circ$ se factorise à travers $\oU^{\star \circ}$.
D'après \eqref{higgs3-RGG17a}, on a un isomorphisme canonique de $\bB_{\pi_1(\oX'^\circ,\oy)}$
\begin{equation}\label{higgs3-coh13h}
\nu_{\oy}(\varphi_\ox(\cC_n^{(r)}))\stackrel{\sim}{\rightarrow}
\underset{\underset{(U,\fp)\in \fV_\ox(\bQ)^\circ}{\longrightarrow}}{\lim}\ \cC^{\oy,(r)}_U/p^n\cC^{\oy,(r)}_U.
\end{equation}
En vertu de (\cite{ag2} 11.10 et (9.7.6)), celui-ci induit un isomorphisme 
\begin{equation}\label{higgs3-coh13i}
\rH^q(\oX'^\circ_\fet, \varphi_\ox(\cC_n^{(r)}))
\stackrel{\sim}{\rightarrow} \underset{\underset{(U,\fp)\in \fV_\ox(\bQ)^\circ}{\longrightarrow}}{\lim}\
\rH^q(\pi_1(\oU^{\star \circ},\oy), \cC^{\oy,(r)}_U/p^n\cC^{\oy,(r)}_U).
\end{equation} 
On déduit de \eqref{higgs3-coh13f} et \eqref{higgs3-coh13i} un isomorphisme  
\begin{equation}\label{higgs3-coh13j}
(\rR^q\sigma_{n*}(\cC^{(r)}_n))_\ox\stackrel{\sim}{\rightarrow} 
\underset{\underset{(U,\fp)\in \fV_\ox(\bQ)^\circ}{\longrightarrow}}{\lim}\
\rH^q(\pi_1(\oU^{\star \circ},\oy), \cC^{\oy,(r)}_U/p^n\cC^{\oy,(r)}_U).
\end{equation}

D'autre part, on démontre, comme dans \eqref{higgs3-coh2gb}, qu'on a un isomorphisme canonique 
\begin{equation}\label{higgs3-coh13k}
(\co_{\oX_n})_\ox\stackrel{\sim}{\rightarrow} \underset{\underset{(U,\fp)\in \fV_\ox^\circ}{\longrightarrow}}{\lim}\
\co_{\oX_n}(\oU^\star), 
\end{equation}
où l'on considère $\co_{\oX_n}$ à gauche comme un faisceau de $X_{s,\et}$ et à droite comme un faisceau de $\oX_\et$. 

La fibre du morphisme \eqref{higgs3-coh13a} en $\ox$ s'identifie à la limite inductive sur la catégorie
$\fV_\ox(\bQ)^\circ$ du morphisme canonique
\begin{equation}\label{higgs3-coh13l}
\co_{\oX_n}(\oU^\star)\rightarrow (\cC^{\oy,(r)}_U/p^n\cC^{\oy,(r)}_U)^{\pi_1(\oU^{\star \circ},\oy)}.
\end{equation}
Comme les limites inductives filtrantes sont exactes, la proposition résulte alors de (\cite{ag1} 12.7),
compte tenu de \ref{higgs3-RGG65} et de la preuve de \ref{higgs3-RGG105}. 

\begin{cor}\label{higgs3-coh14}
Soient $r,r'$ deux nombres rationnels tels que $r>r'>0$. Alors~:
\begin{itemize}
\item[{\rm (i)}] L'homomorphisme canonique de $X_{s,\et}^{\mN^\circ}$ \eqref{higgs3-RGG21}
\begin{equation}\label{higgs3-coh14a}
\co_{\bvoX}\rightarrow \bvsigma_*(\bvcC^{(r)})
\end{equation}
est presque-injectif. Notons $\bvcH^{(r)}$ son conoyau. 
\item[{\rm (ii)}] Il existe un nombre rationnel $\alpha>0$ tel que le morphisme 
\begin{equation}\label{higgs3-coh14b}
\bvcH^{(r)}\rightarrow \bvcH^{(r')}
\end{equation} 
induit par l'homomorphisme canonique $\bvalpha^{r,r'}\colon \bvcC^{(r)}\rightarrow \bvcC^{(r')}$ \eqref{higgs3-RGG180c}
soit annulé par $p^{\alpha}$. 
\item[{\rm (iii)}] Il existe un nombre rationnel $\beta>0$ tel que pour tout entier $q\geq 1$, 
le morphisme canonique de $X_{s,\et}^{\mN^\circ}$
\begin{equation}\label{higgs3-coh14c}
\rR^q\bvsigma_*(\bvcC^{(r)})\rightarrow \rR^q\bvsigma_*(\bvcC^{(r')})
\end{equation}
soit annulé par $p^{\beta}$.
\end{itemize}
\end{cor}

Cela résulte de \ref{higgs3-coh13}, \ref{higgs3-spsa99}(i) et \eqref{higgs3-spsa1h}.

\begin{lem}\label{higgs3-coh15}
Soient $r,r'$ deux nombres rationnels tels que $r>r'>0$. Alors~:
\begin{itemize}
\item[{\rm (i)}] L'homomorphisme canonique de $X_{s,\zar}^{\mN^\circ}$ \eqref{higgs3-RGG21c}
\begin{equation}\label{higgs3-coh15a}
\co_{\bvoX}\rightarrow \bvtau_*(\bvcC^{(r)})
\end{equation}
est presque-injectif. Notons $\bvcK^{(r)}$ son conoyau. 
\item[{\rm (ii)}] Il existe un nombre rationnel $\alpha>0$ tel que le morphisme 
\begin{equation}\label{higgs3-coh15b}
\bvcK^{(r)}\rightarrow \bvcK^{(r')}
\end{equation} 
induit par l'homomorphisme canonique $\bvalpha^{r,r'}\colon \bvcC^{(r)}\rightarrow \bvcC^{(r')}$ \eqref{higgs3-RGG180c} 
soit annulé par $p^{\alpha}$. 
\item[{\rm (iii)}] Pour tout entier $q\geq 1$, il existe un nombre rationnel $\beta>0$ tel que  
le morphisme canonique de $X_{s,\zar}^{\mN^\circ}$
\begin{equation}\label{higgs3-coh15c}
\rR^q\bvtau_*(\bvcC^{(r)})\rightarrow \rR^q\bvtau_*(\bvcC^{(r')})
\end{equation}
soit annulé par $p^{\beta}$.
\end{itemize}
\end{lem}

Reprenons les notations de \ref{higgs3-coh14}, de plus,
notons $\bvcN^{(r)}$ et $\bvcM^{(r)}$ les noyaux des morphismes \eqref{higgs3-coh14a} et \eqref{higgs3-coh15a}, respectivement.

(i) Comme $\bvcM^{(r)}=\bvu_*(\bvcN^{(r)})$, la proposition résulte de \ref{higgs3-coh14}(i).

(ii) Comme $\rR^1\bvu_*(\co_\bvoX)=0$ d'après \eqref{higgs3-spsa1h} et (\cite{sga4} VII 4.3), on a une suite exacte
\begin{equation}\label{higgs3-coh15d}
0\rightarrow \rR^1\bvu_*(\bvcN^{(r)})\rightarrow \bvcK^{(r)}\rightarrow \bvu_*(\bvcH^{(r)}).
\end{equation}
La proposition résulte alors de \ref{higgs3-coh14}(i)-(ii).  

(iii) Considérons la suite spectrale de Cartan-Leray 
\begin{equation}
{^r \rE}_2^{i,j}=\rR^i\bvu_*(\rR^j\bvsigma_*(\bvcC^{(r)}))\Rightarrow \rR^{i+j}\bvtau_*(\bvcC^{(r)}), 
\end{equation}
et notons $({^r \rE}^q_i)_{0\leq i\leq q}$ la filtration aboutissement sur $\rR^q\bvtau_*(\bvcC^{r)})$, de sorte que l'on a 
\begin{equation}
{^r \rE}_i^q/{^r \rE}_{i+1}^q={^r \rE}_\infty^{i,q-i}.
\end{equation}
Pour tout entier $0\leq i\leq q+1$, posons $r_i=r'+(q+1-i)(r-r')/(q+1)$. D'après \ref{higgs3-coh14}(iii), pour tout entier $0\leq i\leq q-1$, 
il existe un nombre rationnel $\beta_i>0$ tel que le morphisme canonique 
\begin{equation}
{^{r_{i}} \rE}_2^{i,q-i}\rightarrow {^{r_{i+1}} \rE}_2^{i,q-i}
\end{equation} 
soit annulé par $p^{\beta_i}$. Il en est alors de même du morphisme 
${^{r_{i}} \rE}_\infty^{i,q-i}\rightarrow {^{r_{i+1}} \rE}_\infty^{i,q-i}$. Par ailleurs, $\rR^q\bvu_*(\co_\bvoX)=0$
d'après \eqref{higgs3-spsa1h} et (\cite{sga4} VII 4.3). On en déduit que le morphisme canonique 
\begin{equation}
\rR^q\bvu_*(\bvsigma_*(\bvcC^{(r')}))\rightarrow \rR^q\bvu_*(\bvcH^{(r')})
\end{equation}
est presque-injectif d'après \ref{higgs3-coh14}(i). Par suite, en vertu de \ref{higgs3-coh14}(ii), il existe un nombre rationnel $\beta_q>0$ 
tel que le morphisme canonique 
\begin{equation}
{^{r_q} \rE}_2^{q,0}\rightarrow {^{r_{q+1}} \rE}_2^{q,0}
\end{equation} 
soit annulé par $p^{\beta_q}$. Il en est alors de même du morphisme 
${^{r_q} \rE}_\infty^{q,0}\rightarrow {^{r_{q+1}} \rE}_\infty^{q,0}$.
La proposition s'ensuit en prenant $\beta=\sum_{i=0}^q\beta_i$.

\begin{prop}\label{higgs3-coh16}
Soient $r,r'$ deux nombres rationnels tels que $r>r'>0$. Alors~:
\begin{itemize}
\item[{\rm (i)}] L'homomorphisme canonique \eqref{higgs3-RGG21c}
\begin{equation}\label{higgs3-coh16a}
\co_{\fX}\rightarrow \top_*(\bvcC^{(r)})
\end{equation}
est injectif. Notons $\cL^{(r)}$ son conoyau. 
\item[{\rm (ii)}] Il existe un nombre rationnel $\alpha>0$ tel que le morphisme 
\begin{equation}\label{higgs3-coh16b}
\cL^{(r)}\rightarrow \cL^{(r')}
\end{equation} 
induit par l'homomorphisme canonique $\bvalpha^{r,r'}\colon \bvcC^{(r)}\rightarrow \bvcC^{(r')}$ \eqref{higgs3-RGG180c}
soit annulé par $p^{\alpha}$. 
\item[{\rm (iii)}] Pour tout entier $q\geq 1$, il existe un nombre rationnel $\beta>0$ tel que  
le morphisme canonique 
\begin{equation}\label{higgs3-coh16c}
\rR^q\top_*(\bvcC^{(r)})\rightarrow \rR^q\top_*(\bvcC^{(r')})
\end{equation}
soit annulé par $p^{\beta}$.
\end{itemize}
\end{prop}

Reprenons les notations de \ref{higgs3-coh15}, de plus, notons $\bvcM^{(r)}$ le noyau du morphisme \eqref{higgs3-coh15a}.

(i) Le noyau du morphisme \eqref{higgs3-coh16a} est canoniquement isomorphe à $\lambda_*(\bvcM^{(r)})$ \eqref{higgs3-coh0d}.
Il est donc presque-nul en vertu de \ref{higgs3-coh15}(i). Comme $\co_\fX$ est $\co_C$-plat  d'après \ref{higgs3-slad4}(i)
(rig-pur dans la terminologie de \cite{egr1} 2.10.1.4), le morphisme \eqref{higgs3-coh16a} est injectif.

(ii) Comme $\rR^1\lambda_*(\co_\bvoX)=0$ d'après \ref{higgs3-coh5}, on a une suite exacte
\begin{equation}\label{higgs3-coh16d}
0\rightarrow \rR^1\lambda_*(\bvcM^{(r)})\rightarrow \cL^{(r)}\rightarrow \lambda_*(\bvcK^{(r)}).
\end{equation}
La proposition résulte alors de \ref{higgs3-coh15}(i)-(ii).  

(iii) La preuve est similaire à celle de \ref{higgs3-coh15}(iii). Considérons la suite spectrale de Cartan-Leray 
\begin{equation}
{^r \rE}_2^{i,j}=\rR^i\lambda_*(\rR^j\bvtau_*(\bvcC^{(r)}))\Rightarrow \rR^{i+j}\top_*(\bvcC^{(r)}), 
\end{equation}
et notons $({^r \rE}^q_i)_{0\leq i\leq q}$ la filtration aboutissement sur $\rR^q\top_*(\bvcC^{r)})$, de sorte que l'on a 
\begin{equation}
{^r \rE}_i^q/{^r \rE}_{i+1}^q={^r \rE}_\infty^{i,q-i}.
\end{equation}
Pour tout entier $0\leq i\leq q+1$, posons $r_i=r'+(q+1-i)(r-r')/(q+1)$. D'après \ref{higgs3-coh15}(iii), pour tout entier $0\leq i\leq q-1$, 
il existe un nombre rationnel $\beta_i>0$ tel que le morphisme canonique 
\begin{equation}
{^{r_{i}} \rE}_2^{i,q-i}\rightarrow {^{r_{i+1}} \rE}_2^{i,q-i}
\end{equation} 
soit annulé par $p^{\beta_i}$. Il en est alors de même du morphisme 
${^{r_{i}} \rE}_\infty^{i,q-i}\rightarrow {^{r_{i+1}} \rE}_\infty^{i,q-i}$. Par ailleurs, $\rR^q\bvu_*(\co_\bvoX)=0$
d'après \ref{higgs3-coh5}. On en déduit que le morphisme canonique 
\begin{equation}
\rR^q\lambda_*(\bvtau_*(\bvcC^{(r')}))\rightarrow \rR^q\lambda_*(\bvcK^{(r')})
\end{equation}
est presque-injectif d'après \ref{higgs3-coh15}(i). Par suite, en vertu de \ref{higgs3-coh15}(ii), il existe un nombre rationnel $\beta_q>0$ 
tel que le morphisme canonique 
\begin{equation}
{^{r_q} \rE}_2^{q,0}\rightarrow {^{r_{q+1}} \rE}_2^{q,0}
\end{equation} 
soit annulé par $p^{\beta_q}$. Il en est alors de même du morphisme 
${^{r_q} \rE}_\infty^{q,0}\rightarrow {^{r_{q+1}} \rE}_\infty^{q,0}$.
La proposition s'ensuit en prenant $\beta=\sum_{i=0}^q\beta_i$.

\begin{cor}\label{higgs3-coh25}
Soient $r$, $r'$ deux nombres rationnels tels que $r>r'>0$. 
\begin{itemize}
\item[{\rm (i)}] L'homomorphisme canonique 
\begin{equation}\label{higgs3-coh25a}
u^r\colon \co_{\fX}[\frac 1 p]\rightarrow \top_*(\bvcC^{(r)})[\frac 1 p]
\end{equation}
admet (en tant que morphisme $\co_\fX[\frac 1 p]$-linéaire) un inverse à gauche canonique
\begin{equation}
v^r\colon \top_*(\bvcC^{(r)})[\frac 1 p]\rightarrow \co_{\fX}[\frac 1 p].
\end{equation}
\item[{\rm (ii)}] Le composé 
\begin{equation}
\top_*(\bvcC^{(r)})[\frac 1 p]\stackrel{v^r}{\longrightarrow} \co_{\fX}[\frac 1 p]
\stackrel{u^{r'}}{\longrightarrow} \top_*(\bvcC^{(r')})[\frac 1 p]
\end{equation}
est l'homomorphisme canonique.
\item[{\rm (iii)}] Pour tout entier $q\geq 1$, le morphisme canonique 
\begin{equation}\label{higgs3-coh25c}
\rR^q\top_*(\bvcC^{(r)})[\frac 1 p]\rightarrow \rR^q\top_*(\bvcC^{(r')})[\frac 1 p]
\end{equation}
est nul. 
\end{itemize}
\end{cor}

En effet, d'après \ref{higgs3-coh16}(i)-(ii), $u^r$ est injectif et il existe un et un unique morphisme $\co_\fX[\frac 1 p]$-linéaire
\begin{equation}
v^{r,r'}\colon  \top_*(\bvcC^{(r)})[\frac 1 p]\rightarrow  \co_\fX[\frac 1 p]
\end{equation}
tel que $u^{r'}\circ v^{r,r'}$ soit l'homomorphisme canonique $\top_*(\bvcC^{(r)})[\frac 1 p]\rightarrow 
\top_*(\bvcC^{(r')})[\frac 1 p]$. Comme on a 
$u^{r'}\circ v^{r,r'}\circ u^r=u^{r'}$, on en déduit que $v^{r,r'}$ est un inverse à gauche de $u^r$.  
On vérifie aussitôt qu'il ne dépend pas de $r'$; d'où les propositions (i) et (ii).  
La proposition (iii) résulte aussitôt de \ref{higgs3-coh16}(iii).

\begin{cor}\label{higgs3-coh17}
L'homomorphisme canonique 
\begin{equation}\label{higgs3-coh17a}
\co_{\fX}[\frac 1 p]\rightarrow \underset{\underset{r\in \mQ_{>0}}{\longrightarrow}}{\lim}\  \top_*(\bvcC^{(r)})[\frac 1 p]
\end{equation}
est un isomorphisme, et pour tout entier $q\geq 1$,
\begin{equation}\label{higgs3-coh17b}
\underset{\underset{r\in \mQ_{>0}}{\longrightarrow}}{\lim}\ \rR^q\top_*(\bvcC^{(r)})[\frac 1 p] =0.
\end{equation}
\end{cor}

\subsection{}\label{higgs3-coh18}
Pour tout nombre rationnel $r\geq 0$ et tout entier $n\geq 1$, on rappelle que 
la $\ocB_n$-dérivation universelle de $\cC_n^{(r)}$ \eqref{higgs3-RGG18e} 
\begin{equation}
d_n^{(r)}\colon \cC_n^{(r)}\rightarrow \sigma_n^*(\xi^{-1}\tOmega^1_{\oX_n/\oS_n})\otimes_{\ocB_n}\cC_n^{(r)}
\end{equation}
est un $\ocB_n$-champ de Higgs à coefficients dans $\sigma_n^*(\xi^{-1}\tOmega^1_{\oX_n/\oS_n})$ \eqref{higgs3-RGG181}. 
On note $\mK^\bullet(\cC^{(r)}_n,p^rd^{(r)}_n)$ le complexe de Dolbeault
du $\ocB_n$-module de Higgs $(\cC^{(r)}_n,p^rd^{(r)}_n)$ (\cite{ag1} 2.8.2)
et $\tmK^\bullet(\cC^{(r)}_n,p^rd^{(r)}_n)$ le complexe de Dolbeault augmenté
\begin{equation}\label{higgs3-coh18a}
\ocB_n\rightarrow \mK^0(\cC^{(r)}_n,p^rd^{(r)}_n)\rightarrow \mK^1(\cC^{(r)}_n,p^rd^{(r)}_n)\rightarrow 
\mK^2(\cC^{(r)}_n,p^rd^{(r)}_n)\rightarrow \dots,
\end{equation}
où $\ocB_n$ est placé en degré $-1$ et la différentielle $\ocB_n\rightarrow \cC^{(r)}_n$ est l'homomorphisme canonique. 

Pour tous nombres rationnels $r\geq r'\geq 0$, on a \eqref{higgs3-RGG18i}
\begin{equation}\label{higgs3-coh18b}
p^r(\id \otimes \alpha^{r,r'}_n) \circ d^{(r)}_n=p^{r'}d^{(r')}_n\circ \alpha^{r,r'}_n,
\end{equation}
où $\alpha_n^{r,r'}\colon \cC_n^{(r)}\rightarrow \cC_n^{(r')}$ est l'homomorphisme \eqref{higgs3-RGG24d}.
Par suite, $\alpha_n^{r,r'}$ induit un morphisme 
\begin{equation}\label{higgs3-coh18c}
\nu_n^{r,r'}\colon \tmK^\bullet(\cC^{(r)}_n,p^rd^{(r)}_n)\rightarrow \tmK^\bullet(\cC^{(r')}_n,p^{r'}d^{(r')}_n).
\end{equation}

\begin{prop}\label{higgs3-coh19}
Pour tous nombres rationnels $r>r'>0$, 
il existe un nombre rationnel $\alpha\geq 0$ tel que pour tous entiers $n$ et $q$ avec $n\geq 1$,   
le morphisme 
\begin{equation}\label{higgs3-coh19a}
\rH^q(\nu_n^{r,r'})\colon \rH^q(\tmK^\bullet(\cC^{(r)}_n,p^rd^{(r)}_n))\rightarrow 
\rH^q(\tmK^\bullet(\cC^{(r')}_n,p^{r'}d^{(r')}_n))
\end{equation}
soit annulé par $p^\alpha$. 
\end{prop}

Soient $n$ un entier $\geq 1$, $\ox$ un point géométrique de $X$ au-dessus de $s$, $X'$ le localisé strict de $X$ en $\ox$,
\begin{equation}\label{higgs3-coh19b}
\varphi_\ox\colon \tE\rightarrow \oX'^\circ_\fet
\end{equation}
le foncteur \eqref{higgs3-tfa6d}. 
En vertu de \ref{higgs3-sli6}, $\oX'$ est normal et strictement local (et en particulier intègre). 
Comme $\oX'$ est $\oS$-plat, $\oX'^\circ$ est intègre et non-vide. 
Soient $v\colon \oy\rightarrow \oX'^\circ$ un point géométrique générique,  
\begin{equation}\label{higgs3-coh19d}
\nu_{\oy}\colon \oX'^\circ_\fet \stackrel{\sim}{\rightarrow}\bB_{\pi_1(\oX'^\circ,\oy)}
\end{equation}
le foncteur fibre associé  \eqref{higgs3-not6c}. 
On note encore $\oy$ le point géométrique de $\oX^\circ$ et 
$u\colon \oy\rightarrow X'$ le morphisme induits par $v$. 
On obtient ainsi un point $(\oy\rightsquigarrow \ox)$ de $X_\et\gtimes_{X_\et}\oX^\circ_\et$. 

On désigne par  $\fV_\ox$ (resp. $\fV_\ox(\bQ)$) la catégorie des voisinages 
du point de $X_\et$ associé à $\ox$ dans le site $\Et_{/X}$ (resp. $\bQ$ \eqref{higgs3-RGG85}).
Pour tout $(U,\fp\colon \ox\rightarrow U)\in \ob(\fV_\ox)$, on note encore $\fp\colon X'\rightarrow U$ le morphisme déduit de 
$\fp$, et on pose
\begin{equation}\label{higgs3-coh19c}
\ofp^\circ=\fp\times_X\oX^\circ\colon \oX'^\circ \rightarrow \oU^\circ.
\end{equation}
On note aussi (abusivement) $\oy$ le point géométrique $\ofp^\circ(v(\oy))$ de $\oU^\circ$.
On observera que $\oy$ est localisé en un point générique de $\oU^\circ$ car $\ofp^\circ$ est plat. 
Comme $\oU$ est localement irréductible \eqref{higgs3-sli3}, 
il est la somme des schémas induits sur ses composantes irréductibles. 
Notons $\oU^\star$  la composante irréductible de $\oU$ contenant $\oy$. De même, $\oU^\circ$
est la somme des schémas induits sur ses composantes irréductibles, et $\oU^{\star \circ}=\oU^\star\times_XX^\circ$ 
est la composante irréductible de $\oU^\circ$ contenant $\oy$. 
Le morphisme $\ofp^\circ$ se factorise à travers $\oU^{\star \circ}$.

D'après \eqref{higgs3-RGG3e} et \eqref{higgs3-RGG17a}, on a des isomorphismes canoniques de $\bB_{\pi_1(\oX'^\circ,\oy)}$
\begin{eqnarray}
\nu_{\oy}(\varphi_\ox(\ocB_n))&\stackrel{\sim}{\rightarrow}& 
\underset{\underset{(U,\fp)\in \fV_\ox^\circ}{\longrightarrow}}{\lim}\ \oR^{\oy}_U/p^n\oR^{\oy}_U,\label{higgs3-coh19f}\\
\nu_{\oy}(\varphi_\ox(\cC_n^{(r)}))&\stackrel{\sim}{\rightarrow}&
\underset{\underset{(U,\fp)\in \fV_\ox(\bQ)^\circ}{\longrightarrow}}{\lim}\ \cC^{\oy,(r)}_U/p^n\cC^{\oy,(r)}_U. \label{higgs3-coh19e}
\end{eqnarray}
Les anneaux sous-jacents à ces représentations sont canoniquement isomorphes aux 
fibres de $\ocB_n$ et $\cC_n^{(r)}$ en $\rho(\oy \rightsquigarrow \ox)$ (\cite{ag2} 10.31 et 9.9).
D'autre part, on a des isomorphismes canoniques \eqref{higgs3-coh13k} et \eqref{higgs3-coh2gb} 
\begin{eqnarray}
(\co_{\oX_n})_\ox\stackrel{\sim}{\rightarrow} \underset{\underset{(U,\fp)\in \fV_\ox^\circ}{\longrightarrow}}{\lim}\
\co_{\oX_n}(\oU^\star), \label{higgs3-coh19hh}\\
(\tOmega^1_{\oX_n/\oS_n})_\ox\stackrel{\sim}{\rightarrow} \underset{\underset{(U,\fp)\in \fV_\ox^\circ}{\longrightarrow}}{\lim}\
\tOmega^1_{\oX_n/\oS_n}(\oU^\star),\label{higgs3-coh19h}
\end{eqnarray}
où l'on considère $\co_{\oX_n}$ et $\tOmega^1_{\oX_n/\oS_n}$
à gauche comme des faisceaux de $X_{s,\et}$ et à droite comme des faisceaux de $\oX_\et$. 
Ces modules sont canoniquement isomorphes aux fibres de $\sigma_s^*(\co_{\oX_n})$ et 
$\sigma_s^*(\tOmega^1_{\oX_n/\oS_n})$ en $\rho(\oy \rightsquigarrow \ox)$ (\cite{ag2} (10.18.1)).
Il résulte de \ref{higgs3-RGG171} que 
la fibre de la dérivation $d^{(r)}_n$ \eqref{higgs3-RGG18e} en $\rho(\oy \rightsquigarrow \ox)$ s'identifie à la limite inductive 
sur la catégorie $\fV_\ox(\bQ)^\circ$ des $(\oR_U^\oy/p^n\oR_U^\oy)$-dérivations universelles 
\begin{equation}
\cC^{\oy,(r)}_U/p^n\cC^{\oy,(r)}_U\rightarrow \xi^{-1}\tOmega^1_{X/S}(U)\otimes_{\co_X(U)} 
(\cC^{\oy,(r)}_U/p^n\cC^{\oy,(r)}_U).
\end{equation}

D'après \ref{higgs3-TFT5}, la famille des points $\rho(\oy\rightsquigarrow \ox)$  de $\tE_s$ est conservative. 
Comme les limites inductives filtrantes sont exactes, 
la proposition résulte alors de (\cite{ag1} 12.3(i)), compte tenu de \ref{higgs3-RGG65} et de la preuve de \ref{higgs3-RGG105}.

\subsection{}\label{higgs3-coh20}\index{10001075@$\bvd^{(r)}$|textbf}
Pour tout nombre rationnel $r\geq 0$, on note encore 
\begin{equation}\label{higgs3-coh20a}
\bvd^{(r)}\colon \bvcC^{(r)}\rightarrow \top^*(\xi^{-1}\tOmega^1_{\fX/\cS})\otimes_{\bvocB}\bvcC^{(r)}
\end{equation}
la $\bvocB$-dérivation induite par $\bvd^{(r)}$ \eqref{higgs3-RGG180b} et l'isomorphisme \eqref{higgs3-coh1cc},
que l'on identifie à la $\bvocB$-dérivation universelle de $\bvcC^{(r)}$.
C'est un $\bvocB$-champ de Higgs à coefficients dans $\top^*(\xi^{-1}\tOmega^1_{\fX/\cS})$ \eqref{higgs3-RGG181}.
On désigne par $\mK^\bullet(\bvcC^{(r)},p^r\bvd^{(r)})$ le complexe de Dolbeault
du $\bvocB$-module de Higgs $(\bvcC^{(r)},p^r\bvd^{(r)})$
et par $\tmK^\bullet(\bvcC^{(r)},p^r\bvd^{(r)})$ le complexe de Dolbeault augmenté
\begin{equation}\label{higgs3-coh20d}
\bvocB\rightarrow \mK^0(\bvcC^{(r)},p^r\bvd^{(r)})\rightarrow \mK^1(\bvcC^{(r)},p^r\bvd^{(r)})\rightarrow \dots 
\rightarrow \mK^n(\bvcC^{(r)},p^r\bvd^{(r)})\rightarrow \dots,
\end{equation}
où $\bvocB$ est placé en degré $-1$ et la différentielle $\bvocB\rightarrow \bvcC^{(r)}$ est l'homomorphisme canonique. 

Pour tous nombres rationnels $r\geq r'\geq 0$, on a \eqref{higgs3-RGG180d} 
\begin{equation}\label{higgs3-coh20e}
p^r(\id \otimes \bvalpha^{r,r'}) \circ \bvd^{(r)}=p^{r'}\bvd^{(r')}\circ \bvalpha^{r,r'},
\end{equation}
où $\bvalpha^{r,r'}\colon \bvcC^{(r)}\rightarrow \bvcC^{(r')}$ est l'homomorphisme \eqref{higgs3-RGG180c}. Par suite, 
$\bvalpha^{r,r'}$ induit un morphisme de complexes 
\begin{equation}\label{higgs3-coh20f}
\bvnu^{r,r'}\colon \tmK^\bullet(\bvcC^{(r)},p^r\bvd^{(r)})\rightarrow \tmK^\bullet(\bvcC^{(r')},p^{r'}\bvd^{(r')}).
\end{equation}

On note $\bMod_{\mQ}(\bvocB)$ la catégorie des $\bvocB$-modules de $\tE_s^{\mN^\circ}$ à isogénie près \eqref{higgs3-caip4},
et on désigne par $\mK^\bullet_\mQ(\bvcC^{(r)},p^r\bvd^{(r)})$ et $\tmK^\bullet_\mQ(\bvcC^{(r)},p^r\bvd^{(r)})$ 
les images des complexes 
$\mK^\bullet(\bvcC^{(r)},p^r\bvd^{(r)})$ et $\tmK^\bullet(\bvcC^{(r)},p^r\bvd^{(r)})$ dans $\bMod_{\mQ}(\bvocB)$.

\begin{prop}\label{higgs3-coh21}
Pour tous nombres rationnels  $r> r'> 0$ et tout entier $q$, le morphisme canonique \eqref{higgs3-coh20f} 
\begin{equation}\label{higgs3-coh21a}
\rH^q(\bvnu^{r,r'}_\mQ)\colon \rH^q(\tmK^\bullet_\mQ(\bvcC^{(r)},p^r\bvd^{(r)}))\rightarrow 
\rH^q(\tmK^\bullet_\mQ(\bvcC^{(r')},p^{r'}\bvd^{(r')}))
\end{equation}
est nul.
\end{prop}

Cela résulte de \ref{higgs3-coh19} et \ref{higgs3-spsa99}(i).

\begin{cor}\label{higgs3-coh24}
Soient $r$, $r'$ deux nombres rationnels tels que $r>r'>0$. 
\begin{itemize}
\item[{\rm (i)}] Le morphisme canonique 
\begin{equation}\label{higgs3-coh24a}
u^r\colon \bvocB_\mQ\rightarrow \rH^0(\mK^\bullet_\mQ(\bvcC^{(r)},p^r\bvd^{(r)}))
\end{equation}
admet un inverse à gauche canonique 
\begin{equation}
v^r\colon  \rH^0(\mK^\bullet_\mQ(\bvcC^{(r)},p^r\bvd^{(r)})) \rightarrow \bvocB_\mQ.
\end{equation}
\item[{\rm (ii)}] Le composé 
\begin{equation}
\rH^0(\mK^\bullet_\mQ(\bvcC^{(r)},p^r\bvd^{(r)})) \stackrel{v^r}{\longrightarrow} \bvocB_\mQ
\stackrel{u^{r'}}{\longrightarrow} \rH^0(\mK^\bullet_\mQ(\bvcC^{(r')},p^r\bvd^{(r')}))
\end{equation}
est le morphisme canonique.
\item[{\rm (iii)}] Pour tout entier $q\geq 1$, le morphisme canonique 
\begin{equation}\label{higgs3-coh24c}
\rH^q(\mK^\bullet_\mQ(\bvcC^{(r)},p^r\bvd^{(r)}))\rightarrow 
\rH^q(\mK^\bullet_\mQ(\bvcC^{(r')},p^{r'}\bvd^{(r')}))
\end{equation}
est nul.
\end{itemize}
\end{cor}

En effet, considérons le diagramme commutatif canonique (sans la flèche pointillée)
\begin{equation}
\xymatrix{
{\bvocB_\mQ}\ar@{=}[d]\ar[r]^-(0.5){u^r}&{\rH^0(\mK^\bullet_\mQ(\bvcC^{(r)},p^r\bvd^{(r)}))}\ar@{.>}[ld]_-(0.5){v^{r,r'}}
\ar[d]^-(0.5){\varpi^{r,r'}}
\ar@{->>}[r]&{\rH^0(\tmK^\bullet_\mQ(\bvcC^{(r)},p^r\bvd^{(r)}))}\ar[d]^{\rH^0(\bvnu^{r,r'}_\mQ)}\\
{\bvocB_\mQ}\ar[r]_-(0.5){u^{r'}}& {\rH^0(\mK^\bullet_\mQ(\bvcC^{(r')},p^{r'}\bvd^{(r')}))}\ar@{->>}[r]&
{\rH^0(\tmK^\bullet_\mQ(\bvcC^{(r')},p^{r'}\bvd^{(r')}))}}
\end{equation}
Il résulte de \ref{higgs3-coh21} que $u^r$ et par suite $u^{r'}$ sont injectifs, et qu'il existe un et un unique morphisme
$v^{r,r'}$ comme ci-dessus tel que $\varpi^{r,r'}=u^{r'}\circ v^{r,r'}$. Comme on a 
$u^{r'}\circ v^{r,r'}\circ u^r=u^{r'}$, on en déduit que $v^{r,r'}$ est un inverse à gauche de $u^r$.  
On vérifie aussitôt qu'il ne dépend pas de $r'$; d'où les propositions (i) et (ii).  
La proposition (iii) résulte aussitôt de \ref{higgs3-coh21}. 

\begin{cor}\label{higgs3-coh22}
Le morphisme canonique 
\begin{equation}\label{higgs3-coh22a}
\bvocB_\mQ\rightarrow \underset{\underset{r\in \mQ_{>0}}{\longrightarrow}}{\lim}\ 
\rH^0(\mK^\bullet_\mQ(\bvcC^{(r)},p^r\bvd^{(r)}))
\end{equation}
est un isomorphisme, et pour tout entier $q\geq 1$, 
\begin{equation}\label{higgs3-coh22b}
\underset{\underset{r\in \mQ_{>0}}{\longrightarrow}}{\lim}\ \rH^q(\mK^\bullet_\mQ(\bvcC^{(r)},p^r\bvd^{(r)}))=0.
\end{equation}
\end{cor}

Cela résulte de  \ref{higgs3-coh24}.

\begin{rema}\label{higgs3-coh23}
Les limites inductives filtrantes ne sont pas a priori représentables dans la catégorie $\bMod_\mQ(\bvocB)$. 
\end{rema}

\section{Modules de Dolbeault}\label{higgs3-MF}

\subsection{}\label{higgs3-MF1}\index{10001203@$\bMod^\atf_\mQ(\bvocB)$}
\index{1000929@$\top\colon (\tE_s^{\mN^\circ},\bvocB)\rightarrow (X_{s,\zar},\co_{\fX})$}\index{1000928@$\fX$}
Les hypothèses et notations générales de § \ref{higgs3-RGG} sont en vigueur dans cette section.
On désigne, de plus, par $\bMod(\bvocB)$ la catégorie des $\bvocB$-modules de $\tE_s^{\mN^\circ}$ \eqref{higgs3-RGG21},
par $\bMod^\ad(\bvocB)$ (resp. $\bMod^\atf(\bvocB)$)
la sous-catégorie pleine formée des $\bvocB$-modules adiques (resp. $\bvocB$-modules adiques de type fini)  
\eqref{higgs3-spad3} et par $\bMod_{\mQ}(\bvocB)$ 
(resp. $\bMod^\ad_{\mQ}(\bvocB)$, resp. $\bMod^\atf_{\mQ}(\bvocB)$)
la catégorie des objets de $\bMod(\bvocB)$ (resp. $\bMod^\ad(\bvocB)$, resp. $\bMod^\atf(\bvocB)$) 
à isogénie près \eqref{higgs3-caip1a}. 
La catégorie $\bMod_\mQ(\bvocB)$ est alors abélienne et les foncteurs canoniques 
\begin{equation}\label{higgs3-MF1a}
\bMod^\atf_{\mQ}(\bvocB)\rightarrow \bMod^\ad_{\mQ}(\bvocB)\rightarrow\bMod_{\mQ}(\bvocB)
\end{equation} 
sont pleinement fidèles. 

On désigne par $\fX$ le $\cS$-schéma formel complété $p$-adique de $\oX$ et par 
$\xi^{-1}\tOmega^1_{\fX/\cS}$ le complété $p$-adique du $\co_\oX$-module 
$\xi^{-1}\tOmega^1_{\oX/\oS}=\xi^{-1}\tOmega^1_{X/S}\otimes_{\co_X}\co_{\oX}$  (cf. \ref{higgs3-not1}).
On note $\bMod^\coh(\co_{\fX})$ (resp. $\bMod^{\coh}(\co_{\fX}[\frac 1 p])$)
la catégorie des $\co_{\fX}$-modules (resp. $\co_{\fX}[\frac 1 p]$-modules) cohérents de $X_{s,\zar}$ \eqref{higgs3-formel1}.

Soit
\begin{equation}\label{higgs3-MF1ab}
\top\colon (\tE_s^{\mN^\circ},\bvocB)\rightarrow (X_{s,\zar},\co_{\fX})
\end{equation}
le morphisme de topos annelés défini dans \eqref{higgs3-coh0e}.
Le foncteur $\top_*$ induit un foncteur additif et exact à gauche que l'on note encore 
\begin{equation}\label{higgs3-MF1b}
\top_*\colon \bMod_{\mQ}(\bvocB) \rightarrow \bMod(\co_{\fX}[\frac 1 p]).
\end{equation}
D'après \ref{higgs3-formel2}, le foncteur $\top^*$ induit un foncteur additif que l'on note encore
\begin{equation}\label{higgs3-MF1c}
\top^*\colon \bMod^\coh(\co_{\fX}[\frac 1 p]) \rightarrow \bMod_{\mQ}^\atf(\bvocB). 
\end{equation}
Pour tout $\co_\fX[\frac 1 p]$-module cohérent $\cF$ et tout $\bvocB_\mQ$-module $\cG$, 
on a un homomorphisme canonique bifonctoriel 
\begin{equation}\label{higgs3-MF1d}
\Hom_{\bvocB_\mQ}(\top^*(\cF),\cG)\rightarrow\Hom_{\co_\fX[\frac 1 p]}(\cF,\top_*(\cG)),
\end{equation}
qui est injectif d'après \eqref{higgs3-formel1c} et \ref{higgs3-formel2}. 
On appelle abusivement l'{\em adjoint} d'un morphisme $\bvocB_\mQ$-linéaire $\top^*(\cF)\rightarrow \cG$
son image par l'homomorphisme \eqref{higgs3-MF1d}. 

On note 
\begin{eqnarray}
\rR\top_*\colon \bD^+(\bMod_\mQ(\bvocB))&\rightarrow& \bD^+(\bMod(\co_{\fX}[\frac 1 p])), \label{higgs3-MF1ee}\\
\rR^q \top_*\colon \bMod_\mQ(\bvocB)&\rightarrow& \bMod(\co_{\fX}[\frac 1 p]), \ \ \ \ (q\in \mN), \label{higgs3-MF1e}
\end{eqnarray}
les foncteurs dérivés droits du foncteur $\top_*$ \eqref{higgs3-MF1b}. 
Ces notations n'induisent aucune confusion avec celles des foncteurs dérivés 
droits du foncteur $\top_*\colon \bMod(\bvocB)\rightarrow \bMod(\co_{\fX})$, 
puisque le foncteur de localisation $\bMod(\bvocB)\rightarrow \bMod_\mQ(\bvocB)$ 
est exact et transforme les objets injectifs en des objets injectifs.

\subsection{}\label{higgs3-MF2}
Soient $\cM$ un $\co_\fX$-module, $\cN$ un $\bvocB$-module, $q$ un entier $\geq 0$. 
Le morphisme d'adjonction $\cM\rightarrow \top_*(\top^*(\cM))$ et le cup-produit induisent un morphisme bifonctoriel
\begin{equation}\label{higgs3-MF2a}
\cM\otimes_{\co_\fX}\rR^q\top_*(\cN)\rightarrow \rR^q\top_*(\top^*(\cM)\otimes_{\bvocB}\cN). 
\end{equation}
On peut faire les remarques suivantes~:
\begin{itemize}
\item[(i)] Pour tout $\co_\fX$-module $\cM'$, le composé 
\begin{equation}
\xymatrix{
{\cM\otimes_{\co_\fX}\cM'\otimes_{\co_\fX}\rR^q\top_*(\cN)}\ar[r]\ar[rd]&
{\cM\otimes_{\co_\fX}\rR^q\top_*(\top^*(\cM')\otimes_{\bvocB}\cN)}\ar[d]\\ 
&{\rR^q\top_*(\top^*(\cM\otimes_{\co_\fX}\cM')\otimes_{\bvocB}\cN)}}
\end{equation}
des morphismes induits par les morphismes \eqref{higgs3-MF2a} relatifs à $\cM$ et $\cM'$, 
n'est autre que le morphisme \eqref{higgs3-MF2a} relatif à $\cM\otimes_{\co_\fX}\cM'$. 
\item[(ii)] Lorsque $q=0$, le morphisme \eqref{higgs3-MF2a} est le composé 
\begin{equation}\label{higgs3-MF2b}
\cM\otimes_{\co_\fX}\top_*(\cN)\rightarrow \top_*(\top^*(\cM\otimes_{\co_\fX}\top_*(\cN)))\rightarrow  
\top_*(\top^*(\cM)\otimes_{\bvocB}\cN),
\end{equation}
où la première flèche est le morphisme d'adjonction et 
la seconde flèche est induite par le morphisme canonique $\top^*(\top_*(\cN))\rightarrow \cN$.
Son adjoint  
\begin{equation}\label{higgs3-MF2c}
\top^*(\cM\otimes_{\co_\fX}\top_*(\cN))\rightarrow\top^*(\cM)\otimes_{\bvocB}\cN
\end{equation}
est donc induit par le morphisme canonique $\top^*(\top_*(\cN))\rightarrow \cN$.  
\end{itemize}

\subsection{}\label{higgs3-MF4}
Soient $\cF$ un $\co_\fX[\frac 1 p]$-module cohérent, $\cG$ un $\bvocB_\mQ$-module, $q$ un entier $\geq 0$. 
Compte tenu de \ref{higgs3-formel2}, le morphisme \eqref{higgs3-MF2a} induit un morphisme bifonctoriel
\begin{equation}\label{higgs3-MF4a}
\cF\otimes_{\co_\fX[\frac 1 p]}\rR^q\top_*(\cG)\rightarrow\rR^q\top_*(\top^*(\cF)\otimes_{\bvocB_\mQ}\cG).
\end{equation}
On peut faire les remarques suivantes~:
\begin{itemize}
\item[(i)] Soit $\cF'$ un $\co_\fX[\frac 1 p]$-module cohérent. Il résulte de \ref{higgs3-MF2}(i) 
que le composé 
\begin{equation}
\xymatrix{
{\cF\otimes_{\co_\fX[\frac 1 p]}\cF'\otimes_{\co_\fX[\frac 1 p]}\rR^q\top_*(\cG)}\ar[r]\ar[rd]&
{\cF\otimes_{\co_\fX[\frac 1 p]}\rR^q\top_*(\top^*(\cF')\otimes_{\bvocB_\mQ}\cG)}\ar[d]\\
&{\rR^q\top_*(\top^*(\cF\otimes_{\co_\fX[\frac 1 p]}\cF')\otimes_{\bvocB_\mQ}\cG)}}
\end{equation}
des morphismes induits par les morphismes \eqref{higgs3-MF4a} relatifs à $\cF$ et $\cF'$,
n'est autre que le morphisme \eqref{higgs3-MF4a} relatif à $\cF\otimes_{\co_\fX[\frac 1 p]}\cF'$.
\item[(ii)] Soient $\cL$ un $\co_\fX[\frac 1 p]$-module cohérent, 
$u\colon \top^*(\cL)\rightarrow \cG$ un morphisme $\bvocB_\mQ$-linéaire, 
$v\colon \cL\rightarrow \top_*(\cG)$ le morphisme adjoint \eqref{higgs3-MF1d}. Il résulte alors de \ref{higgs3-MF2}(ii) et \ref{higgs3-formel2}
que le morphisme 
\begin{equation}\label{higgs3-MF4b}
\cF\otimes_{\co_\fX[\frac 1 p]}\cL \rightarrow \top_*(\top^*(\cF)\otimes_{\bvocB_\mQ} \cG)
\end{equation}
induit par \eqref{higgs3-MF4a} et $v$, est l'adjoint du morphisme 
\begin{equation}\label{higgs3-MF4c}
\top^*(\cF\otimes_{\co_\fX[\frac 1 p]}\cL)\rightarrow\top^*(\cF)\otimes_{\bvocB_\mQ}\cG
\end{equation}
induit par $u$.
\end{itemize}

\begin{lem}\label{higgs3-MF6}
{\rm (i)}\ Soient $\cM$ un $\co_\fX$-module localement libre de type fini, 
$\cN$ un $\bvocB$-module, $q$ un entier $\geq 0$. 
Alors le morphisme canonique \eqref{higgs3-MF2a} 
\begin{equation}\label{higgs3-MF6d}
\cM\otimes_{\co_\fX}\rR^q\top_*(\cN)\rightarrow\rR^q\top_*(\top^*(\cM)\otimes_{\bvocB}\cN)
\end{equation}
est un isomorphisme.

{\rm (ii)}\ Soient $\cF$ un $\co_\fX[\frac 1 p]$-module localement projectif de type fini \eqref{higgs3-not3},
$\cG$ un $\bvocB_\mQ$-module, $q$ un entier $\geq 0$. 
Alors le morphisme canonique \eqref{higgs3-MF4a} 
\begin{equation}\label{higgs3-MF6a}
\cF\otimes_{\co_\fX[\frac 1 p]}\rR^q\top_*(\cG)\rightarrow\rR^q\top_*(\top^*(\cF)\otimes_{\bvocB_\mQ}\cG)
\end{equation}
est un isomorphisme.
\end{lem}
On se limite à démontrer (ii); la preuve de (i) est similaire et plus simple. 
Il existe un recouvrement ouvert de Zariski $(U_i)_{i\in I}$ de $X$ tel que pour tout $i\in I$, la restriction 
de $\cF$ à $(U_i)_s$ soit un facteur direct d'un $(\co_\fX|U_i)[\frac 1 p]$-module libre de type fini. 
Compte tenu de \ref{higgs3-TFT15}, on peut alors se borner au cas où $\cF$ est un facteur direct
d'un $\co_\fX[\frac 1 p]$-module libre de type fini, et même au cas où $\cF$ est un 
$\co_\fX[\frac 1 p]$-module libre de type fini, auquel cas l'assertion est évidente.

\subsection{}\label{higgs3-MF12}\index{10001210@$\bIH^\coh_\mQ(\co_\fX,\xi^{-1}\tOmega^1_{\fX/\cS})$}
\index{10001211@$\bMH^\coh(\co_\fX[\frac 1 p], \xi^{-1}\tOmega^1_{\fX/\cS})$}
On désigne par $\bIH(\co_\fX,\xi^{-1}\tOmega^1_{\fX/\cS})$ la catégorie des $\co_\fX$-isogénies  
de Higgs à coefficients dans $\xi^{-1}\tOmega^1_{\fX/\cS}$ \eqref{higgs3-imh1} et
par $\bIH^\coh(\co_\fX,\xi^{-1}\tOmega^1_{\fX/\cS})$ la sous-catégorie pleine 
formée des quadruplets $(\cM,\cN,u,\theta)$ tels que $\cM$ et $\cN$ soient des $\co_\fX$-modules cohérents.
Ce sont des catégories additives. 
On note $\bIH_\mQ(\co_\fX,\xi^{-1}\tOmega^1_{\fX/\cS})$ (resp. $\bIH^\coh_\mQ(\co_\fX,\xi^{-1}\tOmega^1_{\fX/\cS})$)
la catégorie des objets de $\bIH(\co_\fX,\xi^{-1}\tOmega^1_{\fX/\cS})$ 
(resp. $\bIH^\coh(\co_\fX,\xi^{-1}\tOmega^1_{\fX/\cS})$) à isogénie près  \eqref{higgs3-caip1a}. 

On sous-entend par $\co_\fX[\frac 1 p]$-module de Higgs à coefficients dans $\xi^{-1}\tOmega^1_{\fX/\cS}$, 
un $\co_\fX[\frac 1 p]$-module de Higgs à coefficients dans $\xi^{-1}\tOmega^1_{\fX/\cS}[\frac 1 p]$ (\cite{ag1} 2.8).
Dans la suite, on omettra le champ de Higgs de la notation d'un module de Higgs 
lorsqu'on en n'a pas explicitement besoin. 
On désigne par $\bMH(\co_\fX[\frac 1 p], \xi^{-1}\tOmega^1_{\fX/\cS})$ la catégorie des 
$\co_\fX[\frac 1 p]$-modules de Higgs à coefficients dans $\xi^{-1}\tOmega^1_{\fX/\cS}$
et par $\bMH^\coh(\co_\fX[\frac 1 p], \xi^{-1}\tOmega^1_{\fX/\cS})$ la sous-catégorie pleine formée des modules
de Higgs dont le $\co_\fX[\frac 1 p]$-module sous-jacent est cohérent.
Le foncteur \eqref{higgs3-formel3a}
\begin{equation}\label{higgs3-MF12a}
\begin{array}[t]{clcr}
\bIH(\co_\fX,\xi^{-1}\tOmega^1_{\fX/\cS})&\rightarrow& \bMH(\co_\fX[\frac 1 p], \xi^{-1}\tOmega^1_{\fX/\cS})\\
(\cM,\cN,u,\theta)&\mapsto& (\cM_{\mQ_p}, (\id \otimes u_{\mQ_p}^{-1})\circ\theta_{\mQ_p})
\end{array}
\end{equation}
induit un foncteur 
\begin{equation}\label{higgs3-MF12b}
\bIH_\mQ(\co_\fX,\xi^{-1}\tOmega^1_{\fX/\cS})\rightarrow \bMH(\co_\fX[\frac 1 p], \xi^{-1}\tOmega^1_{\fX/\cS}).
\end{equation}
D'après \ref{higgs3-formel5}, celui-ci induit une équivalence de catégories 
\begin{equation}\label{higgs3-MF12c}
\bIH^\coh_\mQ(\co_\fX,\xi^{-1}\tOmega^1_{\fX/\cS})\stackrel{\sim}{\rightarrow} 
\bMH^\coh(\co_\fX[\frac 1 p], \xi^{-1}\tOmega^1_{\fX/\cS}).
\end{equation}

\begin{defi}\label{higgs3-MF123}\index{Fibre de Higgs@$\co_\fX[\frac 1 p]$-fibré de Higgs}
On appelle {\em $\co_\fX[\frac 1 p]$-fibré de Higgs à coefficients dans $\xi^{-1}\tOmega^1_{\fX/\cS}$} 
tout $\co_\fX[\frac 1 p]$-module de Higgs à coefficients dans $\xi^{-1}\tOmega^1_{\fX/\cS}$ 
dont le $\co_\fX[\frac 1 p]$-module sous-jacent est localement projectif de type fini \eqref{higgs3-not3}. 
\end{defi}

\subsection{}\label{higgs3-MF15}\index{10001220@$\Xi^r_\mQ$}\index{10001221@$\fS^r$, $\cK^r$, $\top^{r+}$, $\top^r_+$ ($r\in \mQ_{\geq 0}$)}
Soit $r$ un nombre rationnel $\geq 0$. On note encore 
\begin{equation}\label{higgs3-MF15h}
\bvd^{(r)}\colon \bvcC^{(r)}\rightarrow \top^*(\xi^{-1}\tOmega^1_{\fX/\cS})\otimes_{\bvocB}\bvcC^{(r)}
\end{equation}
la $\bvocB$-dérivation induite par $\bvd^{(r)}$ \eqref{higgs3-RGG180b} et l'isomorphisme \eqref{higgs3-coh1cc},
que l'on identifie à la $\bvocB$-dérivation universelle de $\bvcC^{(r)}$.
C'est un $\bvocB$-champ de Higgs à coefficients dans $\top^*(\xi^{-1}\tOmega^1_{\fX/\cS})$ \eqref{higgs3-RGG181}.
On désigne par $\Xi^r$ la catégorie des $p^r$-isoconnexions intégrables 
relativement à l'extension $\bvcC^{(r)}/\bvocB$ \eqref{higgs3-isoco1}. C'est une catégorie additive.
On note $\Xi^r_\mQ$ la catégorie des objets de $\Xi^r$ à isogénie près  \eqref{higgs3-caip1a}. 
D'après \ref{higgs3-isoco3} et \ref{higgs3-RGG181}(iii), tout objet de $\Xi^r$ est une $\bvocB$-isogénie
de Higgs à coefficients dans $\top^*(\xi^{-1}\tOmega^1_{\fX/\cS})$ \eqref{higgs3-imh1}. 
En particulier, on peut associer fonctoriellement à tout objet de $\Xi^r_\mQ$ un complexe de Dolbeault dans 
$\bMod_{\mQ}(\bvocB)$ (cf. \ref{higgs3-imh2}).

Considérons le foncteur 
\begin{equation}\label{higgs3-MF15b}
\fS^r\colon\bMod(\bvocB)\rightarrow \Xi^r, \ \ \ \cM\mapsto (\bvcC^{(r)}\otimes_{\bvocB}\cM,\bvcC^{(r)}\otimes_{\bvocB}\cM,
\id,p^r\bvd^{(r)}\otimes \id),
\end{equation}
et notons encore
\begin{equation}\label{higgs3-MF15bb}
\fS^r\colon \bMod_\mQ(\bvocB)\rightarrow \Xi^r_\mQ
\end{equation}
le foncteur induit. Considérons par ailleurs le foncteur 
\begin{equation}\label{higgs3-MF15a}
\cK^r\colon \Xi^r\rightarrow \bMod(\bvocB),\ \ \ (\cF,\cG,u,\nabla)\mapsto \ker(\nabla),
\end{equation}
et notons encore
\begin{equation}\label{higgs3-MF15aa}
\cK^r\colon \Xi^r_\mQ\rightarrow \bMod_\mQ(\bvocB)
\end{equation}
le foncteur induit. 
Il est clair que le foncteur \eqref{higgs3-MF15b} est un adjoint à gauche du foncteur \eqref{higgs3-MF15a}. 
Par suite, le foncteur \eqref{higgs3-MF15bb} est un adjoint à gauche du foncteur \eqref{higgs3-MF15aa}.

D'après \ref{higgs3-isoco3}, si $(\cN,\cN',v,\theta)$ est une $\co_\fX$-isogénie de Higgs
à coefficients dans $\xi^{-1}\tOmega^1_{\fX/\cS}$,
\begin{equation}\label{higgs3-MF15e}
(\bvcC^{(r)}\otimes_{\bvocB}\top^*(\cN),\bvcC^{(r)}\otimes_{\bvocB}\top^*(\cN'),\id \otimes_{\bvocB}\top^*(v),
p^r \bvd^{(r)} \otimes\top^*(v)+\id \otimes \top^*(\theta))
\end{equation}
est un objet de $\Xi^r$. On obtient ainsi un foncteur 
\begin{equation}\label{higgs3-MF15f}
\top^{r+}\colon\bIH(\co_\fX,\xi^{-1}\tOmega^1_{\fX/\cS})\rightarrow \Xi^r.
\end{equation}
D'après \eqref{higgs3-MF12c}, celui-ci induit un foncteur que l'on note encore
\begin{equation}\label{higgs3-MF15ff}
\top^{r+}\colon\bMH^\coh(\co_\fX[\frac 1 p], \xi^{-1}\tOmega^1_{\fX/\cS})\rightarrow \Xi^r_\mQ.
\end{equation}

Soit $(\cF,\cG,u,\nabla)$ un objet de $\Xi^r$.
Compte tenu de \ref{higgs3-MF6}(i), $\nabla$ induit un morphisme $\co_\fX$-linéaire~: 
\begin{equation}\label{higgs3-MF15c}
\top_*(\nabla)\colon \top_*(\cF)\rightarrow \xi^{-1}\tOmega^1_{\fX/\cS}\otimes_{\co_\fX}\top_*(\cG).
\end{equation}
On déduit facilement de \ref{higgs3-MF4}(i) que $(\top_*(\cF),\top_*(\cG),\top_*(u),\top_*(\nabla))$ est une $\co_\fX$-isogénie 
de Higgs à coefficients dans $\xi^{-1}\tOmega^1_{\fX/\cS}$. On obtient ainsi un foncteur
\begin{equation}\label{higgs3-MF15d}
\top^r_+\colon \Xi^r\rightarrow \bIH(\co_\fX,\xi^{-1}\tOmega^1_{\fX/\cS}).
\end{equation}
Le composé des foncteurs \eqref{higgs3-MF15d} et \eqref{higgs3-MF12a} induit un foncteur  que l'on note encore
\begin{equation}\label{higgs3-MF15dd}
\top^r_+\colon \Xi^r_\mQ\rightarrow \bMH(\co_\fX[\frac 1 p], \xi^{-1}\tOmega^1_{\fX/\cS}).
\end{equation}

Il est clair que  le foncteur \eqref{higgs3-MF15f} est un adjoint à gauche du foncteur \eqref{higgs3-MF15d}. 
On en déduit que pour tous $\cN\in \ob(\bMH^\coh(\co_\fX[\frac 1 p], \xi^{-1}\tOmega^1_{\fX/\cS}))$ 
et $\cA\in \ob(\Xi^r_\mQ)$, on a un homomorphisme canonique bifonctoriel 
\begin{equation}\label{higgs3-MF15g}
\Hom_{\Xi^r_\mQ}(\top^{r+}(\cN),\cA)\rightarrow
\Hom_{\bMH(\co_\fX[\frac 1 p], \xi^{-1}\tOmega^1_{\fX/\cS})}(\cN,\top^r_+(\cA)),
\end{equation}
qui est injectif d'après \ref{higgs3-formel4} et \ref{higgs3-formel5}. 
On appelle abusivement l'{\em adjoint} d'un morphisme $\top^{r+}(\cN)\rightarrow \cA$ de $\Xi^r_\mQ$,
son image par l'homomorphisme \eqref{higgs3-MF15g}.

\subsection{}\label{higgs3-MF17}\index{10001225@$\epsilon^{r,r'}$}
Soient $r$, $r'$ deux nombres rationnels tels que $r\geq r'\geq 0$, $(\cF,\cG,u,\nabla)$ une $p^r$-isoconnexion 
intégrable relativement à l'extension $\bvcC^{(r)}/\bvocB$. 
D'après \eqref{higgs3-coh20e}, il existe un et un unique morphisme $\bvocB$-linéaire
\begin{equation}\label{higgs3-MF17b}
\nabla'\colon \bvcC^{(r')}\otimes_{\bvcC^{(r)}}\cF\rightarrow \top^*(\xi^{-1}\tOmega^1_{\fX/\cS}) \otimes_{\bvocB}
\bvcC^{(r')}\otimes_{\bvcC^{(r)}} \cG
\end{equation}
tel que pour toutes sections locales $x'$ de $\bvcC^{(r')}$ et $s$ de $\cF$, on ait 
\begin{equation}\label{higgs3-MF17c}
\nabla'(x'\otimes_{\bvcC^{(r)}} s)=p^{r'} \bvd^{(r')}(x')\otimes_{\bvcC^{(r)}}u(s) + x'\otimes_{\bvcC^{(r)}}\nabla(s).
\end{equation}
Le quadruplet $(\bvcC^{(r')}\otimes_{\bvcC^{(r)}}\cF,\bvcC^{(r')}\otimes_{\bvcC^{(r)}}\cG,
\id\otimes_{\bvcC^{(r)}}u,\nabla')$ est une $p^{r'}$-isoconnexion intégrable relativement à l'extension 
$\bvcC^{(r')}/\bvocB$. On obtient ainsi un foncteur 
\begin{equation}\label{higgs3-MF17d}
\epsilon^{r,r'}\colon \Xi^r\rightarrow \Xi^{r'}.
\end{equation}
Celui-ci induit un foncteur que l'on note encore 
\begin{equation}\label{higgs3-MF17e}
\epsilon^{r,r'}\colon \Xi^r_\mQ\rightarrow \Xi^{r'}_\mQ.
\end{equation}

On a un isomorphisme canonique de foncteurs de $\bMod(\bvocB)$ dans $\Xi^{r'}$ 
(resp.  de $\bMod_\mQ(\bvocB)$ dans $\Xi^{r'}_\mQ$) 
\begin{equation}\label{higgs3-MF17g}
\epsilon^{r,r'}\circ \fS^r\stackrel{\sim}{\longrightarrow} \fS^{r'}.
\end{equation}
Par ailleurs, on a un isomorphisme canonique de foncteurs de $\bIH(\co_\fX,\xi^{-1}\tOmega^1_{\fX/\cS})$ dans $\Xi^{r'}$ 
(resp.  de $\bMH^\coh(\co_\fX[\frac 1 p], \xi^{-1}\tOmega^1_{\fX/\cS})$ dans $\Xi^{r'}_\mQ$ ) 
\begin{equation}\label{higgs3-MF17h}
\epsilon^{r,r'}\circ \top^{r+}\stackrel{\sim}{\longrightarrow} \top^{r'+}.
\end{equation}

Le diagramme 
\begin{equation}\label{higgs3-MF17i}
\xymatrix{
{\cF}\ar[r]^-(0.5){\nabla}\ar[d]_{\bvalpha^{r,r'}\otimes_{\bvcC^{(r)}} \id}&
{\top^*(\xi^{-1}\tOmega^1_{\fX/\cS})\otimes_{\bvocB}\cG}\ar[d]^{\id\otimes_{\bvocB}\bvalpha^{r,r'}\otimes_{\bvcC^{(r)}} \id}\\
{\bvcC^{(r')}\otimes_{\bvcC^{(r)}}\cF}\ar[r]^-(0.5){\nabla'}&
{\top^*(\xi^{-1}\tOmega^1_{\fX/\cS})\otimes_{\bvocB}\bvcC^{(r')}\otimes_{\bvcC^{(r)}}\cG}}
\end{equation}
est clairement commutatif. On en déduit un morphisme canonique de foncteurs de $\Xi^r$ dans $\bMod(\bvocB)$ (resp. 
de $\Xi^r_\mQ$ dans $\bMod_\mQ(\bvocB)$)
\begin{equation}\label{higgs3-MF17k}
\cK^r\rightarrow \cK^{r'}\circ \epsilon^{r,r'}.
\end{equation}
On en déduit aussi un morphisme canonique de foncteurs de $\Xi^r$ dans $\bIH(\co_\fX,\xi^{-1}\tOmega^1_{\fX/\cS})$ 
(resp.  de $\Xi^r_\mQ$ dans $\bMH(\co_\fX[\frac 1 p], \xi^{-1}\tOmega^1_{\fX/\cS})$) 
\begin{equation}\label{higgs3-MF17j}
\top^r_+\longrightarrow \top^{r'}_+\circ \epsilon^{r,r'}.
\end{equation}

Pour tout nombre rationnel $r''$ tel $r'\geq r''\geq 0$, on a un isomorphisme canonique de foncteurs 
de $\Xi^r$ dans $\Xi^{r''}$ (resp. de $\Xi^r_\mQ$ dans $\Xi^{r''}_\mQ$)
\begin{equation}
\epsilon^{r',r''}\circ \epsilon^{r,r'}\stackrel{\sim}{\rightarrow}\epsilon^{r,r''}.
\end{equation}

\begin{rema}
Sous les hypothèses de \eqref{higgs3-MF17}, $(\bvcC^{(r')}\otimes_{\bvcC^{(r)}}\cF,\bvcC^{(r')}\otimes_{\bvcC^{(r)}}\cG,
\id\otimes_{\bvcC^{(r)}}u,p^{r-r'}\nabla')$ est la $p^r$-isoconnexion intégrable relativement à l'extension
$\bvcC^{(r')}/\bvocB$ déduite de $(\cF,\cG,u,\nabla)$
par extension des scalaires par $\bvalpha^{r,r'}$, définie dans \ref{higgs3-isoco2}. Ce décalage s'explique par le fait que
l'homomorphisme canonique $\Omega^1_{\bvcC^{(r)}/\bvocB}\rightarrow \Omega^1_{\bvcC^{(r')}/\bvocB}$
s'identifie à 
\begin{equation}
p^{r-r'}\id\otimes \bvalpha^{r,r'}\colon   \top^*(\xi^{-1}\tOmega^1_{\fX/\cS})\otimes_{\bvocB}\bvcC^{(r)}\rightarrow
\top^*(\xi^{-1}\tOmega^1_{\fX/\cS})\otimes_{\bvocB}\bvcC^{(r')}.
\end{equation}
\end{rema}

\begin{defi}\label{higgs3-MF8}\index{associes@associés, $r$-associés (des modules ---)}\index{admissible@$r$-admissible (triplet ---)}
Soient $\cM$ un objet de $\bMod^\atf_{\mQ}(\bvocB)$ \eqref{higgs3-MF1}, 
$\cN$ un $\co_\fX[\frac 1 p]$-fibré de Higgs à coefficients dans $\xi^{-1}\tOmega^1_{\fX/\cS}$ \eqref{higgs3-MF123}.
\begin{itemize}
\item[(i)] Soit $r$ un nombre rationnel $>0$.  On dit que $\cM$ et  $\cN$ sont {\em $r$-associés} 
s'il existe un isomorphisme de $\Xi^r_\mQ$ 
\begin{equation}\label{higgs3-MF8a}
\alpha\colon \top^{r+}(\cN) \stackrel{\sim}{\rightarrow}\fS^r(\cM).
\end{equation}
On dit alors aussi que le triplet $(\cM,\cN,\alpha)$ est {\em $r$-admissible}. 
\item[(ii)] On dit que $\cM$ et  $\cN$ sont {\em associés} s'il existe un nombre rationnel $r>0$ tel que 
$\cM$ et  $\cN$ soient $r$-associés.
\end{itemize}
\end{defi}

On notera que pour tous nombres rationnels $r\geq r'>0$, 
si $\cM$ et  $\cN$ sont $r$-associés, ils sont $r'$-associés, compte tenu de \eqref{higgs3-MF17g} et \eqref{higgs3-MF17h}. 

\begin{defi}\label{higgs3-MF14}\index{Module de Dolbeault@$\bvocB_\mQ$-module de Dolbeault}
\index{Fibre de Higgs@$\co_\fX[\frac 1 p]$-fibré de Higgs!soluble}
(i)\ On appelle {\em $\bvocB_\mQ$-module de Dolbeault} tout objet de $\bMod^\atf_{\mQ}(\bvocB)$ pour lequel il existe 
un $\co_\fX[\frac 1 p]$-fibré de Higgs associé, à coefficients dans $\xi^{-1}\tOmega^1_{\fX/\cS}$.

(ii)\ On dit qu'un $\co_\fX[\frac 1 p]$-fibré de Higgs à coefficients dans $\xi^{-1}\tOmega^1_{\fX/\cS}$ est {\em soluble}
s'il admet un module de Dolbeault associé. 
\end{defi}

On désigne par $\bMod_\mQ^\Dolb(\bvocB)$ la sous-catégorie pleine de $\bMod^\atf_\mQ(\bvocB)$ 
formée des $\bvocB_\mQ$-modules de Dolbeault, et par $\bMH^\sol(\co_\fX[\frac 1 p], \xi^{-1}\tOmega^1_{\fX/\cS})$ 
la sous-catégorie pleine de  $\bMH(\co_\fX[\frac 1 p], \xi^{-1}\tOmega^1_{\fX/\cS})$
formée des $\co_\fX[\frac 1 p]$-fibrés de Higgs solubles à coefficients dans $\xi^{-1}\tOmega^1_{\fX/\cS}$.

\begin{prop}\label{higgs3-plat3}
Tout $\bvocB_\mQ$-module de Dolbeault est $\bvocB_\mQ$-plat \eqref{higgs3-caip5}.
\end{prop}

Soient $\cM$ un $\bvocB_\mQ$-module de Dolbeault, 
$\cN$ un $\co_\fX[\frac 1 p]$-fibré de Higgs à coefficients dans 
$\xi^{-1}\tOmega^1_{\fX/\cS}$, $r$ un nombre rationnel $>0$, 
$\alpha\colon \top^{r+}(\cN) \stackrel{\sim}{\rightarrow}\fS^r(\cM)$
un isomorphisme de $\Xi^r_\mQ$. Comme le $\co_\fX[\frac 1 p]$-module $\cN$ est localement libre de type fini, 
le $\bvcC^{(r)}_\mQ$-module $\top^*(\cN)\otimes_{\bvocB_\mQ}\bvcC^{(r)}_\mQ$ est plat d'après \ref{higgs3-caip7}(iii). 
On en déduit que $\cM\otimes_{\bvocB_\mQ}\bvcC^{(r)}_\mQ$ est $\bvcC^{(r)}_\mQ$-plat. 
Par suite,  $\cM$ est $\bvocB_\mQ$-plat en vertu de \ref{higgs3-caip5d} et \ref{higgs3-plat2}. 

\subsection{}\label{higgs3-MF16}\index{10001230@$\cH$}
Pour tout $\bvocB_\mQ$-module $\cM$ et tous nombres rationnels $r\geq r'\geq 0$, 
le morphisme \eqref{higgs3-MF17j} et l'isomorphisme \eqref{higgs3-MF17g} induisent un morphisme de 
$\bMH(\co_\fX[\frac 1 p], \xi^{-1}\tOmega^1_{\fX/\cS})$
\begin{equation}\label{higgs3-MF16b}
\top^r_+(\fS^r(\cM))\rightarrow \top^{r'}_+(\fS^{r'}(\cM)).
\end{equation}
On obtient ainsi un système inductif filtrant $(\top^r_+(\fS^r(\cM)))_{r\in \mQ_{\geq 0}}$. 
On désigne par $\cH$ le foncteur 
\begin{equation}\label{higgs3-MF16a}
\cH\colon \bMod_\mQ(\bvocB)\rightarrow \bMH(\co_\fX[\frac 1 p], \xi^{-1}\tOmega^1_{\fX/\cS}), \ \ \ \cM\mapsto 
\underset{\underset{r\in \mQ_{>0}}{\longrightarrow}}{\lim}\ \top^r_+(\fS^r(\cM)). 
\end{equation}

Pour tout objet $\cN$ de $\bMH(\co_\fX[\frac 1 p], \xi^{-1}\tOmega^1_{\fX/\cS})$ 
et tous nombres rationnels $r\geq r'\geq 0$, 
le morphisme \eqref{higgs3-MF17k} et l'isomorphisme \eqref{higgs3-MF17h} induisent un morphisme de $\bMod_{\mQ}(\bvocB)$
\begin{equation}\label{higgs3-MF16c}
\cK^r(\top^{r+}(\cN))\rightarrow \cK^{r'}(\top^{r'+}(\cN)).
\end{equation}
On obtient ainsi un système inductif filtrant $(\cK^r(\top^{r+}(\cN)))_{r\geq  0}$. 
On rappelle \eqref{higgs3-coh23} que les limites inductives filtrantes ne sont pas a priori représentables dans la catégorie 
$\bMod_\mQ(\bvocB)$.

\begin{lem}\label{higgs3-MF10}
On a un isomorphisme canonique de $\bMH(\co_\fX[\frac 1 p], \xi^{-1}\tOmega^1_{\fX/\cS})$
\begin{equation}\label{higgs3-MF10a}
(\co_\fX[\frac 1 p],0)\stackrel{\sim}{\rightarrow}\cH(\bvocB_\mQ).
\end{equation}
\end{lem}

Cela résulte de \ref{higgs3-coh17}. 

\begin{lem}\label{higgs3-MF101}
Soient $\cN$ un $\co_\fX[\frac 1 p]$-fibré de Higgs à coefficients dans $\xi^{-1}\tOmega^1_{\fX/\cS}$ \eqref{higgs3-MF123},
$r$ un nombre rationnel $\geq 0$.   
On a alors un isomorphisme canonique de $\bMH(\co_\fX[\frac 1 p], \xi^{-1}\tOmega^1_{\fX/\cS})$
\begin{equation}\label{higgs3-MF101a}
\gamma^r\colon \cN\otimes_{\co_\fX[\frac 1 p]}\top^r_+(\fS^r(\bvocB_\mQ))\stackrel{\sim}{\rightarrow}
\top^r_+(\top^{r+}(\cN)),
\end{equation}
où le membre de gauche est le produit tensoriel des modules de Higgs {\rm (\cite{ag1} 2.8.8)}. 
De plus, on a les propriétés suivantes~:
\begin{itemize}
\item[{\rm (i)}] Le morphisme 
\begin{equation}\label{higgs3-MF101b}
\cN\rightarrow \top^r_+(\top^{r+}(\cN))
\end{equation}
induit par $\gamma^r$ et le morphisme canonique $\co_\fX[\frac 1 p]\rightarrow \top^r_+(\fS^r(\bvocB_\mQ))$,
est l'adjoint du morphisme identité de $\top^{r+}(\cN)$ \eqref{higgs3-MF15g}.
\item[{\rm (ii)}] Pour tout nombre rationnel $r'$ tel que $r\geq r'\geq 0$,  le diagramme 
\begin{equation}\label{higgs3-MF101c}
\xymatrix{
{\cN\otimes_{\co_\fX[\frac 1 p]}\top^r_+(\fS^r(\bvocB_\mQ))}\ar[r]^-(0.5){\gamma^r}\ar[d]&
{\top^r_+(\top^{r+}(\cN))}\ar[d]\\
{\cN\otimes_{\co_\fX[\frac 1 p]}\top^{r'}_+(\fS^{r'}(\bvocB_\mQ))}\ar[r]^-(0.5){\gamma^{r'}}&
{\top^{r'}_+(\top^{r'+}(\cN))}}
\end{equation}
où les flèches verticales sont induites par le morphisme \eqref{higgs3-MF17j} et par les isomorphismes \eqref{higgs3-MF17g} 
et \eqref{higgs3-MF17h}, est commutatif. 
\end{itemize}
\end{lem}

En effet, d'après \ref{higgs3-MF6}(ii),
on a des isomorphismes canoniques de $\co_\fX[\frac 1 p]$-modules 
\begin{eqnarray}
\cN\otimes_{\co_\fX[\frac 1 p]}\top_*(\bvcC^{(r)}_\mQ)&\stackrel{\sim}{\rightarrow}&\top_*(\top^*(\cN)\otimes_{\bvocB_\mQ}
\bvcC^{(r)}_\mQ),\label{higgs3-MF101d}\\
\ \ \ \ \ \ \xi^{-1}\tOmega^1_{\fX/\cS}\otimes_{\co_\fX}\top_*(\top^*(\cN)\otimes_{\bvocB_\mQ} \bvcC^{(r)}_\mQ)
&\stackrel{\sim}{\rightarrow}&\top_*(\top^*(\xi^{-1}\tOmega^1_{\fX/\cS}\otimes_{\co_\fX}\cN)\otimes_{\bvocB_\mQ}
\bvcC^{(r)}_\mQ),\label{higgs3-MF101e}\\
\xi^{-1}\tOmega^1_{\fX/\cS}\otimes_{\co_\fX}\cN\otimes_{\co_\fX[\frac 1 p]}\top_*(\bvcC^{(r)}_\mQ)
&\stackrel{\sim}{\rightarrow}&\top_*(\top^*(\xi^{-1}\tOmega^1_{\fX/\cS}\otimes_{\co_\fX}\cN)\otimes_{\bvocB_\mQ}
\bvcC^{(r)}_\mQ).\label{higgs3-MF101f}
\end{eqnarray}
Le troisième isomorphisme est induit par les deux premiers d'après \ref{higgs3-MF4}(i). 
De plus, compte tenu du caractère bifonctoriel de l'isomorphisme \eqref{higgs3-MF6a}, le diagramme 
\begin{equation}\label{higgs3-MF101g}
\xymatrix{
{\cN\otimes_{\co_\fX[\frac 1 p]}\top_*(\bvcC^{(r)}_\mQ)}\ar[r]\ar[d]_{\theta\otimes \id+p^r\id\otimes \top_*(\bvd^{(r)}_\mQ)}&{\top_*(\top^*(\cN)\otimes_{\bvocB_\mQ} \bvcC^{(r)}_\mQ)}\ar[d]^{\top_*(\top^*(\theta)\otimes \id+p^r\id\otimes \bvd^{(r)})}\\
{\xi^{-1}\tOmega^1_{\fX/\cS}\otimes_{\co_\fX}\cN\otimes_{\co_\fX[\frac 1 p]}\top_*(\bvcC^{(r)}_\mQ)}\ar[r]&
{\top_*(\top^*(\xi^{-1}\tOmega^1_{\fX/\cS}\otimes_{\co_\fX}\cN)\otimes_{\bvocB_\mQ}
\bvcC^{(r)}_\mQ)}}
\end{equation}
où $\theta$ est le champs de Higgs de $\cN$, est commutatif. 
On prend alors pour  $\gamma^r$ \eqref{higgs3-MF101a} l'isomorphisme \eqref{higgs3-MF101d}. 
La proposition (i) résulte de \ref{higgs3-MF2}(ii) et \ref{higgs3-formel5}.
La proposition (ii) est une conséquence du caractère bifonctoriel de l'isomorphisme \eqref{higgs3-MF6a}.

\subsection{}\label{higgs3-MF91}
Soient $r$ un nombre rationnel $>0$, $(\cM,\cN,\alpha)$ un triplet $r$-admissible. 
Pour tout nombre rationnel $r'$ tel que $0< r'\leq r$, on désigne par 
\begin{equation}\label{higgs3-MF91a}
\alpha^{r'}\colon \top^{r'+}(\cN) \stackrel{\sim}{\rightarrow}\fS^{r'}(\cM)
\end{equation}
l'isomorphisme de $\Xi^{r'}_\mQ$ induit par $\epsilon^{r,r'}(\alpha)$ et les isomorphismes \eqref{higgs3-MF17g} et \eqref{higgs3-MF17h}, 
et par 
\begin{equation}\label{higgs3-MF91b}
\beta^{r'}\colon \cN \rightarrow\top^{r'}_+(\fS^{r'}(\cM))
\end{equation}
son adjoint \eqref{higgs3-MF15g}.

\begin{prop}\label{higgs3-MF9}
Les hypothèses étant celles de \eqref{higgs3-MF91}, soient, de plus, $r'$, $r''$ deux nombres rationnels 
tels que $0<r''< r'\leq r$. Alors~:
\begin{itemize}
\item[{\rm (i)}] Le morphisme composé 
\begin{equation}\label{higgs3-MF9a}
\cN\stackrel{\beta^{r'}}{\longrightarrow} \top^{r'}_+(\fS^{r'}(\cM))\longrightarrow \cH(\cM),
\end{equation}
où la seconde flèche est le morphisme canonique \eqref{higgs3-MF16a}, 
est un isomorphisme, indépendant de $r'$. 
\item[{\rm (ii)}] Le morphisme composé 
\begin{equation}\label{higgs3-MF9b}
\top^{r'}_+(\fS^{r'}(\cM))\longrightarrow \cH(\cM)\stackrel{\sim}{\longrightarrow}\cN\stackrel{\beta^{r''}}{\longrightarrow} 
\top^{r''}_+(\fS^{r''}(\cM))
\end{equation}
où la première flèche est le morphisme canonique \eqref{higgs3-MF16a} et la deuxième flèche est l'isomorphisme inverse de 
\eqref{higgs3-MF9a}, est le morphisme canonique \eqref{higgs3-MF16b}.
\end{itemize}
\end{prop}

(i) Pour tout nombre rationnel $0<t\leq r$, on désigne par  
\begin{equation}\label{higgs3-MF9c}
\gamma^t\colon \cN\otimes_{\co_\fX[\frac 1 p]}\top^t_+(\fS^t(\bvocB_\mQ)) 
\stackrel{\sim}{\rightarrow} \top^t_+(\top^{t+}(\cN))
\end{equation} 
l'isomorphisme \eqref{higgs3-MF101a}  de $\bMH(\co_\fX[\frac 1 p], \xi^{-1}\tOmega^1_{\fX/\cS})$, et par 
\begin{equation}\label{higgs3-MF9d}
\delta^t\colon \cN\otimes_{\co_\fX[\frac 1 p]}\top^t_+(\fS^t(\bvocB_\mQ)) 
\stackrel{\sim}{\rightarrow} \top^t_+(\fS^t(\cM))
\end{equation} 
le composé $\top^t_+(\alpha^t)\circ \gamma^t$. Le diagramme 
\begin{equation}\label{higgs3-MF9e}
\xymatrix{
{\cN\otimes_{\co_\fX[\frac 1 p]}\top^{r'}_+(\fS^{r'}(\bvocB_\mQ))}\ar[r]^-(0.5){\delta^{r'}}\ar[d]&{\top^{r'}_+(\fS^{r'}(\cM))}\ar[d]\\
{\cN\otimes_{\co_\fX[\frac 1 p]}\top^{r''}_+(\fS^{r''}(\bvocB_\mQ))}\ar[r]^-(0.5){\delta^{r''}}&{\top^{r''}_+(\fS^{r''}(\cM))}}
\end{equation}
où les flèches verticales sont les morphismes canoniques \eqref{higgs3-MF16b}, 
est commutatif en vertu de \ref{higgs3-MF101}(ii). Les isomorphismes $(\delta^t)_{0<t\leq r}$ 
induisent par passage à la limite inductive un isomorphisme de $\bMH(\co_\fX[\frac 1 p], \xi^{-1}\tOmega^1_{\fX/\cS})$
\begin{equation}\label{higgs3-MF9f}
\delta\colon \cN\otimes_{\co_\fX[\frac 1 p]}\cH(\bvocB_\mQ) \stackrel{\sim}{\rightarrow} \cH(\cM).
\end{equation}
Considérons le diagramme commutatif
\begin{equation}\label{higgs3-MF9g}
\xymatrix{
{\cN}\ar[r]^-(0.5){\iota^{r'}}\ar[rd]&
{\cN\otimes_{\co_\fX[\frac 1 p]}\top^{r'}_+(\fS^{r'}(\bvocB_\mQ))}\ar[r]^-(0.5){\delta^{r'}}
\ar[d]&{\top^{r'}_+(\fS^{r'}(\cM))}\ar[d]\\
&{\cN\otimes_{\co_\fX[\frac 1 p]}\cH(\bvocB_\mQ)}\ar[r]^-(0.5){\delta}&{\cH(\cM)}}
\end{equation}
où $\iota^{r'}$ est induit par le morphisme canonique $\co_\fX[\frac 1 p]\rightarrow \top^{r'}_+(\fS^{r'}(\bvocB_\mQ))$ 
et les flèches verticales sont les morphismes canoniques. D'après \ref{higgs3-MF101}(i), on a 
\begin{equation}\label{higgs3-MF9h}
\delta^{r'}\circ \iota^{r'}=\top^{r'}_+(\alpha^{r'})\circ\gamma^{r'}\circ \iota^{r'}=\beta^{r'}.
\end{equation} 
La proposition s'ensuit en vertu de \ref{higgs3-MF10}. 

(ii) Cela résulte de \eqref{higgs3-MF9e}, \eqref{higgs3-MF9g} et \ref{higgs3-coh25}(ii).

\begin{cor}\label{higgs3-MF20}
Pour tout $\bvocB_\mQ$-module de Dolbeault $\cM$, 
$\cH(\cM)$ \eqref{higgs3-MF16a} est un $\co_\fX[\frac 1 p]$-fibré de Higgs soluble associé à $\cM$.  
En particulier, $\cH$ induit un foncteur que l'on note encore
\begin{equation}\label{higgs3-MF20a}
\cH\colon \bMod^\Dolb_\mQ(\bvocB)\rightarrow \bMH^\sol(\co_\fX[\frac 1 p], \xi^{-1}\tOmega^1_{\fX/\cS}), 
\ \ \ \cM\mapsto \cH(\cM).
\end{equation}
\end{cor}

\begin{cor}\label{higgs3-MF99}
Pour tout $\bvocB_\mQ$-module de Dolbeault $\cM$, il existe un nombre rationnel $r>0$
et un isomorphisme de $\Xi^{r}_\mQ$
\begin{equation}\label{higgs3-MF99a}
\alpha\colon \top^{r+}(\cH(\cM))\stackrel{\sim}{\rightarrow} \fS^r(\cM)
\end{equation}
vérifiant les propriétés suivantes. Pour tout nombre rationnel $r'$ tel que $0< r'\leq r$, notons
\begin{equation}\label{higgs3-MF99b}
\alpha^{r'}\colon \top^{r'+}(\cH(\cM)) \stackrel{\sim}{\rightarrow}\fS^{r'}(\cM)
\end{equation}
l'isomorphisme de $\Xi^{r'}_\mQ$ induit par $\epsilon^{r,r'}(\alpha)$ et 
les isomorphismes \eqref{higgs3-MF17g} et \eqref{higgs3-MF17h}, et 
\begin{equation}\label{higgs3-MF99c}
\beta^{r'}\colon \cH(\cM)\rightarrow\top^{r'}_+(\fS^{r'}(\cM))
\end{equation}
son adjoint \eqref{higgs3-MF15g}. Alors~:
\begin{itemize}
\item[{\rm (i)}] Pour tout nombre rationnel $r'$ tel que $0< r'\leq r$, 
le morphisme $\beta^{r'}$ est un inverse à droite du morphisme canonique 
$\varpi^{r'}\colon\top^{r'}_+(\fS^{r'}(\cM))\rightarrow \cH(\cM)$. 
\item[{\rm (ii)}] Pour tous nombres rationnels $r'$ et $r''$ tels que $0<r''<r'\leq r$, le composé
\begin{equation}\label{higgs3-MF99d}
\top^{r'}_+(\fS^{r'}(\cM))\stackrel{\varpi^{r'}}{\longrightarrow} 
\cH(\cM)\stackrel{\beta^{r''}}{\longrightarrow} \top^{r''}_+(\fS^{r''}(\cM))
\end{equation}
est le morphisme canonique. 
\end{itemize}
\end{cor}

\begin{rema}\label{higgs3-MF92}
Sous les hypothèses de \ref{higgs3-MF99}, l'isomorphisme $\alpha$ n'est a priori pas uniquement déterminé par $(\cM,r)$, 
mais pour tout nombre rationnel $0<r'<r$, le morphisme $\alpha^{r'}$ \eqref{higgs3-MF99b} ne dépend que de $\cM$,
et il en dépend fonctoriellement (cf. la preuve de \ref{higgs3-MF21}).
\end{rema}

\subsection{}\label{higgs3-MF220}
Soient $r$ un nombre rationnel $>0$, $(\cM,\cN,\alpha)$ un triplet $r$-admissible. 
Pour éviter toute ambiguïté avec \eqref{higgs3-MF91a}, notons 
\begin{equation}\label{higgs3-MF220a}
\calpha\colon \fS^r(\cM)\rightarrow \top^{r+}(\cN)
\end{equation}
l'inverse de $\alpha$ dans $\Xi^r_\mQ$. 
Pour tout nombre rationnel $r'$ tel que $0< r'\leq r$, on désigne par 
\begin{equation}\label{higgs3-MF220b}
\calpha^{r'}\colon \fS^{r'}(\cM) \stackrel{\sim}{\rightarrow}\top^{r'+}(\cN)
\end{equation}
l'isomorphisme de $\Xi^{r'}_\mQ$ induit par $\epsilon^{r,r'}(\calpha)$ et les isomorphismes \eqref{higgs3-MF17g} et \eqref{higgs3-MF17h}, et par
\begin{equation}\label{higgs3-MF220c}
\cbeta^{r'}\colon \cM \rightarrow\cK^{r'}(\top^{r'+}(\cN))
\end{equation}
le morphisme adjoint. 

\begin{prop}\label{higgs3-MF22}
Les hypothèses étant celles de \eqref{higgs3-MF220}, soient, de plus, $r'$, $r''$ deux nombres rationnels 
tels que $0<r''< r'\leq r$. Alors~:
\begin{itemize} 
\item[{\rm (i)}] La limite inductive $\cV(\cN)$ du système inductif $(\cK^t(\top^{t+}(\cN)))_{t\in \mQ_{>0}}$ \eqref{higgs3-MF16c}
est représentable dans $\bMod_\mQ(\bvocB)$.
\item[{\rm (ii)}] Le morphisme composé 
\begin{equation}\label{higgs3-MF22a}
\cM\stackrel{\cbeta^{r'}}{\longrightarrow} \cK^{r'}(\top^{r'+}(\cN))\longrightarrow \cV(\cN),
\end{equation}
où la seconde flèche est le morphisme canonique, est un isomorphisme, indépendant de $r'$. 
\item[{\rm (iii)}] Le morphisme composé 
\begin{equation}\label{higgs3-MF22b}
\cK^{r'}(\top^{r'+}(\cN))\longrightarrow \cV(\cN)\stackrel{\sim}{\longrightarrow} \cM 
\stackrel{\cbeta^{r''}}{\longrightarrow} \cK^{r''}(\top^{r''+}(\cN))
\end{equation}
où la première flèche est le morphisme canonique et la seconde flèche est l'isomorphisme inverse de \eqref{higgs3-MF22a},
est le morphisme canonique \eqref{higgs3-MF16c}.
\end{itemize}
\end{prop}

(i) Comme $\cM$ est $\bvocB_\mQ$-plat d'après \ref{higgs3-plat3}, pour tout nombre rationnel $t\geq 0$, 
on a un isomorphisme canonique de $\bMod_\mQ(\bvocB)$
\begin{equation}\label{higgs3-MF22c}
\gamma^t\colon \cM\otimes_{\bvocB_\mQ}\cK^t(\fS^t(\bvocB_\mQ))
\stackrel{\sim}{\rightarrow}\cK^t(\fS^t(\cM)).
\end{equation}
On désigne par 
\begin{equation}\label{higgs3-MF22d}
\delta^t\colon \cM\otimes_{\bvocB_\mQ}\cK^t(\fS^t(\bvocB_\mQ))
\stackrel{\sim}{\rightarrow}\cK^t(\top^{t+}(\cN))
\end{equation}
le composé $\cK^t(\calpha^t)\circ \gamma^t$. Le diagramme 
\begin{equation}\label{higgs3-MF22e}
\xymatrix{
{\cM\otimes_{\bvocB_\mQ}\cK^{r'}(\fS^{r'}(\bvocB_\mQ))}\ar[d]\ar[r]^-(0.6){\delta^{r'}}&
{\cK^{r'}(\top^{r'+}(\cN))}\ar[d]\\
{\cM\otimes_{\bvocB_\mQ}\cK^{r''}(\fS^{r''}(\bvocB_\mQ))}\ar[r]^-(0.6){\delta^{r''}}&
{\cK^{r''}(\top^{r''+}(\cN))}}
\end{equation}
où les flèches verticales sont induites par le morphisme \eqref{higgs3-MF17k} et les isomorphismes 
\eqref{higgs3-MF17g} et \eqref{higgs3-MF17h}, est clairement commutatif. La proposition résulte alors de \ref{higgs3-coh24}. 

(ii) D'après \ref{higgs3-coh24}, le morphisme canonique 
\begin{equation}\label{higgs3-MF22f}
\cM\rightarrow \underset{\underset{t\in \mQ_{>0}}{\longrightarrow}}{\lim}\ 
\cM\otimes_{\bvocB_\mQ} \cK^t(\fS^t(\bvocB_\mQ))
\end{equation}
est un isomorphisme. Les isomorphismes $(\delta^t)_{0<t\leq r}$ 
induisent alors par passage à la limite inductive un isomorphisme
\begin{equation}\label{higgs3-MF22g}
\delta\colon \cM\stackrel{\sim}{\rightarrow} \cV(\cN).
\end{equation}
Il résulte aussitôt des définitions que le diagramme 
\begin{equation}\label{higgs3-MF22h}
\xymatrix{
{\cM\otimes_{\bvocB_\mQ} \cK^{r'}(\fS^{r'}(\bvocB_\mQ))}
\ar[r]^-(0.5){\delta^{r'}}&{\cK^{r'}(\top^{r'+}(\cN))}\ar[d]\\
{\cM}\ar[r]^-(0.5){\delta}\ar[u]^-(0.5){\iota^{r'}}&{\cV(\cN)}}
\end{equation}
où $\iota^{r'}$ est induit par le morphisme canonique $\bvocB_\mQ\rightarrow \cK^{r'}(\fS^{r'}(\bvocB_\mQ))$
et la flèche non-libellée est le morphisme canonique, est commutatif. On vérifie aussitôt qu'on a 
\begin{equation}\label{higgs3-MF22i}
\delta^{r'}\circ \iota^{r'}=\cK^{r'}(\calpha^{r'})\circ \gamma^{r'}\circ \iota^{r'}=\cbeta^{r'}.
\end{equation}
La proposition s'ensuit. 

(iii) Cela résulte de \eqref{higgs3-MF22e} et \ref{higgs3-coh24}(ii).

\begin{cor}\label{higgs3-MF23}\index{10001240@$\cV$}
On a un foncteur
\begin{equation}\label{higgs3-MF23a}
\cV\colon \bMH^\sol(\co_\fX[\frac 1 p], \xi^{-1}\tOmega^1_{\fX/\cS})\rightarrow \bMod^\Dolb_\mQ(\bvocB), 
\ \ \ \cN\mapsto \underset{\underset{r\in \mQ_{>0}}{\longrightarrow}}{\lim}\ \cK^{r}(\top^{r+}(\cN))).
\end{equation}
De plus, pour tout objet $\cN$ de $\bMH^\sol(\co_\fX[\frac 1 p], \xi^{-1}\tOmega^1_{\fX/\cS})$, 
$\cV(\cN)$ est associé à $\cN$.
\end{cor}

\begin{cor}\label{higgs3-MF24}
Pour tout $\co_\fX[\frac 1 p]$-fibré de Higgs soluble $\cN$ à coefficients dans $\xi^{-1}\tOmega^1_{\fX/\cS}$,
il existe un nombre rationnel $r>0$ et un isomorphisme de $\Xi^r_\mQ$
\begin{equation}\label{higgs3-MF24a}
\calpha\colon \fS^r(\cV(\cN))\stackrel{\sim}{\rightarrow}\top^{r+}(\cN)
\end{equation}
vérifiant les propriétés suivantes. Pour tout nombre rationnel $r'$ tel que $0< r'\leq r$, notons
\begin{equation}\label{higgs3-MF24b}
\calpha^{r'}\colon\fS^{r'}(\cV(\cN))\stackrel{\sim}{\rightarrow}\top^{r'+}(\cN)
\end{equation}
l'isomorphisme de $\Xi^{r'}_\mQ$ induit par $\epsilon^{r,r'}(\calpha)$ et 
les isomorphismes \eqref{higgs3-MF17g} et \eqref{higgs3-MF17h}, et 
\begin{equation}\label{higgs3-MF24c}
\cbeta^{r'}\colon \cV(\cN)\rightarrow\cK^{r'}(\top^{r'+}(\cN))
\end{equation}
son adjoint. Alors~:
\begin{itemize}
\item[{\rm (i)}] Pour tout nombre rationnel $r'$ tel que $0<r'\leq r$,
le morphisme $\cbeta^{r'}$ est un inverse à droite du morphisme canonique
$\varpi^{r'}\colon \cK^{r'}(\top^{r'+}(\cN))\rightarrow \cV(\cN)$.  
\item[{\rm (ii)}] Pour tous nombres rationnels $r'$ et $r''$ tels que $0<r''<r'\leq r$, le composé
\begin{equation}\label{higgs3-MF24d}
\cK^{r'}(\top^{r'+}(\cN))\stackrel{\varpi^{r'}}{\longrightarrow} 
\cV(\cN)\stackrel{\cbeta^{r''}}{\longrightarrow} \cK^{r''}(\top^{r''+}(\cN))
\end{equation}
est le morphisme canonique. 
\end{itemize}
\end{cor}

\begin{rema}\label{higgs3-MF241}
Sous les hypothèses de \ref{higgs3-MF24}, l'isomorphisme $\calpha$ n'est a priori pas uniquement déterminé par $(\cN,r)$, 
mais pour tout nombre rationnel $0<r'<r$, 
le morphisme $\calpha^{r'}$ \eqref{higgs3-MF24b} ne dépend que de $\cN$ et il en dépend fonctoriellement (cf. la preuve
de \ref{higgs3-MF21}). 
\end{rema}

\begin{teo}\label{higgs3-MF21}
Les foncteurs \eqref{higgs3-MF20a} et \eqref{higgs3-MF23a} 
\begin{equation}
\xymatrix{
{\bMod^\Dolb_\mQ(\bvocB)}\ar@<1ex>[r]^-(0.5){\cH}&{\bMH^\sol(\co_\fX[\frac 1 p], \xi^{-1}\tOmega^1_{\fX/\cS})}
\ar@<1ex>[l]^-(0.5){\cV}}
\end{equation}
sont des équivalences de catégories quasi-inverses l'une de l'autre. 
\end{teo}

Pour tout objet $\cM$ de $\bMod^\Dolb_\mQ(\bvocB)$,
$\cH(\cM)$ est un $\co_\fX[\frac 1 p]$-fibré de Higgs soluble associé à $\cM$, en vertu de \ref{higgs3-MF20}. 
Choisissons un nombre rationnel $r_\cM>0$ et un isomorphisme de $\Xi^{r_\cM}_\mQ$ 
\begin{equation}\label{higgs3-MF21a}
\alpha_\cM\colon \top^{r_\cM+}(\cH(\cM))\stackrel{\sim}{\rightarrow} \fS^{r_\cM}(\cM)
\end{equation}
vérifiant les propriétés du \ref{higgs3-MF99}. Pour tout nombre rationnel $r$ tel que $0< r\leq r_\cM$, on désigne par
\begin{equation}\label{higgs3-MF21b}
\alpha^r_\cM\colon \top^{r+}(\cH(\cM)) \stackrel{\sim}{\rightarrow}\fS^r(\cM)
\end{equation}
l'isomorphisme de $\Xi^r_\mQ$ induit par $\epsilon^{r_\cM,r}(\alpha_\cM)$ et 
les isomorphismes \eqref{higgs3-MF17g} et \eqref{higgs3-MF17h}, par 
\begin{eqnarray}
\calpha_\cM\colon \fS^{r_\cM}(\cM)&\stackrel{\sim}{\rightarrow}& \top^{r_\cM+}(\cH(\cM)),\label{higgs3-MF21c}\\
\calpha^r_\cM\colon \fS^r(\cM) &\stackrel{\sim}{\rightarrow}&\top^{r+}(\cH(\cM)),\label{higgs3-MF21d}
\end{eqnarray}
les inverses de $\alpha_\cM$ et $\alpha_\cM^r$, respectivement, et par
\begin{eqnarray}
\beta^r_\cM\colon \cH(\cM)&\rightarrow&\top^r_+(\fS^r(\cM)),\label{higgs3-MF21e}\\
\cbeta^r_\cM\colon \cM&\rightarrow&\cK^r(\top^{r+}(\cH(\cM))),\label{higgs3-MF21f}
\end{eqnarray}
les morphismes adjoints de $\alpha^r_\cM$ et $\calpha^r_\cM$ respectivement. On notera que $\calpha^r_\cM$ est
induit par $\epsilon^{r_\cM,r}(\calpha_\cM)$ et les isomorphismes \eqref{higgs3-MF17g} et \eqref{higgs3-MF17h}.
D'après \ref{higgs3-MF22}(ii), le morphisme composé
\begin{equation}\label{higgs3-MF21g}
\cM\stackrel{\cbeta^{r}_\cM}{\longrightarrow}\cK^r(\top^{r+}(\cH(\cM)))\longrightarrow \cV(\cH(\cM)),
\end{equation}
où la seconde flèche est le morphisme canonique, 
est un isomorphisme, qui dépend a priori de $\alpha_\cM$ mais pas de $r$.  Montrons que cet isomorphisme
ne dépend que de $\cM$ (mais pas du choix de $\alpha_\cM$) et qu'il en dépend fonctoriellement. 
Il suffit de montrer que pour tout morphisme $u\colon \cM\rightarrow \cM'$ de $\bMod^\Dolb_\mQ(\bvocB)$
et tout nombre rationnel $0<r<\inf(r_\cM,r_{\cM'})$, le diagramme de $\Xi^r_\mQ$
\begin{equation}
\xymatrix{
{\top^{r+}(\cH(\cM))}\ar[r]^-(0.5){\alpha^r_\cM}\ar[d]_{\top^{r+}(\cH(u))}&{\fS^r(\cM)}\ar[d]^{\fS^r(u)}\\
{\top^{r+}(\cH(\cM'))}\ar[r]^-(0.5){\alpha^r_{\cM'}}&{\fS^r(\cM')}}
\end{equation}
est commutatif. Soient $r$, $r'$ deux nombres rationnels tels que $0<r<r'<\inf(r_\cM,r_{\cM'})$. 
Considérons le diagramme 
\begin{equation}
\xymatrix{
{\top^{r'}_+(\fS^{r'}(\cM))}\ar[r]^-(0.4){\varpi^{r'}_\cM}\ar[d]_{\top^{r'}_+(\fS^{r'}(u))}&
{\cH(\cM)}\ar[r]^-(0.5){\beta^{r}_\cM}\ar[d]^{\cH(u)}&{\top^{r}_+(\fS^{r}(\cM))}\ar[d]^{\top^{r}_+(\fS^{r}(u))}\\
{\top^{r'}_+(\fS^{r'}(\cM'))}\ar[r]^-(0.4){\varpi^{r'}_{\cM'}}&{\cH(\cM')}\ar[r]^-(0.5){\beta^{r}_{\cM'}}&{\top^{r}_+(\fS^{r}(\cM'))}}
\end{equation}
où $\varpi^{r'}_\cM$ et $\varpi^{r'}_{\cM'}$ sont les morphismes canoniques. 
Il résulte de \ref{higgs3-MF99}(ii) que le grand rectangle est commutatif. 
Comme le carré de gauche est commutatif et que $\varpi^{r'}_{\cM'}$ est surjectif d'après \ref{higgs3-MF99}(i),
le carré de droite est aussi commutatif. L'assertion recherchée s'ensuit compte tenu de l'injectivité de \eqref{higgs3-MF15g}. 

De même, pour tout objet $\cN$ de $\bMH^\sol(\co_\fX[\frac 1 p], \xi^{-1}\tOmega^1_{\fX/\cS})$,
$\cV(\cN)$ est un $\bvocB_\mQ$-module de Dolbeault associé à $\cN$, en vertu de \ref{higgs3-MF23}.
Choisissons un nombre rationnel $r_\cN>0$ et un isomorphisme de $\Xi^{r_\cN}_\mQ$ 
\begin{equation}\label{higgs3-MF21aa}
\calpha_\cN\colon \fS^{r_\cN}(\cV(\cN))\stackrel{\sim}{\rightarrow} \top^{r_\cN+}(\cN)
\end{equation}
vérifiant les propriétés du \ref{higgs3-MF24}. Pour tout nombre rationnel $r$ tel que $0< r\leq r_\cN$, on désigne par
\begin{equation}\label{higgs3-MF21bb}
\calpha^r_\cN\colon \fS^r(\cV(\cN))\stackrel{\sim}{\rightarrow} \top^{r+}(\cN)
\end{equation}
l'isomorphisme de $\Xi^r_\mQ$ induit par $\epsilon^{r_\cN,r}(\calpha_\cN)$ et 
les isomorphismes \eqref{higgs3-MF17g} et \eqref{higgs3-MF17h}, par 
\begin{eqnarray}
\alpha_\cN\colon \top^{r_\cN+}(\cN) &\stackrel{\sim}{\rightarrow}&\fS^{r_\cN}(\cV(\cN)),\label{higgs3-MF21cc}\\
\alpha^r_\cN\colon \top^{r+}(\cN)&\stackrel{\sim}{\rightarrow} &\fS^r(\cV(\cN)),\label{higgs3-MF21dd}
\end{eqnarray}
les inverses de $\calpha_\cM$ et $\calpha_\cN^r$, respectivement, et par
\begin{eqnarray}
\cbeta^r_\cN\colon \cV(\cN)&\rightarrow&\cK^r(\top^{r+}(\cN)),\label{higgs3-MF21ee}\\
\beta^r_\cN\colon \cN&\rightarrow&\top^r_+(\fS^r(\cV(\cN))),\label{higgs3-MF21ff}
\end{eqnarray}
les morphismes adjoints de $\calpha^r_\cN$ et $\alpha^r_\cN$, respectivement. D'après \ref{higgs3-MF9}(i), 
le morphisme composé 
\begin{equation}
\cN\stackrel{\beta^r}{\longrightarrow} \top^r_+(\fS^r(\cV(\cN)))\longrightarrow \cH(\cV(\cN)),
\end{equation}
où la seconde flèche est le morphisme canonique, est un isomorphisme, qui dépend a priori de $\calpha_\cN$
mais pas de $r$. Montrons que cet isomorphisme
ne dépend que de $\cN$ (mais pas du choix de $\calpha_\cN$) et qu'il en dépend fonctoriellement. 
Il suffit de montrer que pour tout morphisme $v\colon \cN\rightarrow \cN'$ de 
$\bMH^\sol(\co_\fX[\frac 1 p], \xi^{-1}\tOmega^1_{\fX/\cS})$ et tout nombre rationnel 
$0<r<\inf(r_\cN,r_{\cN'})$, le diagramme de $\Xi^r_\mQ$
\begin{equation}
\xymatrix{
{\fS^r(\cV(\cN))}\ar[r]^-(0.5){\calpha^r_\cN}\ar[d]_{\fS^r(\cV(v))}&{\top^{r+}(\cN)}\ar[d]^{\top^{r+}(v)}\\
{\fS^r(\cV(\cN'))}\ar[r]^-(0.5){\calpha^r_{\cN'}}&{\top^{r+}(\cN')}}
\end{equation}
est commutatif. Soient $r$, $r'$ deux nombres rationnels tels que $0<r<r'<\inf(r_\cN,r_{\cN'})$. 
Considérons le diagramme de $\bMod_\mQ(\bvocB)$
\begin{equation}
\xymatrix{
{\cK^{r'}(\top^{r'+}(\cN))}\ar[r]^-(0.5){\varpi^{r'}_\cN}\ar[d]_{\cK^{r'}(\top^{r'+}(v))}&
{\cV(\cN)}\ar[r]^-(0.5){\cbeta^{r}_\cM}\ar[d]^{\cV(v)}&{\cK^{r}(\top^{r+}(\cN))}\ar[d]^{\cK^{r}(\top^{r+}(v))}\\
{\cK^{r'}(\top^{r'+}(\cN'))}\ar[r]^-(0.5){\varpi^{r'}_{\cN'}}&{\cV(\cN')}\ar[r]^-(0.5){\cbeta^{r}_{\cN'}}&{\cK^{r}(\top^{r+}(\cN'))}}
\end{equation}
où $\varpi^{r'}_\cN$ et $\varpi^{r'}_{\cN'}$ sont les morphismes canoniques. 
Il résulte de \ref{higgs3-MF24}(ii) que le grand rectangle est commutatif. 
Comme le carré de gauche est commutatif et que $\varpi_\cN^{r'}$ est inversible à droite d'après \ref{higgs3-MF24}(i),
le carré de droite est aussi commutatif~; d'où l'assertion recherchée.

\subsection{}\label{higgs3-MF25}
Soient $(\oy\rightsquigarrow \ox)$ un point de $X_\et\gtimes_{X_\et}\oX^\circ_\et$ \eqref{higgs3-tfa41} tel que $\ox$ soit
au-dessus de $s$, $X'$ le localisé strict de $X$ en $\ox$, $R'_1=\Gamma(\oX',\co_{\oX'})$, 
$\hR'_1$ son séparé complété $p$-adique. Pour tout $\co_\fX$-module cohérent $\cF$, on pose (abusivement)
\begin{equation}\label{higgs3-MF25a}
\cF_\ox=\underset{\underset{n\in \mN^\circ}{\longleftarrow}}{\lim}\ \Gamma(\oX'_n,\cF\otimes_{\co_\fX}\co_{\oX'_n}).
\end{equation}
D'après \ref{higgs3-formel2}, le foncteur $\bMod^\coh(\co_\fX)\rightarrow \bMod(\hR'_1)$ ainsi défini,  
induit un foncteur que l'on note encore
\begin{equation}\label{higgs3-MF25b}
\bMod^\coh(\co_\fX[\frac 1 p])\rightarrow \bMod(\hR'_1[\frac 1 p]), \ \ \ \cF\mapsto \cF_\ox. 
\end{equation}

D'après \ref{higgs3-TFT5}(i), $\rho(\oy\rightsquigarrow \ox)$ est un point de $\tE_s$. 
Pour tout objet $\fF=(\fF_n)_{n\in \mN}$ de  $\tE_s^{\mN^\circ}$, on pose (abusivement)
\begin{equation}\label{higgs3-MF25c}
\fF_{\rho(\oy\rightsquigarrow \ox)}=\underset{\underset{n\in \mN^\circ}{\longleftarrow}}{\lim}\ (\fF_n)_{\rho(\oy\rightsquigarrow \ox)}.
\end{equation}
On prendra garde que le foncteur $\tE_s^{\mN^\circ}\rightarrow \Ens$ ainsi défini n'est pas a priori un foncteur fibre.
Reprenons les notations de \ref{higgs3-RGG16}. D'après \eqref{higgs3-RGG16m} et (\cite{ag2} 10.31 et 9.8), 
on a un isomorphisme canonique
\begin{equation}\label{higgs3-MF25d}
\bvocB_{\rho(\oy\rightsquigarrow \ox)}\stackrel{\sim}{\rightarrow} \hoR^\oy_{X'}.
\end{equation}
Le foncteur \eqref{higgs3-MF25c} induit des foncteurs que l'on note encore
\begin{eqnarray}
\bMod(\bvocB)&\rightarrow& \bMod(\hoR^\oy_{X'}), \ \ \ \fM\mapsto \fM_{\rho(\oy\rightsquigarrow \ox)},\label{higgs3-MF25e}\\
\bMod_\mQ(\bvocB)&\rightarrow& \bMod(\hoR^\oy_{X'}[\frac 1 p]), \ \ \ \cM\mapsto \cM_{\rho(\oy\rightsquigarrow \ox)}.\label{higgs3-MF25f}
\end{eqnarray}

D'après \ref{higgs3-RGG172}, \ref{higgs3-RGG17} et (\cite{ag2} 10.31 et 9.8), pour tout nombre rationnel $r\geq 0$, 
on a un isomorphisme canonique
\begin{equation}\label{higgs3-MF25g}
\bvcC^{(r)}_{\rho(\oy\rightsquigarrow \ox)} \stackrel{\sim}{\rightarrow}\hcC^{\oy,(r)}_{X'}.
\end{equation}

\begin{lem}\label{higgs3-MF251}
Les hypothèses étant celles de \eqref{higgs3-MF25}, soient, de plus, $\cF$ un $\co_\fX$-module cohérent,
$(U,\fp\colon \ox\rightarrow U)$ un objet de $\fV_\ox$ \eqref{higgs3-RGG3} tel que $U$ soit affine, 
$\cU$ le schéma formel complété $p$-adique de $\oU$. Alors,
\begin{itemize}
\item[{\rm (i)}]  Le $\hR'_1$-module $\cF_\ox$ est de type fini, et est complet et séparé pour la topologie $p$-adique.
\item[{\rm (ii)}] On a un isomorphisme canonique et fonctoriel
\begin{equation}\label{higgs3-MF251b}
\cF_\ox\stackrel{\sim}{\rightarrow} \Gamma(\cU,\cF\otimes_{\co_\fX}\co_\cU)\hotimes_{\co_\cU(\cU)}\hR'_1.
\end{equation}
\item[{\rm (iii)}] On a un isomorphisme canonique et fonctoriel
\begin{equation}\label{higgs3-MF251c}
(\top^*(\cF))_{\rho(\oy\rightsquigarrow \ox)} \stackrel{\sim}{\rightarrow} \cF_\ox\hotimes_{\hR'_1} \hoR^\oy_{X'}.
\end{equation}
\end{itemize}
\end{lem}

(i) Cela résulte de (\cite{ac} chap.~III § 2.11 prop.~14).

(ii) Cela résulte immédiatement de la définition. 

(iii) D'après (\cite{sga4} VIII 5.2 et VII 5.8), on a un isomorphisme canonique 
\begin{equation}\label{higgs3-MF251d}
(\cF\otimes_{\co_\fX}\co_{\oX_n})_\ox\stackrel{\sim}{\rightarrow} \Gamma(\oX'_n,\cF\otimes_{\co_\fX}\co_{\oX'_n}),
\end{equation}
où l'on considère à gauche $\cF\otimes_{\co_\fX}\co_{\oX_n}$ comme un faisceau de $X_{s,\et}$ \eqref{higgs3-TFT9}. 
Compte tenu de \eqref{higgs3-not2e}, \eqref{higgs3-coh0g} et (\cite{ag2} (10.18.1)), on en déduit un isomorphisme canonique et fonctoriel
\begin{equation}\label{higgs3-MF251e}
(\top^*(\cF))_{\rho(\oy\rightsquigarrow \ox)} \stackrel{\sim}{\rightarrow} 
\underset{\underset{n\in \mN^\circ}{\longleftarrow}}{\lim}\ 
(\Gamma(\oX'_n,\cF\otimes_{\co_\fX}\co_{\oX'_n})\otimes_{R'_1} \oR^\oy_{X'}).
\end{equation}
La proposition s'ensuit compte tenu de (\cite{ac} chap.~III § 2.11 cor.~1 de prop.~14). 

\subsection{}\label{higgs3-MF250}
Conservons les hypothèses et notations de \ref{higgs3-MF25}; posons, de plus, $R'=\Gamma(X',\co_{X'})$ et
\begin{equation}\label{higgs3-MF250a}
\tOmega^1_{X/S}(X') = \Gamma(X',\tOmega^1_{X/S}\otimes_{\co_X}\co_{X'}). 
\end{equation}
Pour toute $R'_1$-algèbre $A$, on sous-entend par $A$-module de Higgs 
à coefficients dans $\xi^{-1}\tOmega^1_{X/S}(X')$ un $A$-module de Higgs à coefficients dans 
$\xi^{-1}\tOmega^1_{X/S}(X')\otimes_{R'}A$ (\cite{ag1} 2.8). 
Compte tenu de \eqref{higgs3-MF12c} et du fait que le $R'$-module $\tOmega^1_{X/S}(X')$ est libre de type fini, 
le foncteur \eqref{higgs3-MF25b} induit un foncteur 
\begin{equation}\label{higgs3-MF250c}
\bMH^\coh(\co_\fX[\frac 1 p],\xi^{-1}\tOmega^1_{\fX/\cS}) \rightarrow \bMH(\hR'_1[\frac 1 p],\xi^{-1}\tOmega^1_{X/S}(X')),
\ \ \ (\cN,\theta)\mapsto (\cN_\ox,\theta_\ox).
\end{equation}

Soit $r$ un nombre rationnel $\geq 0$. D'après \eqref{higgs3-RGG16l}, on a un isomorphisme canonique
\begin{equation}
\Omega^1_{\cC^{\oy,(r)}_{X'}/\hoR^\oy_{X'}}\stackrel{\sim}{\rightarrow} \xi^{-1}\tOmega^1_{X/S}(X')\otimes_{R'}\cC^{\oy,(r)}_{X'}.
\end{equation}
On désigne par
\begin{equation}\label{higgs3-MF250d}
d_{\cC^{\oy,(r)}_{X'}}\colon \cC^{\oy,(r)}_{X'}\rightarrow \xi^{-1}\tOmega^1_{X/S}(X')\otimes_{R'}\cC^{\oy,(r)}_{X'}
\end{equation}
la $\hoR^\oy_{X'}$-dérivation universelle de $\cC^{\oy,(r)}_{X'}$ (cf. \ref{higgs3-RGG16}) et par 
\begin{equation}\label{higgs3-MF250e}
d_{\hcC^{\oy,(r)}_{X'}}\colon \hcC^{\oy,(r)}_{X'}\rightarrow \xi^{-1}\tOmega^1_{X/S}(X')\otimes_{R'}\hcC^{\oy,(r)}_{X'}
\end{equation}
son prolongement aux complétés $p$-adiques (le $R'$-module $\tOmega^1_{X/S}(X')$ étant libre de type fini).
Ce sont des $\hoR^\oy_{X'}$-champs de Higgs à coefficients dans $\xi^{-1}\tOmega^1_{X/S}(X')$ d'après (\cite{ag1} 2.12 et 2.16).

On désigne par $\Xi^{\atf,r}$ la sous-catégorie pleine de $\Xi^r$ \eqref{higgs3-MF15} formée des $p^r$-isoconnexions intégrables 
$(\cF,\cG,u,\nabla)$ relativement à l'extension $\bvcC^{(r)}/\bvocB$ telles que les $\bvcC^{(r)}$-modules 
$\cF$ et $\cG$ soient adiques de type fini \eqref{higgs3-spad3}, et par  $\Xi'^r$ la catégorie des $p^r$-isoconnexions $p$-adiques intégrables 
relativement à l'extension $\hcC_{X'}^{\oy,(r)}/\hoR_{X'}^\oy$ \eqref{higgs3-isoco4}. Ce sont des catégories additives. 
On note $\Xi^{\atf,r}_\mQ$ et $\Xi'^r_\mQ$ les catégories des objets de $\Xi^{\atf,r}$ et $\Xi'^r$ à isogénie près \eqref{higgs3-caip1a}.
Compte tenu de (\cite{ac} chap.~III § 2.11 prop.~14 et cor.~1), le foncteur \eqref{higgs3-MF25c} induit un foncteur additif 
\begin{equation}\label{higgs3-MF250f}
\Xi^{\atf,r}\rightarrow \Xi'^r.
\end{equation}
D'après \ref{higgs3-isoco5}, tout objet de $\Xi'^r$ est une $\hoR^\oy_{X'}$-isogénie
de Higgs à coefficients dans $\xi^{-1}\tOmega^1_{X/S}(X')$. 
On dispose donc du foncteur 
\begin{equation}\label{higgs3-MF250g}
\Xi'^r_{\mQ}\rightarrow \bMH(\hoR^\oy_{X'}[\frac 1 p],\xi^{-1}\tOmega^1_{X/S}(X')), \ \ \ (M,N,u,\nabla)\mapsto (M_{\mQ_p},
(\id \otimes u_{\mQ_p}^{-1})\circ\nabla_{\mQ_p}).
\end{equation}
On en déduit un foncteur
\begin{equation}\label{higgs3-MF250h}
\phi_{\rho(\oy\rightsquigarrow \ox)} \colon \Xi^{\atf,r}_\mQ\rightarrow  \bMH(\hoR^\oy_{X'}[\frac 1 p],\xi^{-1}\tOmega^1_{X/S}(X')).
\end{equation}

Soit $\fM$ un $\bvocB$-module adique de type fini. D'après (\cite{ac} chap.~III § 2.11 prop.~14 et cor.~1),
le $\hoR_{X'}^\oy$-module $\fM_{\rho(\oy\rightsquigarrow \ox)}$ est de type fini, et est complet et séparé pour 
la topologie $p$-adique. On a un isomorphisme canonique et fonctoriel \eqref{higgs3-MF15bb}
\begin{equation}\label{higgs3-MF250i}
\phi_{\rho(\oy\rightsquigarrow \ox)}(\fS^r(\fM))\stackrel{\sim}{\rightarrow}
((\hcC^{(r),\oy}_{X'}\hotimes_{\hoR^\oy_{X'}}\fM_{\rho(\oy\rightsquigarrow \ox)})_{\mQ_p},
(p^rd_{\hcC^{\oy,(r)}_{X'}}\hotimes \id)_{\mQ_p}).
\end{equation}

\begin{lem}\label{higgs3-MF252}
Les hypothèses étant celles de \eqref{higgs3-MF25} et \eqref{higgs3-MF250}, 
soient, de plus, $(\cN,\theta)$ un $\co_\fX[\frac 1 p]$-fibré de Higgs à coefficients 
dans $\xi^{-1}\tOmega^1_{\fX/\cS}$. Alors, 
\begin{itemize}
\item[{\rm (i)}] Le $\hR'_1[\frac 1 p]$-module $\cN_\ox$ est projectif de type fini.
\item[{\rm (ii)}] Pour tout nombre rationnel $r\geq 0$, on a un isomorphisme canonique et fonctoriel
\begin{equation}\label{higgs3-MF252a}
\phi_{\rho(\oy\rightsquigarrow \ox)}(\top^{r+}(\cN,\theta))\stackrel{\sim}{\rightarrow}
(\hcC^{(r),\oy}_{X'}\otimes_{\hR'_1}\cN_\ox, p^rd_{\hcC^{\oy,(r)}_{X'}}\otimes \id+\id\otimes \theta_\ox).
\end{equation}
\end{itemize}
\end{lem}

En effet, soit $\fN$ un $\co_\fX$-module cohérent tel que $\fN_{\mQ_p}=\cN$ \eqref{higgs3-formel2}.

(i) Soient $(U,\fp\colon \ox\rightarrow U)$ un objet de $\fV_\ox(\bP)$, $\cU$ le schéma formel complété $p$-adique de $\oU$. 
D'après \ref{higgs3-MF251}(ii), on a un isomorphisme canonique 
\begin{equation}\label{higgs3-MF252b}
\Gamma(\cU,\fN\otimes_{\co_\fX}\co_\cU)\hotimes_{\co_\cU(\cU)}\hR'_1\stackrel{\sim}{\rightarrow} \fN_\ox.
\end{equation}
On a $\Gamma(\cU,\fN\otimes_{\co_\fX}\co_\cU)\otimes_{\mZ_p}\mQ_p=\Gamma(\cU,\cN\otimes_{\co_\fX}\co_\cU)$
(\cite{egr1} (2.10.5.1)). D'après \ref{higgs3-formel20}, le $\co_\cU(\cU)[\frac 1 p]$-module 
$\Gamma(\cU,\cN\otimes_{\co_\fX}\co_\cU)$ est projectif de type fini. 
On en déduit, en vertu de  \ref{higgs3-formel21} et \ref{higgs3-RGG160}(ii), que le morphisme 
\begin{equation}\label{higgs3-MF252c}
\Gamma(\cU,\fN\otimes_{\co_\fX}\co_\cU)\otimes_{\co_\cU(\cU)}\hR'_1[\frac 1 p]\rightarrow \cN_\ox
\end{equation}
induit par \eqref{higgs3-MF252b} est un isomorphisme~; d'où la proposition.

(ii) D'après \ref{higgs3-MF251}(iii), on a un isomorphisme canonique et fonctoriel
\begin{equation}\label{higgs3-MF252d}
(\top^*(\fN)\otimes_{\bvocB}\bvcC^{(r)})_{\rho(\oy\rightsquigarrow \ox)}
\stackrel{\sim}{\rightarrow} \fN_\ox\hotimes_{\hR'_1} \hcC^{(r),\oy}_{X'}.
\end{equation}
La proposition résulte alors de (i), \ref{higgs3-formel21} et \ref{higgs3-RGG160}(ii).

\begin{lem}\label{higgs3-MF253}
Les hypothèses étant celles de \eqref{higgs3-MF25} et \eqref{higgs3-MF250}, 
soient, de plus, $\cM$ un objet de $\bMod^\atf_\mQ(\bvocB)$, $(\cN,\theta)$ un $\co_\fX[\frac 1 p]$-fibré de Higgs à coefficients dans 
$\xi^{-1}\tOmega^1_{\fX/\cS}$, $r$ un nombre rationnel $>0$, 
\begin{equation}\label{higgs3-MF253a}
\alpha\colon \top^{r+}(\cN,\theta) \stackrel{\sim}{\rightarrow}\fS^r(\cM)
\end{equation}
un isomorphisme de $\Xi^r_\mQ$ \eqref{higgs3-MF15}. Alors,
\begin{itemize}
\item[{\rm (i)}] Le $\hoR_{X'}^\oy[\frac 1 p]$-module $\cM_{\rho(\oy\rightsquigarrow \ox)}$ est projectif de type fini.
\item[{\rm (ii)}] L'isomorphisme $\phi_{\rho(\oy\rightsquigarrow \ox)}(\alpha)$ \eqref{higgs3-MF250h} induit 
un isomorphisme $\hcC^{\oy,(r)}_{X'}$-linéaire 
\begin{equation}\label{higgs3-MF253b}
\alpha_{\rho(\oy\rightsquigarrow \ox)}\colon \hcC^{\oy,(r)}_{X'}\otimes_{\hR'_1}\cN_\ox\stackrel{\sim}{\rightarrow}
\hcC^{\oy,(r)}_{X'}\otimes_{\hoR^\oy_{X'}}\cM_{\rho(\oy\rightsquigarrow \ox)}
\end{equation}
de $\hoR^\oy_{X'}$-modules de Higgs à coefficients dans $\xi^{-1}\tOmega^1_{X/S}(X')$,
où $\cN_\ox$ est muni du champ de Higgs $\theta_\ox$ \eqref{higgs3-MF250c},
$\hcC^{\oy,(r)}_{X'}$ est muni du champ de Higgs $p^rd_{\hcC^{\oy,(r)}_{X'}}$ \eqref{higgs3-MF250e}
et $\cM_{\rho(\oy\rightsquigarrow \ox)}$ est muni du champ de Higgs nul.
\end{itemize}
\end{lem}

En effet, soit $\fM$ un objet de $\bMod^\atf(\bvocB)$ tel que $\cM=\fM_{\mQ_p}$.
D'après \eqref{higgs3-MF250i} et \ref{higgs3-MF252}, $\phi_{\rho(\oy\rightsquigarrow \ox)}(\alpha)$ 
est un isomorphisme $\hcC^{\oy,(r)}_{X'}$-linéaire 
\begin{equation}\label{higgs3-MF253c}
\hcC^{\oy,(r)}_{X'}\otimes_{\hR'_1}\cN_\ox\stackrel{\sim}{\rightarrow}
(\hcC^{\oy,(r)}_{X'}\hotimes_{\hoR^\oy_{X'}}\fM_{\rho(\oy\rightsquigarrow \ox)})\otimes_{\mZ_p}\mQ_p
\end{equation}
de $\hoR^\oy_{X'}$-modules de Higgs à coefficients dans $\xi^{-1}\tOmega^1_{X/S}(X')$, 
où $\cN_\ox$ est muni du champ de Higgs $\theta_\ox$,
$\hcC^{\oy,(r)}_{X'}$ est muni du champ de Higgs $p^rd_{\hcC^{\oy,(r)}_{X'}}$
et $\fM_{\rho(\oy\rightsquigarrow \ox)}$ est muni du champ de Higgs nul.
Considérons une $\hoR_{X'}^\oy$-augmentation $u\colon \hcC^{\oy,(r)}_{X'}\rightarrow  \hoR_{X'}^\oy$, 
qui existe d'après \eqref{higgs3-RGG16l}.  
En vertu de \ref{higgs3-formel21}, \ref{higgs3-RGG160}(ii), \ref{higgs3-MF252}(i) et \eqref{higgs3-MF253c}, le morphisme canonique 
\[
(\hoR^\oy_{X'}\otimes_{\hcC^{\oy,(r)}_{X'}}(\hcC^{\oy,(r)}_{X'}\hotimes_{\hoR^\oy_{X'}}\fM_{\rho(\oy\rightsquigarrow \ox)}))
\otimes_{\mZ_p}\mQ_p\rightarrow
(\hoR^\oy_{X'}\hotimes_{\hcC^{\oy,(r)}_{X'}}(\hcC^{\oy,(r)}_{X'}\hotimes_{\hoR^\oy_{X'}}\fM_{\rho(\oy\rightsquigarrow \ox)}))
\otimes_{\mZ_p}\mQ_p
\]
est un isomorphisme. Le membre de droite s'identifie canoniquement à $\fM_{\rho(\oy\rightsquigarrow \ox)}\otimes_{\mZ_p}\mQ_p$. 
L'isomorphisme \eqref{higgs3-MF253c} induit donc par changement de base par $u$ un isomorphisme 
\begin{equation}\label{higgs3-MF253e}
\hoR^\oy_{X'}\otimes_{\hR'_1}\cN_\ox\stackrel{\sim}{\rightarrow}
\fM_{\rho(\oy\rightsquigarrow \ox)}\otimes_{\mZ_p}\mQ_p= \cM_{\rho(\oy\rightsquigarrow \ox)}.
\end{equation}
Par suite, le $\hoR_{X'}^\oy[\frac 1 p]$-module $\cM_{\rho(\oy\rightsquigarrow \ox)}$ 
est projectif de type fini en vertu de \ref{higgs3-MF252}(i). 
D'après \ref{higgs3-formel21} et \ref{higgs3-RGG160}(ii), le morphisme canonique 
\begin{equation}\label{higgs3-MF253f}
(\hcC^{\oy,(r)}_{X'}\otimes_{\hoR^\oy_{X'}}\fM_{\rho(\oy\rightsquigarrow \ox)})\otimes_{\mZ_p}\mQ_p\rightarrow
(\hcC^{\oy,(r)}_{X'}\hotimes_{\hoR^\oy_{X'}}\fM_{\rho(\oy\rightsquigarrow \ox)})\otimes_{\mZ_p}\mQ_p
\end{equation}
est donc un isomorphisme~; d'où la proposition.

\begin{lem}\label{higgs3-MF26}
Pour tout $\bvocB_\mQ$-module plat $\cM$ et tout entier $q\geq 0$, 
on a un isomorphisme canonique fonctoriel \eqref{higgs3-coh20}
\begin{equation}
\underset{\underset{r\in \mQ_{>0}}{\longrightarrow}}{\lim}\
\rR^q\top_*(\cM\otimes_{\bvocB_\mQ}\mK^\bullet_\mQ(\bvcC^{(r)},p^r\bvd^{(r)}))\stackrel{\sim}{\rightarrow} \rR^q\top_*(\cM),
\end{equation}
où $\rR^q\top_*(\cM\otimes_{\bvocB_\mQ}\mK^\bullet_\mQ(\bvcC^{(r)},p^r\bvd^{(r)}))$ désigne l'hypercohomologie du foncteur
$\top_*$ \eqref{higgs3-MF1b} par rapport au complexe $\cM\otimes_{\bvocB_\mQ}\mK^\bullet_\mQ(\bvcC^{(r)},p^r\bvd^{(r)})$. 
\end{lem}

En effet, la suite spectrale d'hypercohomologie du foncteur $\top_*$ induit,
pour tout nombre rationnel $r\geq 0$, une suite spectrale fonctorielle
\begin{equation}
{^r\rE}_2^{i,j}=
\rR^{i}\top_*(\cM\otimes_{\bvocB_\mQ}\rH^{j}(\mK^\bullet_\mQ(\bvcC^{(r)},p^r\bvd^{(r)})))\Rightarrow \rR^{i+j}\top_*(\cM\otimes_{\bvocB_\mQ}\mK^\bullet_\mQ(\bvcC^{(r)},p^r\bvd^{(r)})).
\end{equation}
D'après \ref{higgs3-coh24}(iii), pour tous entiers $i\geq 0$ et $j\geq 1$ et tous nombres rationnels $r>r'>0$, le morphisme canonique 
\begin{equation}
{^r\rE}_2^{i,j}\rightarrow {^{r'}\rE}_2^{i,j}
\end{equation}
est nul. On a donc 
\begin{equation}
\underset{\underset{r\in \mQ_{>0}}{\longrightarrow}}{\lim}\ {^r\rE}_2^{i,j}=0.
\end{equation} 
Par ailleurs, il résulte de \ref{higgs3-coh24}(ii) que les morphismes canoniques 
\begin{equation}
\bvocB_\mQ\rightarrow 
\rH^0(\mK^\bullet_\mQ(\bvcC^{(r)},p^r\bvd^{(r)})),
\end{equation}
pour $r\in \mQ_{>0}$, induisent un isomorphisme
\begin{equation}
\rR^i\top_*(\cM)\stackrel{\sim}{\rightarrow} \underset{\underset{r\in \mQ_{>0}}{\longrightarrow}}{\lim}\ {^r\rE}_2^{i,0}.
\end{equation}
Comme les limites inductives filtrantes sont représentables dans $\bMod(\co_\fX[\frac 1 p])$ et qu'elles commutent 
aux limites projectives finies (\cite{sga4} II 4.3), la proposition s'ensuit.

\begin{lem}\label{higgs3-MF28}
Soient $\cN$ un $\co_\fX[\frac 1 p]$-fibré de Higgs à coefficients dans $\xi^{-1}\tOmega_{\fX/\cS}$, $q$ un entier $\geq 0$. 
Notons $\mK^\bullet(\cN)$ le complexe de Dolbeault de $\cN$ {\rm (\cite{ag1} 2.8.2)} et 
pour tout nombre rationnel $r\geq 0$, $\mK^\bullet(\top^{r+}(\cN))$
le complexe de Dolbeault de $\top^{r+}(\cN)$ \eqref{higgs3-MF15}. On a alors un isomorphisme canonique fonctoriel 
\begin{equation}\label{higgs3-MF28a}
\underset{\underset{r\in \mQ_{>0}}{\longrightarrow}}{\lim}\
\rR^q\top_*(\mK^\bullet(\top^{r+}(\cN)))\stackrel{\sim}{\rightarrow} \rH^q(\mK^\bullet(\cN)),
\end{equation}
où $\rR^q\top_*(\mK^\bullet(\top^{r+}(\cN)))$ désigne l'hypercohomologie du foncteur
$\top_*$ \eqref{higgs3-MF1b} par rapport au complexe $\mK^\bullet(\top^{r+}(\cN))$.
\end{lem}

Soient $r$ un nombre rationnel $>0$, $i$ et $j$ deux entiers $\geq 0$. 
Compte tenu de \ref{higgs3-MF6}(i), $\bvd^{(r)}$ \eqref{higgs3-coh20a}
induit un morphisme $\co_\fX$-linéaire
\begin{equation}\label{higgs3-MF28b}
\delta^{j,(r)}\colon \rR^j\top_*(\bvcC^{(r)})\rightarrow \xi^{-1}\tOmega^1_{\fX/\cS}\otimes_{\co_\fX}\rR^j\top_*(\bvcC^{(r)}),
\end{equation}
qui est clairement un $\co_\fX$-champ de Higgs sur $\rR^j\top_*(\bvcC^{(r)})$ à coefficients dans 
$\xi^{-1}\tOmega^1_{\fX/\cS}$. On note $\theta$ le $\co_\fX[\frac 1 p]$-champ de Higgs sur $\cN$
et $\vartheta^{j,(r)}_\tot=\theta\otimes \id+p^r\id \otimes \delta^{j,(r)}$ le 
$\co_\fX[\frac 1 p]$-champ de Higgs total sur $\cN\otimes_{\co_\fX}\rR^j\top_*(\hcC^{(r)})$ (\cite{ag1} 2.8.8).
D'après \ref{higgs3-MF6}(ii), on a un isomorphisme canonique $\co_\fX[\frac 1 p]$-linéaire
\begin{equation}\label{higgs3-MF28c}
\rR^j\top_*(\mK^i(\top^{r+}(\cN)))\stackrel{\sim}{\rightarrow} 
\mK^i(\cN\otimes_{\co_\fX}\rR^j\top_*(\bvcC^{(r)}),\vartheta_\tot^{j,(r)}),
\end{equation} 
compatible avec les différentielles des deux complexes de Dolbeault. 

Par ailleurs, on a une suite spectrale canonique fonctorielle
\begin{equation}\label{higgs3-MF28d}
{^r\rE}_1^{i,j}=\rR^j\top_*(\mK^i(\top^{r+}(\cN)))\Rightarrow \rR^{i+j}\top_*(\mK^\bullet(\top^{r+}(\cN))).
\end{equation}
D'après \ref{higgs3-coh17} et \eqref{higgs3-MF28c}, pour tout $i\geq 0$, on a un isomorphisme canonique 
\begin{equation}\label{higgs3-MF28e}
\underset{\underset{r\in \mQ_{>0}}{\longrightarrow}}{\lim}\ {^r\rE}_1^{i,0}\stackrel{\sim}{\rightarrow}\mK^i(\cN,\theta),
\end{equation}
et pour tout $j\geq 1$, on a 
\begin{equation}\label{higgs3-MF28f}
\underset{\underset{r\in \mQ_{>0}}{\longrightarrow}}{\lim}\ {^r\rE}_1^{i,j}=0.
\end{equation}
De plus, les isomorphismes \eqref{higgs3-MF28e}  (pour $i\in \mN$) forment un isomorphisme de complexes. 
La proposition s'ensuit (\cite{sga4} II 4.3).

\begin{teo}\label{higgs3-MF27}
Soient $\cM$ un $\bvocB_\mQ$-module de Dolbeault, $q$ un entier $\geq 0$. 
Notons $\mK^\bullet(\cH(\cM))$ le complexe de Dolbeault du $\co_\fX[\frac 1 p]$-fibré de Higgs 
$\cH(\cM)$ {\rm (\cite{ag1} 2.8.2)}. On a alors un isomorphisme canonique fonctoriel de $\co_\fX[\frac 1 p]$-modules
\begin{equation}\label{higgs3-MF27a}
\rR^q\top_*(\cM)\stackrel{\sim}{\rightarrow}\rH^q(\mK^\bullet(\cH(\cM))).
\end{equation}
\end{teo}

En effet, $\cH(\cM)$ est un $\co_\fX[\frac 1 p]$-fibré de Higgs soluble associé à $\cM$ en vertu de \ref{higgs3-MF20}.
Choisissons un nombre rationnel $r_\cM>0$ et un isomorphisme de $\Xi^{r_\cM}_\mQ$ 
\begin{equation}\label{higgs3-MF27b}
\alpha_\cM\colon \top^{r_\cM+}(\cH(\cM))\stackrel{\sim}{\rightarrow} \fS^{r_\cM}(\cM)
\end{equation}
vérifiant les propriétés du \ref{higgs3-MF99}. Pour tout nombre rationnel $r$ tel que $0< r< r_\cM$, on désigne par
\begin{equation}\label{higgs3-MF27c}
\alpha^{r}_\cM\colon \top^{r+}(\cH(\cM)) \stackrel{\sim}{\rightarrow}\fS^r(\cM)
\end{equation}
l'isomorphisme de $\Xi^{r}_\mQ$ induit par $\epsilon^{r_\cM,r}(\alpha_\cM)$ et 
les isomorphismes \eqref{higgs3-MF17g} et \eqref{higgs3-MF17h}. D'après la preuve de \ref{higgs3-MF21}, 
$\alpha^r_\cM$ ne dépend que de $\cM$ (mais pas de $\alpha_\cM$)
et il en dépend fonctoriellement. On note $\mK^\bullet(\top^{r+}(\cH(\cM)))$ 
le complexe de Dolbeault de $\top^{r+}(\cH(\cM))$ dans $\bMod_{\mQ}(\bvocB)$ (cf. \ref{higgs3-MF15}). 
Comme $\cM$ est $\bvocB_\mQ$-plat d'après \ref{higgs3-plat3}, $\alpha^{r}_\cM$ induit un isomorphisme \eqref{higgs3-coh20}
\begin{equation}
\mK^\bullet(\top^{r+}(\cH(\cM))) \stackrel{\sim}{\rightarrow}
\cM\otimes_{\bvocB_\mQ}\mK^\bullet_\mQ(\bvcC^{(r)},p^{r}\bvd^{(r)}).
\end{equation}
On en déduit un isomorphisme canonique fonctoriel de $\co_\fX[\frac 1 p]$-modules
\begin{equation}
\underset{\underset{r\in \mQ_{>0}}{\longrightarrow}}{\lim}\ 
\rR^q\top_*(\mK^\bullet(\top^{r+}(\cH(\cM)))) \stackrel{\sim}{\rightarrow}
\underset{\underset{r\in \mQ_{>0}}{\longrightarrow}}{\lim}\ 
\rR^q\top_*(\cM\otimes_{\bvocB_\mQ}\mK^\bullet_\mQ(\bvcC^{(r)},p^{r}\bvd^{(r)})),
\end{equation}
où $\rR^q\top_*(-)$ désigne l'hypercohomologie du foncteur $\top_*$. 
Le théorème s'ensuit compte tenu de \ref{higgs3-MF26} et \ref{higgs3-MF28}.

\section{Modules de Dolbeault sur un petit schéma affine}\label{higgs3-AFE}

\subsection{}\label{higgs3-AFE1}
Les hypothèses et notations générales de § \ref{higgs3-RGG} et § \ref{higgs3-MF} sont en vigueur dans cette section~;
on suppose de plus que $X$ est un objet de $\bQ$ \eqref{higgs3-RGG85}, autrement dit, 
que les conditions suivantes sont satisfaites~:
\begin{itemize}
\item[{\rm (i)}] $X$ est affine et connexe~;
\item[{\rm (ii)}] $f\colon (X,\cM_X)\rightarrow (S,\cM_S)$ admet une carte adéquate \eqref{higgs3-slad3}~;
\item[{\rm (iii)}] il existe une carte fine et saturée $M\rightarrow \Gamma(X,\cM_X)$ 
pour $(X,\cM_X)$ induisant un isomorphisme 
\begin{equation}
M\stackrel{\sim}{\rightarrow}\Gamma(X,\cM_X)/(X,\co_X^\times).
\end{equation}
\end{itemize}
On pose $R=\Gamma(X,\co_X)$, $R_1=R\otimes_{\co_K}\co_\oK$ et 
\begin{equation}
\tOmega^1_{R/\co_K}=\Gamma(X,\tOmega^1_{X/S}).
\end{equation}
On désigne par $\hRun$ le séparé complété $p$-adique de $R_1$,
par $\delta\colon \tE_s\rightarrow \tE$ le plongement canonique \eqref{higgs3-RGG2bc} et par 
\begin{equation}\label{higgs3-AFE1a}
\beta\colon \tE \rightarrow \oX^\circ_\fet
\end{equation}
le morphisme canonique \eqref{higgs3-tfa3a}. Pour tout entier $n\geq 1$, on note  
\begin{equation}\label{higgs3-AFE1b}
\beta_n\colon (\tE_s,\ocB_n) \rightarrow (\oX^\circ_\fet,\ocB_{X,n})
\end{equation}
le morphisme de topos annelés défini par le morphisme de topos $\beta\circ \delta$ et par l'homomorphisme canonique 
$\ocB_{X,n}\rightarrow \beta_*(\ocB_n)$ (cf. \ref{higgs3-RGG2}). On rappelle que ce dernier n'est pas en général un 
isomorphisme \eqref{higgs3-TFT1g}. 
On désigne par $\bvocB_X$ l'anneau $(\ocB_{X,n+1})_{n\in \mN}$ de $(\oX^\circ_\fet)^{\mN^\circ}$ 
et par  
\begin{equation}\label{higgs3-AFE1c}
\bvbeta\colon (\tE_s^{\mN^\circ},\bvocB) \rightarrow ((\oX^\circ_\fet)^{\mN^\circ},\bvocB_X)
\end{equation}
le morphisme de topos annelés  induit par les $(\beta_{n+1})_{n\in \mN}$.

Si $A$ est un anneau et $M$ un $A$-module, on note encore $A$ (resp. $M$) le faisceau constant 
de valeur $A$ (resp. $M$) de $\oX^\circ_\fet$ ou $(\oX^\circ_\fet)^{\mN^\circ}$, selon le contexte.

\begin{prop}\label{higgs3-AFE6}
Pour tout $\co_\fX$-module cohérent $\cN$, on a un isomorphisme $\bvocB$-linéaire canonique et fonctoriel 
\begin{equation}\label{higgs3-AFE6a}
\bvbeta^*(\cN(\fX)\otimes_{\hRun}\bvocB_X)\stackrel{\sim}{\rightarrow}\top^*(\cN).
\end{equation}
\end{prop}

Pour tout entier $n\geq 1$, on pose $\cN_n=\cN/p^n\cN$, que 
l'on considère comme un $\co_{\oX_n}$-module de $X_{s,\et}$ ou $X_\et$, selon le contexte (cf. \ref{higgs3-not2} et \ref{higgs3-TFT9}). 
D'après  \eqref{higgs3-coh0g} et la remarque suivant \eqref{higgs3-not2e}, on a un isomorphisme canonique
\begin{equation}\label{higgs3-AFE6b}
\top^*(\cN)\stackrel{\sim}{\rightarrow}(\sigma^*_{n+1}(\cN_{n+1}))_{n\in \mN}.
\end{equation}
Pour tout entier $n\geq 1$, on a un isomorphisme canonique \eqref{higgs3-TFT8d}
\begin{equation}\label{higgs3-AFE6c}
\sigma_n^*(\cN_n)\stackrel{\sim}{\rightarrow}\sigma^{-1}(\cN_n)\otimes_{\sigma^{-1}(\co_{\oX_n})}\ocB_n.
\end{equation}
Donc en vertu de (\cite{ag2} 5.34(ii), 8.9 et 5.17), le $\ocB_n$-module $\sigma^*_n(\cN_n)$ est 
le faisceau de $\tE$ associé au préfaisceau 
\begin{equation}\label{higgs3-AFE6d}
\{U\mapsto \cN_n(U_s)\otimes_{\co_{\oX_n}(U_s)}\ocB_{U,n}\}, \ \ \ (U\in \ob(\Et_{/X})).
\end{equation}
Pour tout objet affine $U$ de $\Et_{/X}$, on a des isomorphismes canoniques 
\begin{eqnarray}
\cN(U_s)&\stackrel{\sim}{\rightarrow}&\cN(\fX)\otimes_{\hRun}\co_\fX(U_s),\\
\cN_n(U_s)&\stackrel{\sim}{\rightarrow}&\cN(U_s)/p^n\cN(U_s),\label{higgs3-AFE6e}\\
\co_{\oX_n}(U_s)&\stackrel{\sim}{\rightarrow}&\co_\fX(U_s)/p^n\co_\fX(U_s).\label{higgs3-AFE6f}
\end{eqnarray}
Par ailleurs, en vertu de (\cite{ag2} 5.34(i), 8.9 et 5.17), le $\ocB_n$-module 
$\beta_n^*(\cN(\fX)\otimes_{\hRun}\bvocB_{X,n})$ 
est le faisceau de $\tE$ associé au préfaisceau 
\begin{equation}\label{higgs3-AFE6g}
\{U\mapsto \cN(\fX)\otimes_{\hRun}\ocB_{U,n}\}, \ \ \ (U\in \ob(\Et_{/X})).
\end{equation}
D'après \eqref{higgs3-RGG86e}, on en déduit un isomorphisme $\ocB_n$-linéaire canonique et fonctoriel
\begin{equation}\label{higgs3-AFE6h}
\beta^*_n(\cN(\fX)\otimes_{\hRun}\ocB_{X,n})\stackrel{\sim}{\rightarrow}\sigma_n^*(\cN_n).
\end{equation}
La proposition s'ensuit compte tenu de \eqref{higgs3-spsa1g} et \eqref{higgs3-AFE6b}.

\subsection{}\label{higgs3-AFE2}
Soit $r$ un nombre rationnel $\geq 0$. 
On désigne par $\bvcF_X^{(r)}$ le $\bvocB_X$-module $(\cF_{X,n+1}^{(r)})_{n\in \mN}$ 
et par $\bvcC_X^{(r)}$ la $\bvocB_X$-algèbre $(\cC_{X,n+1}^{(r)})_{n\in \mN}$ (cf. \ref{higgs3-RGG22}).
D'après \ref{higgs3-spsa99}(i) et \eqref{higgs3-spsa2a}, on a une suite exacte de $\bvocB_X$-modules 
\begin{equation}\label{higgs3-AFE2a}
0\rightarrow \bvocB_X\rightarrow \bvcF^{(r)}_X\rightarrow 
\xi^{-1}\tOmega^1_{R/\co_K}\otimes_{R}\bvocB_X\rightarrow 0.
\end{equation} 
Compte tenu de \ref{higgs3-spsa99}(i) et \eqref{higgs3-spsa2c}, on a un isomorphisme canonique de $\bvocB_X$-algèbres
\begin{equation}\label{higgs3-AFE2b}
\bvcC^{(r)}_X\stackrel{\sim}{\rightarrow} 
\underset{\underset{m\geq 0}{\longrightarrow}}\lim\ \rS^m_{\bvocB_X}(\bvcF^{(r)}_X).
\end{equation}

Pour tous nombres rationnels $r\geq r'\geq 0$, les morphismes 
$(\tta_{X,n+1}^{r,r'})_{n\in \mN}$ \eqref{higgs3-RGG22g} induisent un morphisme $\bvocB_X$-linéaire 
\begin{equation}\label{higgs3-AFE2c}
\bvtta^{r,r'}_X\colon \bvcF^{(r)}_X\rightarrow \bvcF^{(r')}_X.
\end{equation}
Les homomorphismes 
$(\alpha_{X,n+1}^{r,r'})_{n\in \mN}$ \eqref{higgs3-RGG22h} induisent un homomorphisme de $\bvocB_X$-algèbres 
\begin{equation}\label{higgs3-AFE2d}
\bvalpha^{r,r'}_X\colon \bvcC^{(r)}_X\rightarrow \bvcC^{(r')}_X.
\end{equation}
Pour tous nombres rationnels $r\geq r'\geq r''\geq 0$, on a
\begin{equation}\label{higgs3-AFE2e}
\bvtta^{r,r''}_X=\bvtta^{r',r''}_X \circ \bvtta^{r,r'}_X\ \ \ {\rm et}\ \ \  
\bvalpha^{r,r''}_X=\bvalpha^{r',r''}_X \circ \bvalpha^{r,r'}_X.
\end{equation}

On a un isomorphisme canonique $\bvcC_X^{(r)}$-linéaire
\begin{equation}\label{higgs3-AFE2f}
\Omega^1_{\bvcC_X^{(r)}/\bvocB_X}\stackrel{\sim}{\rightarrow} 
\xi^{-1}\tOmega^1_{R/\co_K}\otimes_R\bvcC_X^{(r)}.
\end{equation}
La $\bvocB_X$-dérivation universelle de $\bvcC_X^{(r)}$ correspond via cet isomorphisme à l'unique 
$\bvocB_X$-dérivation 
\begin{equation}\label{higgs3-AFE2g}
\bvd_X^{(r)}\colon \bvcC_X^{(r)}\rightarrow \xi^{-1}\tOmega^1_{R/\co_K}\otimes_R\bvcC_X^{(r)}
\end{equation}
qui prolonge le morphisme canonique $\bvcF^{(r)}_X\rightarrow \xi^{-1}\tOmega^1_{R/\co_K}\otimes_{R}\bvocB_X$ 
\eqref{higgs3-AFE2a}. Comme 
\begin{equation}
\xi^{-1}\tOmega^1_{R/\co_K}\otimes_R\bvocB_X=\bvd^{(r)}_X(\bvcF_X^{(r)}) 
\subset \bvd_X^{(r)}(\bvcC_X^{(r)}),
\end{equation} 
la dérivation $\bvd_X^{(r)}$ est un $\bvocB_X$-champ de Higgs à coefficients dans 
$\xi^{-1}\tOmega^1_{R/\co_K}$ d'après (\cite{ag1} 2.12). 
Pour tous nombres rationnels $r\geq r'\geq 0$, on a 
\begin{equation}\label{higgs3-AFE2h}
p^{r-r'}(\id \otimes \bvalpha^{r,r'}_X) \circ \bvd_X^{(r)}=\bvd_X^{(r')}\circ \bvalpha^{r,r'}_X.
\end{equation}

\begin{prop}\label{higgs3-AFE3}
Pour tout nombre rationnel $r\geq 0$, les morphismes canoniques
\begin{eqnarray}
\bvbeta^*(\bvcF^{(r)}_X)&\stackrel{\sim}{\rightarrow}&\bvcF^{(r)},\label{higgs3-AFE3a}\\
\bvbeta^*(\bvcC^{(r)}_X)&\stackrel{\sim}{\rightarrow}&\bvcC^{(r)},\label{higgs3-AFE3b}
\end{eqnarray}
sont des isomorphismes. De plus, pour tous nombres rationnels 
$r\geq r'\geq 0$, les morphismes $\bvbeta^*(\bvtta_X^{r,r'})$ et $\bvbeta^*(\bvalpha_X^{r,r'})$ 
s'identifient aux morphismes $\bvtta^{r,r'}$ \eqref{higgs3-RGG180cc} et $\bvalpha^{r,r'}$ \eqref{higgs3-RGG180c}, respectivement.
\end{prop}

Pour tout $U\in \ob(\Et_{/X})$, on note $g_U\colon U\rightarrow X$ le morphisme canonique. 
En vertu de (\cite{ag2} 5.34(i), 8.9 et 5.17), pour tout entier $n\geq 1$, le $\ocB_n$-module $\beta_n^*(\cF^{(r)}_{X,n})$ 
est canoniquement isomorphe au faisceau associé au préfaisceau sur $E$ défini par la correspondance 
\begin{equation}
\{U\mapsto (\ogg^\circ_U)^*_\fet(\cF^{(r)}_{X,n})\otimes_{(\ogg^\circ_U)^*_\fet(\ocB_{X,n})}\ocB_{U,n}\}.
\end{equation}
Pour tout $Y\in \ob(\bQ)$, l'homomorphisme canonique 
\begin{equation}
(\ogg^\circ_Y)^*_\fet(\cF^{(r)}_{X,n})\otimes_{(\ogg^\circ_Y)^*_\fet(\ocB_{X,n})}\ocB_{Y,n}\rightarrow 
\cF^{(r)}_{Y,n}
\end{equation} 
est un isomorphisme en vertu de \ref{higgs3-RGG23}.  Par suite, le morphisme 
\begin{equation}
\beta^*_n(\cF^{(r)}_{X,n})\rightarrow\cF^{(r)}_n,
\end{equation}
adjoint du morphisme canonique $\cF^{(r)}_{X,n}\rightarrow \beta_{n*}(\cF^{(r)}_n)$, est un isomorphisme d'après \eqref{higgs3-RGG86e}.
On en déduit, compte tenu de \eqref{higgs3-spsa1g} et \eqref{higgs3-spsa2a}, que le morphisme canonique \eqref{higgs3-AFE3a}
est un isomorphisme. On démontre de même que l'homomorphisme canonique \eqref{higgs3-AFE3b} est un
isomorphisme~; on peut aussi le déduire de \eqref{higgs3-AFE3a}. La dernière assertion est évidente par adjonction.

\begin{rema}\label{higgs3-AFE4}
Il résulte de \ref{higgs3-AFE3} que pour tout nombre rationnel $r\geq 0$, 
$\bvbeta^*(\bvd_X^{(r)})$ s'identifie à la dérivation $\bvd^{(r)}$ \eqref{higgs3-RGG180b}.
On peut construire l'identification explicitement comme suit. 
D'après la preuve de \ref{higgs3-RGG18}(ii), pour tout entier $n\geq 1$, le diagramme
\begin{equation}\label{higgs3-AFE4c}
\xymatrix{
{\cF^{(r)}_{X,n}}\ar[r]\ar[d]&{\xi^{-1}\tOmega^1_{R/\co_K}\otimes_R\ocB_{X,n}}\ar[d]\\
{\beta_*(\cF^{(r)}_n)}\ar[r]&{\beta_*(\sigma_n^*(\xi^{-1}\tOmega^1_{\oX_n/\oS_n}))}}
\end{equation}
où les flèches verticales sont les morphismes canoniques et  
les flèches horizontales proviennent des suites exactes \eqref{higgs3-RGG22a} et \eqref{higgs3-RGG18a}, est commutatif. 
On en déduit par adjonction un diagramme commutatif
\begin{equation}\label{higgs3-AFE4d}
\xymatrix{
{\beta_n^*(\cF^{(r)}_{X,n})}\ar[r]\ar[d]&{\beta_n^*(\xi^{-1}\tOmega^1_{R/\co_K}\otimes_R\ocB_{X,n})}\ar[d]\\
{\cF^{(r)}_n}\ar[r]&{\sigma_n^*(\xi^{-1}\tOmega^1_{\oX_n/\oS_n})}}
\end{equation}
dont les flèches verticales sont des isomorphismes, d'après les preuves de \ref{higgs3-AFE6} et \ref{higgs3-AFE3}.
Par suite, le diagramme 
\begin{equation}\label{higgs3-AFE4f}
\xymatrix{
{\bvbeta^*(\bvcC_X^{(r)})}\ar[rr]^-(0.5){\bvbeta^*(\bvd^{(r)}_X)}\ar[d]&&
{\bvbeta^*(\xi^{-1}\tOmega^1_{R/\co_K}\otimes_R\bvcC^{(r)}_X)}\ar[d]\\
{\bvcC^{(r)}}\ar[rr]^-(0.5){\bvd^{(r)}}&&{\bvsigma^*(\xi^{-1}\tOmega^1_{\bvoX/\bvoS})\otimes_{\bvocB}\bvcC^{(r)}}}
\end{equation}
où les flèches verticales sont les isomorphismes induits par \eqref{higgs3-AFE6a} et \eqref{higgs3-AFE3b} est commutatif. 
\end{rema}

\subsection{}\label{higgs3-AFE7}
Soit $r$ un nombre rationnel $\geq 0$. 
On désigne par $\bMod(\bvocB_X)$ la catégorie des $\bvocB_X$-modules,
par $\Theta^r$ la catégorie des $p^r$-isoconnexions intégrables 
relativement à l'extension $\bvcC_X^{(r)}/\bvocB_X$ \eqref{higgs3-isoco1} et par $\fS^r_X$ le foncteur 
\begin{equation}\label{higgs3-AFE7a}
\fS^r_X\colon\bMod(\bvocB_X)\rightarrow \Theta^r, \ \ \ 
\cM\mapsto (\bvcC^{(r)}_X\otimes_{\bvocB_X}\cM,\bvcC^{(r)}_X\otimes_{\bvocB_X}\cM,
\id,p^r\bvd^{(r)}_X\otimes \id).
\end{equation}
D'après \ref{higgs3-isoco3}, si $(\cN,\cN',v,\theta)$ est une $\co_\fX$-isogénie de Higgs
à coefficients dans $\xi^{-1}\tOmega^1_{\fX/\cS}$ \eqref{higgs3-MF12},
\begin{equation}\label{higgs3-AFE7b}
(\bvcC^{(r)}_X\otimes_{\hRun}\cN(\fX),\bvcC^{(r)}_X\otimes_{\hRun}\cN'(\fX),\id \otimes_{\hRun}v,
p^r \bvd^{(r)}_X \otimes v+\id \otimes \theta)
\end{equation}
est un objet de $\Theta^r$. On obtient ainsi un foncteur 
\begin{equation}\label{higgs3-AFE7c}
\top^{r+}_X\colon\bIH(\co_\fX,\xi^{-1}\tOmega^1_{\fX/\cS})\rightarrow \Theta^r.
\end{equation}

\begin{prop}\label{higgs3-AFE5}
Pour tout nombre rationnel $r\geq 0$, les diagrammes de foncteurs
\begin{equation}\label{higgs3-AFE5a}
\xymatrix{
{\bMod(\bvocB_X)}\ar[r]^-(0.5){\fS^r_X}\ar[d]_{\bvbeta^*}&{\Theta^r}\ar[d]^{\bvbeta^*}\\
{\bMod(\bvocB)}\ar[r]^-(0.5){\fS^r}&{\Xi^r}}
\end{equation}
\begin{equation}\label{higgs3-AFE5b}
\xymatrix{
{\bIH^\coh(\co_\fX,\xi^{-1}\tOmega^1_{\fX/\cS})}\ar[r]^-(0.5){\top^{r+}_X}\ar[rd]_-(0.5){\top^{r+}}&{\Theta^r}\ar[d]^{\bvbeta^*}\\
&{\Xi^r}}
\end{equation}
où le foncteur image inverse par $\bvbeta$
pour les $p^r$-isoconnexions est défini dans \eqref{higgs3-isoco2}, sont commutatifs à isomorphismes
canoniques près.  
\end{prop}

Cela résulte de \ref{higgs3-AFE6}, \ref{higgs3-AFE3} et \ref{higgs3-AFE4}. 

\begin{prop}\label{higgs3-AFE9}
Soient $\epsilon, r$ deux nombres rationnels tels que $\epsilon>r+\frac{1}{p-1}$ et $r>0$, 
$\cN$ un $\co_\fX$-module cohérent et $\cS$-plat, 
$\theta$ un $\co_\fX$-champ de Higgs sur $\cN$ à coefficients dans $\xi^{-1}\tOmega^1_{\fX/\cS}$ tel que 
\begin{equation}
\theta(\cN)\subset  p^\epsilon \xi^{-1}\tOmega^1_{\fX/\cS}\otimes_{\co_\fX} \cN.
\end{equation}
On note encore $\cN$ la $\co_\fX$-isogénie de Higgs $(\cN,\cN,\id,\theta)$
à coefficients dans $\xi^{-1}\tOmega^1_{\fX/\cS}$. 
Il existe alors un $\bvocB_X$-module adique de type fini
$\cM$ de $(\oX^\circ_\fet)^{\mN^\circ}$ et un isomorphisme de $\Theta^r$
\begin{equation}
\top^{r+}_X(\cN)\stackrel{\sim}{\rightarrow}
\fS^r_X(\cM).
\end{equation}
\end{prop}

On peut clairement supposer $X_s$ non-vide, de sorte que $(X,\cM_X)$ satisfait les hypothèses de (\cite{ag1} 6.2). 
Soit $\oy$ un point géométrique générique de $\oX^\circ$. Comme $\oX$ est localement irréductible \eqref{higgs3-sli1}, 
il est la somme des schémas induits sur ses composantes irréductibles. On note  
$\oX_{\langle \oy\rangle}$ la composante irréductible de $\oX$ contenant $\oy$. De même, $\oX^\circ$
est la somme des schémas induits sur ses composantes irréductibles, 
et $\oX^{\circ}_{\langle \oy\rangle}=\oX_{\langle \oy\rangle}\times_XX^\circ$ 
est la composante irréductible de $\oX^\circ$ contenant $\oy$.
On pose $R^\oy_1=\Gamma(\oX_{\langle \oy\rangle},\co_{\oX})$ et 
$\Delta_\oy=\pi_1(\oX^{\circ}_{\langle \oy\rangle},\oy)$.
On note $\hR^\oy_1$ le séparé complété $p$-adique de $R^\oy_1$, 
$\bB_{\Delta_\oy}$ le topos classifiant de $\Delta_\oy$, 
\begin{equation}\label{higgs3-PMH7b}
\nu_\oy\colon \oX^{\circ}_{\langle \oy\rangle,\fet} \stackrel{\sim}{\rightarrow} \bB_{\Delta_\oy}
\end{equation}
le foncteur fibre en $\oy$ \eqref{higgs3-not6c}, $\oR^\oy_X$ l'anneau défini dans \eqref{higgs3-tfa112c}
et $\hoR^\oy_X$ son séparé complété $p$-adique. On désigne par
$\cC_X^{\oy,(r)}$ la $\hoR^\oy_X$-algèbre définie dans \eqref{higgs3-RGG19b}
par $\hcC_X^{\oy,(r)}$ son séparé complété $p$-adique, par
\begin{equation}
d_{\cC_X^{\oy,(r)}} \colon \cC_X^{\oy,(r)}\rightarrow \xi^{-1}\tOmega^1_{R/\co_K}\otimes_R\cC_X^{\oy,(r)}
\end{equation}
la $\hoR^\oy_X$-dérivation universelle de $\cC_X^{\oy,(r)}$ et par
\begin{equation}
d_{\hcC_X^{\oy,(r)}}\colon \hcC_X^{\oy,(r)}\rightarrow \xi^{-1}\tOmega^1_{R/\co_K}\otimes_R\hcC_X^{\oy,(r)}
\end{equation}
son prolongement aux complétés (on notera que le $R$-module $\tOmega^1_{R/\co_K}$ est libre de type fini). 
Il résulte de \ref{higgs3-tfa113} et des définitions qu'on a des isomorphismes canoniques
\begin{eqnarray}
\nu_\oy(\ocB_X|\oX^{\circ}_{\langle \oy\rangle})&\stackrel{\sim}{\rightarrow}&\oR_X^\oy,\label{higgs3-PMH7c}\\
\nu_\oy(\cC^{(r)}_{X,n}|\oX^{\circ}_{\langle \oy\rangle})&\stackrel{\sim}{\rightarrow}&\cC_X^{\oy,(r)}/p^n \cC_X^{\oy,(r)}.\label{higgs3-PMH7d}
\end{eqnarray}

Les objets $\Delta_\oy$, $R^\oy_1$, $\oR^\oy_X$
et $\cC_X^{\oy,(r)}$ correspondent aux objets $\Delta$, $R_1$, $\oR$ et $\cC^{(r)}$ définis dans \cite{ag1},
en prenant $\tkappa=\oy$ dans (\cite{ag1} 6.7). 
Nous utiliserons dans la suite les constructions de (\cite{ag1} § 13).
Posons $N_\oy=\Gamma(\oX_{\langle \oy\rangle,s},\cN)$ et notons
\begin{equation}
\theta_\oy\colon N_\oy\rightarrow \xi^{-1}\tOmega^1_{R/\co_K}\otimes_RN_\oy
\end{equation}
le $\hR^\oy_1$-champ de Higgs à coefficients dans $\xi^{-1}\tOmega^1_{R/\co_K}$ induit par $\theta$,
qui est $\epsilon$-quasi-petit dans le sens de (\cite{ag1} 13.4).
On lui associe par le foncteur (\cite{ag1} (13.10.10))
une $\hR^\oy_1$-représentation quasi-petite $\varphi_\oy$ de $\Delta_\oy$ sur $N_\oy$.
D'après (\cite{ag1} 13.17), on a un $\hcC^{\oy,(r)}_X$-isomorphisme 
$\Delta_\oy$-équivariant de modules à $p^r$-connexions $p$-adiques relativement à l'extension 
$\hcC^{\oy,(r)}_X/\hoR^\oy$,
\begin{equation}\label{higgs3-higgs12a}
u_\oy\colon N_\oy\hotimes_{\hR^\oy_1} \hcC_X^{\oy,(r)}
\stackrel{\sim}{\rightarrow} N_\oy\hotimes_{\hR^\oy_1} \hcC_X^{\oy,(r)},
\end{equation}
où $\hcC_X^{\oy,(r)}$ est muni de l'action canonique de $\Delta_\oy$ 
et de la $p^r$-connexion $p$-adique $p^r d_{\hcC_X^{\oy,(r)}}$,
le module $N_\oy$ de la source est muni de l'action triviale de $\Delta_\oy$ 
et du $\hR^\oy_1$-champ de Higgs $\theta_\oy$, et le module $N_\oy$ 
du but est muni de l'action $\varphi_\oy$ de $\Delta_\oy$ et du $\hR^\oy_1$-champ de Higgs nul.

Il existe un système projectif $(\cL_{n+1})_{n\in \mN}$ de $R_1$-modules de $\oX^\circ_\fet$ tel que 
pour tout point géométrique générique $\oy$ de $\oX^\circ$, 
on ait un isomorphisme de systèmes projectifs de $R^\oy_1$-représentations de $\Delta_\oy$, 
\begin{equation}
(\nu_\oy(\cL_{n+1}|\oX^{\circ}_{\langle \oy\rangle}))_{n\in \mN}\stackrel{\sim}{\rightarrow} 
(N_\oy/p^{n+1}N_\oy,\varphi_\oy)_{n\in \mN}.
\end{equation}
Ceci résulte aussitôt de la définition du foncteur (\cite{ag1} (13.10.10)) (cf. \cite{ag2} 9.8). 
D'après \ref{higgs3-not60}, $\cL_n$ est de type fini sur $R_1$.  
Pour tous entiers $m\geq n\geq 1$, le morphisme $\cL_m/p^n\cL_m\rightarrow \cL_n$ induit par le
morphisme de transition $\cL_m\rightarrow \cL_n$ est un isomorphisme. 
On pose
\begin{equation}
\cM=(\cL_{n+1}\otimes_{R_1}\ocB_{X,n+1})_{n\in \mN}.
\end{equation} 
C'est un $\bvocB_X$-module adique de type fini de $(\oX^\circ_\fet)^{\mN^\circ}$ en vertu de \ref{higgs3-spad15}. 
Les isomorphismes \eqref{higgs3-higgs12a} induisent un isomorphisme $\bvcC_X^{(r)}$-linéaire 
\begin{equation}
u\colon \cN(\fX)\otimes_{\hRun}\bvcC^{(r)}_X\stackrel{\sim}{\rightarrow}\cM\otimes_{\bvocB_X}\bvcC^{(r)}_X
\end{equation}
tel que le diagramme 
\begin{equation}
\xymatrix{
{\cN(\fX)\otimes_{\hRun}\bvcC^{(r)}_X}\ar[r]^u\ar[d]_{\theta\otimes \id +p^r\id\otimes \bvd_X^{(r)}} 
&{\cM\otimes_{\bvocB_X}\bvcC^{(r)}_X}\ar[d]^{p^r\id\otimes \bvd_X^{(r)}}\\
{\xi^{-1}\tOmega^1_{R/\co_K}\otimes _R\cN(\fX)\otimes_{\hRun}\bvcC^{(r)}_X}\ar[r]^{\id\otimes u}&
{\xi^{-1}\tOmega^1_{R/\co_K}\otimes _R\cM\otimes_{\bvocB_X}\bvcC^{(r)}_X}}
\end{equation}
soit commutatif~; d'où la proposition. 

\begin{cor}\label{higgs3-AFE10}
Sous les hypothèses de \eqref{higgs3-AFE9}, si $\cN_{\mQ_p}$ est un $\co_\fX[\frac 1 p]$-module localement projectif de type fini, 
$(\cN_{\mQ_p},\theta_{\mQ_p})$ est un $\co_\fX[\frac 1 p]$-fibré de Higgs soluble. 
\end{cor}
Cela résulte de \ref{higgs3-AFE5} et \ref{higgs3-AFE9}.

\section{Image inverse d'un module de Dolbeault par un morphisme étale}\label{higgs3-FHT}

\subsection{}\label{higgs3-FHT1}
Les hypothèses et notations générales de § \ref{higgs3-RGG} et § \ref{higgs3-MF} sont en vigueur dans cette section.
Soit, de plus, $g\colon X'\rightarrow X$ un morphisme étale de type fini. 
On munit $X'$ de la structure logarithmique $\cM_{X'}$ image inverse de $\cM_X$  
et on note $f'\colon (X',\cM_{X'})\rightarrow (S,\cM_S)$ le morphisme induit par $f$ et $g$. 
On observera que $f'$ est adéquat \eqref{higgs3-slad6} et que $X'^\circ=X^\circ\times_XX'$ 
est le sous-schéma ouvert maximal de $X'$ où la structure logarithmique $\cM_{X'}$ est triviale.
On munit $\oX'$ et $\coX'$ des structures logarithmiques $\cM_{\oX'}$ et $\cM_{\coX'}$ 
images inverses de $\cM_{X'}$. 
Il existe essentiellement un unique morphisme étale $\tg\colon \tX'\rightarrow \tX$ 
qui s'insère dans un diagramme cartésien~\eqref{higgs3-RGG1}
\begin{equation}\label{higgs3-FHT1c}
\xymatrix{
{\coX'}\ar[r]\ar[d]_{\cog}&{\tX'}\ar[d]^{\tg}\\
{\coX}\ar[r]&{\tX}}
\end{equation}
On munit $\tX'$ de la structure logarithmique $\cM_{\tX'}$ image inverse de $\cM_{\tX}$, de sorte que 
$(\tX',\cM_{\tX'})$ est une $(\cA_2(\oS),\cM_{\cA_2(\oS)})$-déformation lisse de $(\coX',\cM_{\coX'})$.

On associe à $(f',\tX',\cM_{\tX'})$ des objets analogues à ceux définis dans § \ref{higgs3-RGG} et § \ref{higgs3-MF} pour 
$(f,\tX,\cM_{\tX})$, 
qu'on note par les mêmes symboles affectés d'un exposant $^\prime$.  
On désigne par
\begin{eqnarray}
\Phi\colon \tE'&\rightarrow& \tE,\label{higgs3-FHT1d}\\
\Phi_s\colon \tE'_s&\rightarrow& \tE_s,\label{higgs3-FHT1e}
\end{eqnarray}
les morphismes de topos \eqref{higgs3-tfa4b} et \eqref{higgs3-TFT12h} induits par fonctorialité par $g$.
D'après (\cite{ag2} 10.14), $\Phi$ s'identifie au morphisme de localisation de $\tE$ en $\sigma^*(X')$. 
De plus, l'homomorphisme canonique $\Phi^{-1}(\ocB)\rightarrow \ocB'$ est un isomorphisme
en vertu de \ref{higgs3-tfa14}(i).
Pour tout entier $n\geq 1$, $\Phi_s$ est sous-jacent à un morphisme canonique de topos annelés \eqref{higgs3-TFT12j} 
\begin{equation}\label{higgs3-FHT1f}
\Phi_n\colon (\tE'_s,\ocB'_n)\rightarrow (\tE_s,\ocB_n).
\end{equation}
L'homomorphisme $\Phi_s^*(\ocB_n)\rightarrow \ocB'_n$ étant un isomorphisme \eqref{higgs3-TFT13},
il n'y a pas de différence pour les $\ocB_n$-modules entre l'image
inverse par $\Phi_s$ au sens des faisceaux abéliens et l'image inverse par $\Phi_n$ au sens des modules.
Le diagramme de morphismes de topos annelés  \eqref{higgs3-TFT12jj}
\begin{equation}\label{higgs3-FHT1g}
\xymatrix{
{(\tE'_s,\ocB'_n)}\ar[r]^{\Phi_n}\ar[d]_{\sigma'_n}&{(\tE_s,\ocB_n)}\ar[d]^{\sigma_n}\\
{(X'_{s,\et},\co_{\oX'_n})}\ar[r]^{\ogg_n}&{(X_{s,\et},\co_{\oX_n})}}
\end{equation}
où $\ogg_n$ est le morphisme induit par $g$,
est commutatif à isomorphisme canonique près.

\subsection{}\label{higgs3-FHT2}
Tout objet de $E'$ étant naturellement un objet de $E$, on note $\jmath\colon E'\rightarrow E$ le foncteur canonique. 
Celui-ci se factorise à travers une équivalence de catégories 
\begin{equation}\label{higgs3-FHT2a}
E'\stackrel{\sim}{\rightarrow}E_{/(\oX'^\circ\rightarrow X')},
\end{equation}
qui est même une équivalence de catégories au-dessus de $\Et_{/X'}$, 
où l'on considère $E_{/(\oX'^\circ\rightarrow X')}$ comme une ($\Et_{/X'}$)-catégorie par changement 
de base du foncteur fibrant canonique $\pi\colon E\rightarrow \Et_{/X}$ \eqref{higgs3-RGG2d}. 
D'après (\cite{ag2} 5.38), la topologie co-évanescente de $E'$ est induite par celle de $E$ au moyen de $\jmath$. 
Par suite, $\jmath$ est continu et cocontinu (\cite{sga4} III 5.2). De plus, $\Phi$ s'identifie au morphisme
de localisation de $\tE$ en $\sigma^*(X')=(\oX'^\circ\rightarrow X')^a$ (\cite{ag2} 10.14). 
En particulier, $\Phi^*$ n'est autre que le foncteur de restriction par $\jmath$. 
On désigne par $\bQ'$ la sous-catégorie pleine de $\Et_{/X'}$ des objets qui sont dans $\bQ$ \eqref{higgs3-RGG85} et par
\begin{equation}\label{higgs3-FHT2b}
\pi'_{\bQ'}\colon E'_{\bQ'}\rightarrow \bQ'
\end{equation}
le catégorie fibrée déduite par changement de base  du foncteur fibrant canonique $\pi'\colon E'\rightarrow \Et_{/X'}$.
Le foncteur $\jmath$ induit donc un foncteur $\jmath_\bQ\colon E'_{\bQ'}\rightarrow E_\bQ$
qui s'insère dans un diagramme commutatif à isomorphisme canonique près
\begin{equation}\label{higgs3-FHT2c}
\xymatrix{
{E'_{\bQ'}}\ar[r]^{\jmath_\bQ}\ar[d]_{u'}&{E_\bQ}\ar[d]^u\\
E'\ar[r]^\jmath&E}
\end{equation}
où $u$ et $u'$ sont les foncteurs de projection canoniques. 
Les foncteurs $u$ et $u'$ sont pleinement fidèles, et 
la catégorie $E_\bQ$ (resp. $E'_{\bQ'}$) est $\mU$-petite et topologiquement 
génératrice du site $E$ (resp. $E'$) d'après \ref{higgs3-RGG85}. 
Il résulte aussitôt de \eqref{higgs3-FHT2a} que $\jmath_\bQ$ se factorise à travers une équivalence de catégories
\begin{equation}\label{higgs3-FHT2d}
E'_{\bQ'}\stackrel{\sim}{\rightarrow}(E_{\bQ})_{/\hu^*(\oX'^\circ\rightarrow X')},
\end{equation}
où $\hu^*(\oX'^\circ\rightarrow X')$ est le préfaisceau sur $E_\bQ$ déduit de $(\oX'^\circ\rightarrow X')$
par restriction par $u$. On munit $E_\bQ$ (resp. $E'_{\bQ'}$) de la topologie induite par 
celle de $E$ (resp. $E'$). Le foncteur 
$\Phi^*\colon \tE\rightarrow \tE'$ étant essentiellement surjectif, la topologie de $E'_{\bQ'}$ est induite par celle de $E_\bQ$
par $\jmath_\bQ$ (\cite{sga4} II 2.2). Par suite, $\jmath_\bQ$ est continu et cocontinu (\cite{sga4} III 5.2).

\subsection{}\label{higgs3-FHT3}
Soient $Y$ un objet de $\bQ'$ \eqref{higgs3-FHT2} tel que $Y_s\not=\emptyset$, $\oy$ un point géométrique de $\oY^\circ$,
$\oY^\star$ la composante irréductible de $\oY$ contenant $\oy$ \eqref{higgs3-sli3}. 
Considérons les objets associés à $Y$ dans \ref{higgs3-RGG10} et \ref{higgs3-RGG11} relatifs à $(f,\tX,\cM_\tX)$. 
On notera que l'anneau $\oR^\oy_Y$ \eqref{higgs3-tfa112c},
la structure logarithmique $\cM_{\cA_2(\oY^\oy)}$ sur $\cA_2(\oY^\oy)$ \eqref{higgs3-RGG10e}
et le $\hoR^\oy_Y$-module $\rT^\oy_Y$ \eqref{higgs3-RGG10h} ne changent pas que l'on utilise $f$ ou $f'$. 
Remplaçant $(f,\tX,\cM_\tX)$ par $(f',\tX',\cM_{\tX'})$, 
on désigne par $\cL'^{\oy}_{Y}$ le $\trT^\oy_Y$-torseur 
de $\hoY^\oy_\zar$ défini dans \eqref{higgs3-RGG10g}, par $\cF'^\oy_Y$ le 
$\hoR^\oy_Y$-module défini dans \eqref{higgs3-RGG10m} et par $\cC'^\oy_Y$
la $\hoR^\oy_Y$-algèbre définie dans \eqref{higgs3-RGG10n}. 
Soient $U$ un ouvert de Zariski de $\hoY^{\oy}$, $\tU$ l'ouvert de $\cA_2(\oY^\oy)$ défini par $U$. 
Considérons le diagramme commutatif (sans la flèche pointillée)
\begin{equation}\label{higgs3-FHT3a}
\xymatrix{
{(U,\cM_{\hoY^{\oy}}|U)}\ar[r]\ar[d]_{i_Y|U}&{(\coX',\cM_{\coX'})}\ar[r]^{\cog}\ar[d]&{(\coX,\cM_{\coX})}\ar[d]\\
{(\tU,\cM_{\cA_2(\oY^{\oy})}|\tU)}\ar@{.>}[r]^-(0.5)\psi&{(\tX',\cM_{\tX'})}\ar[r]^{\tg}&{(\tX,\cM_{\tX})}}
\end{equation}
Comme $\tg$ est étale, l'application 
\begin{equation}\label{higgs3-FHT3b}
\cL'^\oy_Y(U)\rightarrow \cL^{\oy}_{Y}(U), \ \ \ \psi\mapsto \tg\circ \psi,
\end{equation}
est bijective. On en déduit un isomorphisme $\hoR^\oy_Y$-linéaire et $\pi_1(\oY^{\star \circ},\oy)$-équivariant
\begin{equation}\label{higgs3-FHT3c}
\cF^\oy_Y\stackrel{\sim}{\rightarrow} \cF'^\oy_Y,
\end{equation}
qui s'insère dans un diagramme commutatif 
\begin{equation}\label{higgs3-FHT3d}
\xymatrix{
0\ar[r]&{\hoR^{\oy}_Y}\ar[r]\ar@{=}[d]&{\cF^{\oy}_Y}\ar[r]\ar[d]&
{\xi^{-1}\tOmega^1_{X/S}(Y)\otimes_{\co_X(Y)}\hoR^{\oy}_Y}\ar[r]\ar[d]&0\\
0\ar[r]&{\hoR^{\oy}_{Y}}\ar[r]&{\cF'^{\oy}_{Y}}\ar[r]&
{\xi^{-1}\tOmega^1_{X'/S}(Y)\otimes_{\co_{X'}(Y)}\hoR^{\oy}_{Y}}\ar[r]&0}
\end{equation} 
où les lignes horizontales sont les suites exactes \eqref{higgs3-RGG10m}. 
L'isomorphisme \eqref{higgs3-FHT3c} induit un isomorphisme $\pi_1(\oY^{\star \circ},\oy)$-équivariant de $\hoR^\oy_Y$-algèbres
\begin{equation}\label{higgs3-FHT3e}
\cC^\oy_Y\stackrel{\sim}{\rightarrow} \cC'^\oy_Y.
\end{equation}

Soit $n$ un entier $\geq 1$. Remplaçant $(f,\tX,\cM_\tX)$ par $(f',\tX',\cM_{\tX'})$, 
on désigne par $\cF'_{Y,n}$ le $\ocB'_Y$-module de $\oY^\circ_\fet$ défini dans 
\eqref{higgs3-RGG14a} et par $\cC'_{Y,n}$ la $\ocB'_Y$-algèbre de $\oY^\circ_\fet$ définie dans \eqref{higgs3-RGG14b}.
D'après \ref{higgs3-tfa14}(ii), on a  un isomorphisme canonique d'anneaux de $\oY^\circ_\fet$
\begin{equation}\label{higgs3-FHT3f}
\ocB_Y\stackrel{\sim}{\rightarrow}\ocB'_Y.
\end{equation}
Compte tenu de ce qui précède, on a un isomorphisme $\ocB_Y$-linéaire canonique
\begin{equation}\label{higgs3-FHT3g}
\cF_{Y,n}\stackrel{\sim}{\rightarrow}\cF'_{Y,n}.
\end{equation}
On en déduit un isomorphisme de $\ocB_Y$-algèbres 
\begin{equation}\label{higgs3-FHT3h}
\cC_{Y,n}\stackrel{\sim}{\rightarrow}\cC'_{Y,n}. 
\end{equation}

\subsection{}\label{higgs3-FHT4}
Soient $n$ un entier $\geq 1$, $r$ un nombre rationnel $\geq 0$. 
D'après (\cite{sga4} III 2.3(2)), comme le foncteur canonique $\jmath_\bQ\colon E'_{\bQ'}\rightarrow E_\bQ$ est cocontinu \eqref{higgs3-FHT2}, 
les isomorphismes \eqref{higgs3-FHT3g} induisent un isomorphisme de $\ocB'_n$-modules 
\begin{equation}\label{higgs3-FHT4d}
\rho_n\colon \Phi_n^*(\cF_n)\stackrel{\sim}{\rightarrow}\cF'_n
\end{equation}
qui s'insère dans un diagramme commutatif
\begin{equation}\label{higgs3-FHT4dd}
\xymatrix{
0\ar[r]&{\Phi_n^*(\ocB_n)}\ar[r]\ar@{=}[d]&{\Phi_n^*(\cF_n)}\ar[r]
\ar[d]^{\rho_n}&{\Phi_n^*(\sigma_n^*(\xi^{-1}\tOmega^1_{\oX_n/\oS_n}))}\ar[r]\ar[d]&0\\
0\ar[r]&{\ocB'_n}\ar[r]&{\cF'_n}\ar[r]&{\sigma'^*_n(\xi^{-1}\tOmega^1_{\oX'_n/\oS_n})}\ar[r]&0}
\end{equation}
où les lignes horizontales sont les suites exactes déduites de \eqref{higgs3-RGG18a} (cf. la preuve de \ref{higgs3-RGG18}(ii)). 
De même, les isomorphismes \eqref{higgs3-FHT3h} induisent un isomorphisme de $\ocB'_n$-algèbres
\begin{equation}\label{higgs3-FHT4e}
\gamma_n\colon \Phi_n^*(\cC_n)\stackrel{\sim}{\rightarrow}\cC'_n,
\end{equation}
compatible avec $\rho_n$ via les isomorphismes \eqref{higgs3-RGG18b}. 

D'après \eqref{higgs3-FHT4dd}, $\rho_n$ induit un isomorphisme $\ocB'_n$-linéaire 
\begin{equation}\label{higgs3-FHT4f}
\rho_n^{(r)}\colon \Phi_n^*(\cF^{(r)}_n)\rightarrow \cF'^{(r)}_n
\end{equation}
qui s'insère dans un diagramme commutatif
\begin{equation}\label{higgs3-FHT4ff}
\xymatrix{
0\ar[r]&{\Phi_n^*(\ocB_n)}\ar[r]\ar@{=}[d]&{\Phi_n^*(\cF^{(r)}_n)}\ar[r]
\ar[d]^{\rho_n^{(r)}}&{\Phi_n^*(\sigma_n^*(\xi^{-1}\tOmega^1_{\oX_n/\oS_n}))}\ar[r]\ar[d]&0\\
0\ar[r]&{\ocB'_n}\ar[r]&{\cF'^{(r)}_n}\ar[r]&{\sigma'^*_n(\xi^{-1}\tOmega^1_{\oX'_n/\oS_n})}\ar[r]&0}
\end{equation}
où les lignes horizontales sont les suites exactes déduites de \eqref{higgs3-RGG18a}. 
On en déduit un isomorphisme de $\ocB'_n$-algèbres 
\begin{equation}\label{higgs3-FHT4g}
\gamma_n^{(r)}\colon \Phi_n^*(\cC^{(r)}_n)\stackrel{\sim}{\rightarrow} \cC'^{(r)}_n.
\end{equation}

\subsection{}\label{higgs3-FHT5}
On désigne par $\fX$ (resp. $\fX'$) le schéma formel complété $p$-adique de $\oX$ (resp. $\oX'$), par
\begin{equation}\label{higgs3-FHT5a}
\fgg\colon \fX'\rightarrow \fX
\end{equation} 
le prolongement de $\ogg\colon \oX'\rightarrow \oX$ aux complétés et par
\begin{equation}\label{higgs3-FHT5b}
\bvPhi\colon (\tE'^{\mN^\circ}_s,\bvocB')\rightarrow (\tE^{\mN^\circ}_s,\bvocB)
\end{equation} 
le morphisme de topos annelés induit par les morphismes $(\Phi_n)_{n\geq 1}$ \eqref{higgs3-FHT1f} (cf. \ref{higgs3-spsa35}).  
On note encore 
\begin{equation}
\bvPhi^*\colon \bMod_\mQ(\bvocB)\rightarrow \bMod_\mQ(\bvocB')
\end{equation}
le foncteur induit par l'image inverse par $\bvPhi$. 
D'après \ref{higgs3-TFT17}, $\bvPhi$ est canoniquement isomorphe 
au morphisme de localisation du topos annelé $(\tE_s^{\mN^\circ},\bvocB)$ en $\lambda^*(\sigma^*_s(X'_s))$,
où $\lambda\colon \tE_s^{\mN^\circ}\rightarrow \tE_s$ est le morphisme de topos défini dans \eqref{higgs3-spsa3c}.
Par suite, il n'y a pas de différence pour les $\bvocB$-modules entre l'image
inverse par $\bvPhi$ au sens des faisceaux abéliens et l'image inverse au sens des modules.
Le diagramme de morphismes de topos annelés 
\begin{equation}\label{higgs3-FHT5d}
\xymatrix{
{(\tE'^{\mN^\circ}_s,\bvocB')}\ar[r]^{\bvPhi}\ar[d]_{\top'}&{(\tE^{\mN^\circ}_s,\bvocB)}\ar[d]^\top\\
{(X'_{s,\zar},\co_{\fX'})}\ar[r]^{\fgg}&{(X_{s,\zar},\co_{\fX})}}
\end{equation}
où $\top$ et $\top'$ sont les morphismes de topos annelés définis dans \eqref{higgs3-coh0e},
est commutatif à isomorphisme canonique près \eqref{higgs3-TFT12jc}. Le morphisme canonique 
\begin{equation}\label{higgs3-FHT5da}
\fgg^*(\xi^{-1}\tOmega^1_{\fX/\cS})\rightarrow \xi^{-1}\tOmega^1_{\fX'/\cS}
\end{equation}
étant un isomorphisme, il induit par \eqref{higgs3-FHT5d} un isomorphisme 
\begin{equation}\label{higgs3-FHT5db}
\delta\colon \bvPhi^*(\top^*(\xi^{-1}\tOmega^1_{\fX/\cS}))\stackrel{\sim}{\rightarrow}\top'^*(\xi^{-1}\tOmega^1_{\fX'/\cS}).
\end{equation}

Soit $r$ un nombre rationnel $\geq 0$. Compte tenu de \eqref{higgs3-spsa1g} et \eqref{higgs3-spsa2a}, les isomorphismes 
$(\rho_n^{(r)})_{n\geq 1}$ \eqref{higgs3-FHT4f} induisent un isomorphisme $\bvocB'$-linéaire 
\begin{equation}\label{higgs3-FHT5e}
\bvrho^{(r)}\colon \bvPhi^*(\bvcF^{(r)})\stackrel{\sim}{\rightarrow} \bvcF'^{(r)}.
\end{equation}
De même, les isomorphismes $(\gamma_n^{(r)})_{n\geq 1}$ \eqref{higgs3-FHT4g} 
induisent un isomorphisme de $\bvocB'$-algèbres 
\begin{equation}\label{higgs3-FHT5f}
\bvgamma^{(r)}\colon \bvPhi^*(\bvcC^{(r)})\stackrel{\sim}{\rightarrow} \bvcC'^{(r)}.
\end{equation}
Il résulte aussitôt de \eqref{higgs3-FHT4ff} que le diagramme 
\begin{equation}\label{higgs3-FHT5g}
\xymatrix{
{\bvPhi^*(\bvcC^{(r)})}\ar[r]^-(0.5){\bvgamma^{(r)}}\ar[d]_{\bvPhi^*(\bvd^{(r)})}&{\bvcC'^{(r)}}\ar[d]^{\bvd'^{(r)}}\\
{\bvPhi^*(\top^*(\xi^{-1}\tOmega^1_{\fX/\cS})\otimes_{\bvocB}\bvcC^{(r)})}\ar[r]^-(0.5){\delta\otimes \bvgamma^{(r)}}&
{\top'^*(\xi^{-1}\tOmega^1_{\fX'/\cS})\otimes_{\bvocB'}\bvcC'^{(r)}}}
\end{equation}
où $\bvd^{(r)}$ et $\bvd'^{(r)}$ sont les dérivations \eqref{higgs3-MF15h}, est commutatif. 
Pour tous nombres rationnels $r\geq r'\geq 0$, le diagramme 
\begin{equation}\label{higgs3-FHT5h}
\xymatrix{
{\bvPhi^*(\bvcC^{(r)})}\ar[r]^-(0.5){\bvgamma^{(r)}}\ar[d]_{\bvPhi^*(\bvalpha^{r,r'})}&{\bvcC'^{(r)}}\ar[d]^{\bvalpha'^{r,r'}}\\
{\bvPhi^*(\bvcC^{(r')})}\ar[r]^-(0.5){\bvgamma^{(r')}}&{\bvcC'^{(r')}}}
\end{equation}
où $\bvalpha^{r,r'}$ et  $\bvalpha'^{r,r'}$ sont les homomorphismes canoniques \eqref{higgs3-RGG180c}, est commutatif.

\subsection{}\label{higgs3-FHT6}
On désigne (abusivement) par 
\begin{equation}\label{higgs3-FHT6a}
\fgg^*\colon \bMH(\co_\fX[\frac 1 p],\xi^{-1}\tOmega^1_{\fX/\cS})\rightarrow \bMH(\co_{\fX'}[\frac 1 p],
\xi^{-1}\tOmega^1_{\fX'/\cS})
\end{equation}
le foncteur image inverse pour les modules de Higgs (\cite{ag1} 2.9) induit par $\fgg$ et le morphisme canonique
\eqref{higgs3-FHT5da}. On définit de même un foncteur image inverse \eqref{higgs3-MF12}
\begin{equation}\label{higgs3-FHT6c}
\fgg^*\colon \bIH(\co_\fX,\xi^{-1}\tOmega^1_{\fX/\cS})\rightarrow \bIH(\co_{\fX'},\xi^{-1}\tOmega^1_{\fX'/\cS}).
\end{equation}
Celui-ci induit un foncteur que l'on note encore 
\begin{equation}\label{higgs3-FHT6d}
\fgg^*\colon \bIH_\mQ(\co_\fX,\xi^{-1}\tOmega^1_{\fX/\cS})\rightarrow \bIH_\mQ(\co_{\fX'},\xi^{-1}\tOmega^1_{\fX'/\cS}).
\end{equation}
Le diagramme de foncteurs  
\begin{equation}\label{higgs3-FHT6e}
\xymatrix{
{\bIH_\mQ(\co_\fX,\xi^{-1}\tOmega^1_{\fX/\cS})}\ar[r]\ar[d]_{\fgg^*}&
{\bMH(\co_\fX[\frac 1 p],\xi^{-1}\tOmega^1_{\fX/\cS})}\ar[d]^{\fgg^*}\\
{\bIH_\mQ(\co_{\fX'},\xi^{-1}\tOmega^1_{\fX'/\cS})}\ar[r]&
{\bMH(\co_{\fX'}[\frac 1 p],\xi^{-1}\tOmega^1_{\fX'/\cS})}}
\end{equation}
où les flèches horizontales sont les foncteurs \eqref{higgs3-MF12b} est commutatif à isomorphisme canonique près. 

\subsection{}\label{higgs3-FHT7}
Soit $r$ un nombre rationnel $\geq 0$. 
D'après \ref{higgs3-isoco2} et \eqref{higgs3-FHT5g}, pour tout objet $(\cF,\cG,u,\nabla)$ 
de $\Xi^r$ \eqref{higgs3-MF15}, $(\bvPhi^*(\cF),\bvPhi^*(\cG),\bvPhi^*(u),\bvPhi^*(\nabla))$ s'identifie à un
objet de $\Xi'^r$ au moyen des isomorphismes $\bvgamma^{(r)}$ \eqref{higgs3-FHT5f} et $\delta$ \eqref{higgs3-FHT5db}. On en déduit un foncteur qu'on note encore
\begin{equation}\label{higgs3-FHT7b}
\bvPhi^*\colon \Xi^r\rightarrow \Xi'^r.
\end{equation}
Celui-ci induit un foncteur que l'on note encore
\begin{equation}\label{higgs3-FHT7c}
\bvPhi^*\colon \Xi^{r}_\mQ\rightarrow \Xi'^{r}_\mQ.
\end{equation}

Les diagrammes de foncteurs
\begin{equation}\label{higgs3-FHT7d}
\xymatrix{
{\bMod(\bvocB)}\ar[r]^-(0.5){\fS^r}\ar[d]_{\bvPhi^*}&{\Xi^r}\ar[d]^{\bvPhi^*}\\
{\bMod(\bvocB')}\ar[r]^-(0.5){\fS'^r}&{\Xi'^r}}
\end{equation}
où les flèches horizontales sont les foncteurs \eqref{higgs3-MF15b}, et 
\begin{equation}\label{higgs3-FHT7e}
\xymatrix{
{\bIH(\co_\fX,\xi^{-1}\tOmega^1_{\fX/\cS})}\ar[r]^-(0.5){\top^{r+}}\ar[d]_{\fgg^*}&{\Xi^r}\ar[d]^{\bvPhi^*}\\
{\bIH(\co_{\fX'},\xi^{-1}\tOmega^1_{\fX'/\cS})}\ar[r]^-(0.5){\top'^{r+}}&{\Xi'^r}}
\end{equation}
où les flèches horizontales sont les foncteurs \eqref{higgs3-MF15f}, sont
clairement commutatifs à isomorphismes canoniques près.  
D'après \ref{higgs3-TFT17}, le diagramme de foncteurs
\begin{equation}\label{higgs3-FHT7g}
\xymatrix{
{\Xi^r}\ar[r]^-(0.5){\cK^r}\ar[d]_{\Phi^*}&{\bMod(\bvocB)}\ar[d]^{\bvPhi^*}\\
{\Xi'^r}\ar[r]^-(0.5){\cK'^r}&{\bMod(\bvocB')}}
\end{equation}
où $\cK^r$ et $\cK'^r$ sont les foncteurs \eqref{higgs3-MF15a}, est commutatif à isomorphisme canonique près.

Le morphisme de changement de base relatif au diagramme \eqref{higgs3-FHT5d} induit un morphisme de foncteurs
de $\Xi^r$ dans $\bIH(\co_{\fX'},\xi^{-1}\tOmega^1_{\fX'/\cS})$ 
\begin{equation}\label{higgs3-FHT7h}
\fgg^*\circ \top^r_+\rightarrow \top'^r_+\circ \bvPhi^*,
\end{equation}
où $\top^r_+$ et $\top'^r_+$ sont les foncteurs \eqref{higgs3-MF15d}.
D'après (\cite{sga4} XVII 2.1.3), celui-ci est l'adjoint du morphisme 
\begin{equation}
\top'^{r+}\circ \fgg^*\circ \top^r_+\stackrel{\sim}{\rightarrow} \bvPhi^*\circ \top^{r+} \circ \top^r_+ \rightarrow \bvPhi^*,
\end{equation} 
où la première flèche est l'isomorphisme sous-jacent au diagramme \eqref{higgs3-FHT7e} 
et la seconde flèche est la flèche d'adjonction.
Par suite, pour tout objet $\cN$ de $\bIH(\co_\fX,\xi^{-1}\tOmega^1_{\fX/\cS})$ et tout objet $\cF$ de $\Xi^r$, 
le diagramme d'applications d'ensembles
\begin{equation}\label{higgs3-FHT7hh}
\xymatrix{
{\Hom_{\Xi^r}(\top^{r+}(\cN),\cF)}\ar[r]\ar[d]_a&
{\Hom_{\bIH(\co_\fX,\xi^{-1}\tOmega^1_{\fX/\cS})}(\cN,\top^r_+(\cF))}\ar[d]^b\\
{\Hom_{\Xi'^r}(\top'^{r+}(\fgg^*(\cN)),\bvPhi^*(\cF))}\ar[r]&{\Hom_{\bIH(\co_{\fX'},\xi^{-1}\tOmega^1_{\fX'/\cS})}(\fgg^*(\cN),\top'^r_+(\bvPhi^*(\cF)))}}
\end{equation}
où les flèches horizontales sont les isomorphismes d'adjonction, $a$ est induit par 
le foncteur $\bvPhi^*$ et l'isomorphisme sous-jacent au diagramme \eqref{higgs3-FHT7e}  
et $b$ est induit par le foncteur $\fgg^*$ et le morphisme \eqref{higgs3-FHT7h}, est commutatif. 

Pour tous nombres rationnels $r\geq r'\geq 0$, le diagramme de foncteurs 
\begin{equation}\label{higgs3-FHT7f}
\xymatrix{
{\Xi^r}\ar[r]^-(0.5){\epsilon^{r,r'}}\ar[d]_{\bvPhi^*}&{\Xi^{r'}}\ar[d]^{\bvPhi^*}\\
{\Xi'^r}\ar[r]^-(0.5){\epsilon'^{r,r'}}&{\Xi'^{r'}}}
\end{equation}
où les flèches horizontales sont les foncteurs \eqref{higgs3-MF17d}, est  
commutatif à isomorphisme canonique près.
Il résulte aussitôt de \eqref{higgs3-FHT5h} que le diagramme de morphismes de foncteurs 
\begin{equation}\label{higgs3-FHT7ff}
\xymatrix{
{\fgg^*\circ \top^r_+}\ar[rr]\ar[d]&&{\fgg^*\circ \top^{r'}_+\circ \epsilon^{r,r'}}\ar[d]\\
{\top'^r_+\circ \bvPhi^*}\ar[r]&{\top'^{r'}_+\circ \epsilon^{r,r'}\circ \bvPhi^*}\ar@{=}[r]&{\top'^{r'}_+\circ \bvPhi^*\circ \epsilon^{r,r'}}}
\end{equation}
où  les flèches horizontales sont induites par le morphisme \eqref{higgs3-MF17j}, 
les flèches verticales sont induites par le morphisme \eqref{higgs3-FHT7h}
et l'identification notée avec un symbole $=$ provient du diagramme \eqref{higgs3-FHT7f}, est commutatif.
Par suite, le morphisme composé 
\begin{equation}\label{higgs3-FHT7ii}
\fgg^*\circ \top^r_+\circ \fS^r\rightarrow \top'^r_+\circ \bvPhi^*\circ \fS \stackrel{\sim}{\rightarrow} \top'^r_+\circ \fS'^r\circ \bvPhi^*,
\end{equation}
où la première flèche est induite par \eqref{higgs3-FHT7h} et la seconde flèche est l'isomorphisme sous-jacent au 
diagramme \eqref{higgs3-FHT7d},
induit par passage à la limite inductive, pour $r\in \mQ_{>0}$, un morphisme de foncteurs de $\bMod_\mQ(\bvocB)$
dans $\bMH(\co_{\fX'}[\frac 1 p],\xi^{-1}\tOmega^1_{\fX'/\cS})$
\begin{equation}\label{higgs3-FHT7i}
\fgg^*\circ \cH\rightarrow \cH'\circ \bvPhi^*,
\end{equation} 
où $\cH$ et $\cH'$ sont les foncteurs \eqref{higgs3-MF16a}.

\begin{prop}\label{higgs3-FHT11}
Supposons que $g$ soit une immersion ouverte. Alors~:
\begin{itemize}
\item[{\rm (i)}] Pour tout nombre rationnel $r\geq 0$, le morphisme \eqref{higgs3-FHT7h} est un isomorphisme. 
Il rend commutatif le diagramme de foncteurs
\begin{equation}\label{higgs3-FHT11b}
\xymatrix{
{\Xi^r}\ar[d]_{\bvPhi^*}\ar[r]^-(0.5){\top^r_{+}}&{\bIH(\co_\fX,\xi^{-1}\tOmega^1_{\fX/\cS})}\ar[d]^{\fgg^*}\\
{\Xi'^r}\ar[r]^-(0.5){\top'^r_{+}}&{\bIH(\co_{\fX'},\xi^{-1}\tOmega^1_{\fX'/\cS})}}
\end{equation}
\item[{\rm (ii)}]  Le morphisme \eqref{higgs3-FHT7i} est un isomorphisme. Il rend commutatif le diagramme de foncteurs  
\begin{equation}\label{higgs3-FHT11a}
\xymatrix{
{\bMod_\mQ(\bvocB)}\ar[r]^-(0.5)\cH\ar[d]_{\bvPhi^*}&
{\bMH(\co_\fX[\frac 1 p],\xi^{-1}\tOmega^1_{\fX/\cS})}\ar[d]^{\fgg^*}\\
{\bMod_\mQ(\bvocB')}\ar[r]^-(0.5){\cH'}&{\bMH(\co_{\fX'}[\frac 1 p],\xi^{-1}\tOmega^1_{\fX'/\cS})}}
\end{equation}
\end{itemize}
\end{prop}

(i) Cela résulte de \ref{higgs3-TFT15}. 

(ii) Cela résulte de (i) et des définitions.

\begin{prop}\label{higgs3-FHT8}
Soient $\cM$ un $\bvocB_\mQ$-module de Dolbeault, 
$\cN$ un $\co_\fX[\frac 1 p]$-fibré de Higgs soluble à coefficients dans $\xi^{-1}\tOmega^1_{\fX/\cS}$.
Alors $\bvPhi^*(\cM)$ est un $\bvocB'_\mQ$-module de Dolbeault et $\fgg^*(\cN)$ est un 
$\co_{\fX'}[\frac 1 p]$-fibré de Higgs soluble à coefficients dans $\xi^{-1}\tOmega^1_{\fX'/\cS}$.
Si de plus, $\cM$ et $\cN$ sont associés, $\bvPhi^*(\cM)$  et $\fgg^*(\cN)$ sont associés. 
\end{prop}

En effet, $\fgg^*(\cN)$ est un $\co_{\fX'}[\frac 1 p]$-fibré de Higgs à coefficients dans $\xi^{-1}\tOmega^1_{\fX'/\cS}$
et $\bvPhi^*(\cM)$ est un objet de $\bMod^\atf_\mQ(\bvocB')$. 
Supposons qu'il existe un nombre rationnel $r>0$ et un isomorphisme de $\Xi^r_\mQ$
\begin{equation}\label{higgs3-FHT8a}
\alpha\colon \top^{r+}(\cN)\stackrel{\sim}{\rightarrow}\fS^r(\cM).
\end{equation}
Compte tenu de \eqref{higgs3-FHT7d} et \eqref{higgs3-FHT7e}, $\bvPhi^*(\alpha)$ induit un isomorphisme de $\Xi'^r_\mQ$
\begin{equation}\label{higgs3-FHT8b}
\alpha'\colon \top'^{r+}(\fgg^*(\cN))\stackrel{\sim}{\rightarrow}\fS'^r(\bvPhi^*(\cM));
\end{equation}
d'où la proposition.

\subsection{}\label{higgs3-FHT9}
D'après \ref{higgs3-FHT8}, $\bvPhi^*$ induit un foncteur 
\begin{equation}\label{higgs3-FHT9a}
\bvPhi^*\colon \bMod_\mQ^\Dolb(\bvocB)\rightarrow \bMod_\mQ^\Dolb(\bvocB'),
\end{equation}
et $\fgg^*$ induit un foncteur 
\begin{equation}\label{higgs3-FHT9b}
\fgg^*\colon \bMH^\sol(\co_\fX[\frac 1 p],\xi^{-1}\tOmega^1_{\fX/\cS})\rightarrow \bMH^\sol(\co_{\fX'}[\frac 1 p],
\xi^{-1}\tOmega^1_{\fX'/\cS}).
\end{equation}

\begin{prop}\label{higgs3-FHT10}
{\rm (i)}\ Le diagramme de foncteurs  
\begin{equation}\label{higgs3-FHT10a}
\xymatrix{
{\bMod_\mQ^\Dolb(\bvocB)}\ar[r]^-(0.5)\cH\ar[d]_{\bvPhi^*}&
{\bMH^\sol(\co_\fX[\frac 1 p],\xi^{-1}\tOmega^1_{\fX/\cS})}\ar[d]^{\fgg^*}\\
{\bMod_\mQ^\Dolb(\bvocB')}\ar[r]^-(0.5){\cH'}&{\bMH^\sol(\co_{\fX'}[\frac 1 p],\xi^{-1}\tOmega^1_{\fX'/\cS})}}
\end{equation}
où $\cH$ et $\cH'$ sont les foncteurs \eqref{higgs3-MF20a} est commutatif à isomorphisme canonique près.

{\rm (ii)}\ Le diagramme de foncteurs  
\begin{equation}\label{higgs3-FHT10b}
\xymatrix{
{\bMH^\sol(\co_\fX[\frac 1 p],\xi^{-1}\tOmega^1_{\fX/\cS})}\ar[d]_{\fgg^*}
\ar[r]^-(0.5)\cV&{\bMod_\mQ^\Dolb(\bvocB)}\ar[d]^{\bvPhi^*}\\
{\bMH^\sol(\co_{\fX'}[\frac 1 p],\xi^{-1}\tOmega^1_{\fX'/\cS})}\ar[r]^-(0.5){\cV'}&{\bMod_\mQ^\Dolb(\bvocB')}}
\end{equation}
où $\cV$ et $\cV'$ sont les foncteurs \eqref{higgs3-MF23a} est commutatif à isomorphisme canonique près.
\end{prop}

(i) Pour tout objet $\cM$ de $\bMod^\Dolb_\mQ(\bvocB)$, $\cM$ et $\cH(\cM)$ sont associés 
en vertu de \ref{higgs3-MF20}. Choisissons un nombre rationnel $r_\cM>0$ et un isomorphisme de $\Xi^{r_\cM}_\mQ$ 
\begin{equation}\label{higgs3-FHT10c}
\alpha_\cM\colon \top^{r_\cM+}(\cH(\cM))\stackrel{\sim}{\rightarrow} \fS^{r_\cM}(\cM)
\end{equation}
vérifiant les propriétés du \ref{higgs3-MF99}. Pour tout nombre rationnel $r$ tel que $0< r\leq r_\cM$, on désigne par
\begin{equation}\label{higgs3-FHT10d}
\alpha^r_\cM\colon \top^{r+}(\cH(\cM)) \stackrel{\sim}{\rightarrow}\fS^r(\cM)
\end{equation}
l'isomorphisme de $\Xi^r_\mQ$ induit par $\epsilon^{r_\cM,r}(\alpha_\cM)$ \eqref{higgs3-MF17e} et 
les isomorphismes \eqref{higgs3-MF17g} et \eqref{higgs3-MF17h}. Compte tenu de \eqref{higgs3-FHT7d} et \eqref{higgs3-FHT7e}, 
$\bvPhi^*(\alpha_\cM)$ induit un isomorphisme de $\Xi'^{r_\cM}_\mQ$
\begin{equation}\label{higgs3-FHT10e}
\alpha'_\cM\colon \top'^{r_\cM+}(\fgg^*(\cH(\cM)))\stackrel{\sim}{\rightarrow}\fS'^{r_\cM}(\bvPhi^*(\cM)).
\end{equation}
De même, $\bvPhi^*(\alpha^r_\cM)$ induit un isomorphisme de $\Xi'^r_\mQ$
\begin{equation}\label{higgs3-FHT8ee}
\alpha'^r_\cM\colon \top'^{r+}(\fgg^*(\cH(\cM)))\stackrel{\sim}{\rightarrow}\fS'^r(\bvPhi^*(\cM)),
\end{equation}
qu'on peut aussi déduire de $\epsilon'^{r_\cM,r}(\alpha'_\cM)$ par \eqref{higgs3-FHT7f}.
On désigne par  
\begin{equation}\label{higgs3-FHT10f}
\beta'^r_\cM\colon \fgg^*(\cH(\cM))\rightarrow\top'^{r}_+(\fS'^r(\bvPhi^*(\cM))
\end{equation}
son adjoint \eqref{higgs3-MF15g}. 
D'après \ref{higgs3-FHT8}, $\bvPhi^*(\cM)$ est un $\bvocB'_\mQ$-module de Dolbeault et 
$\fgg^*(\cH(\cM))$ est un $\co_{\fX'}[\frac 1 p]$-fibré de Higgs soluble à coefficients 
dans $\xi^{-1}\tOmega^1_{\fX'/\cS}$, associé à $\bvPhi^*(\cM)$. 
Par suite, en vertu de \ref{higgs3-MF9}(i), le morphisme composé 
\begin{equation}\label{higgs3-FHT10g}
\fgg^*(\cH(\cM))\stackrel{\beta'^r_\cM}{\longrightarrow} \top'^{r}_+(\fS'^r(\bvPhi^*(\cM))\longrightarrow \cH'(\bvPhi^*(\cM)),
\end{equation}
où la seconde flèche est le morphisme canonique \eqref{higgs3-MF16a}, est un isomorphisme 
qui dépend a priori de $\alpha_\cM$ mais pas de $r$.   
D'après la preuve de \ref{higgs3-MF21}, pour tout morphisme $u\colon \cM\rightarrow \cM'$ de $\bMod^\Dolb_\mQ(\bvocB)$ 
et tout nombre rationnel $r$ tel que $0<r<\inf(r_\cM,r_{\cM'})$, le diagramme de $\Xi^r_\mQ$
\begin{equation}
\xymatrix{
{\top^{r+}(\cH(\cM))}\ar[r]^-(0.5){\alpha^r_\cM}\ar[d]_{\top^{r+}(\cH(u))}&{\fS^r(\cM)}\ar[d]^{\fS^r(u)}\\
{\top^{r+}(\cH(\cM'))}\ar[r]^-(0.5){\alpha^r_{\cM'}}&{\fS^r(\cM')}}
\end{equation}
est commutatif. On en déduit que l'isomorphisme composé \eqref{higgs3-FHT10g}
\begin{equation}\label{higgs3-FHT10h}
\fgg^*(\cH(\cM))\stackrel{\sim}{\rightarrow} \cH'(\bvPhi^*(\cM))
\end{equation}
ne dépend que de $\cM$ (mais pas du choix de $\alpha_\cM$) et qu'il en dépend fonctoriellement~; d'où la proposition. 
 
(ii) La preuve est similaire à celle de (i) et est laissée au lecteur.

\begin{remas}\label{higgs3-FHT12}
Soit $\cM$ un $\bvocB_\mQ$-module de Dolbeault.
\begin{itemize}
\item[(i)] Le morphisme canonique \eqref{higgs3-FHT7i}
\begin{equation}\label{higgs3-FHT12a}
\fgg^*(\cH(\cM))\rightarrow \cH'(\bvPhi^*(\cM))
\end{equation}
est un isomorphisme~; c'est l'isomorphisme sous-jacent au diagramme commutatif \eqref{higgs3-FHT10a}.
En effet, reprenons les notations de la preuve de \ref{higgs3-FHT10}(i) et notons, de plus,
\begin{equation}\label{higgs3-FHT12b}
\beta^r_\cM\colon \cH(\cM)\rightarrow \top^r_+(\fS^r(\cM))
\end{equation}
le morphisme adjoint de $\alpha^r_\cM$. Il résulte de \eqref{higgs3-FHT7hh}
que le morphisme $\beta'^r_\cM$ \eqref{higgs3-FHT10f} est égal au composé 
\begin{equation}\label{higgs3-FHT12c}
\xymatrix{
{\fgg^*(\cH(\cM))}\ar[r]^-(0.5){\fgg^*(\beta^r_\cM)}&{\fgg^*(\top^r_+(\fS^r(\cM)))}\ar[r]&
{\top'^r_+(\bvPhi^*(\fS^r(\cM)))}\ar[r]^-(0.5)\sim&{\top'^r_+(\fS'^r(\bvPhi^*(\cM)))}},
\end{equation}
où la deuxième flèche est le morphisme \eqref{higgs3-FHT7h} et la dernière flèche est l'isomorphisme sous-jacent 
au diagramme \eqref{higgs3-FHT7d}. Par ailleurs, la limite inductive des morphismes $\beta^r_\cM$, pour $r\in \mQ_{>0}$, 
est l'identité, et la limite inductive des morphismes $\beta'^r_\cM$, pour $r\in \mQ_{>0}$,
est égale à l'isomorphisme composé \eqref{higgs3-FHT10g}, sous-jacent au diagramme commutatif \eqref{higgs3-FHT10a}. 
\item[(ii)] Soient $r$ un nombre rationnel $>0$, 
\begin{equation}\label{higgs3-FHT12d}
\alpha\colon \top^{r+}(\cH(\cM))\stackrel{\sim}{\rightarrow}\fS^r(\cM)
\end{equation}
un isomorphisme de $\Xi^r_\mQ$ vérifiant les propriétés de \ref{higgs3-MF99}. 
Compte tenu de (i), \eqref{higgs3-FHT7d} et \eqref{higgs3-FHT7e}, on peut identifier $\bvPhi^*(\alpha)$ à un isomorphisme
\begin{equation}\label{higgs3-FHT12e}
\alpha'\colon \top'^{r+}(\cH'(\bvPhi^*(\cM)))\stackrel{\sim}{\rightarrow}\fS'^r(\bvPhi^*(\cM)).
\end{equation}
Par ailleurs, $\bvPhi^*(\cM)$ est un $\bvocB'_\mQ$-module de Dolbeault d'après \ref{higgs3-FHT8}.
Il résulte aussitôt de \ref{higgs3-MF9} que $\alpha'$ vérifie les propriétés de \ref{higgs3-MF99}. 
\end{itemize}
\end{remas}

\section{Catégorie fibrée des modules de Dolbeault}\label{higgs3-PMH}

\subsection{}\label{higgs3-CMD1}
Les hypothèses et notations générales de § \ref{higgs3-RGG} et § \ref{higgs3-MF} sont en vigueur dans cette section.
On note $\psi$ le morphisme composé
\begin{equation}\label{higgs3-CMD1a}
\psi\colon \tE_s^{\mN^\circ}\stackrel{\lambda}{\longrightarrow} \tE_s \stackrel{\sigma_s}{\longrightarrow} X_{s,\et}
\stackrel{\iota_\et}{\longrightarrow} X_\et,
\end{equation}
où $\lambda$ est le morphisme de topos défini dans \eqref{higgs3-spsa3c}, $\sigma_s$ est le morphisme canonique de topos \eqref{higgs3-RGG2e} 
et $\iota\colon X_s\rightarrow X$ est l'injection canonique.
Pour tout objet $U$ de $\Et_{/X}$, on désigne par $f_U\colon (U,\cM_X|U)\rightarrow (S,\cM_S)$ 
le morphisme induit par $f$, et par $\tU\rightarrow \tX$ l'unique morphisme étale 
qui relève $\coU\rightarrow \coX$, de sorte que $(\tU,\cM_{\tX}|\tU)$ est une 
$(\cA_2(\oS),\cM_{\cA_2(\oS)})$-déformation lisse de $(\coU,\cM_{\coX}|\coU)$. 
Le localisé du topos annelé $(\tE_s^{\mN^\circ},\bvocB)$ en $\psi^*(U)$  
est canoniquement équivalent au topos annelé analogue associé à $f_U$ en vertu de \ref{higgs3-TFT17}.
Pour tout nombre rationnel $r\geq 0$, la restriction de la $\bvocB$-algèbre $\bvcC^{(r)}$ 
à $\psi^*(U)$ est canoniquement isomorphe à la $(\bvocB|\psi^*(U))$-algèbre analogue 
associée à la déformation $(\tU,\cM_\tX|\tU)$, d'après \eqref{higgs3-FHT5f}.
On désigne par $\bMod(\bvocB|\psi^*(U))$ la catégorie des $(\bvocB|\psi^*(U))$-modules de 
$(\tE_s^{\mN^\circ})_{/\psi^*(U)}$, par $\bMod_\mQ(\bvocB|\psi^*(U))$ la catégorie des 
$(\bvocB|\psi^*(U))$-modules à isogénie près, par $\Xi^r_U$ la catégorie des $p^r$-isoconnexions intégrables 
relativement à l'extension $(\bvcC^{(r)}|\psi^*(U))/(\bvocB|\psi^*(U))$ (cf. \ref{higgs3-MF15}), par $\Xi^r_{U,\mQ}$ 
la catégorie des objets de $\Xi^r_U$ à isogénie près et par 
$\bMod_\mQ^\Dolb(\bvocB|\psi^*(U))$ la catégorie des $(\bvocB|\psi^*(U))_\mQ$-modules de 
Dolbeault relativement à la déformation $(\tU,\cM_\tX|\tU)$ (cf. \ref{higgs3-MF14}).
D'après \ref{higgs3-FHT8}, pour tout morphisme $g\colon U'\rightarrow U$ de $\Et_{/X}$, le foncteur de restriction
\begin{equation}\label{higgs3-CMD1b}
\bMod_\mQ(\bvocB|\psi^*(U))\rightarrow \bMod_\mQ(\bvocB|\psi^*(U')), \ \ \ \cM\mapsto \cM|\psi^*(U'),
\end{equation}
induit un foncteur 
\begin{equation}\label{higgs3-CMD1c}
\bMod_\mQ^\Dolb(\bvocB|\psi^*(U))\rightarrow \bMod_\mQ^\Dolb(\bvocB|\psi^*(U')).
\end{equation}
On note
\begin{eqnarray}
\Xi^r_{U}\rightarrow \Xi^r_{U'}, && A\mapsto A|\psi^*(U'),\label{higgs3-CMD1d}\\
\Xi^r_{U,\mQ}\rightarrow \Xi^r_{U',\mQ}, && B\mapsto B|\psi^*(U'), \label{higgs3-CMD1e}
\end{eqnarray}
les foncteurs de restriction définis dans \eqref{higgs3-FHT7b} et \eqref{higgs3-FHT7c}, respectivement. 

\begin{lem}\label{higgs3-CMD3}
Soient $r$ un nombre rationnel $\geq 0$, 
$A,B$ deux objets de $\Xi^r_{X,\mQ}$, $(U_i)_{i\in I}$ un recouvrement étale de $X$. 
Pour tous $(i,j)\in I^2$, on pose $U_{ij}=U_i\times_XU_j$. Alors le diagramme d'applications d'ensembles
\begin{equation}\label{higgs3-CMD3a}
\Hom_{\Xi^r_{X,\mQ}}(A,B)\rightarrow \prod_{i\in I}\Hom_{\Xi^r_{U_i,\mQ}}(A|\psi^*(U_i),B|\psi^*(U_i))\rightrightarrows 
\prod_{(i,j)\in I^2}\Hom_{\Xi^r_{U_{ij},\mQ}}(A|\psi^*(U_{ij}),B|\psi^*(U_{ij}))
\end{equation}
est exact.
\end{lem}
En effet, comme $X$ est quasi-compact, on peut supposer $I$ fini, auquel cas l'assertion résulte facilement de \ref{higgs3-caip7}. 

\subsection{}\label{higgs3-CMD2}\index{10001510@$\MOD_\mQ^\Dolb(\bvocB)$}
On note $\MOD(\bvocB)$ la $(\tE_s^{\mN^\circ})$-catégorie fibrée (et même scindée \cite{sga1} VI § 9) 
des $\bvocB$-modules sur $\tE_s^{\mN^\circ}$ (\cite{giraud2} II 3.4.1). 
C'est un champ au-dessus de $\tE_s^{\mN^\circ}$ d'après (\cite{giraud2} II 3.4.4). On désigne par 
$\Et_{\coh/X}$ la sous-catégorie pleine de $\Et_{/X}$ formée des schémas étales 
de présentation finie sur $X$ et par 
\begin{equation}\label{higgs3-CMD2a}
\MOD'(\bvocB)\rightarrow \Et_{\coh/X}
\end{equation}
le changement de base de $\MOD(\bvocB)$ (\cite{sga1} VI § 3) par $\psi^*\circ \varepsilon$,  où
$\psi$ est le morphisme \eqref{higgs3-CMD1a} et $\varepsilon \colon \Et_{\coh/X}\rightarrow X_{\et}$ est le foncteur canonique. 
C'est aussi un champ d'après (\cite{giraud2} II 3.1.1). 
On en déduit une catégorie fibrée 
\begin{equation}\label{higgs3-CMD2c}
\MOD'_\mQ(\bvocB)\rightarrow \Et_{\coh/X},
\end{equation}
dont la fibre au-dessus d'un objet $U$ de $\Et_{\coh/X}$ est la catégorie $\bMod_\mQ(\bvocB|\psi^*(U))$
et le foncteur image inverse par un morphisme 
$U'\rightarrow U$ de $\Et_{\coh/X}$ est le foncteur de restriction \eqref{higgs3-CMD1b}. 
On notera que ce n'est a priori pas un champ. Elle induit une catégorie fibrée 
\begin{equation}\label{higgs3-CMD2d}
\MOD_\mQ^\Dolb(\bvocB)\rightarrow \Et_{\coh/X}
\end{equation}
dont la fibre au-dessus d'un objet $U$ de $\Et_{\coh/X}$ est la catégorie $\bMod_\mQ^\Dolb(\bvocB|\psi^*(U))$
et le foncteur image inverse par un morphisme $U'\rightarrow U$ de $\Et_{\coh/X}$ 
est le foncteur de restriction \eqref{higgs3-CMD1c}.

\begin{prop}\label{higgs3-CMD4}
Soient $\cM$ un objet de $\bMod^{\atf}_\mQ(\bvocB)$, $(U_i)_{i\in I}$ un recouvrement de $\Et_{\coh/X}$. 
Pour que $\cM$ soit de Dolbeault, il faut et il suffit que pour tout $i\in I$, le $(\bvocB|\psi^*(U_i))_\mQ$-module
$\cM|\psi^*(U_i)$ soit de Dolbeault.
\end{prop}

En effet, la condition est nécessaire en vertu de \ref{higgs3-FHT8}. Supposons que pour tout $i\in I$, $\cM|\psi^*(U_i)$ 
soit de Dolbeault et montrons que $\cM$ est de Dolbeault. 
Comme $X$ est quasi-compact, on peut supposer $I$ fini. Pour tout $i\in I$, notons $\fX_i$ le schéma formel
complété $p$-adique de $\oU_i$. Pour tout $(i,j)\in I^2$, posons $U_{ij}=U_i\times_XU_j$ et notons $\fX_{ij}$
le schéma formel complété $p$-adique de $\oU_{ij}$.  
En vertu de \ref{higgs3-FHT10}(i), on a un isomorphisme canonique de $\co_{\fX_{ij}}[\frac 1 p]$-modules de Higgs
à coefficients dans $\xi^{-1}\tOmega^1_{\fX_{ij}/\cS}$
\begin{equation}
\cH_i(\cM|\psi^*(U_i))\otimes_{\co_{\fX_i}}\co_{\fX_{ij}}\stackrel{\sim}{\rightarrow}\cH_{ij}(\cM|\psi^*(U_{ij})),
\end{equation}
où $\cH_i$ et $\cH_{ij}$ sont les foncteurs \eqref{higgs3-MF20a} associés à $(f_{U_i},\tU_i,\cM_\tX|\tU_i)$ et 
$(f_{U_{ij}},\tU_{ij},\cM_\tX|\tU_{ij})$, respectivement.
On en déduit une donnée de descente $\delta$ sur les modules de Higgs $(\cH_i(\cM|\psi^*(U_i)))_{i\in I}$ relativement 
au recouvrement étale $(\fX_i\rightarrow \fX)_{i\in I}$. Celle-ci étant effective d'après \ref{higgs3-formel6}, 
il existe un $\co_\fX[\frac 1 p]$-fibré de Higgs $\cN$ et pour tout $i\in I$, un isomorphisme de 
$\co_{\fX_i}[\frac 1 p]$-modules de Higgs à coefficients dans $\xi^{-1}\tOmega^1_{\fX_i/\cS}$
\begin{equation}
\cN\otimes_{\co_\fX}\co_{\fX_i}\stackrel{\sim}{\rightarrow} \cH_i(\cM|\psi^*(U_i)),
\end{equation}
qui induisent la donnée de descente $\delta$.

Pour tout $(i,j)\in I^2$ et tout nombre rationnel $r>0$, on note $\top^{r+}_i$ et $\fS^{r}_i$ 
(resp. où $\top^{r+}_{ij}$ et $\fS^r_{ij}$) sont les foncteurs 
\eqref{higgs3-MF15d} et \eqref{higgs3-MF15b} associés à $(f_{U_i},\tU_i,\cM_\tX|\tU_i)$ (resp. $(f_{U_{ij}},\tU_{ij},\cM_\tX|\tU_{ij})$). 
Pour tout $i\in I$, choisissons un nombre rationnel $r_i>0$ et un isomorphisme de $\Xi^{r_i}_{U_i,\mQ}$
\begin{equation}
\alpha_i\colon \top^{r_i+}_i(\cH_i(\cM|\psi^*(U_i)))\stackrel{\sim}{\rightarrow} \fS^{r_i}_i(\cM|\psi^*(U_i))
\end{equation}
vérifiant les propriétés du \ref{higgs3-MF99}. Pour tout $(i,j)\in I^2$, $\cM|\psi^*(U_{ij})$ est de Dolbeault en vertu de \ref{higgs3-FHT8}.  
D'après \eqref{higgs3-FHT7d}, \eqref{higgs3-FHT7e} et \ref{higgs3-FHT10}(i), $\alpha_i|\psi^*(U_{ij})$  s'identifie à un  isomorphisme 
\begin{equation}
\alpha_i|\psi^*(U_{ij})\colon \top^{r_i+}_{ij}(\cH_{ij}(\cM|\psi^*(U_{ij})))\stackrel{\sim}{\rightarrow} \fS^{r_i}_{ij}(\cM|\psi^*(U_{ij})).
\end{equation}
Celui-ci vérifie les propriétés du \ref{higgs3-MF99}, compte tenu de \ref{higgs3-FHT12}.
Pour tout nombre rationnel $r$ tel que $0< r\leq r_i$, on désigne par 
$\epsilon_i^{r_i,r}\colon \Xi^{r_i}_{U_i,\mQ}\rightarrow \Xi^r_{U_i,\mQ}$
le foncteur \eqref{higgs3-MF17e} associé à $(f_{U_i},\tU_i,\cM_\tX|\tU_i)$ et par 
\begin{equation}
\alpha^r_i\colon \top^{r+}_i(\cH_i(\cM|\psi^*(U_i)))\stackrel{\sim}{\rightarrow} \fS^{r}_i(\cM|\psi^*(U_i))
\end{equation}
l'isomorphisme de $\Xi^r_{U_i,\mQ}$ induit par $\epsilon_i^{r_i,r}(\alpha_i)$  et 
les isomorphismes \eqref{higgs3-MF17g} et \eqref{higgs3-MF17h}. 
D'après \eqref{higgs3-FHT7d} et \eqref{higgs3-FHT7e}, on peut identifier $\alpha^r_i$ à un isomorphisme 
\begin{equation}
\alpha^r_i\colon \top^{r+}(\cN)|\psi^*(U_i)\stackrel{\sim}{\rightarrow} \fS^{r}(\cM)|\psi^*(U_i).
\end{equation}
Il résulte de la preuve de \ref{higgs3-MF21} que pour tout nombre rationnel $0<r<\inf(r_i,r_j)$, on a dans $\Xi_{U_{ij},\mQ}^r$
\begin{equation}
\alpha^r_i|\psi^*(U_{ij})=\alpha^r_j|\psi^*(U_{ij}).
\end{equation}
En vertu de \ref{higgs3-CMD3}, pour tout nombre rationnel $0<r<\inf(r_i,i\in I)$, 
les isomorphismes $(\alpha^r_i)_{i\in I}$ se recollent en un isomorphisme de $\Xi^r_\mQ$
\begin{equation}
\alpha^r\colon \top^{r+}(\cN)\stackrel{\sim}{\rightarrow}  \fS^r(\cM).
\end{equation}
Par suite, $\cM$ est de Dolbeault.

\begin{prop}\label{higgs3-CMD5}
Les conditions suivantes sont équivalentes~:
\begin{itemize}
\item[{\rm (i)}] La catégorie fibrée \eqref{higgs3-CMD2d}
\begin{equation}\label{higgs3-CMD5a}
\MOD_\mQ^\Dolb(\bvocB)\rightarrow \Et_{\coh/X}
\end{equation}
est un champ {\rm (\cite{giraud2} II 1.2.1)}.
\item[{\rm (ii)}] Pour tout recouvrement $(U_i\rightarrow U)_{i\in I}$ de $\Et_{\coh/X}$, 
notant $\cU$ (resp. pour tout $i\in I$, $\cU_i$) le schéma formel complété $p$-adique de $\oU$ (resp. $\oU_i$), 
pour qu'un $\co_\cU[\frac 1 p]$-fibré de Higgs $\cN$ à coefficients dans $\xi^{-1}\tOmega^1_{\cU/\cS}$
soit soluble, il faut et il suffit que pour tout $i\in I$, le $\co_{\cU_i}[\frac 1 p]$-fibré de Higgs $\cN\otimes_{\co_\cU}\co_{\cU_i}$ 
à coefficients dans $\xi^{-1}\tOmega^1_{\cU_i/\cS}$ soit soluble.
\end{itemize}
\end{prop}

Soit $(U_i\rightarrow U)_{i\in I}$ un recouvrement de $\Et_{\coh/X}$. Pour tout $(i,j)\in I^2$, posons $U_{ij}=U_i\times_XU_j$. 
Notons $\cU$ le schéma formel complété $p$-adique de $\oU$ et 
$\cH^\star$ et $\cV^\star$ les foncteurs \eqref{higgs3-MF20a} et \eqref{higgs3-MF23a} associés à $(f_{U},\tU,\cM_\tX|\tU)$. 
Pour tout $i\in I$, notons $\cU_i$ le schéma formel complété $p$-adique de $\oU_i$ et 
$\cH_i$ et $\cV_i$ les foncteurs \eqref{higgs3-MF20a} et \eqref{higgs3-MF23a} associés à $(f_{U_i},\tU_i,\cM_\tX|\tU_i)$. 

Montrons d'abord (i)$\Rightarrow$(ii). Soit $\cN$ un $\co_\cU[\frac 1 p]$-fibré de Higgs à coefficients dans 
$\xi^{-1}\tOmega^1_{\cU/\cS}$. 
Si $\cN$ est soluble, pour tout $i\in I$, $\cN\otimes_{\co_\cU}\co_{\cU_i}$ est soluble d'après \ref{higgs3-FHT8}. 
Inversement, supposons que pour tout $i\in I$, $\cN\otimes_{\co_\cU}\co_{\cU_i}$ soit soluble et montrons que $\cN$ est soluble. 
Pour tout $i\in I$, $\cM_i=\cV_i(\cN\otimes_{\co_\cU}\co_{\cU_i})$ est un $\bvocB_\mQ|\psi^*(U_i)$-module de Dolbeault.
D'après \ref{higgs3-FHT10}(ii), la donnée de descente canonique sur les fibrés de Higgs $(\cN\otimes_{\co_\cU}\co_{\cU_i})_{i\in I}$ 
relativement au recouvrement étale $(\cU_i\rightarrow \cU)_{i\in I}$
induit une donnée de descente $\delta$ sur les modules de Dolbeault $(\cM_i)_{i\in I}$
relativement au recouvrement $(U_i\rightarrow U)_{i\in I}$.
Cette dernière étant effective d'après (i), il existe un $(\bvocB_\mQ|\psi^*(U))$-module de Dolbeault 
$\cM$ et pour tout $i\in I$, un isomorphisme de $\bvocB_\mQ|\psi^*(U_i)$-modules 
\begin{equation}
\cM|\psi^*(U_i)\stackrel{\sim}{\rightarrow}\cM_i
\end{equation}
qui induisent la donnée de descente $\delta$. 
En vertu de \ref{higgs3-MF21} et \ref{higgs3-FHT10}(i), on a un isomorphisme canonique de $\co_\cU[\frac 1 p]$-fibrés de Higgs 
$\cH^\star(\cM)\stackrel{\sim}{\rightarrow}\cN$. Par suite, $\cN$ est soluble.

Montrons ensuite (ii)$\Rightarrow$(i).  Pour tous $(\bvocB_\mQ|\psi^*(U))$-modules $\cM$ et $\cM'$,  
le diagramme d'applications d'ensembles
\begin{eqnarray}
\lefteqn{\Hom_{\bvocB_\mQ|\psi^*(U)}(\cM,\cM')\rightarrow \prod_{i\in I}
\Hom_{\bvocB_\mQ|\psi^*(U_i)}(\cM|\psi^*(U_i),\cM'|\psi^*(U_i))}\\
&&\rightrightarrows 
\prod_{(i,j)\in I^2}\Hom_{\bvocB_\mQ|\psi^*(U_{ij})}(\cM|\psi^*(U_{ij}),\cM'|\psi^*(U_{ij}))\nonumber
\end{eqnarray}
est exact. En effet, comme $U$ est quasi-compact, on peut supposer $I$ fini, auquel cas l'assertion résulte de \ref{higgs3-caip7}.

Pour tout $i\in I$, soit $\cM_i$ un $(\bvocB_\mQ|\psi^*(U_i))$-module de Dolbeault et soit $\delta$ une donnée de descente
sur $(\cM_i)_{i\in I}$ relativement au recouvrement $(U_i\rightarrow U)_{i\in I}$. Montrons que $\delta$ est effective. 
Par hypothèse, pour tout $i\in I$, 
$\cN_i=\cH_i(\cM_i)$ est un $\co_{\cU_i}[\frac 1 p]$-fibré de Higgs soluble à coefficients dans $\xi^{-1}\tOmega^1_{\cU_i/\cS}$. 
Compte tenu de \ref{higgs3-FHT10}(i), $\delta$ induit une donnée de descente $\gamma$ 
sur les fibrés de Higgs $(\cN_i)_{i\in I}$ relativement
au recouvrement étale $(\cU_i\rightarrow \cU)_{\in I}$. Celle-ci étant effective d'après \ref{higgs3-formel6}, 
il existe un $\co_\cU[\frac 1 p]$-fibré de Higgs $\cN$ et pour tout $i\in I$, un isomorphisme de 
$\co_{\cU_i}[\frac 1 p]$-modules de Higgs 
\begin{equation}
\cN\otimes_{\co_\cU}\co_{\cU_i}\stackrel{\sim}{\rightarrow} \cN_i,
\end{equation}
qui induisent la donnée de descente $\gamma$. D'après (ii), $\cN$ est soluble. Par suite, 
$\cM=\cV^\star(\cN)$ est un $(\bvocB_\mQ|\psi^*(U))$-module de Dolbeault.
D'après \ref{higgs3-MF21} et \ref{higgs3-FHT10}(ii), pour tout $i\in I$, on a un isomorphisme canonique de $(\bvocB_\mQ|\psi^*(U_i))$-modules
$\cM|\psi^*(U_i)\stackrel{\sim}{\rightarrow}\cM_i$, qui induisent la donnée de descente $\delta$, ce qui prouve l'assertion.

\begin{defi}\label{higgs3-PMH1}\index{Fibre de Higgs@$\co_\fX[\frac 1 p]$-fibré de Higgs!(localement) petit}
Soit $(\cN,\theta)$ un $\co_\fX[\frac 1 p]$-fibré de Higgs  à coefficients dans $\xi^{-1}\tOmega^1_{\fX/\cS}$.
\begin{itemize}
\item[(i)] On dit que $(\cN,\theta)$ est {\em petit} s'il existe un sous-$\co_\fX$-module cohérent $\fN$ de $\cN$
qui l'engendre sur $\co_\fX[\frac 1 p]$ et un nombre rationnel $\varepsilon>\frac{1}{p-1}$ tels que 
\begin{equation}\label{higgs3-PMH1a}
\theta(\fN)\subset  p^\varepsilon \xi^{-1}\tOmega^1_{\fX/\cS}\otimes_{\co_\fX} \fN.
\end{equation}
\item[(ii)] On dit que $(\cN,\theta)$ est {\em localement petit} s'il existe un recouvrement ouvert $(U_i)_{i\in I}$ de $X_s$
tel que pour tout $i\in I$, $(\cN|U_i,\theta|U_i)$ soit petit. 
\end{itemize}
\end{defi}

\begin{rema}\label{higgs3-PMH8}
Supposons $X$ affine et le $\co_X$-module $\tOmega^1_{X/S}$ libre de type fini. 
Posons $R=\Gamma(X,\co_X)$ et $R_1=R\otimes_{\co_K}\co_\oK$
et notons $\hR$ et $\hR_1$ leurs séparés complétés $p$-adiques.
Soit  $(\cN,\theta)$ un $\co_\fX[\frac 1 p]$-fibré de Higgs à coefficients dans $\xi^{-1}\tOmega^1_{\fX/\cS}$.
Posons $N=\Gamma(\fX,\cN)$ qui est un $\hR_1[\frac 1 p]$-module projectif de type fini d'après \ref{higgs3-formel20}, et notons encore 
\begin{equation}\label{higgs3-PMH8a}
\theta\colon N\rightarrow \xi^{-1}\tOmega^1_{X/S}(X)\otimes_RN
\end{equation}
le $\hR_1[\frac 1 p]$-champ de Higgs induit par $\theta$. Pour que $(\cN,\theta)$ soit petit, 
il faut et il suffit qu'il existe un sous-$\hR_1$-module de type fini $N^\circ$ de $N$, qui l'engendre sur $\hR_1[\frac 1 p]$, 
et un nombre rationnel $\varepsilon>\frac{1}{p-1}$ tels que 
\begin{equation}
\theta(N^\circ)\subset  p^\varepsilon \xi^{-1}\tOmega^1_{X/S}(X)\otimes_{R} N^\circ.
\end{equation}
En effet, la condition est nécessaire en vertu de (\cite{egr1} (2.10.5.1)) et elle est suffisante compte tenu de \eqref{higgs3-MF12c}
et (\cite{egr1} 1.10.2).
\end{rema}

\begin{prop}\label{higgs3-PMH10}
Tout $\co_\fX[\frac 1 p]$-fibré de Higgs soluble $(\cN,\theta)$ à coefficients dans $\xi^{-1}\tOmega^1_{\fX/\cS}$ est localement petit.
\end{prop}

En effet, on peut se borner au cas où $X$ est affine et le $\co_X$-module 
$\tOmega^1_{X/S}$ est libre de rang $d$ \eqref{higgs3-FHT8}. Montrons alors que $(\cN,\theta)$ est petit.   
Posons $R=\Gamma(X,\co_X)$ et $R_1=R\otimes_{\co_K}\co_\oK$
et notons $\hR$ et $\hR_1$ leurs séparés complétés $p$-adiques. 
Soit $\omega_1,\dots,\omega_d\in \Gamma(X,\xi^{-1}\tOmega^1_{X/S})$ une $\co_X$-base de $\xi^{-1}\tOmega^1_{X/S}$. 
Posons $N=\Gamma(\fX,\cN)$ qui est un $\hR_1[\frac 1 p]$-module projectif de type fini d'après \ref{higgs3-formel20}, et notons encore 
\begin{equation}\label{higgs3-PMH10a}
\theta\colon N\rightarrow \xi^{-1}\tOmega^1_{X/S}(X)\otimes_RN
\end{equation}
le $\hR_1[\frac 1 p]$-champ de Higgs induit par $\theta$. \'Ecrivons 
\begin{equation}\label{higgs3-PMH10b}
\theta=\sum_{i=1}^d \theta_i\otimes \omega_i,
\end{equation}
où les $\theta_i$ sont des $\hR_1[\frac 1 p]$-endomorphismes de $N$ qui commutent deux à deux.

Soient $r$ un nombre rationnel $>0$, $\cM$ un objet de $\bMod^\atf_\mQ(\bvocB)$, 
\begin{equation}\label{higgs3-PMH10c}
\alpha\colon \top^{r+}(\cN) \stackrel{\sim}{\rightarrow}\fS^r(\cM)
\end{equation}
un isomorphisme de $\Xi^r_\mQ$ \eqref{higgs3-MF8}. 
Soient $(\oy\rightsquigarrow \ox)$ un point de $X_\et\gtimes_{X_\et}\oX^\circ_\et$ \eqref{higgs3-tfa41} tel que $\ox$ soit
au-dessus de $s$, $X'$ le localisé strict de $X$ en $\ox$. Reprenons les notations de \ref{higgs3-MF25} et \ref{higgs3-MF250}; 
toutefois, pour mettre en évidence la dépendance en $\ox$, on pose $R_{\ox}=\Gamma(X',\co_{X'})$ et
$R_{\ox,1}=\Gamma(\oX',\co_{\oX'})$ et on note $\hR_{\ox,1}$ le séparé complété $p$-adique de $R_{\ox,1}$
(au lieu de $R'$, $R'_1$ et $\hR'_1$, respectivement). 
D'après \ref{higgs3-MF253}, $\alpha$ induit un isomorphisme $\hcC^{\oy,(r)}_{X'}$-linéaire
\begin{equation}\label{higgs3-PMH10d}
\hcC^{\oy,(r)}_{X'}\otimes_{\hR_{\ox,1}}\cN_\ox\stackrel{\sim}{\rightarrow}
\hcC^{\oy,(r)}_{X'}\otimes_{\hoR^\oy_{X'}}\cM_{\rho(\oy\rightsquigarrow \ox)}
\end{equation}
de $\hoR^\oy_{X'}$-modules de Higgs à coefficients dans $\xi^{-1}\tOmega^1_{X/S}(X')$ \eqref{higgs3-MF250a},
où $\cN_\ox$ est muni du champ de Higgs $\theta_\ox$ \eqref{higgs3-MF250c},
$\hcC^{\oy,(r)}_{X'}$ est muni du champ de Higgs $p^rd_{\hcC^{\oy,(r)}_{X'}}$ \eqref{higgs3-MF250e}
et $\cM_{\rho(\oy\rightsquigarrow \ox)}$ est muni du champ de Higgs nul \eqref{higgs3-MF25f}.

Soit $\cF$ un $\co_\fX$-module cohérent tel que $\cN=\cF_{\mQ_p}$ \eqref{higgs3-formel2},
de sorte que $N=\Gamma(\fX,\cF)\otimes_{\mZ_p}\mQ_p$ (\cite{egr1}  (2.10.5.1)).
En vertu de \ref{higgs3-formel21}, \ref{higgs3-RGG160}(ii) et \ref{higgs3-MF251}(ii),  on a un isomorphisme canonique 
\begin{equation}\label{higgs3-PMH10e}
\cN_\ox \stackrel{\sim}{\rightarrow}N\otimes_{\hR_1}\hR_{\ox,1}.
\end{equation}
On identifiera dans la suite ces deux modules.
Munissons le $\hR_1[\frac 1 p]$-module $N$ (resp. le $\hR_{\ox,1}[\frac 1 p]$-module $\cN_\ox$, resp. 
le $\hoR^\oy_{X'}[\frac 1 p]$-module $\cN_\ox\otimes_{\hR_{\ox,1}}\hoR^\oy_{X'}$) de la topologie $p$-adique (\cite{ag1} 2.2). 
On observera que la topologie $p$-adique de $\cN_\ox$ est induite par celle de 
$\cN_\ox\otimes_{\hR_{\ox,1}}\hoR^\oy_{X'}$. En effet, comme $\cN_\ox$ est projectif de type fini sur $\hR_{\ox,1}[\frac 1 p]$ 
\eqref{higgs3-MF252}, on peut se réduire au cas où $\cN_\ox$ est libre de type fini, 
et même au cas où $\cN_\ox=\hR_{\ox,1}[\frac 1 p]$, pour lequel l'assertion a été établie dans \ref{higgs3-RGG160}(iv). 

Pour tout $z\in N$ et tout $\un=(n_1,\dots,n_d)\in \mN^d$, 
$p^{-r|\un|}(\prod_{1\leq i\leq d}\frac{1}{n_i!}\theta_i^{n_i})(z)$ 
tend vers $0$ dans $\cN_\ox=N\otimes_{\hR_1}\hR_{\ox,1}$ quand $|\un|$ tend l'infini. 
Cela résulte de l'isomorphisme \eqref{higgs3-PMH10d} par la même preuve que (\cite{ag1} 13.24), en tenant compte de \eqref{higgs3-RGG16l}. 

Soient $\pi$ une uniformisante de $\co_K$,
$(x_i)_{i\in I}$ les points génériques de $X_s$ et pour tout $i\in I$, soit $\ox_i$ un point géométrique de $X$
localisé en $x_i$. Comme $X$ est $S$-plat et que $X_s$ est réduit \eqref{higgs3-slad4}, 
pour tout entier $n\geq 1$, l'homomorphisme canonique 
\begin{equation}
R/\pi^nR\rightarrow \prod_{i\in I}R_{\ox_i}/\pi^nR_{\ox_i}
\end{equation}
est injectif. On en déduit que l'homomorphisme canonique 
\begin{equation}
R_1/p^nR_1\rightarrow \prod_{i\in I}R_{\ox_i,1}/p^nR_{\ox_i,1}
\end{equation}
est injectif. L'homomorphisme $\hR_1\rightarrow \prod_{i\in I}\hR_{\ox_i,1}$ est donc injectif. 
Par ailleurs, d'après (\cite{ac} chap.~III §2.11 prop.~14 et cor.~1), on a 
\begin{equation}
\hR_1/p^n\hR_1\simeq R_1/p^nR_1;
\end{equation}
et de même pour les $R_{\ox_i,1}$. 
On en déduit que $\hR_1\cap p^n(\oplus_{i\in I}\hR_{\ox_i,1})=p^n\hR_1$ et par suite que 
\begin{equation}
\hR_1[\frac 1 p]\cap p^n(\oplus_{i\in I}\hR_{\ox_i,1})=p^n\hR_1.
\end{equation}
Comme le $\hR_1[\frac 1 p]$-module $N$ est projectif de type fini, il s'ensuit que la topologie $p$-adique 
de $N$ est induite par la topologie produit des topologies $p$-adiques sur $\prod_{i\in I}N\otimes_{\hR_1}\hR_{\ox_i,1}$. 

Il résulte de ce qui précède que pour tout $z\in N$ et tout $\un=(n_1,\dots,n_d)\in \mN^d$, 
\begin{equation}
p^{-r|\un|}(\prod_{1\leq i\leq d}\frac{1}{n_i!}\theta_i^{n_i})(z)
\end{equation} 
tend vers $0$ dans $N$ quand $|\un|$ tend l'infini. Soient $N_0$ un sous-$\hR_1$-module de type fini de $N$ 
qui l'engendre sur $\hR_1[\frac 1 p]$, 
$\varepsilon$ un nombre rationnel tel que $\frac{1}{p-1}<\varepsilon<r+\frac{1}{p-1}$. 
Comme la suite $p^{(r-\varepsilon)n}n!$ tend vers $0$ dans $\co_C$ quand $n$ tend vers l'infini, pour tout $z\in N$, 
$p^{-\varepsilon |\un|}(\prod_{1\leq i\leq d}\theta_i^{n_i})(z)$ tend vers $0$ dans $N$ quand $|\un|$ tend l'infini.
On peut donc considérer le sous-$\hR_1$-module 
\begin{equation}
N^\circ=\sum_{\un\in \mN^d}p^{-\varepsilon |\un|}(\prod_{1\leq i\leq d}\theta_i^{n_i})(N_0)
\end{equation}
de $N$. Il est de type fini sur $\hR_1$ et il engendre $N$ sur $\hR_1[\frac 1 p]$. Comme on a 
\begin{equation}
\theta(N^\circ)\subset p^\varepsilon \xi^{-1} \tOmega^1_{X/S}(X) \otimes_R N^\circ,
\end{equation} 
$(\cN,\theta)$ est un petit $\co_\fX[\frac 1 p]$-fibré de Higgs à coefficients dans $\xi^{-1}\tOmega^1_{\fX/\cS}$ \eqref{higgs3-PMH8}.

\begin{prop}\label{higgs3-PMH6}
Supposons que $X$ soit un objet de $\bQ$ \eqref{higgs3-RGG85}, autrement dit, 
que les conditions suivantes soient satisfaites~:
\begin{itemize}
\item[{\rm (i)}] $X$ est affine et connexe~;
\item[{\rm (ii)}] $f\colon (X,\cM_X)\rightarrow (S,\cM_S)$ admet une carte adéquate \eqref{higgs3-slad3}~;
\item[{\rm (iii)}] il existe une carte fine et saturée $M\rightarrow \Gamma(X,\cM_X)$ 
pour $(X,\cM_X)$ induisant un isomorphisme 
\begin{equation}
M\stackrel{\sim}{\rightarrow}\Gamma(X,\cM_X)/\Gamma(X,\co_X^\times).
\end{equation}
\end{itemize}
Alors, tout petit $\co_\fX[\frac 1 p]$-fibré de Higgs à coefficients dans $\xi^{-1}\tOmega^1_{\fX/\cS}$ est soluble.
\end{prop}

Cet énoncé est mentionné à titre de rappel \eqref{higgs3-AFE10}.

\begin{cor}\label{higgs3-PMH7}
Si les conditions de \eqref{higgs3-CMD5} sont remplies, 
tout $\co_\fX[\frac 1 p]$-fibré de Higgs localement petit à coefficients dans $\xi^{-1}\tOmega^1_{\fX/\cS}$ est soluble.
\end{cor}

\printindex

\end{document}